\documentclass[a4paper, twoside, 11pt]{article}

\usepackage[a4paper, margin=1in]{geometry}

\usepackage{amssymb, amsfonts, amsmath, amsthm}
\usepackage{hyperref} 
\usepackage{bm} 
\usepackage{graphicx} 
\usepackage{subcaption} 
\usepackage[title,titletoc]{appendix} 
\usepackage{multirow}
\usepackage{indentfirst} 
\usepackage{cleveref} 

\theoremstyle{plain}

\newtheorem{theorem}{Theorem}[section]
\newtheorem{corollary}{Corollary}[section]
\newtheorem{proposition}{Proposition}[section]

\theoremstyle{definition}
\newtheorem{definition}{Definition}[section]

\theoremstyle{remark}
\newtheorem{remark}{Remark}[section]

\numberwithin{equation}{section}
\numberwithin{figure}{section}
\numberwithin{table}{section}

\def\R{\mathbb{R}}

\def\d{\mathrm{d}}

\def\diag{\mathrm{diag}}

\def\B#1{\left\{#1\right\}} 
\def\abs#1{\left\lvert#1\right\rvert}

\allowdisplaybreaks
\title{Finite difference alternative WENO schemes with Riemann invariant-based local characteristic decompositions 
for compressible Euler equations\footnote{Research is supported in part by NSF grant DMS-2309249.}}
\author{
	Yue~Wu\thanks{Division of Applied Mathematics, Brown University, Providence, RI 02912, USA (yue\_wu3@brown.edu).} \ 
	and Chi-Wang~Shu\thanks{Division of Applied Mathematics, Brown University, Providence, RI 02912, USA (chi-wang\_shu@brown.edu).} 
}
\date{} 

\begin{document}

\maketitle


\textbf{Abstract:} The weighted essentially non-oscillatory (WENO) schemes are widely used for hyperbolic conservation laws due to the ability to resolve discontinuities and maintain high-order accuracy in smooth regions at the same time. For hyperbolic systems, the WENO procedure is usually performed on local characteristic variables that are obtained by local characteristic decompositions to avoid oscillation near shocks. However, such decompositions are often computationally expensive. In this paper, we study a Riemann invariant-based local characteristic decomposition for the compressible Euler equations that reduces the cost. We apply the WENO procedure to the local characteristic fields of the Riemann invariants, where the eigenmatrix is sparse and thus the computational cost can be reduced. It is difficult to obtain the cell averages of Riemann invariants from those of the conserved variables due to the nonlinear relation between them, so we only focus on the finite difference alternative WENO versions. The efficiency and non-oscillatory property of the proposed schemes are well demonstrated by our numerical results. 



\section{Introduction}


It is well known that nonlinear hyperbolic conservation laws, e.g. the compressible Euler equations for gas dynamics, can generate sharp discontinuities in finite time even with smooth initial conditions. Numerically resolving these discontinuities while maintaining high-order accuracy in smooth regions is a challenging task, because spurious oscillations that appear near the discontinuities may mess/blow up the simulation. Inspired by the essentially non-oscillatory (ENO) schemes \cite{10.1016/0021-99918790031-3, 0021999188901775, 10.1016/0021-99918990222-2}, the weighted ENO (WENO) schemes, originally proposed in \cite{10.1006/jcph.1994.1187} and improved to achieve higher order in \cite{10.1006/jcph.1996.0130}, perform well in this regard. Both finite volume and finite difference versions of the WENO schemes have been developed and widely used in practice, and many other variants are proposed as well. For simplicity, we only focus on the WENO-JS version \cite{10.1006/jcph.1996.0130} in this paper. Besides the classical finite difference WENO formulation that is applied to point values of split flux functions, there is also one formulation called alternative WENO (A-WENO) \cite{10.1137/120889885} that works on those of conserved variables, which roots back to early studies of the ENO method \cite{0021999188901775}. 

For scalar hyperbolic conservation laws, the application of WENO methods is straightforward, and the oscillations are well suppressed. However, the situation is more complicated for systems. The component-wise ENO/WENO schemes often produce oscillatory solutions after discontinuities where multiple shocks interact with each other \cite{10.1016/0021-99918790031-3}, such as Riemann problems. To avoid this, the local characteristic decomposition (LCD) procedure should be involved. Both finite volume \cite{10.1016/0021-99918790031-3} and finite difference \cite{10.1016/0021-99918990222-2} versions include the multiplication of the left and right eigenmatrices with corresponding vectors. For a 1D system consisting of \(m\) components, the extra cost at each cell interface of a \(r\)th-order ENO (or a \((2r-1)\)th-order WENO) scheme is \(4r+2\) multiplications of an \(m \times m\)-matrix and an \(m\)-vector, compared to the component-wise version. This costs \((4r+2) m^2\) floating point multiplications in a vanilla implementation. Various experiments have shown that the LCD does take a big portion of the computational time \cite{10.1016/0168-92749290066-M}. 
There have been some attempts to alleviate this extra cost, but none of them works well for the compressible Euler equations in high-order simulations. For example, Jiang and Shu \cite{10.1006/jcph.1996.0130} proposed to use the pressure and the entropy to compute the nonlinear WENO weights, but oscillations can occur near strong shocks. The second-order component-wise central ACM-MUSCL scheme \cite{10.1016/0021-99919090260-8} and the third-order component-wise central PPM method \cite{10.1007/s002110050345} are both non-oscillatory, whereas the LCD is still necessary to suppress oscillations in higher-order central schemes \cite{10.1006/jcph.2002.7191}. Xu and Shu \cite{10.1137/22M1536479} circumvented the LCD procedure by applying the WENO interpolation on the coordinate system of Riemann invariants, if the equations admit such a system. This reduces the per-node cost (other than the WENO part) dependency upon the stencil width \(2r\) from linear to constant, while maintaining non-oscillatory and high-order. Unfortunately, the compressible Euler equations do not belong to this category. 
In this paper, we propose to use the Riemann invariants to improve the sparsity of the eigenmatrix, so as to reduce the cost of the LCD procedure for the compressible Euler equations. Although we cannot completely eliminate the linear dependency upon \(r\) of the per-node LCD cost, a big improvement is still available by bringing the matrix as close to the identity as possible. Because of the nonlinear relation between the conserved variables and the Riemann invariants, it is difficult to obtain the cell averages of the latter from those of the former, so we only concentrate on the finite difference A-WENO versions, which are based on the point values of the conserved variables. 

The rest of the paper is organized as follows. In \Cref{se:pre}, we introduce the Riemann invariants and the compressible Euler equations. In \Cref{se:nm}, we review the A-WENO scheme and propose the Riemann invariant-based LCD. In \Cref{se:nr}, we present our numerical results. Finally, we conclude in \Cref{se:con}.

\section{Preliminaries}\label{se:pre}

\subsection{Riemann invariants}

Consider a system of conservation laws of \(m\) components in one space dimension 
\begin{equation}\label{eq: 1D conservation law}
    \bm u_t + \bm f(\bm u)_x = \bm 0, 
\end{equation}
where \(\bm u = (u_1, u_2, \ldots, u_m)^T \in \R^{m}\) are the conserved variables and \(\bm f(\bm u) = (f_1(\bm u), f_2(\bm u), \ldots, f_m(\bm u))^T \in \R^m\) are the flux functions. We are only interested in the case where \(\bm f(\bm u)\) is a \(C^1\) function of \(\bm u\). Assuming all the functions are differentiable, the system \eqref{eq: 1D conservation law} can be rewritten into a general quasilinear form, 
\begin{equation}\label{eq: 1D quasilinear system}
    \bm u_t + \bm A(\bm u) \bm u_x = \bm 0,
\end{equation}
where \(\bm A(\bm u) = \frac{\partial \bm f}{\partial \bm u}\) is the Jacobian matrix.

\begin{definition}[hyperbolicity in 1D]\label{def:hyperbolicity-in-1D}
    The quasilinear system \eqref{eq: 1D quasilinear system} is hyperbolic if the Jacobian matrix \(\bm A(\bm u)\) is diagonalizable for all \(\bm u\), i.e. there exists an invertible matrix \(\bm R(\bm u)\) such that \(\bm A(\bm u) = \bm R(\bm u) \bm\Lambda(\bm u) \bm R(\bm u)^{-1}\), where \(\bm\Lambda(\bm u) = \diag(\lambda_1(\bm u), \lambda_2(\bm u), \ldots, \lambda_m(\bm u))\) is the diagonal matrix of eigenvalues, \(\bm R(\bm u) = (\bm r_1(\bm u), \bm r_2(\bm u), \ldots, \bm r_m(\bm u))^T\) is the right eigenmatrix, and \(\bm R(\bm u)^{-1} = \bm L(\bm u) = (\bm l_1(\bm u)^T, \bm l_2(\bm u)^T, \cdots, \bm l_m(\bm u)^T)^T\) is the left eigenmatrix. 

    Correspondingly, \eqref{eq: 1D conservation law} is called a system of hyperbolic conservation laws, if the Jacobian matrix \(\frac{\partial \bm f(\bm u)}{\partial \bm u}\) is diagonalizable for all \(\bm u\). 
\end{definition}

\begin{definition}[hyperbolicity in multi-D]
    For a system of conservation laws in \(d\) spatial dimensions, 
    \begin{equation}
        \bm u_t + \sum_{i=1}^d \bm f_{i}(\bm u)_{x_{i}} = \bm 0, 
    \end{equation}
    it is hyperbolic if for any unit vector \(\bm n = (n_1, n_2, \ldots, n_d)\), the following 1D equation 
    \begin{equation}
        \bm u_t + \sum_{i=1}^d n_i \bm f_{i}(\bm u)_{x} = \bm 0 
    \end{equation}
    is hyperbolic. 
\end{definition}

\begin{remark}
    In the following, we will only state the results for 1D systems, but they can be easily extended to multi-D by applying it to each direction. 
\end{remark}

\begin{definition}[Riemann invariant]
    A scalar-valued function \(w(\bm u)\) is called an \(i\)-Riemann invariant for the hyperbolic system \eqref{eq: 1D quasilinear system} if \(\nabla w(\bm u) \cdot \bm r_i(\bm u) = 0\) for any \(\bm u\). 
\end{definition}

\begin{corollary}[closedness under composition of differentiable functions]\label{cor:closedness-under-composition}
    If \(w(\bm u)\) is an \(i\)-Riemann invariant, then for any \(C^1\) function \(g: \R \mapsto \R\), \(g(w(\bm u))\) is also an \(i\)-Riemann invariant. Moreover, a differentiable function of several \(i\)-Riemann invariants is also an \(i\)-Riemann invariant. 
\end{corollary}

\begin{theorem}[jump of Riemann invariants {\cite[Theorem 17.16]{10.1007/978-1-4612-0873-0}}]
    For a system of hyperbolic conservation laws \eqref{eq: 1D conservation law}, the change of an \(i\)-Riemann invariant \(w(\bm u)\) across an isolated \(i\)-rarefaction or \(i\)-contact is \(0\), whereas the change across an isolated \(i\)-shock is \(\mathcal{O}(\epsilon^3)\) as \(\epsilon \to 0\), where \(\epsilon = \abs{\lambda_i(\bm u_l) - \lambda_i(\bm u_r)}\) is the shock strength. 
\end{theorem}

\begin{proposition}[change of variables]
    Consider a hyperbolic system \eqref{eq: 1D quasilinear system} and assume all the functions are differentiable. If there exists a differentiable and invertible transform \(\bm v(\bm u): \R^m \mapsto \R^m\) for all \(\bm u\), then it can be rewritten as another hyperbolic system, 
    \begin{equation}\label{eq: transformed semilinear system}
        \bm v_t + \tilde{\bm A}(\bm v) \bm v_x = \bm 0,
    \end{equation}
    where \(\tilde{\bm A} = \frac{\partial \bm v}{\partial \bm u} \bm A \frac{\partial \bm u}{\partial \bm v}\). Moreover, \(\tilde{\bm A}\) have the same eigenvalues as \(\bm A\), and the corresponding eigenmatrices are \(\tilde{\bm R} = \frac{\partial \bm v}{\partial \bm u} \bm R\) and \(\tilde{\bm L} = \bm L \frac{\partial \bm u}{\partial \bm v}\).
\end{proposition}

\begin{corollary}[sparsity of the eigenmatrix]\label{cor:sparsity-of-eigenmatrix}
    Assume there exists a differentiable and invertible transform \(\bm v(\bm u) = (v_1(\bm u), v_2(\bm u), \ldots, v_m(\bm u))^T\) for all \(\bm u\), where \(v_i(\bm u)\) is a \(j\)-Riemann invariant, then for the transformed system \eqref{eq: transformed semilinear system}, the \((i,j)\)th entry of \(\tilde{\bm R}\) is \(0\). 
\end{corollary}

\begin{remark}[coordinate system of Riemann invariants \cite{10.1137/22M1536479}]
    If \(v_i(\bm u)\) is a \((1,\ldots,i-1,i+1,\ldots,m)\)-Riemann invariant for all \(i\)'s, the system \eqref{eq: 1D conservation law} is endowed with a coordinate system of Riemann invariants which is exactly \(\bm v(\bm u)\). Then for the transformed system \eqref{eq: transformed semilinear system}, the eigenmatrix \(\tilde{\bm R}\) is diagonal and so is \(\tilde{\bm A}\). In this case we can even take \(\tilde{\bm R} = \bm I\). 
\end{remark}

\subsection{The compressible Euler equations}

The 2D compressible Euler equations are 
\begin{subequations}\label{eq: Euler 2D}
    \begin{equation}
        \bm u_t + \bm f(\bm u)_x + \bm g(\bm u)_y = \bm 0, 
    \end{equation}
    where the conserved variables are
    \begin{equation}
        \bm u = \left(\rho, \rho u, \rho v, E\right)^T, 
    \end{equation}
    and the fluxes are 
    \begin{equation}
        \begin{aligned}
            & \bm f(\bm u) = \left(\rho u, \rho u^2 + p, \rho u v, u (E + p)\right)^T, \\
            & \bm g(\bm u) = \left(\rho v, \rho v u, \rho v^2 + p, v (E + p)\right)^T.
        \end{aligned}
    \end{equation}
\end{subequations}
Here, \(\rho\) is the density, \(\vec{u} = (u, v)\) is the velocity, \(E = \frac{\rho}{2} \abs{\vec{u}}^2 + \rho e\) is the total energy, \(e\) is the specific internal energy, and \(p = p(\rho, e)\) is the pressure determined by the equation of state (EOS). Additional derived variables include the total enthalpy \(H = \frac{E + p}{\rho}\), the sonic speed \(c = \sqrt{\left.\frac{\partial p}{\partial \rho}\right|_{s}}\) and the specific entropy \(s\). We consider the simplest \(\gamma\)-law EOS for ideal polytropic gas, 
\begin{equation}\label{eq: EOS}
    p = \rho e (\gamma - 1), 
\end{equation}
with \(\gamma = \frac{c_p}{c_v} > 1\) being the specific heat capacity ratio which is constant. Under this EOS, we have \(c = \sqrt{\frac{\gamma p}{\rho}}\) and \(s = \ln(p \rho^{-\gamma})\). We use \(S = \exp(s) = p \rho^{-\gamma}\) to denote the entropy. 

We list the eigenstructure of the Jacobian matrix in \eqref{eq:Euler-eigen-values}, \eqref{eq:Euler-right-eigen} and \eqref{eq:Euler-left-eigen}. For the 1D version of the equations, the eigenmatrices can be obtained by first restricting \(v = 0\) and then removing the third row and column of the 2D eigenmatrices. By simple algebraic manipulations, it is easy to check that the equations have the following Riemann invariants corresponding to the eigenmatrix in \eqref{eq:Euler-right-eigen}. 
\begin{proposition}[Riemann invariants of the 2D compressible Euler equations]
    The 2D compressible Euler equations \eqref{eq: Euler 2D} with the \(\gamma\)-law EOS \eqref{eq: EOS} has the following Riemann invariants in the \(x\)-direction. 
    \begin{itemize}
        \item 1-Riemann invariants: \(S\), \(u + \frac{2}{\gamma-1}c\), \(v\). 
        \item 2-Riemann invariants: \(u\), \(p\), \(v\). 
        \item 3-Riemann invariants: \(S\), \(u\), \(p\). 
        \item 4-Riemann invariants: \(S\), \(u - \frac{2}{\gamma-1}c\), \(v\). 
    \end{itemize}
\end{proposition}

\section{Numerical methods}\label{se:nm}

\subsection{The A-WENO scheme with transform variables}

Consider the 1D hyperbolic system \eqref{eq: 1D conservation law} on a finite interval \([a, b]\). We use the uniform grid \(I_j = [x_{j-\frac{1}{2}}, x_{j+\frac{1}{2}}]\) for \(j = 1,2,\ldots,N\), where \(x_{j} = a + (j - \frac{1}{2}) \Delta x\) and \(\Delta x = \frac{b - a}{N}\). The semi-discrete \((2r-1)\)th-order finite difference A-WENO scheme is 
\begin{equation}
    \frac{\d}{\d t} \bm u_j + \frac{1}{\Delta x} (\hat{\bm f}_{j+\frac{1}{2}} - \hat{\bm f}_{j-\frac{1}{2}}) = \bm 0. 
\end{equation}
Here, \(\bm u_j\) is the approximation to \(\bm u(t, x_j)\) and \(\hat{\bm f}_{j+\frac{1}{2}}\) is a numerical flux at the cell interface \(x_{j+\frac{1}{2}}\). 

Given a set of transform variables \(\bm v(\bm u)\), the procedure to compute \(\hat{\bm f}_{j+\frac{1}{2}}\)'s is as follows. 
\begin{enumerate}
    \item Obtain ghost point values using the boundary conditions. 
    \item Compute all the flux function values \(\bm f_{j} = \bm f(\bm u_{j})\) and all the transform variables \(\bm v_{j} = \bm v(\bm u_{j})\). 
    \item For each cell interface \(x_{j+\frac{1}{2}}\): \begin{enumerate}
        \item Compute the Roe average \(\bm u_{j+\frac{1}{2}}^{\mathrm{Roe}}\) using \(\bm u_{j}\) and \(\bm u_{j+1}\), then obtain the transformed eigenmatrices \(\tilde{\bm R}_{j+\frac{1}{2}} = \tilde{\bm R}(\bm u_{j+\frac{1}{2}}^{\mathrm{Roe}})\) and \(\tilde{\bm L}_{j+\frac{1}{2}} = \tilde{\bm L}(\bm u_{j+\frac{1}{2}}^{\mathrm{Roe}})\) at the cell interface. 
        \item Compute the local characteristic variables \(\bm w_{j-r+s} = \tilde{\bm L}_{j+\frac{1}{2}} \bm v_{j-r+s}\) in the whole stencil \(s = 1,2,\ldots,2r\). 
        \item Use WENO interpolation (see \Cref{ap:WENO interp} for details) to obtain approximations \(\bm w_{j+\frac{1}{2}}^{-} = \mathrm{WENO}(\bm w_{j-r+1}, \ldots, \bm w_{j+r-1})\) and \(\bm w_{j+\frac{1}{2}}^{+} = \mathrm{WENO}(\bm w_{j+r}, \ldots, \bm w_{j-r+2})\). 
        \item Transform back to \(\bm v_{j+\frac{1}{2}}^{\pm} = \tilde{\bm R}_{j+\frac{1}{2}} \bm w_{j+\frac{1}{2}}^{\pm}\) and then to the conserved variables \(\bm u_{j+\frac{1}{2}}^{\pm} = \bm u(\bm v_{j+\frac{1}{2}}^{\pm})\). 
        \item Compute the lowest-order flux \(\hat{\bm f}_{j+\frac{1}{2}}^{\mathrm{low}} = \hat{\bm f}(\bm u_{j+\frac{1}{2}}^{-}, \bm u_{j+\frac{1}{2}}^{+})\) using some approximate Riemann solver, e.g. the HLL-family flux. 
        \item Compute the high-order flux correction terms (see \Cref{ap:high-order corr} for details) \(\hat{\bm f}_{j+\frac{1}{2}}^{\mathrm{cor}}\) using flux function values \(\bm f_{j-r+s}\) for \(s = 1, 2, \ldots, 2r\). Obtain \(\hat{\bm f}_{j+\frac{1}{2}} = \hat{\bm f}_{j+\frac{1}{2}}^{\mathrm{low}} + \hat{\bm f}_{j+\frac{1}{2}}^{\mathrm{cor}}\). 
    \end{enumerate}
\end{enumerate}

Time is advanced using the strong stability-preserving Runge--Kutta (SSP-RK) methods \cite{10.1137/S003614450036757X}. The time step is determined by the CFL condition, 
\begin{equation}
    \Delta t = \mathrm{CFL} \min_{j} \frac{1}{\frac{\alpha_{j+\frac{1}{2}}}{\Delta x}},
\end{equation}
where \(\alpha_{j+\frac{1}{2}}\) is the estimate of the maximum wave speed at the cell interface \(x_{j+\frac{1}{2}}\) by the selected approximate Riemann solver. 

For the 2D case, the semi-discrete A-WENO scheme is performed in a dimension-by-dimension manner. The time step should be taken as 
\begin{equation}
    \Delta t = \mathrm{CFL} \min_{i,j} \frac{1}{\frac{\alpha_{i\pm\frac{1}{2},j}}{\Delta x} + \frac{\beta_{i,j\pm\frac{1}{2}}}{\Delta y}},
\end{equation} 
where \(\alpha_{i+\frac{1}{2},j}\) is the estimate of the maximum wave speed in the \(x\)-direction at the cell interface \((x_{i+\frac{1}{2}}, y_j)\) and \(\beta_{i,j+\frac{1}{2}}\) is the estimate in the \(y\)-direction at the cell interface \((x_i, y_{j+\frac{1}{2}})\).

\subsection{Reducing the cost of LCD}

Besides the WENO interpolation part, the most costly part that grows with the stencil width \(2r\) or the scheme order \(k = 2r-1\) should be the LCD procedure, which is documented to take a big portion of the total computational time \cite{10.1016/0168-92749290066-M}. For a system of conservation laws consisting of \(m\) components, the A-WENO LCD procedure at \emph{every} cell interface \(x_{j+\frac{1}{2}}\) has the following costs: (1) \(4r\) multiplications between the left eigenmatrix \(\tilde{\bm L}_{j+\frac{1}{2}}\) and vectors of transform variables \(\bm v_{j-r+s}\), (2) \(2\) multiplications between the right eigenmatrix \(\tilde{\bm R}_{j+\frac{1}{2}}\) and vectors of characteristic vectors \(\bm w_{j+\frac{1}{2}}^{\pm}\). A vanilla implementation of the matrix-vector multiplication will cost \(m^2\) floating point multiplications. With these in mind, one could easily see that the major part that grows with the stencil width is the left eigenmatrix-vector multiplication. Therefore, reducing this cost is the key to improve the efficiency of the A-WENO scheme at high-order. 

One attempt for some hyperbolic systems so far is to use the coordinate system of Riemann invariants \cite{10.1137/22M1536479} as the transform variables to completely avoid computing such multiplications, because in this case both the left and right eigenmatrices are diagonal and thus can be replaced by the identity matrix. Unfortunately, the compressible Euler equations \eqref{eq: Euler 2D} with the \(\gamma\)-law EOS \eqref{eq: EOS} is \emph{not} endowed with such a coordinate system. (See detailed proof in \Cref{prop:no-coordinate-Riemann-invariant}.) We have to seek other ways to reduce the cost, and our main freedom lies in choosing ``good'' transform variables \(\bm v(\bm u)\) and designing efficient left eigenmatrix-vector multiplication subroutines. 

If we stick with the trivial transform variables \(\bm v = \bm u\), there are still some possibilities to improve the efficiency. Two main ideas are to suitably normalize the eigenvectors to increase the number of unit elements in the eigenvectors so that to avoid multiplications, and to reuse intermediate results as possible. Take the 2D case as an example. The left eigenvectors are split as \(\bm l_1 = \bm \xi_1 - \bm \xi_2\), \(\bm l_2 = \bm \xi_2 - \bm \xi_3\), \(\bm l_3 = \bm \xi_4\) and \(\bm l_4 = \bm \xi_1 + \bm \xi_2\) as in \eqref{eq:Euler-left-eigen}. When computing \(\bm w = \bm L \bm u\), we first calculate the intermediate results \(\bm\xi_i \cdot \bm u\) and then sum them back to get the final results. Although there are still 4 vectors, cost can be saved because of their sparsity. The inner product between a general vector and the \(\bm \xi_i\)'s require only \(1\), \(4\), \(1\) and \(1\) ``true'' multiplications, so the total number of multiplications can go down from \(16\) to \(7\). 

Inspired by \Cref{cor:sparsity-of-eigenmatrix}, it is possible to improve the sparsity of the eigenmatrix so as to reduce the cost further. We use the following Riemann invariants as transform variables for the 2D compressible Euler equations in \(x\)-direction, \(\bm v(\bm u) = (u - \frac{2}{\gamma-1}c, S^{\frac{1}{2\gamma}}, v, u + \frac{2}{\gamma-1}c)\), whose related eigenmatrices are given in \eqref{eq:transform-right-eigen} and \eqref{eq:transform-left-eigen}. Each of the eigenmatrix has only two non-zero off-diagonals. Moreover, these two off-diagonals are opposite to each other and both lie in the second column. Let us denote the positive off-diagonal element by \(\mu\). We can implement the left eigenmatrix-vector multiplication \(\bm w = \tilde{\bm L} \bm v\) using only \(1\) multiplications: compute the intermediate results \(a = \mu v_2\), then the final result is just \(\bm w = (v_1 + a, v_2, v_3, v_1 - a)\). It's same for both 1D and 2D versions of the equations. 

We list the cost of the original and our new method in \Cref{tab:gemm-cost-comparison}. In terms of the whole scheme, the cost reduction will not be as significant as the reduction of solely the matrix-vector multiplications, because possible additional work is required to compute the forward and backward transforms and the coefficients in the new eigenmatrices. In this particular case, the cost of forward and backward transforms is only constant per-node, but such constant may become relatively large because of the exponential part in \(S^{\frac{1}{2\gamma}} = p^{\frac{1}{2\gamma}} \rho^{-\frac{1}{2}}\) which is expensive for CPUs. It might exceed the reduction of cost gained from our new LCD if not treated properly. Highly optimized SIMD instructions are needed to overcome this issue. 

\begin{table}[htb]
    \centering
    \caption{Number of floating point multiplications of one left eigenmatrix-vector multiplication for the compressible Euler equations of different dimensions.}
    \label{tab:gemm-cost-comparison}
    \begin{tabular}{ccc}
        \hline
        & original method & our new method \\
        \hline
        1D & 9 & 1 \\
        2D & 16 & 1 \\
        \hline
    \end{tabular}
\end{table}

\subsection{Properties of the new scheme}

We emphasize here that with any smooth transform variables, the A-WENO scheme is still conservative, and the convergence rates for smooth solutions will remain the same. Moreover, with our selected Riemann invariant-based transform variables, the scheme has the following property. 

In 1D simulations, the E-property \cite{10.1007/s10915-021-01659-w} refers to the ability to maintain constant velocity \(u\) and pressure \(p\) in the fully discrete solution if they are initially constant, regardless of the initial condition of \(\rho\). Schemes that satisfy the E-property, such as the primitive variable-based and the classical characteristic-wise A-WENO schemes, do not produce oscillations near contact discontinuities. In contrast, schemes that do not satisfy the E-property, such as the component-wise A-WENO scheme, may produce oscillations in these regions. By carefully analyzing the evolution of every component of the variable in a fully discretized scheme, we have the following proposition. 

\begin{proposition}[E-property]
    Assume the WENO procedure is affine-invariant, i.e. the WENO procedure is commutative with affine transformations. Then for the 1D case, both the original and the Riemann invariant-based A-WENO schemes satisfy the E-property. 
    \begin{proof}
        See \cite{10.1007/s10915-021-01659-w} for the proof for the original scheme. The proof for our Riemann invariant-based version is similar. The key point is to notice that when \(u_j\) and \(p_j\) are constants, the WENO interpolation for each component of \(\bm v_{j+\frac{1}{2}}^{-}\) actually gives the same result as that for \((\rho_{j+\frac{1}{2}}^{-})^{-\frac{1}{2}}\), up to affine transformations. 
    \end{proof}
\end{proposition}

\begin{remark}
    Actually, the WENO-JS procedure is not strictly affine-invariant in computation because of the artificial \(\epsilon\) to avoid division by zero. However, such impact is very small because \(\epsilon\) can be taken near machine zero. We refer the readers to the strictly affine-invariant WENO procedure \cite{10.1016/j.apnum.2022.07.007, 10.4208/cicp.OA-2023-0153} for more details. 
\end{remark}

\begin{remark}
    The exponential term \(S^{\frac{1}{2\gamma}}\) may seem arbitrary at first glance, because \Cref{cor:closedness-under-composition} implies that any smooth and invertible function of \(S\) can take the second position in the transform variables while preserving the sparsity of the eigenmatrix. However, it is actually the \emph{only} choice to satisfy the E-property, and it is less oscillatory near contact discontinuities than other choices. 
\end{remark}

\subsection{Positivity-preserving limiters}

During the computation, it is crucial to ensure that certain variables, which are both mathematically and physically required to be positive, do not take negative values so that the computation might break down. This can be achieved by implementing two positivity-preserving (PP) limiters in \cite[Appendix D]{10.1016/j.compfluid.2020.104519}. These limiters only interfere when certain negativity is detected, so in most of the time they contribute little to the computational cost. 

The first one is a PP limiter on the WENO interpolation results. Let us only consider the 1D case, as generalizations to 2D are straightforward. The admissible set of the transform variables \(\bm v = (u - \frac{2}{\gamma-1} c, S^{\frac{1}{2\gamma}}, u + \frac{2}{\gamma-1} c)\), \(\mathcal{G} = \B{\bm v \in \R^3: v_1 < v_3, v_2 > 0}\), is convex, so the convex combination-based limiter \cite{10.1016/j.jcp.2010.08.016} can be directly adopted to the interpolated transform variables \(\bm v_{j+\frac{1}{2}}^{\pm}\) in our framework. 

The second one is the so-called PP flux limiter, which requires a stricter CFL condition to maintain accuracy when it is triggered. The readers are referred to \cite{10.1016/j.jcp.2013.01.024} for details.

\section{Numerical results}\label{se:nr}

In this section, we perform numerical tests using \(5\)th-, \(7\)th- and \(9\)th-order A-WENO schemes. Time is discretized by the explicit third-order SSP-RK method. The classical conserved variable based characteristic decomposition is abbreviated as CH-CON, the component-wise version is abbreviated as CP-CON and the Riemann invariant-based one is abbreviated as CH-RI. We use the HLL numerical flux \cite{10.1137/1025002}. For the lowest-order flux, the positivity-preserving Einfeldt's wave speed estimate \cite{10.1137/0909030, 10.1016/0021-99919190211-3} is used. For the PP flux in the flux limiter, we adopt the more dissipative two-rarefaction wave speed estimate \cite[Section 10.5.2]{10.1007/b79761} that proves to bound the physical wave speeds for \(\gamma \in (1, \frac{5}{3}]\) \cite{10.1016/j.jcp.2016.05.054}. Unless otherwise specified, the CFL number is set to \(0.5\), and the specific heat capacity ratio \(\gamma\) is set to \(1.4\) for air. 

First of all we list the per-time step CPU time cost of the three versions of the A-WENO scheme. Following that, different test examples are presented to verify the accuracy and robustness of our CH-RI version in contrast to the classical CH-CON version. 


\subsection{Efficiency test}

We use \(N = 2000\) nodes for the 1D tests and \((Nx, Ny) = (400,400)\) nodes for the 2D tests. The tests are conducted on \texttt{Intel(R) Core(TM) i7-9750H CPU @ 2.60GHz} with Intel Turbo Boost disabled. The compiler is \texttt{Intel Fortran IFX 2025.0.0 20241008}, and the compiling options are \texttt{-O3 -xHost -r8 -qopenmp-stubs}. The Intel SVML library is triggered by OpenMP SIMD directives to accelerate floating point exponentials. The results are shown in \Cref{tab:CPU-time-cost}. The LCD time saving refers to the percentage of LCD time reduction of the CH-RI version compared to the CH-CON version. 

\begin{table}[htb]
    \centering
    \caption{Efficiency test. Per-time step CPU time cost (in seconds) of the A-WENO scheme.}
    \label{tab:CPU-time-cost}
    \begin{tabular}{cccccc}
        \hline
        & order & CH-RI & CH-CON & CP-CON & LCD time saving \\
        \hline
        \multirow{4}{*}{1D} & 3 & 1.58E-03 & 1.98E-03 & 1.33E-03 & 62\% \\
        & 5 & 1.83E-03 & 2.30E-03 & 1.61E-03 & 68\% \\
        & 7 & 2.23E-03 & 2.72E-03 & 2.01E-03 & 69\% \\
        & 9 & 2.26E-03 & 2.82E-03 & 2.10E-03 & 78\% \\
        \hline
        \multirow{4}{*}{2D} & 3 & 4.98E-01 & 6.36E-01 & 4.69E-01 & 82\% \\
        & 5 & 5.20E-01 & 6.66E-01 & 5.07E-01 & 92\% \\
        & 7 & 5.39E-01 & 6.84E-01 & 5.20E-01 & 88\% \\
        & 9 & 6.21E-01 & 7.67E-01 & 5.99E-01 & 86\% \\
        \hline
    \end{tabular}
\end{table}

\subsection{1D tests}

\subsubsection{Accuracy test 1}

We consider the simple periodic density-transport problem whose exact solution is \((\rho, u, p) = (1+0.2\sin(\pi (x-t)), 1, 1)\). The computation domain is \([0,2]\) with periodic boundaries. Solve up to \(T = 2\). We use \(N_0 = 20\) as the coarsest mesh and use \(h_{0} = \frac{2}{N_0}\) to denote the coarsest mesh size. Using SSP-RK(3,3), the time steps in the following tests are taken to be \(\Delta t = \mathrm{CFL} \frac{\Delta x}{\alpha} \left(\frac{\Delta x}{h_0}\right)^{\frac{k}{3} - 1}\), which makes the time discretization error \(k\)th-order while maintaining stability. Here, \(\alpha\) is the maximum wave speed. The errors are shown by \Cref{tab:1D-accuracy-test-1-k=5,tab:1D-accuracy-test-1-k=7,tab:1D-accuracy-test-1-k=9}. 

\begin{table}[htbp]
    \centering
    \caption{Accuracy test 1. Error table of the density. \(k = 5\).}
    \label{tab:1D-accuracy-test-1-k=5}
    \begin{tabular}{ccccccccc}
        \hline
        & \multicolumn{4}{c}{CH-RI} & \multicolumn{4}{c}{CH-CON} \\
        \(N\) & \(\ell_h^2\) error & order & \(\ell_h^\infty\) error & order & \(\ell_h^2\) error & order & \(\ell_h^\infty\) error & order \\ \hline
        20  & 8.02E-04 & - & 1.10E-03 & - & 4.66E-04 & - & 5.06E-04 & - \\ 
        40  & 3.53E-05 & 4.507  & 5.80E-05 & 4.251  & 1.44E-05 & 5.015  & 1.78E-05 & 4.830  \\ 
        60  & 4.95E-06 & 4.842  & 9.03E-06 & 4.587  & 1.85E-06 & 5.056  & 2.44E-06 & 4.900  \\ 
        80  & 1.19E-06 & 4.961  & 2.29E-06 & 4.766  & 4.34E-07 & 5.050  & 5.78E-07 & 5.000  \\ 
        100  & 3.88E-07 & 5.018  & 7.75E-07 & 4.858  & 1.41E-07 & 5.043  & 1.86E-07 & 5.088  \\ 
        120  & 1.54E-07 & 5.059  & 3.04E-07 & 5.131  & 5.62E-08 & 5.037  & 7.32E-08 & 5.113  \\ 
        140  & 7.03E-08 & 5.098  & 1.39E-07 & 5.094  & 2.59E-08 & 5.034  & 3.36E-08 & 5.045  \\ 
        160  & 3.54E-08 & 5.137  & 6.83E-08 & 5.310  & 1.32E-08 & 5.035  & 1.68E-08 & 5.180 \\ 
        \hline
    \end{tabular}
\end{table}

\begin{table}[htbp]
    \centering
    \caption{Accuracy test 1. Error table of the density. \(k = 7\).}
    \label{tab:1D-accuracy-test-1-k=7}
    \begin{tabular}{ccccccccc}
        \hline
        & \multicolumn{4}{c}{CH-RI} & \multicolumn{4}{c}{CH-CON} \\
        \(N\) & \(\ell_h^2\) error & order & \(\ell_h^\infty\) error & order & \(\ell_h^2\) error & order & \(\ell_h^\infty\) error & order \\ \hline
        20  & 1.28E-04 & - & 2.59E-04 & - & 3.52E-05 & - & 5.74E-05 & - \\
        40  & 2.07E-06 & 5.945  & 5.88E-06 & 5.460  & 4.99E-07 & 6.142  & 1.16E-06 & 5.636  \\ 
        60  & 1.68E-07 & 6.202  & 5.79E-07 & 5.714  & 4.56E-08 & 5.902  & 1.28E-07 & 5.417  \\ 
        80  & 2.54E-08 & 6.562  & 9.97E-08 & 6.118  & 7.91E-09 & 6.088  & 2.60E-08 & 5.548  \\ 
        100  & 5.20E-09 & 7.107  & 2.21E-08 & 6.743  & 1.86E-09 & 6.481  & 6.98E-09 & 5.902  \\ 
        120  & 1.30E-09 & 7.597  & 5.78E-09 & 7.369  & 5.26E-10 & 6.937  & 2.18E-09 & 6.389  \\ 
        140  & 3.84E-10 & 7.923  & 1.72E-09 & 7.876  & 1.69E-10 & 7.361  & 7.50E-10 & 6.911  \\ 
        160  & 1.30E-10 & 8.093  & 5.72E-10 & 8.233  & 6.04E-11 & 7.704  & 2.79E-10 & 7.406 \\ 
        \hline
    \end{tabular}
\end{table}

\begin{table}[htbp]
    \centering
    \caption{Accuracy test 1. Error table of the density. \(k = 9\).}
    \label{tab:1D-accuracy-test-1-k=9}
    \begin{tabular}{ccccccccc}
        \hline
        & \multicolumn{4}{c}{CH-RI} & \multicolumn{4}{c}{CH-CON} \\
        \(N\) & \(\ell_h^2\) error & order & \(\ell_h^\infty\) error & order & \(\ell_h^2\) error & order & \(\ell_h^\infty\) error & order \\ \hline
        20  & 1.95E-05 & - & 3.23E-05 & - & 1.63E-05 & - & 1.60E-05 & - \\
        40  & 4.70E-08 & 8.695  & 8.18E-08 & 8.624  & 3.20E-08 & 8.992  & 3.13E-08 & 8.994  \\ 
        60  & 1.38E-09 & 8.702  & 2.79E-09 & 8.332  & 8.32E-10 & 8.998  & 8.23E-10 & 8.973  \\ 
        80  & 1.09E-10 & 8.839  & 2.32E-10 & 8.651  & 6.24E-11 & 9.002  & 6.22E-11 & 8.976  \\ 
        100  & 1.48E-11 & 8.930  & 3.17E-11 & 8.916  & 8.31E-12 & 9.035  & 8.87E-12 & 8.729  \\ 
        120  & 3.25E-12 & 8.309  & 4.98E-12 & 10.143  & 2.18E-12 & 7.337  & 2.60E-12 & 6.731  \\ 
        \hline
    \end{tabular}
\end{table}

\subsubsection{Accuracy test 2}

Notice that any isentropic solution are completely determined by the initial entropy value \(S_0 = p \rho^{-\gamma}\) and the left- and right- Riemann invariants, \(J^{\pm} = u \pm \frac{2 c}{\gamma - 1}\). When the specific heat ratio is \(\gamma = 3\) and all functions are smooth, \(J^{\pm}\) satisfy two disjoint Burgers' equations \cite{10.1016/j.compfluid.2010.07.018}, so now the exact solution can be constructed from the exact solution of Burgers' equation. Let \(\Omega = [-1, 1]\) with periodic boundaries. Choose \(S_0 = 1\) and \(J^{+} \equiv 2\), and let \(J^{-}\) be the solution with the exact solution of the Burgers' equation with initial condition \(J^{-}|_{t=0} = \sin(\pi x)\). The initial condition in primitive variables is then \(\rho = \frac{1 - \frac{1}{2} \sin(\pi x)}{\sqrt{3}}\), \(u = 1 + \frac{1}{2} \sin(\pi x)\), \(p = \rho^3\). We solve the problem to \(T = 0.2\). The time step is taken to be the same as the previous example. 

The errors are shown in \Cref{tab:1D-accuracy-test-2-k=5,tab:1D-accuracy-test-2-k=7,tab:1D-accuracy-test-2-k=9}. This is a genuinely nontrivial test since all variables have nontrivial variations. We can see from the error tables that the CH-RI version is as accurate as the CH-CON version in terms of the convergence rate. 

\begin{table}[htbp]
    \centering
    \caption{Accuracy test 2. Error table of the density. \(k = 5\).}
    \label{tab:1D-accuracy-test-2-k=5}
    \begin{tabular}{ccccccccc}
        \hline
        & \multicolumn{4}{c}{CH-RI} & \multicolumn{4}{c}{CH-CON} \\
        \(N\) & \(\ell_h^2\) error & order & \(\ell_h^\infty\) error & order & \(\ell_h^2\) error & order & \(\ell_h^\infty\) error & order \\ \hline
        20  & 2.99E-03 & - & 6.84E-03 & - & 9.13E-03 & - & 2.21E-02 & - \\
        40  & 2.98E-04 & 3.324  & 1.05E-03 & 2.701  & 8.91E-04 & 3.358  & 2.45E-03 & 3.172  \\ 
        60  & 5.18E-05 & 4.319  & 2.07E-04 & 4.006  & 2.87E-04 & 2.791  & 6.53E-04 & 3.265  \\ 
        80  & 1.27E-05 & 4.875  & 5.26E-05 & 4.768  & 1.00E-04 & 3.670  & 2.60E-04 & 3.201  \\ 
        100  & 4.00E-06 & 5.192  & 1.66E-05 & 5.168  & 2.90E-05 & 5.547  & 1.03E-04 & 4.158  \\ 
        120  & 1.50E-06 & 5.367  & 6.20E-06 & 5.395  & 1.28E-05 & 4.493  & 4.43E-05 & 4.616  \\ 
        140  & 6.49E-07 & 5.451  & 2.65E-06 & 5.516  & 6.33E-06 & 4.558  & 2.06E-05 & 4.971  \\ 
        160  & 3.12E-07 & 5.481  & 1.32E-06 & 5.212  & 3.38E-06 & 4.707  & 1.07E-05 & 4.871 \\
        \hline
    \end{tabular}
\end{table}

\begin{table}[htbp]
    \centering
    \caption{Accuracy test 2. Error table of the density. \(k = 7\).}
    \label{tab:1D-accuracy-test-2-k=7}
    \begin{tabular}{ccccccccc}
        \hline
        & \multicolumn{4}{c}{CH-RI} & \multicolumn{4}{c}{CH-CON} \\
        \(N\) & \(\ell_h^2\) error & order & \(\ell_h^\infty\) error & order & \(\ell_h^2\) error & order & \(\ell_h^\infty\) error & order \\ \hline
        20  & 2.50E-03 & - & 5.60E-03 & - & 5.24E-03 & - & 1.29E-02 & - \\
        40  & 1.53E-04 & 4.032  & 5.13E-04 & 3.449  & 6.72E-04 & 2.963  & 1.90E-03 & 2.768  \\ 
        60  & 1.82E-05 & 5.242  & 7.43E-05 & 4.766  & 1.14E-04 & 4.376  & 3.94E-04 & 3.878  \\ 
        80  & 3.24E-06 & 6.003  & 1.42E-05 & 5.764  & 1.88E-05 & 6.265  & 8.01E-05 & 5.541  \\ 
        100  & 7.56E-07 & 6.523  & 3.37E-06 & 6.425  & 5.27E-06 & 5.698  & 2.22E-05 & 5.746  \\ 
        120  & 2.15E-07 & 6.890  & 9.79E-07 & 6.789  & 1.79E-06 & 5.930  & 7.87E-06 & 5.692  \\ 
        140  & 7.17E-08 & 7.137  & 3.42E-07 & 6.823  & 6.96E-07 & 6.118  & 3.03E-06 & 6.195  \\ 
        160  & 2.70E-08 & 7.303  & 1.28E-07 & 7.351  & 3.02E-07 & 6.257  & 1.26E-06 & 6.544 \\
        \hline
    \end{tabular}
\end{table}

\begin{table}[htbp]
    \centering
    \caption{Accuracy test 2. Error table of the density. \(k = 9\).}
    \label{tab:1D-accuracy-test-2-k=9}
    \begin{tabular}{ccccccccc}
        \hline
        & \multicolumn{4}{c}{CH-RI} & \multicolumn{4}{c}{CH-CON} \\
        \(N\) & \(\ell_h^2\) error & order & \(\ell_h^\infty\) error & order & \(\ell_h^2\) error & order & \(\ell_h^\infty\) error & order \\ \hline
        20  & 2.00E-03 & - & 4.53E-03 & - & 2.67E-03 & - & 7.13E-03 & - \\
        40  & 1.11E-04 & 4.167  & 3.65E-04 & 3.636  & 3.12E-04 & 3.093  & 9.76E-04 & 2.869  \\ 
        60  & 1.11E-05 & 5.685  & 4.49E-05 & 5.164  & 2.99E-05 & 5.789  & 1.06E-04 & 5.484  \\ 
        80  & 1.60E-06 & 6.721  & 7.21E-06 & 6.358  & 4.28E-06 & 6.754  & 1.63E-05 & 6.506  \\ 
        100  & 3.07E-07 & 7.411  & 1.43E-06 & 7.244  & 7.17E-07 & 8.004  & 3.01E-06 & 7.552  \\ 
        120  & 7.18E-08 & 7.969  & 3.38E-07 & 7.916  & 1.38E-07 & 9.036  & 6.33E-07 & 8.563  \\ 
        140  & 1.97E-08 & 8.390  & 9.35E-08 & 8.343  & 3.31E-08 & 9.271  & 1.61E-07 & 8.873  \\ 
        160  & 6.17E-09 & 8.697  & 3.17E-08 & 8.089  & 1.03E-08 & 8.724  & 4.93E-08 & 8.875 \\
        \hline
    \end{tabular}
\end{table}

\subsubsection{The Sod shock tube problem}

\begin{figure}[htbp]
    \centering

    \begin{subfigure}{0.32\linewidth}
        \centering
        \includegraphics[width=\linewidth]{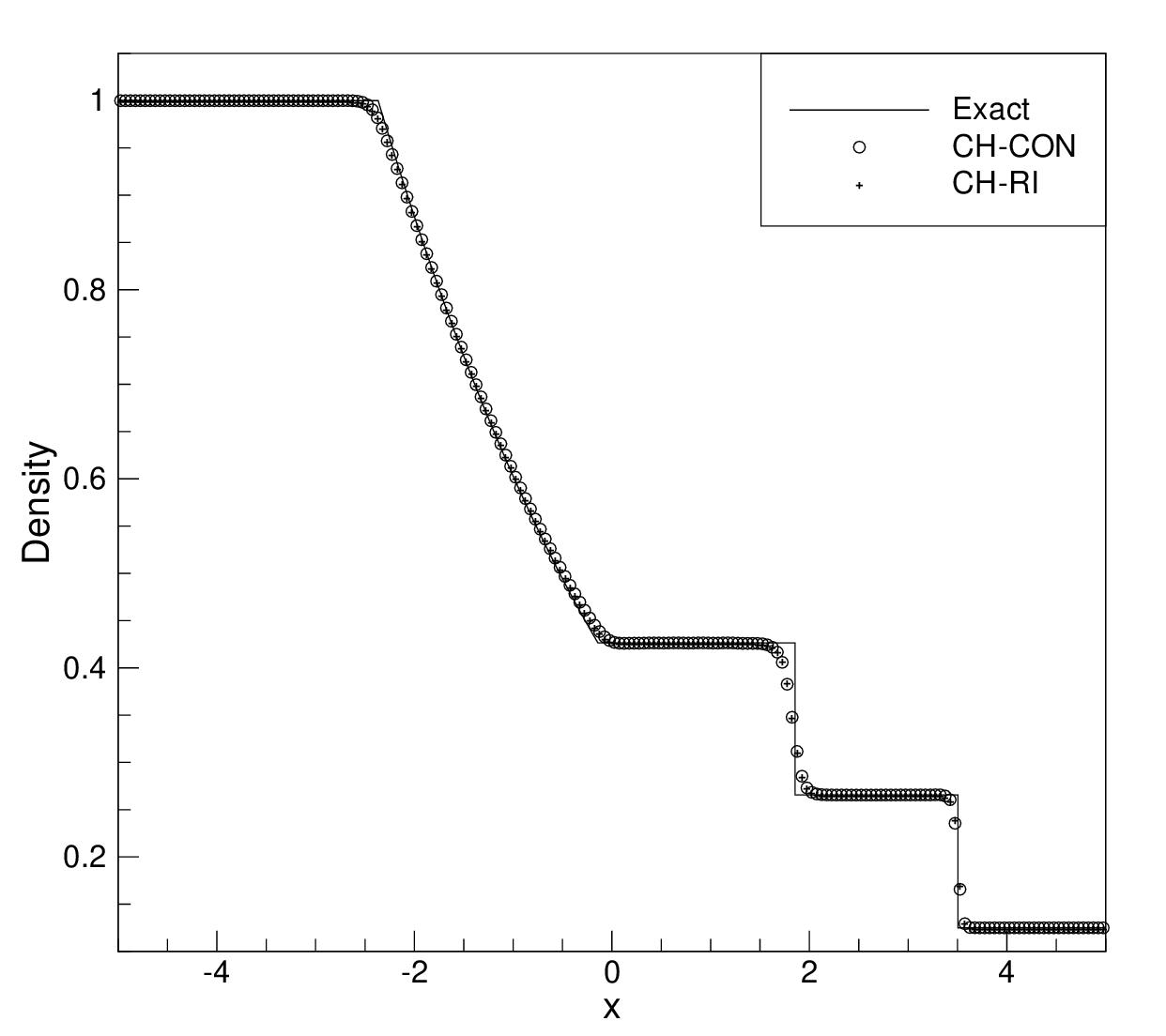}
        \caption{AWENO5}
    \end{subfigure}
    \begin{subfigure}{0.32\linewidth}
        \centering
        \includegraphics[width=\linewidth]{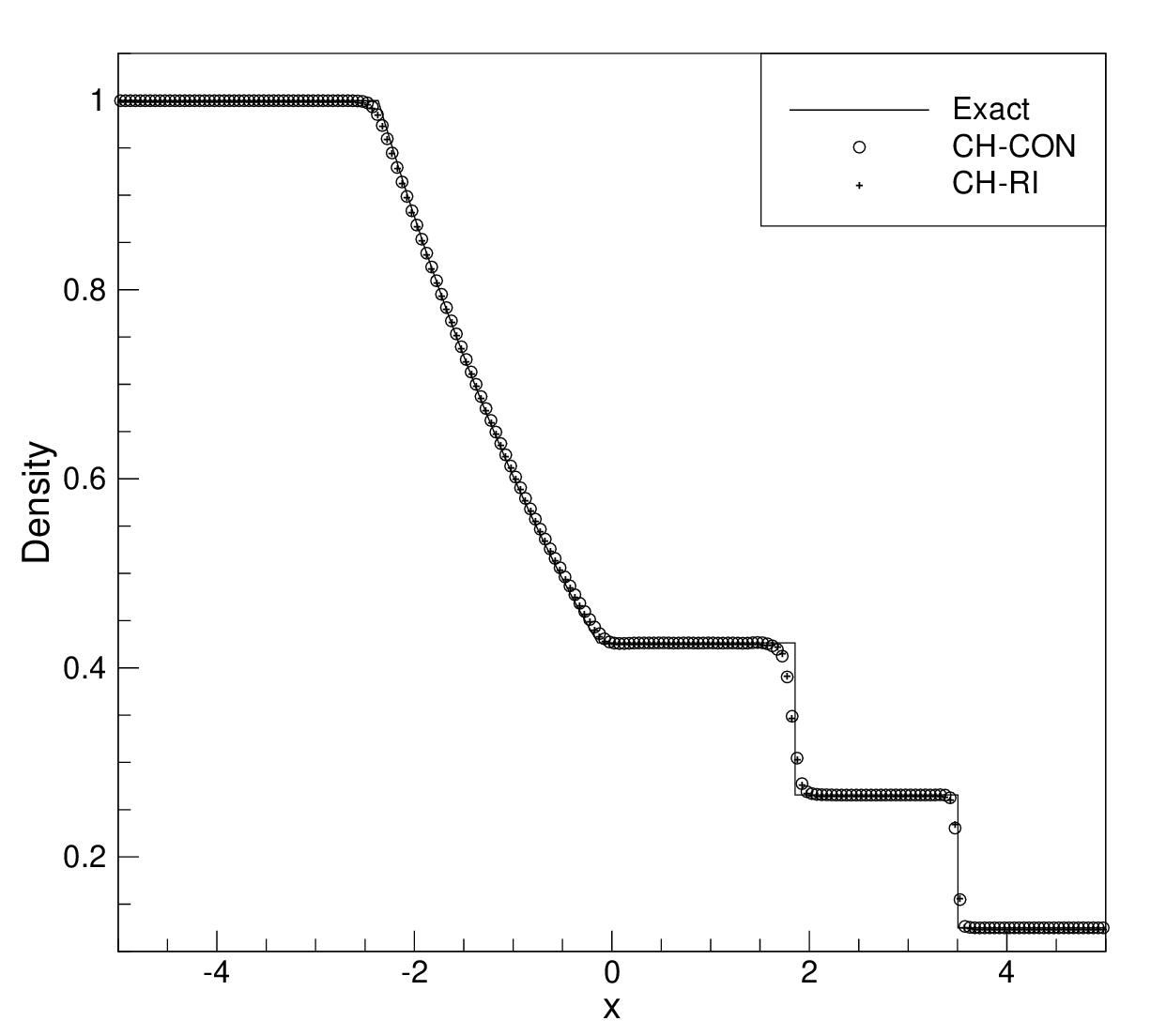}
        \caption{AWENO7}
    \end{subfigure}
    \begin{subfigure}{0.32\linewidth}
        \centering
        \includegraphics[width=\linewidth]{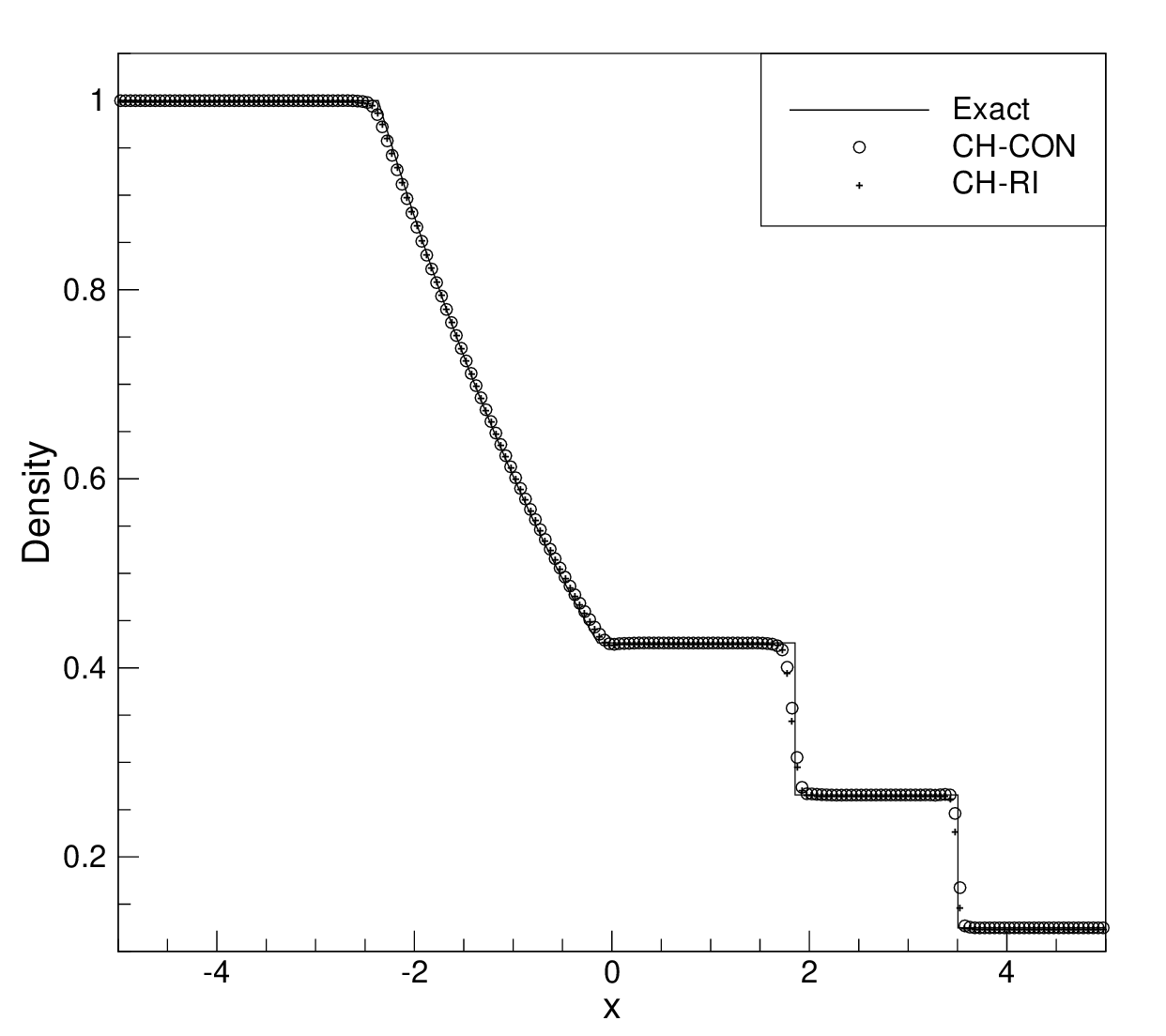}
        \caption{AWENO9}
    \end{subfigure}

    \begin{subfigure}{0.32\linewidth}
        \centering
        \includegraphics[width=\linewidth]{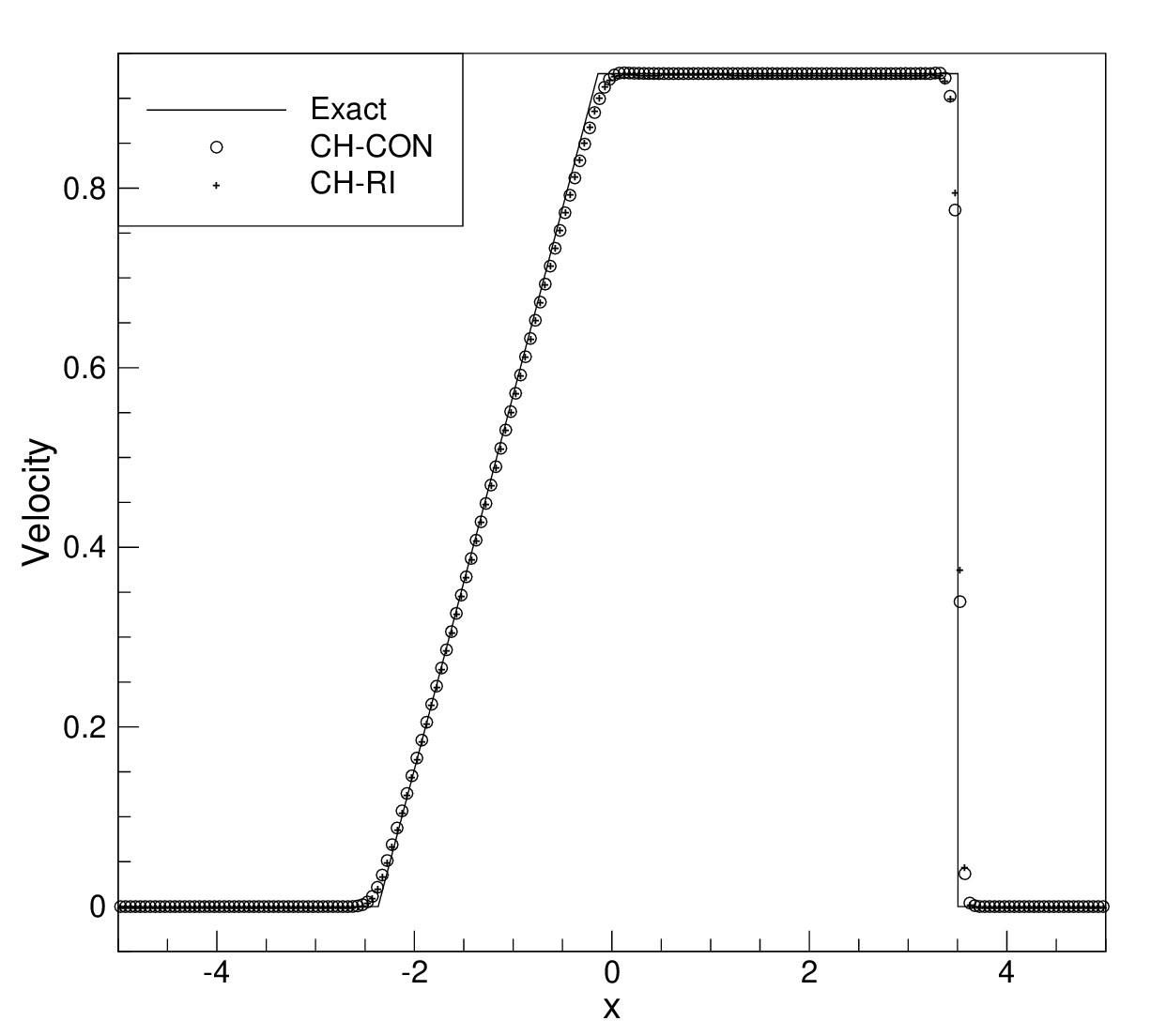}
        \caption{AWENO5}
    \end{subfigure}
    \begin{subfigure}{0.32\linewidth}
        \centering
        \includegraphics[width=\linewidth]{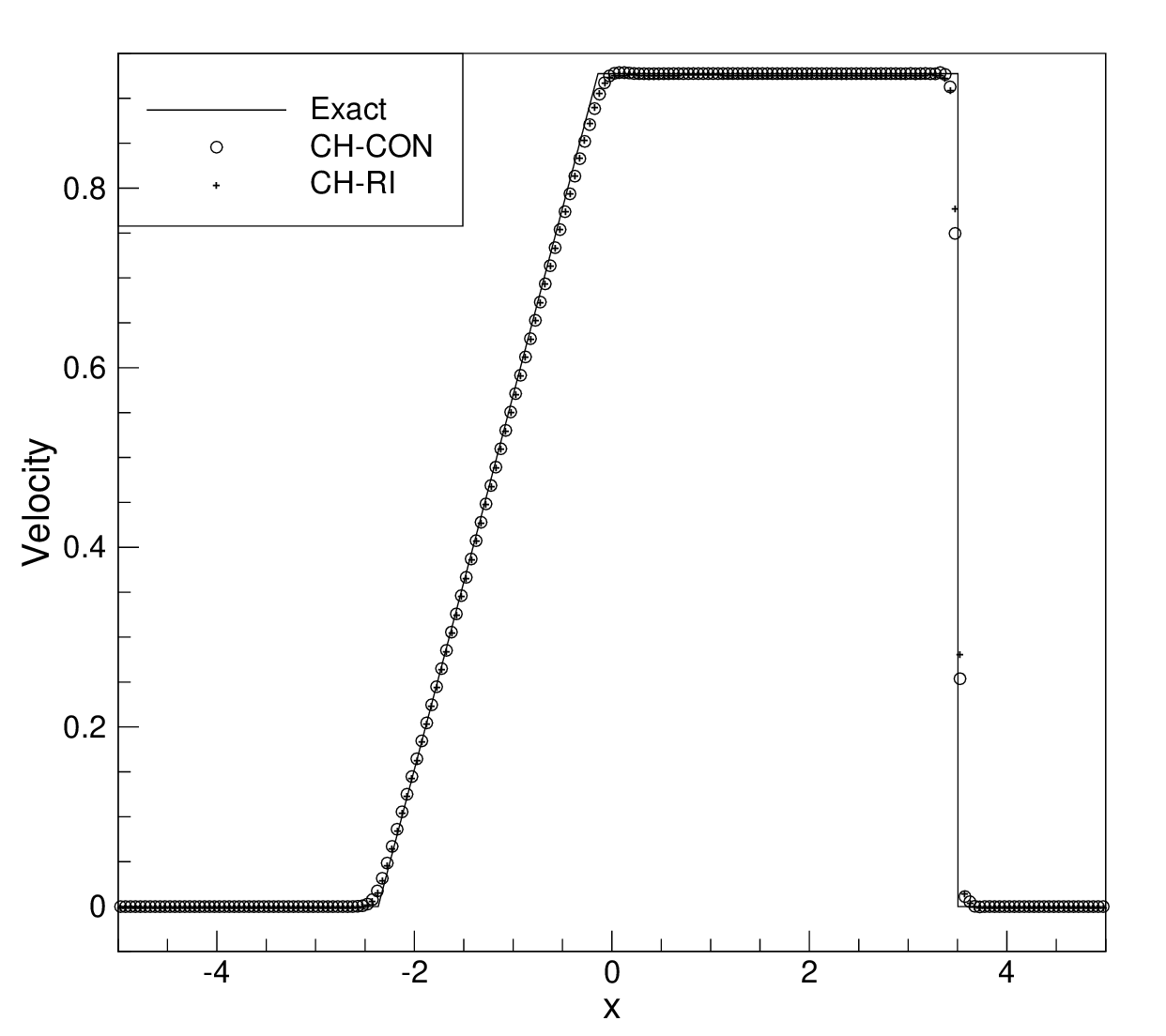}
        \caption{AWENO7}
    \end{subfigure}
    \begin{subfigure}{0.32\linewidth}
        \centering
        \includegraphics[width=\linewidth]{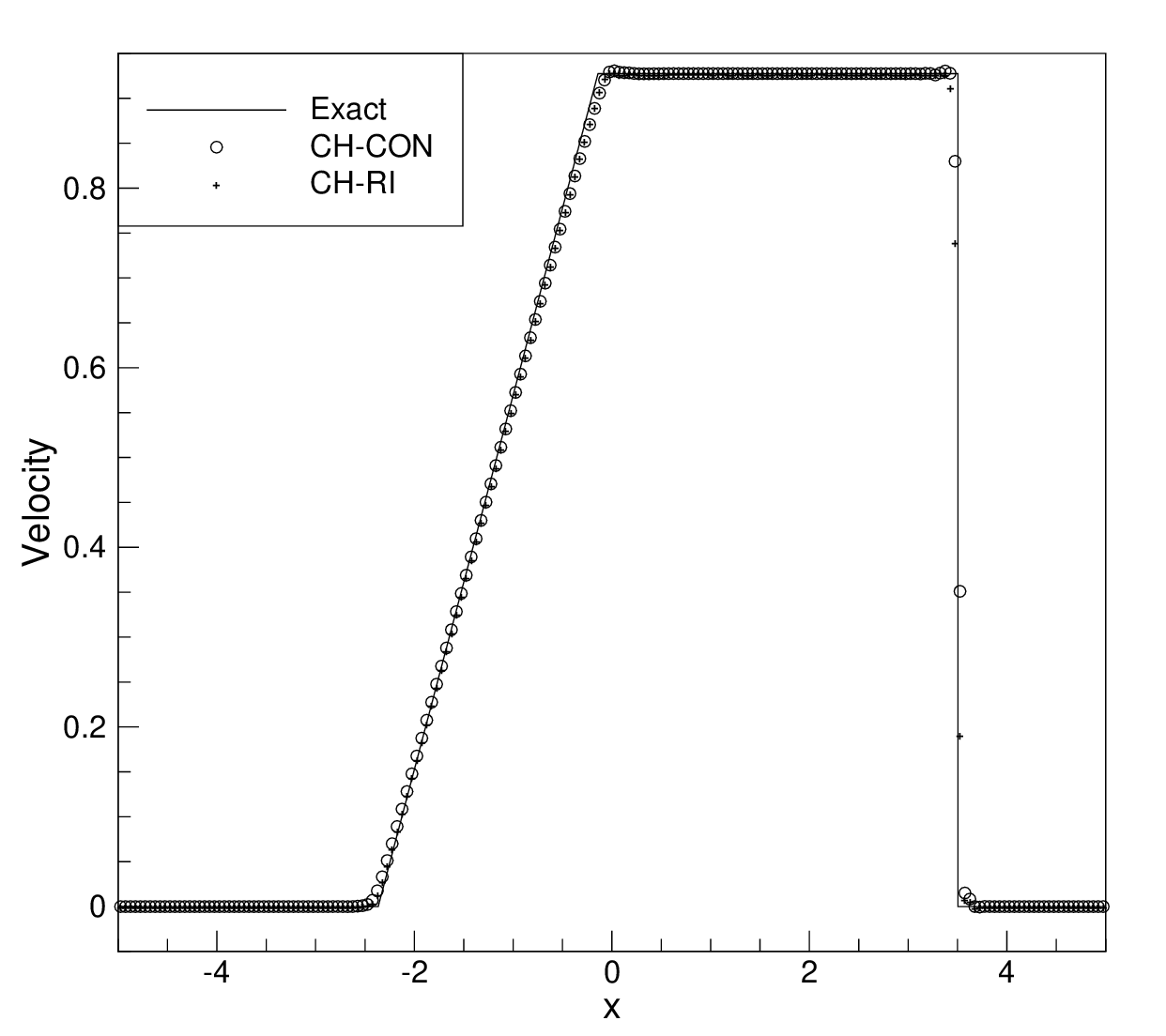}
        \caption{AWENO9}
    \end{subfigure}

    \begin{subfigure}{0.32\linewidth}
        \centering
        \includegraphics[width=\linewidth]{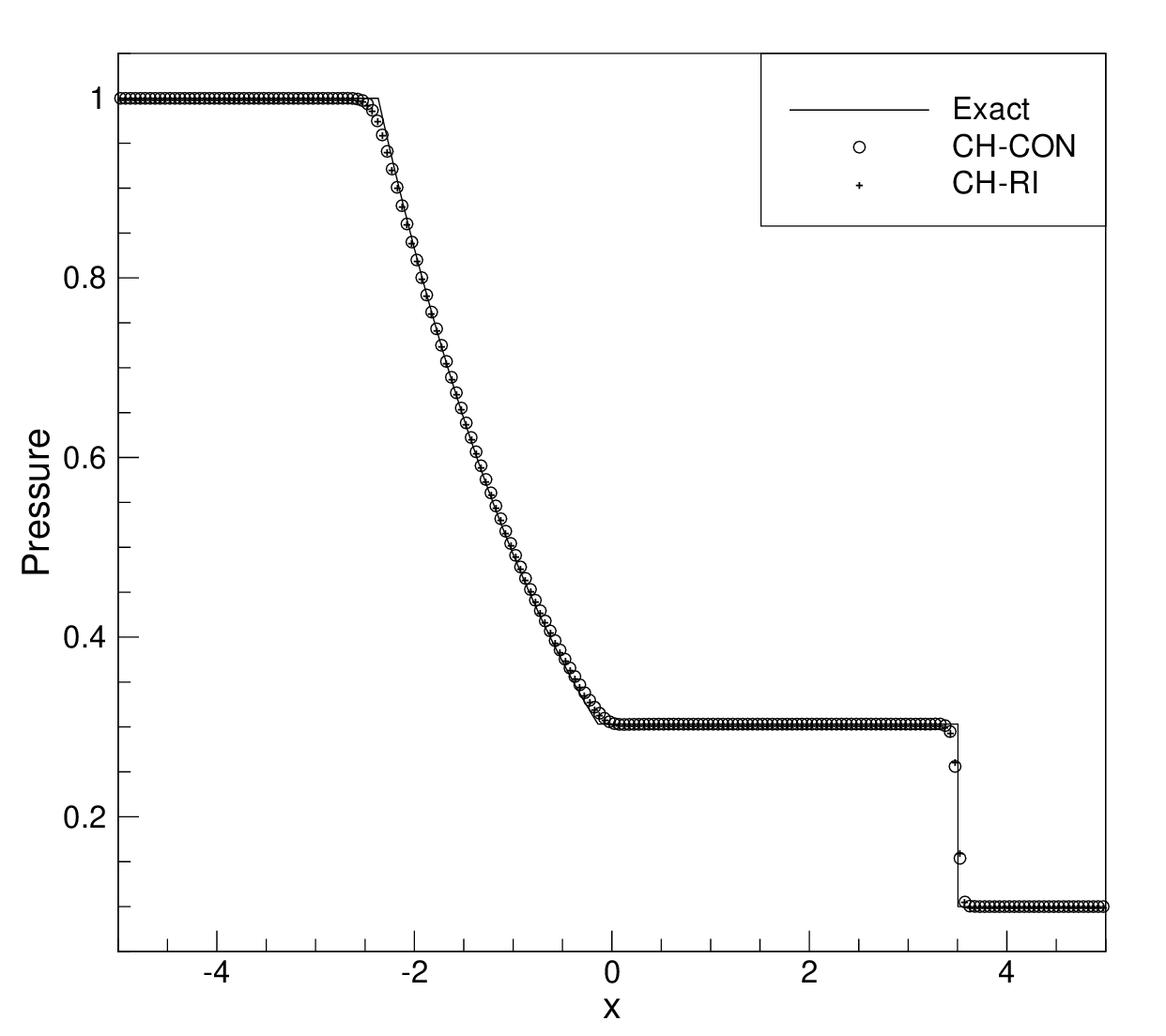}
        \caption{AWENO5}
    \end{subfigure}
    \begin{subfigure}{0.32\linewidth}
        \centering
        \includegraphics[width=\linewidth]{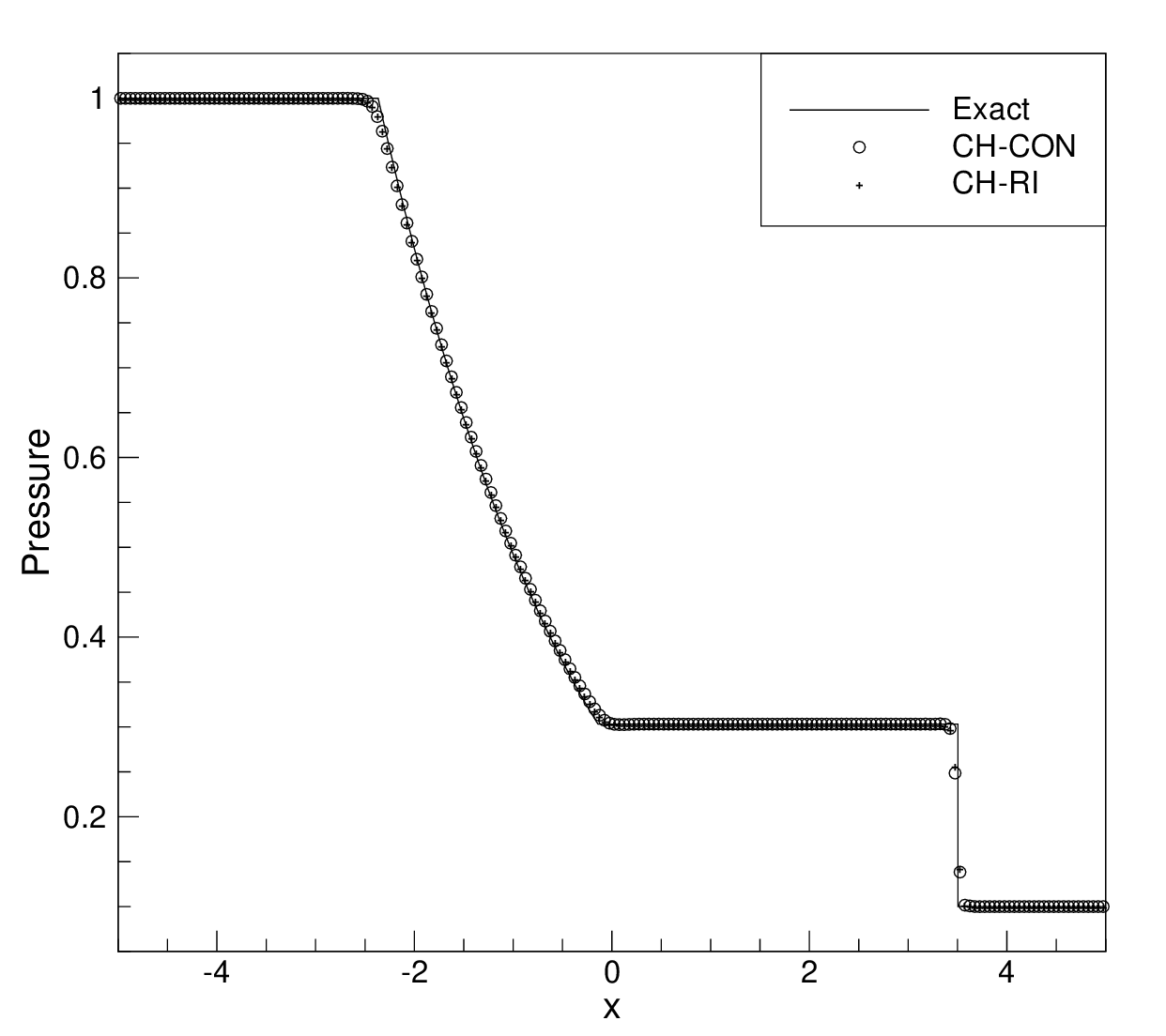}
        \caption{AWENO7}
    \end{subfigure}
    \begin{subfigure}{0.32\linewidth}
        \centering
        \includegraphics[width=\linewidth]{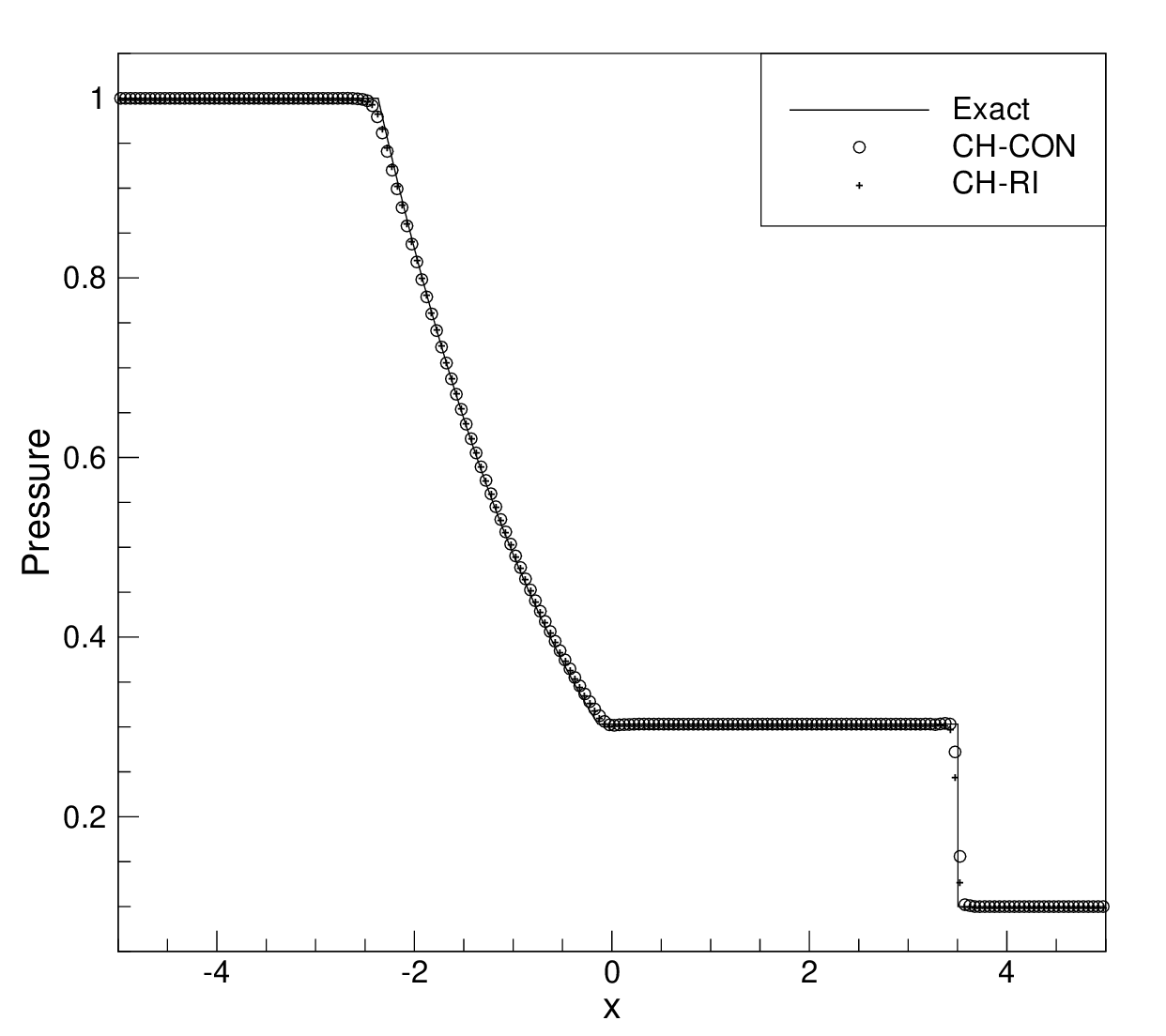}
        \caption{AWENO9}
    \end{subfigure}

    \caption{The Sod shock tube problem at \(T = 2\) with \(N=200\) grids.}
    \label{fig:Sod}
\end{figure}

The initial condition is \((\rho, u, p) = (1, 0, 1)\) for \(x < 0\) and \((\rho, u, p) = (0.125, 0, 0.1)\) for \(x > 0\). The computation domain is \([-5, 5]\). We solve the problem to \(T = 2\). The results are shown in \Cref{fig:Sod}. For this example, both CH-RI and CH-CON have similarly good performance. 

\subsubsection{The Lax shock tube problem}

\begin{figure}[htbp]
    \centering

    \begin{subfigure}{0.32\linewidth}
        \centering
        \includegraphics[width=\linewidth]{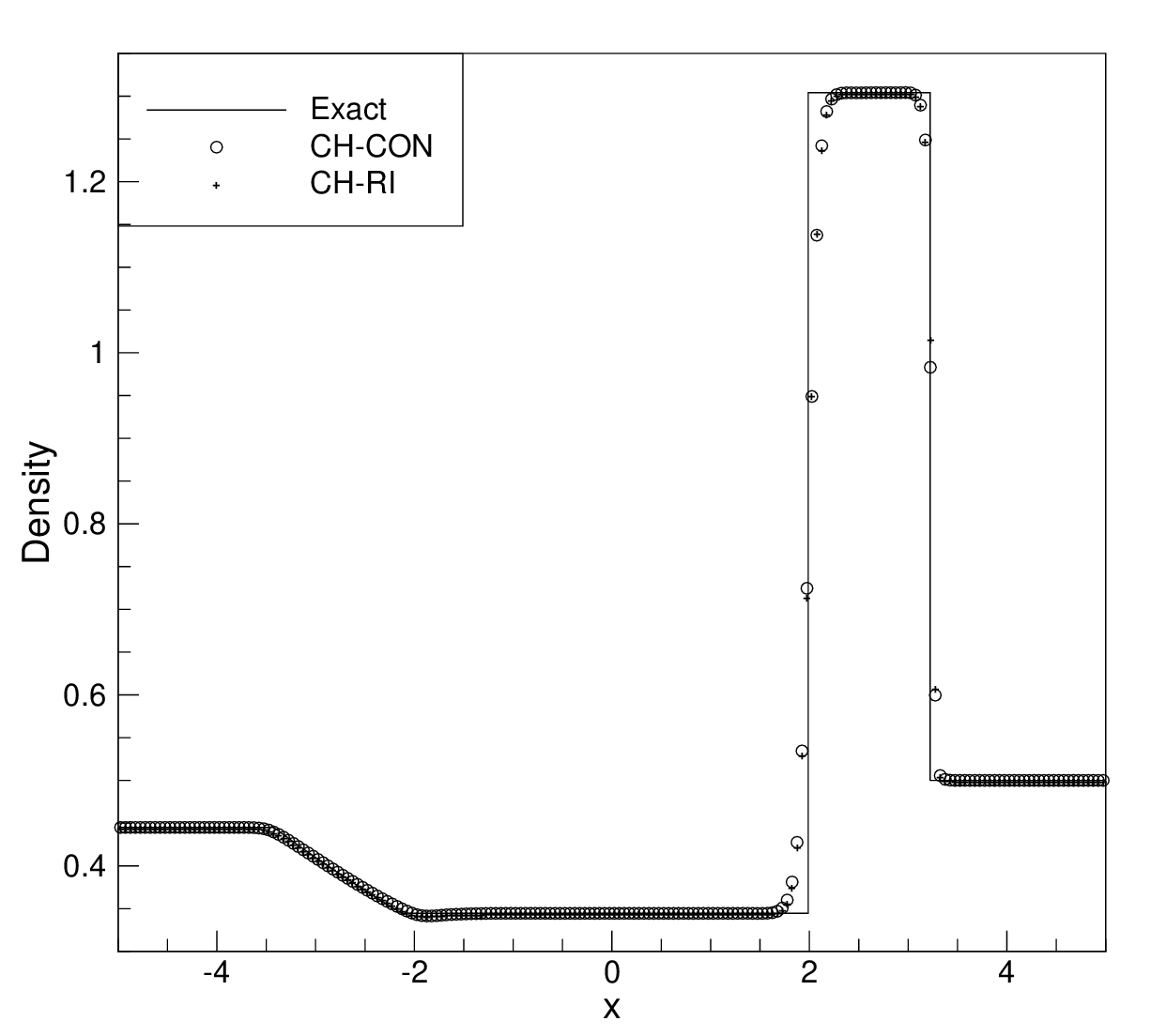}
        \caption{AWENO5}
    \end{subfigure}
    \begin{subfigure}{0.32\linewidth}
        \centering
        \includegraphics[width=\linewidth]{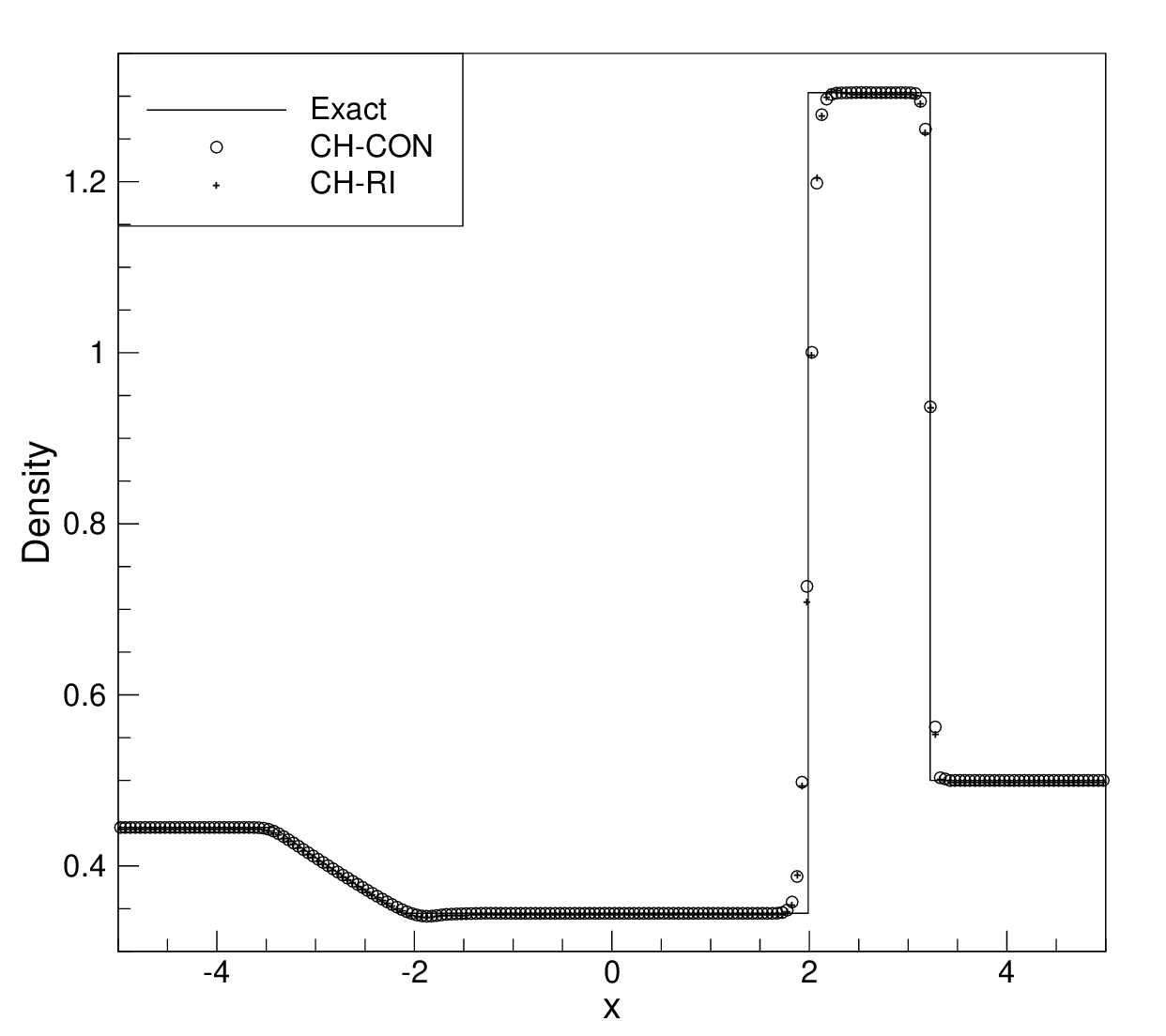}
        \caption{AWENO7}
    \end{subfigure}
    \begin{subfigure}{0.32\linewidth}
        \centering
        \includegraphics[width=\linewidth]{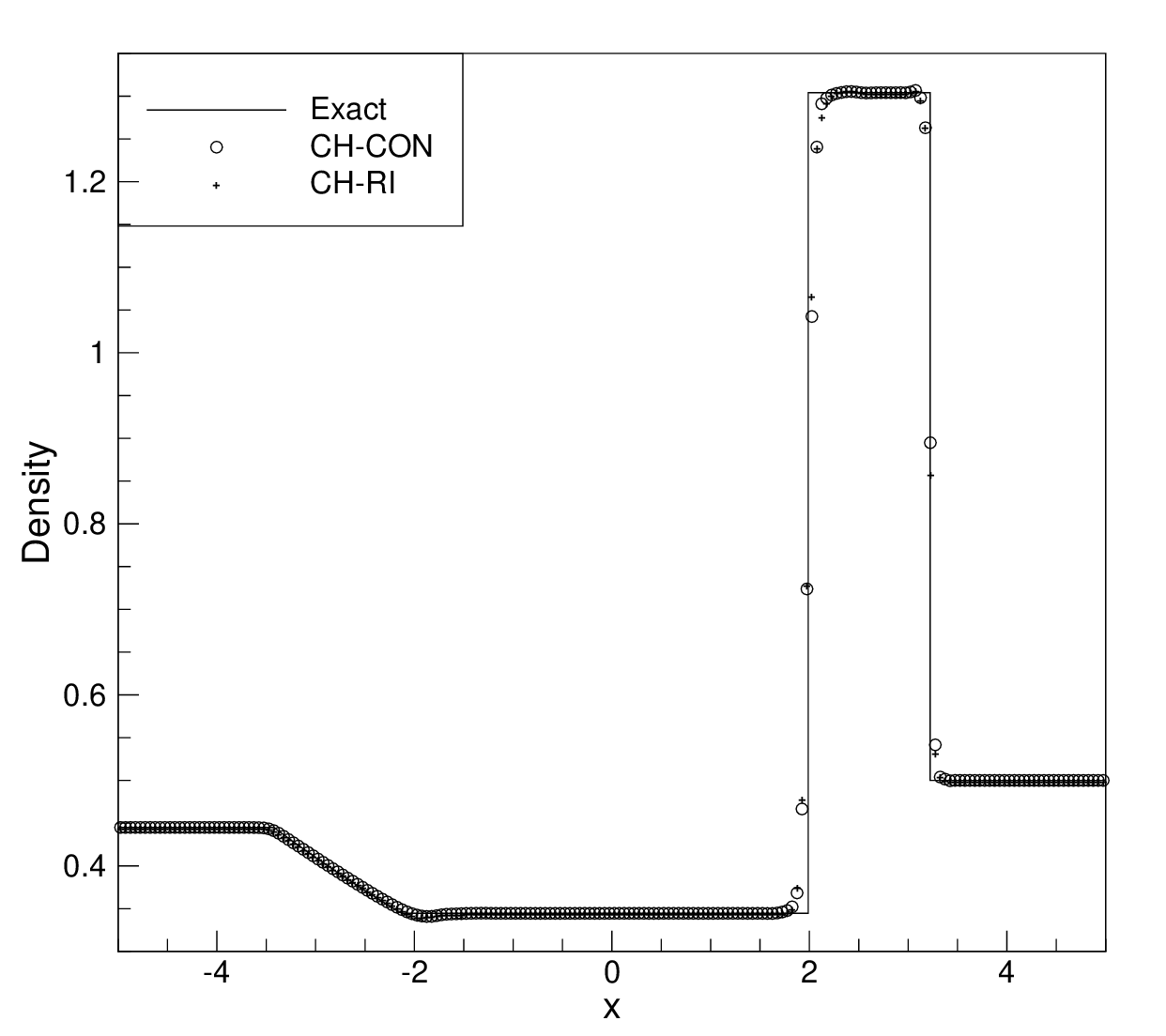}
        \caption{AWENO9}
    \end{subfigure}

    \begin{subfigure}{0.32\linewidth}
        \centering
        \includegraphics[width=\linewidth]{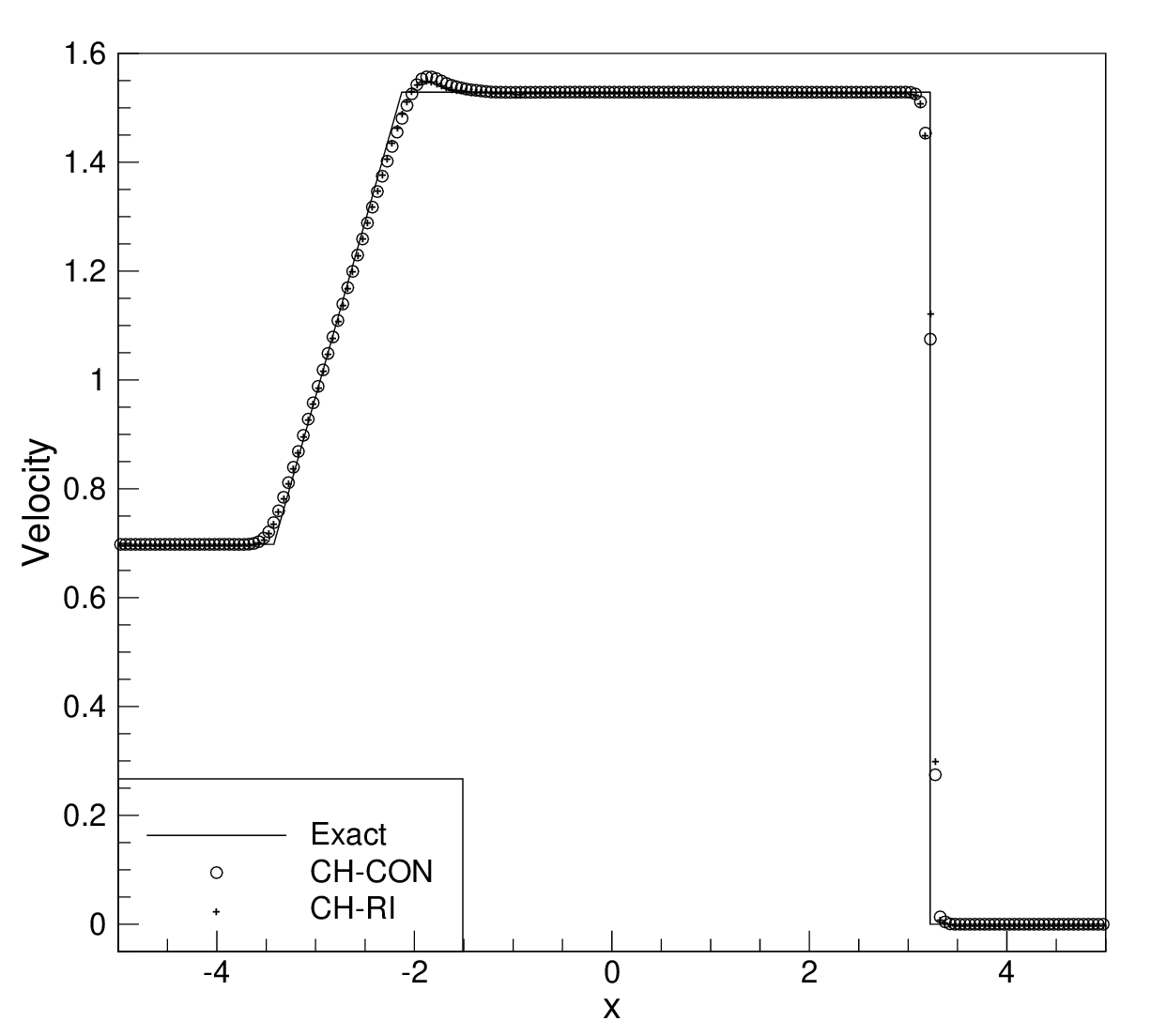}
        \caption{AWENO5}
    \end{subfigure}
    \begin{subfigure}{0.32\linewidth}
        \centering
        \includegraphics[width=\linewidth]{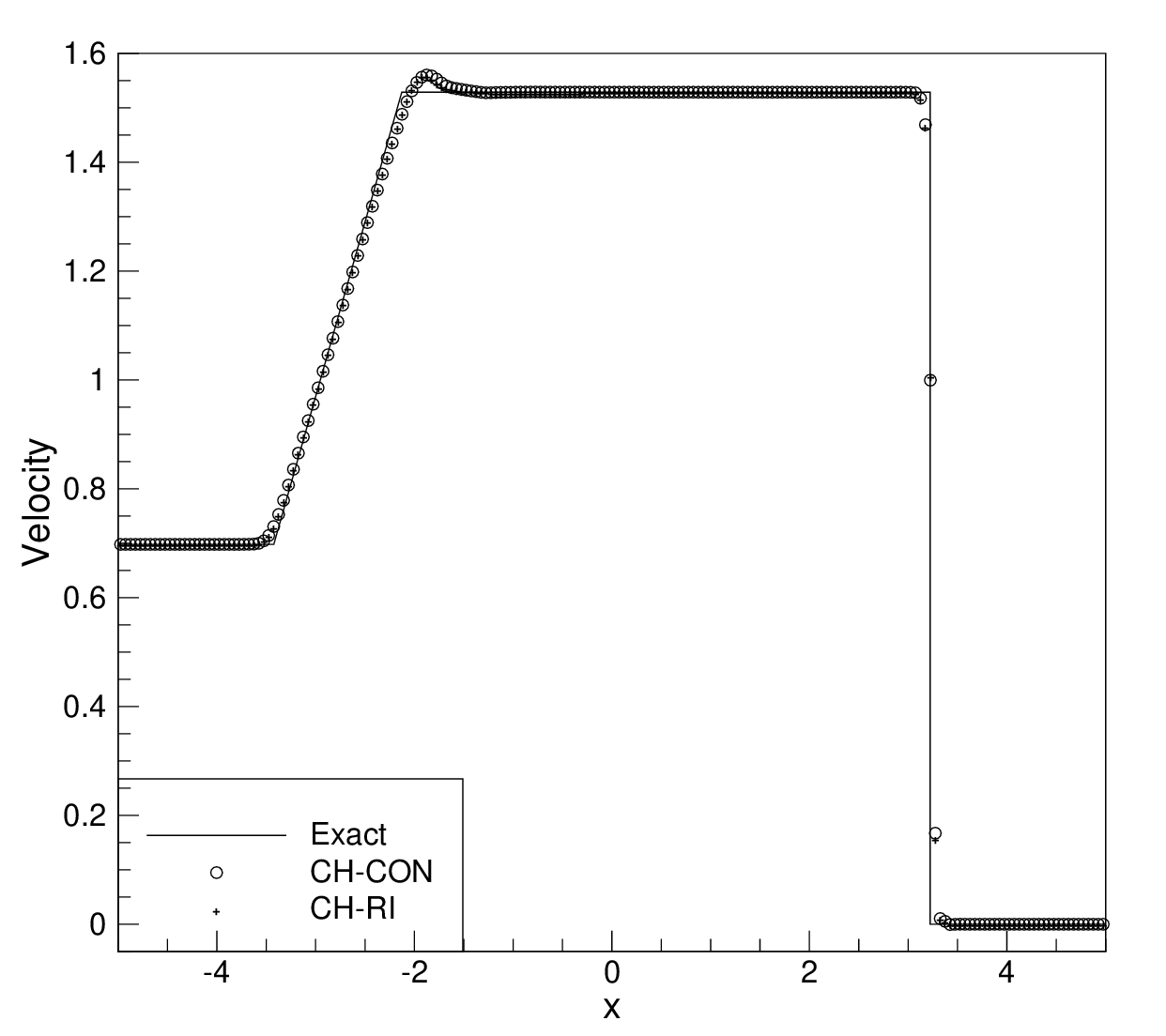}
        \caption{AWENO7}
    \end{subfigure}
    \begin{subfigure}{0.32\linewidth}
        \centering
        \includegraphics[width=\linewidth]{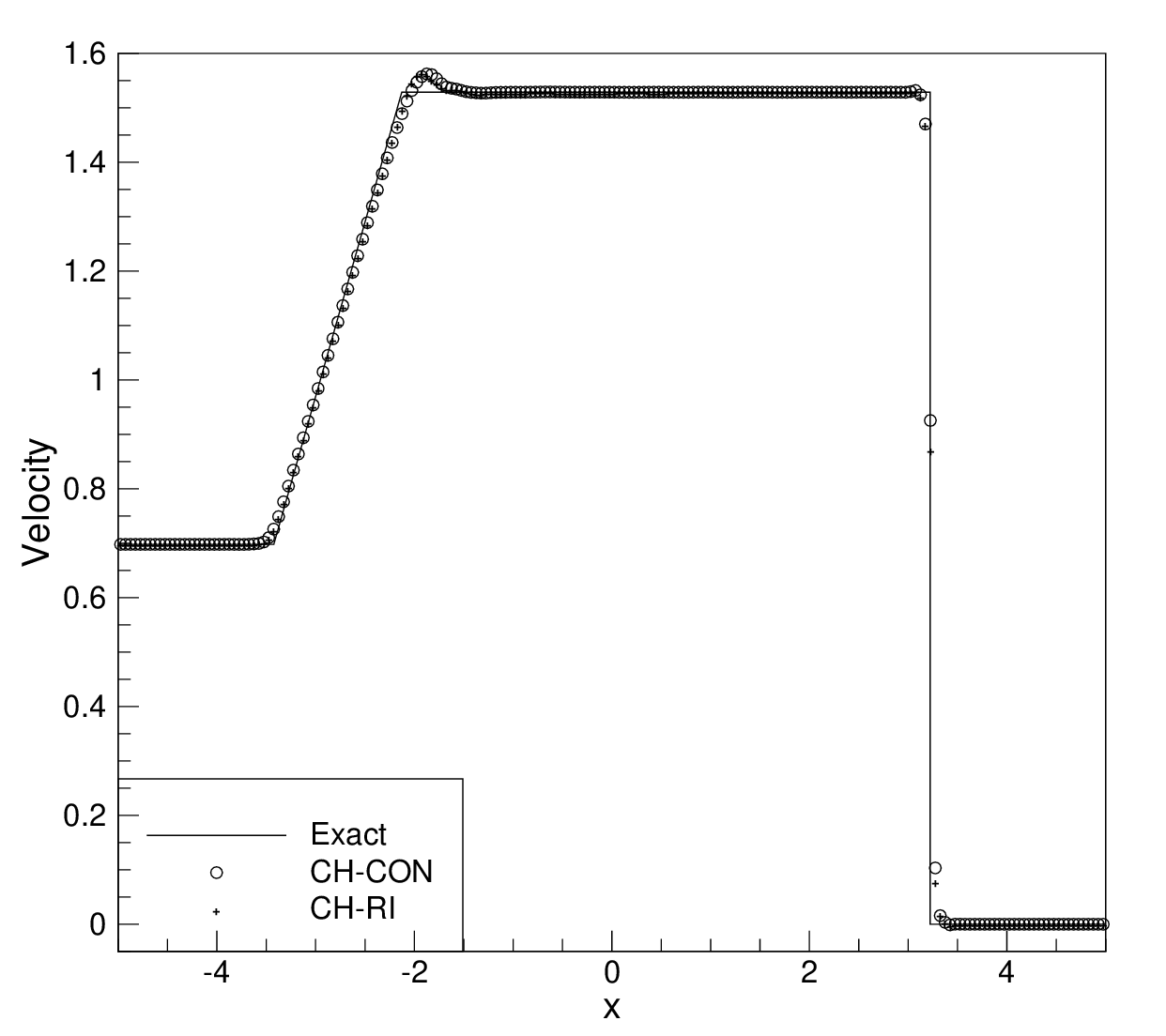}
        \caption{AWENO9}
    \end{subfigure}

    \begin{subfigure}{0.32\linewidth}
        \centering
        \includegraphics[width=\linewidth]{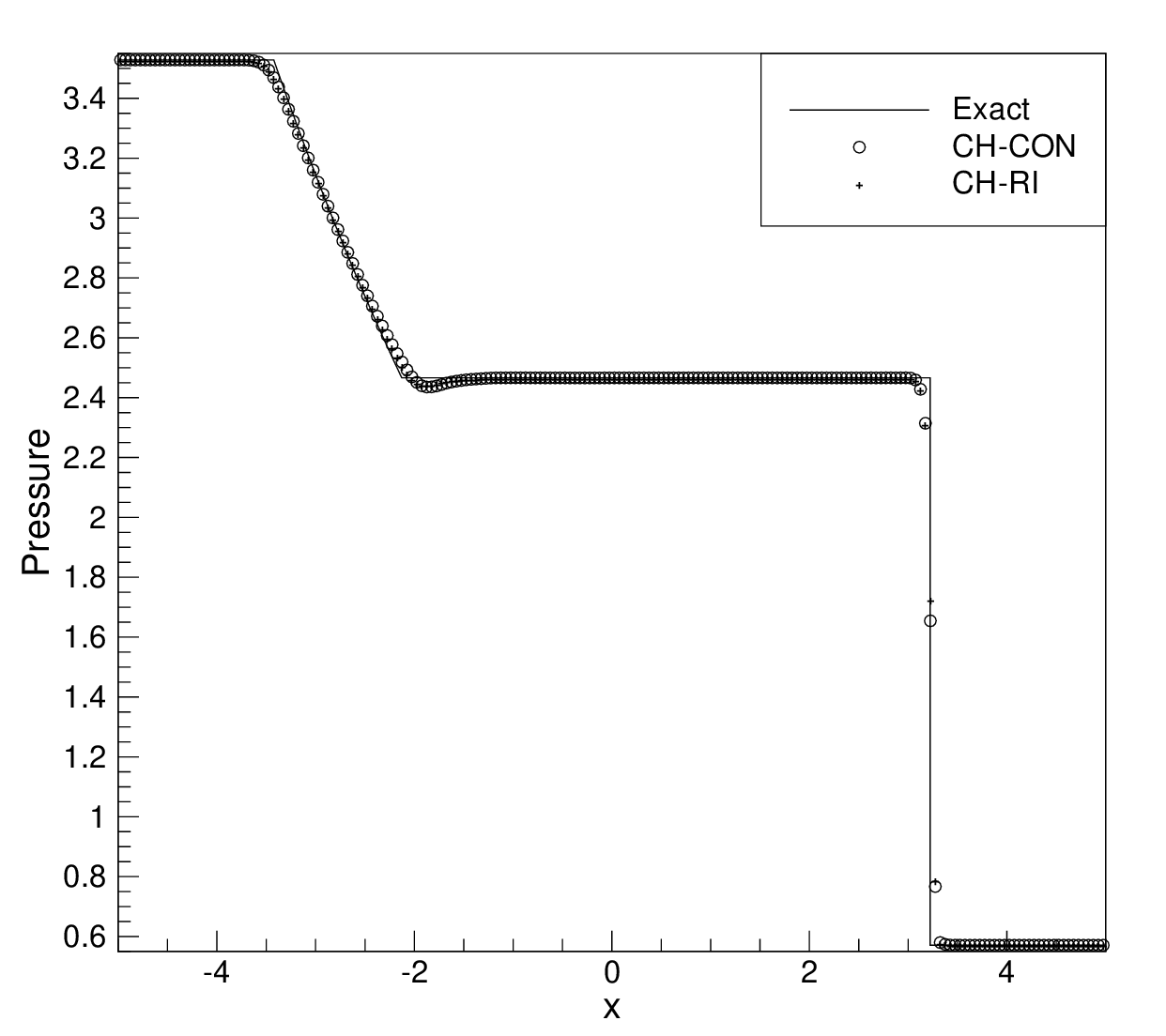}
        \caption{AWENO5}
    \end{subfigure}
    \begin{subfigure}{0.32\linewidth}
        \centering
        \includegraphics[width=\linewidth]{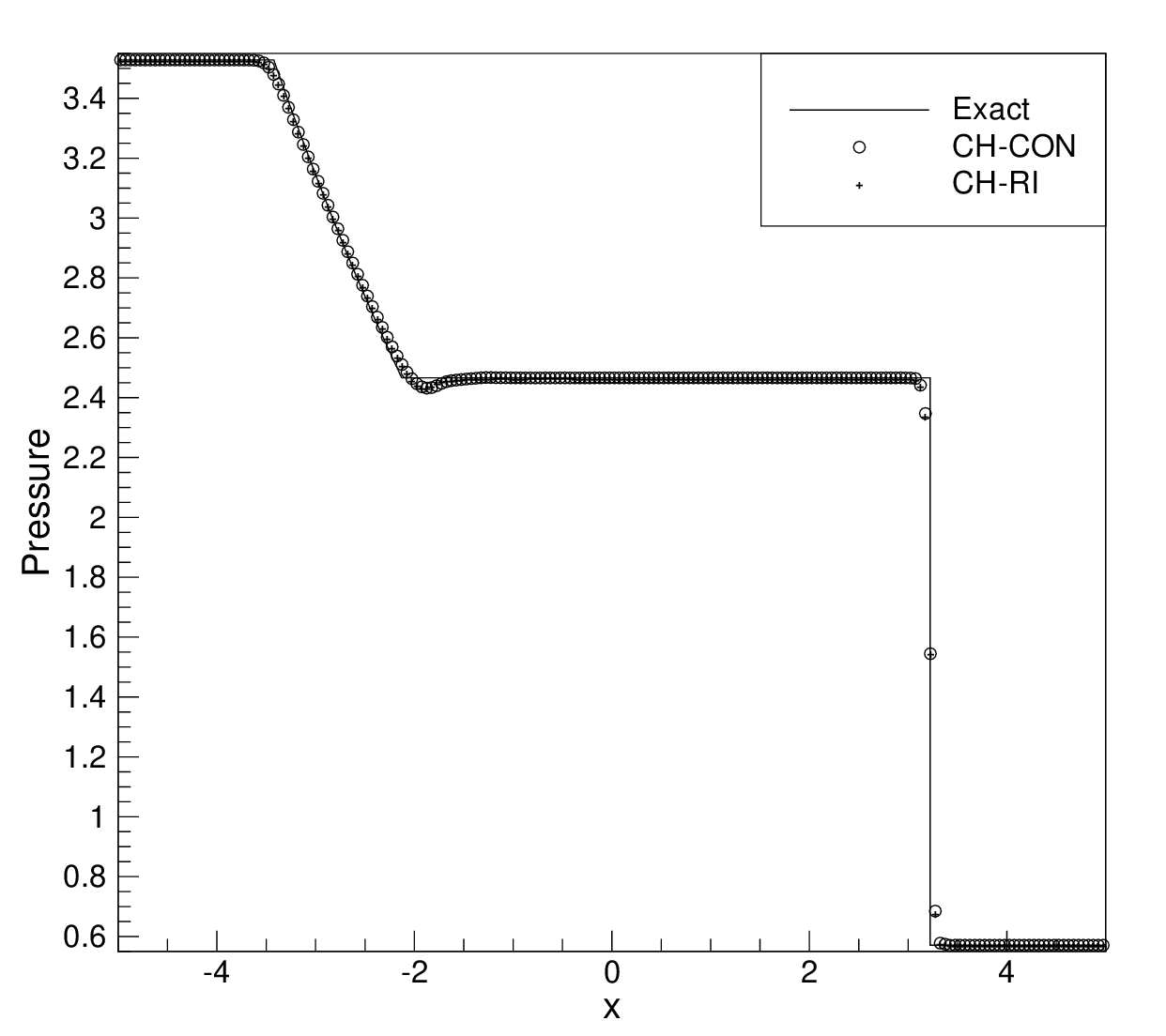}
        \caption{AWENO7}
    \end{subfigure}
    \begin{subfigure}{0.32\linewidth}
        \centering
        \includegraphics[width=\linewidth]{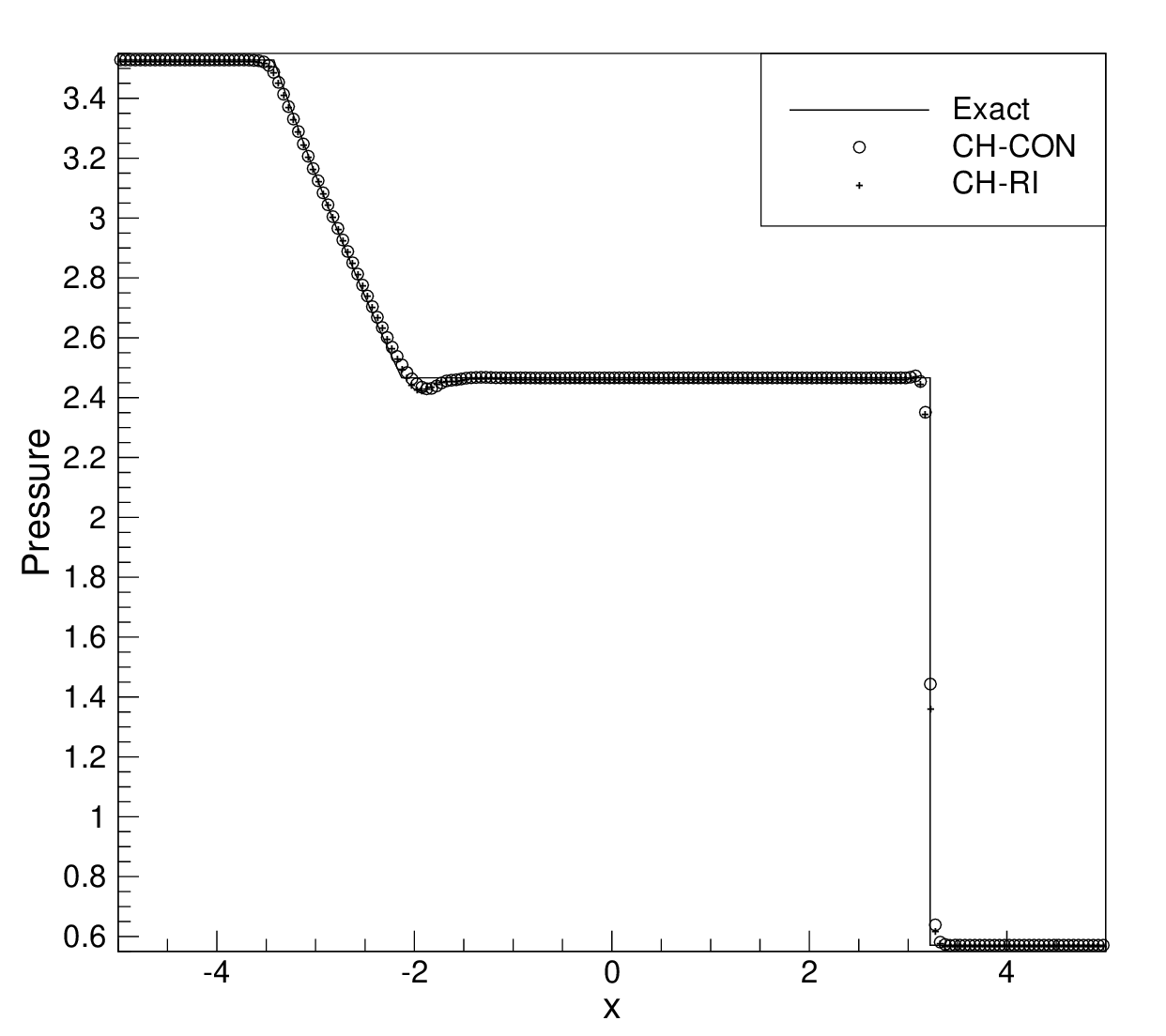}
        \caption{AWENO9}
    \end{subfigure}

    \caption{The Lax shock tube problem at \(T = 1.3\) with \(N=200\) grids.}
    \label{fig:Lax}
\end{figure}

The initial condition is \((\rho, u, p) = (0.445, 0.698, 3.528)\) for \(x < 0\) and \((\rho, u, p) = (0.5, 0, 0.571)\) for \(x > 0\). The computation domain is \([-5, 5]\). We solve the problem to \(T = 1.3\). The results are shown in \Cref{fig:Lax}. The CH-RI and CH-CON agree with each other well. 

\subsubsection{The LeBlanc shock tube problem}

\begin{figure}[htbp]
    \centering

    \begin{subfigure}{0.32\linewidth}
        \centering
        \includegraphics[width=\linewidth]{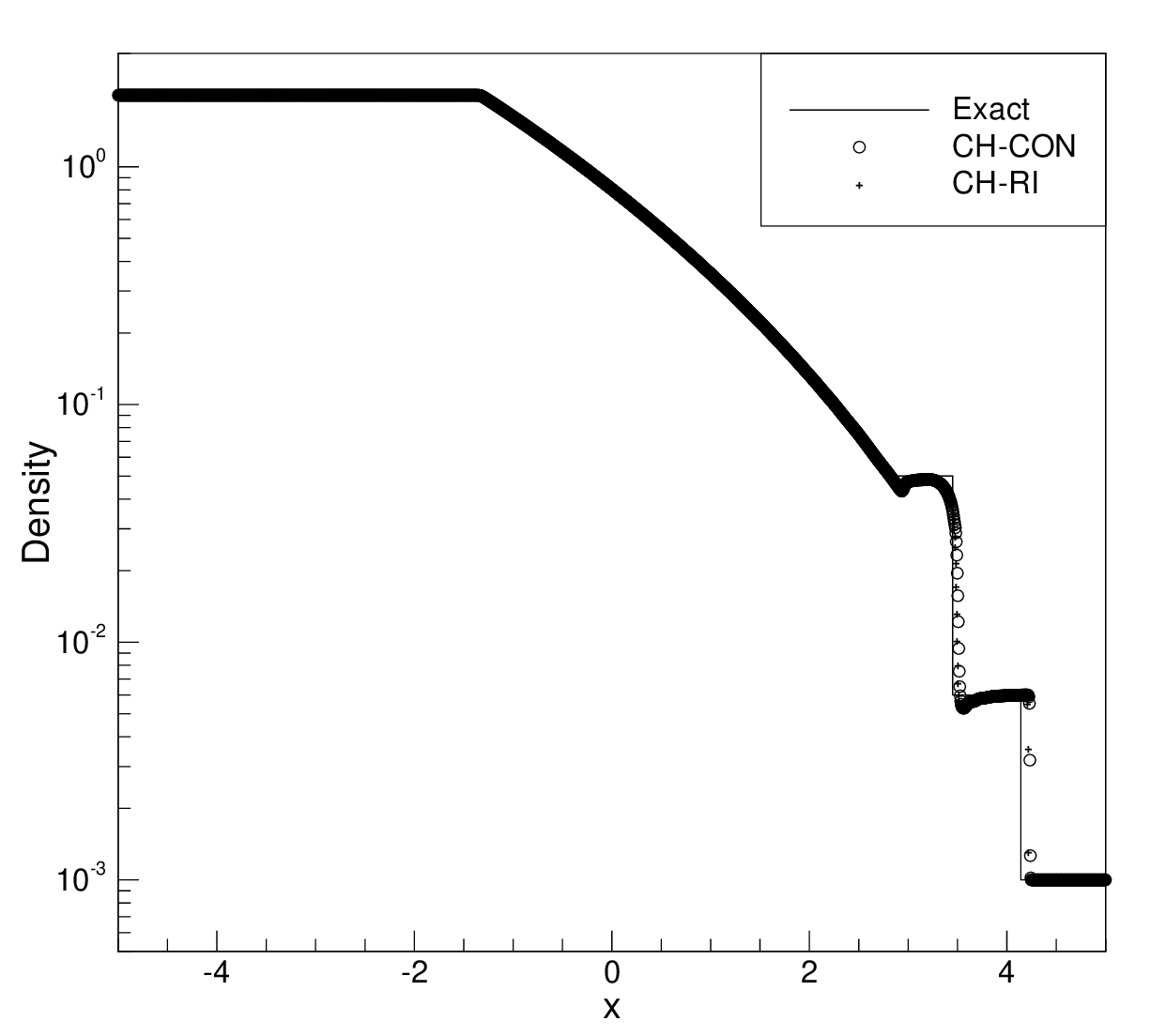}
        \caption{AWENO5}
    \end{subfigure}
    \begin{subfigure}{0.32\linewidth}
        \centering
        \includegraphics[width=\linewidth]{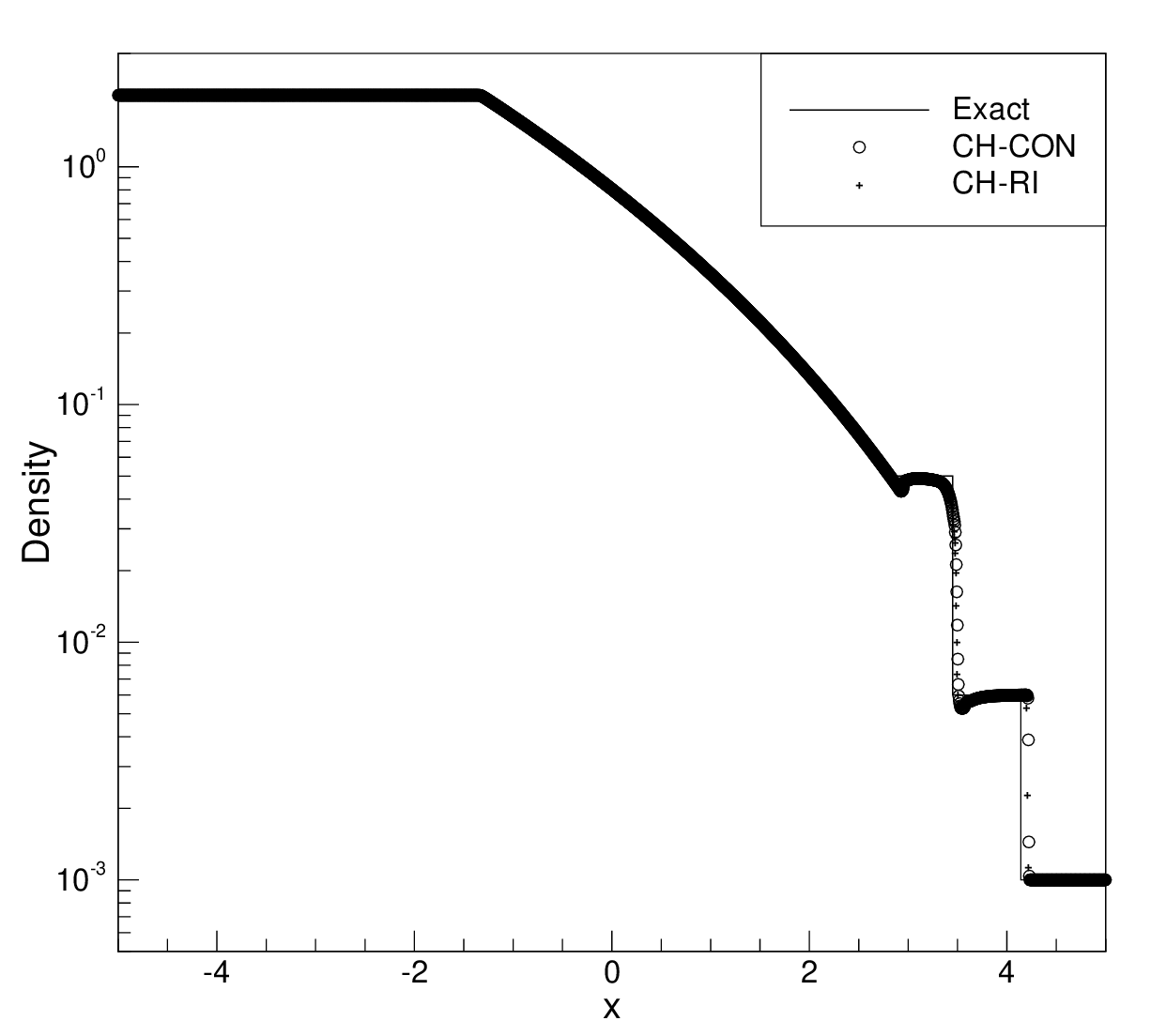}
        \caption{AWENO7}
    \end{subfigure}
    \begin{subfigure}{0.32\linewidth}
        \centering
        \includegraphics[width=\linewidth]{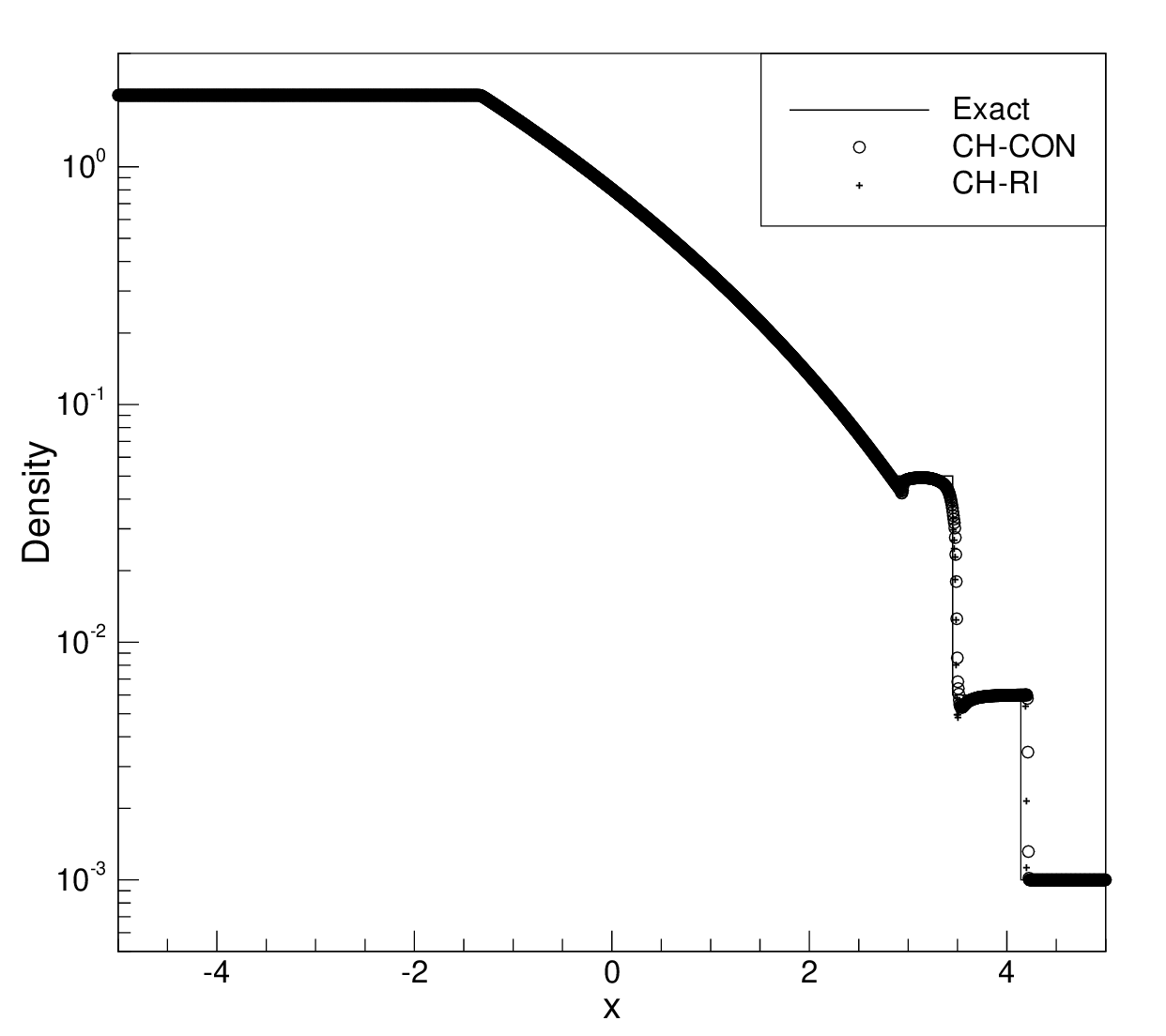}
        \caption{AWENO9}
    \end{subfigure}

    \begin{subfigure}{0.32\linewidth}
        \centering
        \includegraphics[width=\linewidth]{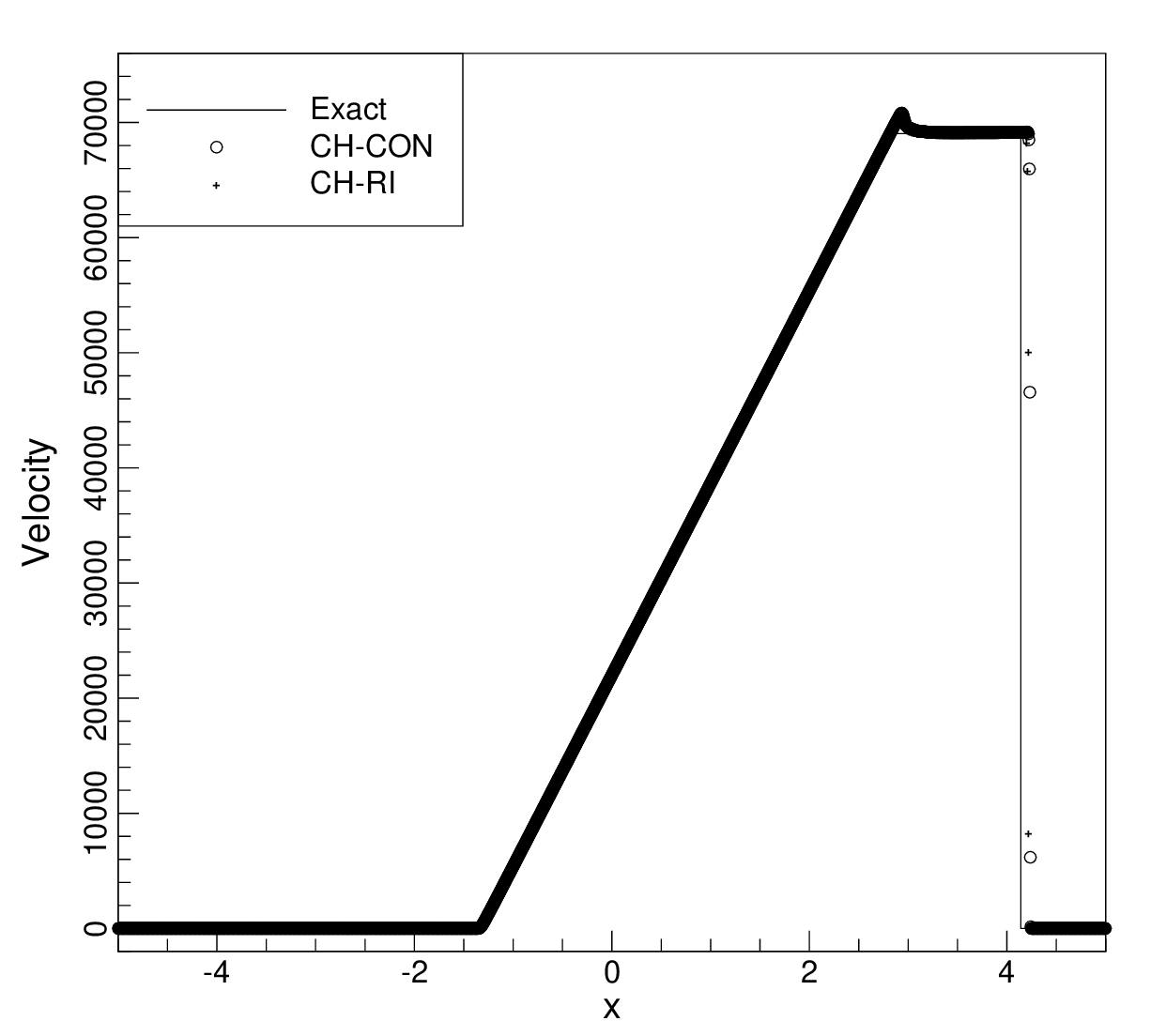}
        \caption{AWENO5}
    \end{subfigure}
    \begin{subfigure}{0.32\linewidth}
        \centering
        \includegraphics[width=\linewidth]{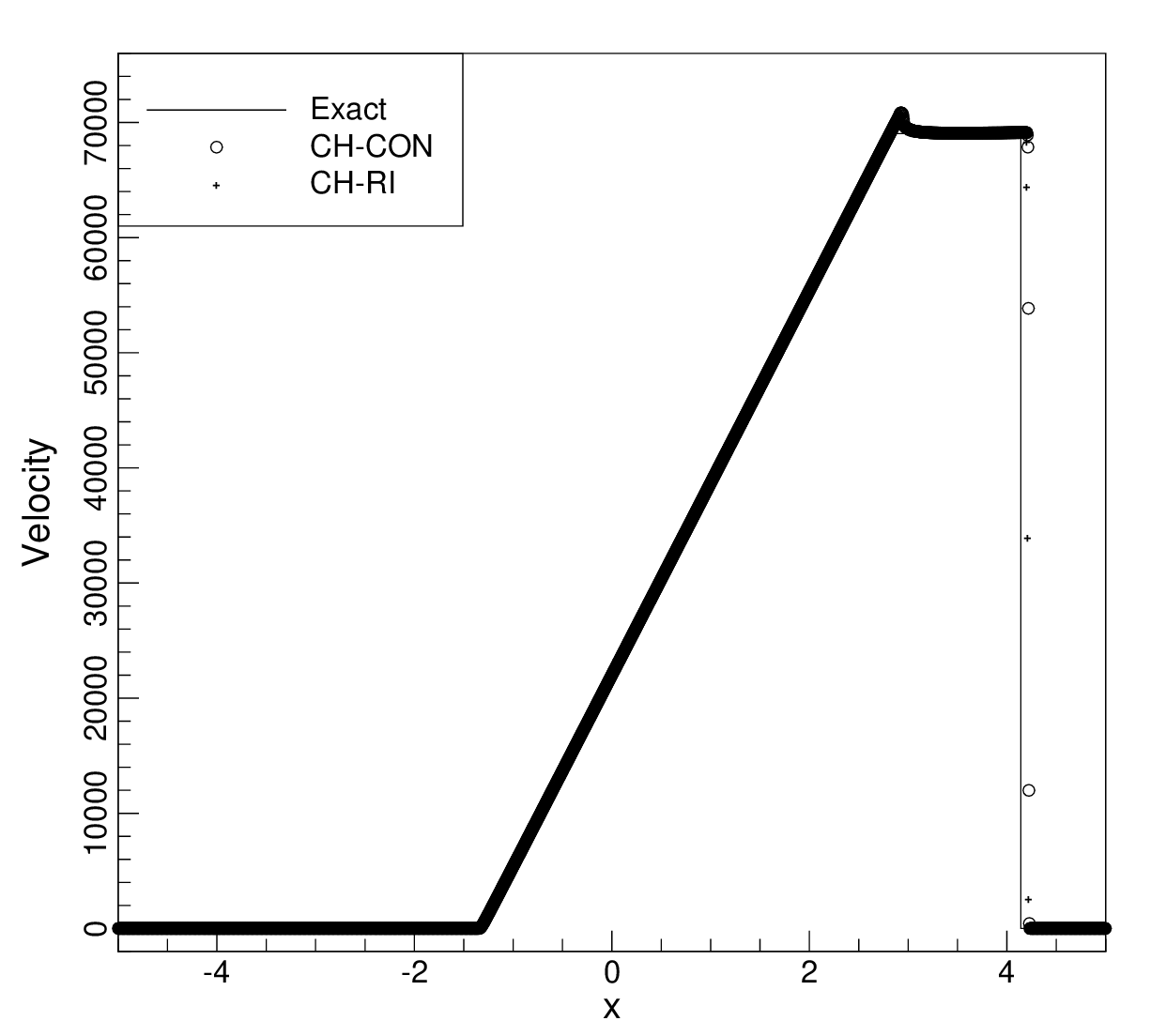}
        \caption{AWENO7}
    \end{subfigure}
    \begin{subfigure}{0.32\linewidth}
        \centering
        \includegraphics[width=\linewidth]{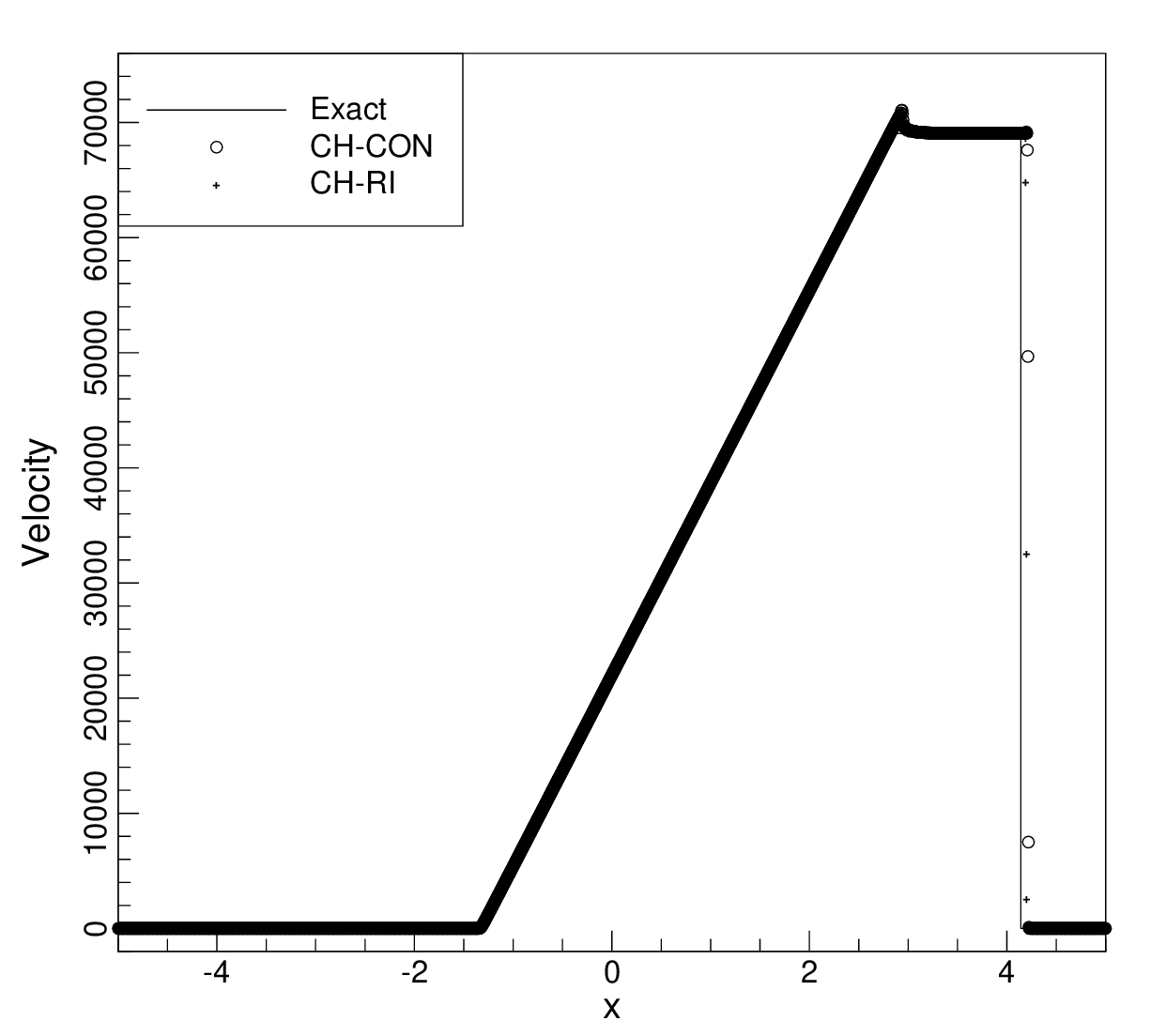}
        \caption{AWENO9}
    \end{subfigure}

    \begin{subfigure}{0.32\linewidth}
        \centering
        \includegraphics[width=\linewidth]{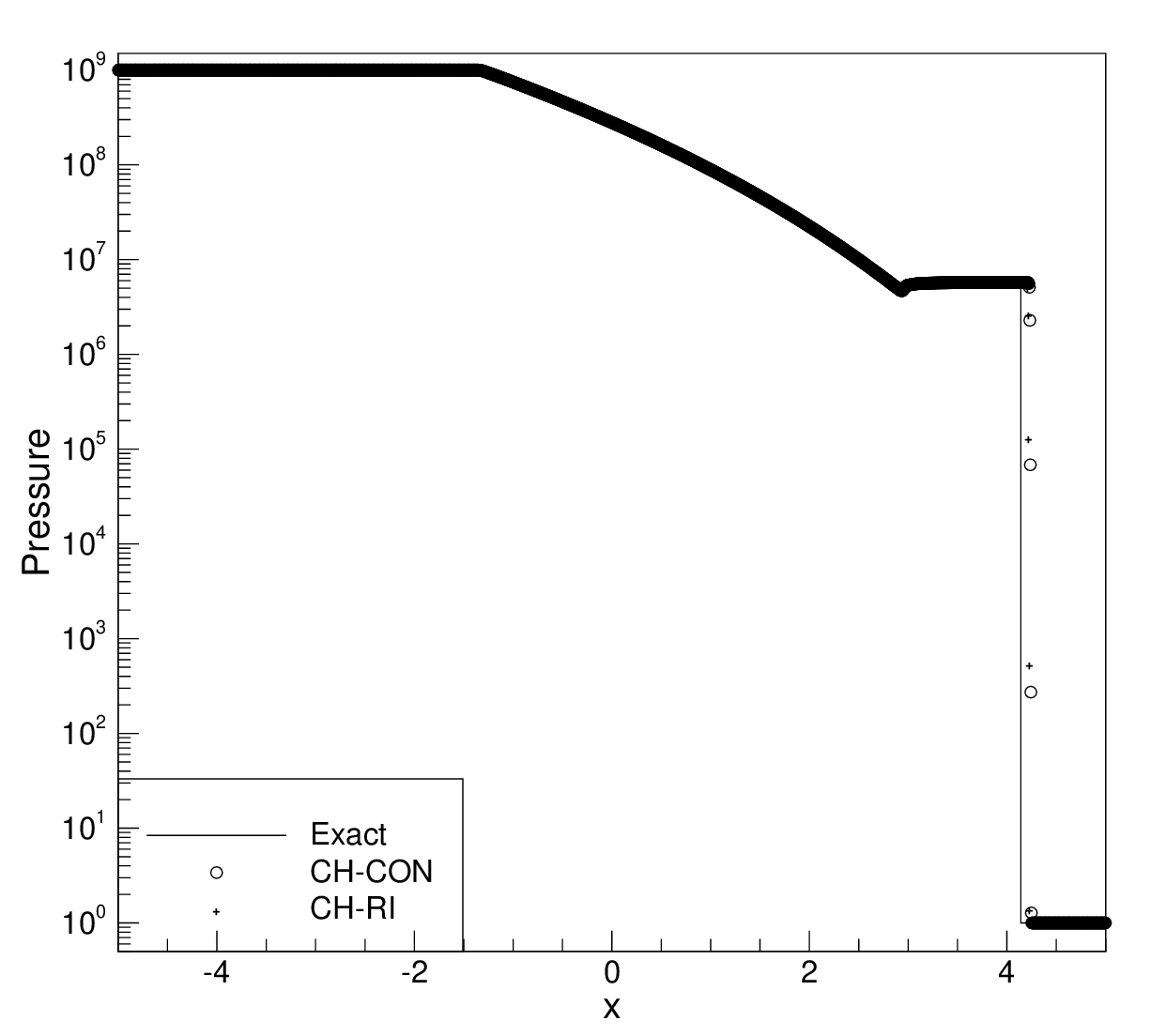}
        \caption{AWENO5}
    \end{subfigure}
    \begin{subfigure}{0.32\linewidth}
        \centering
        \includegraphics[width=\linewidth]{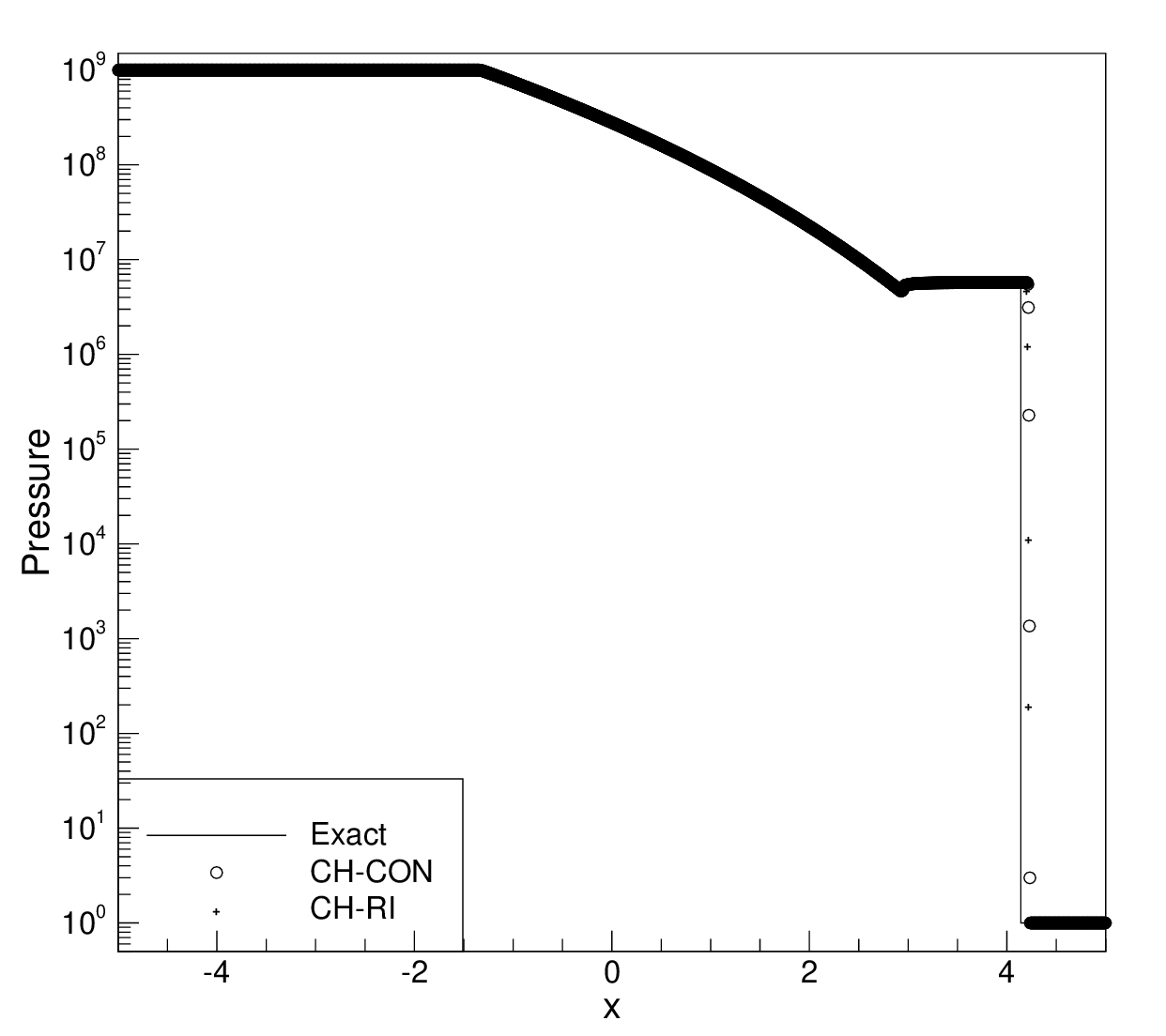}
        \caption{AWENO7}
    \end{subfigure}
    \begin{subfigure}{0.32\linewidth}
        \centering
        \includegraphics[width=\linewidth]{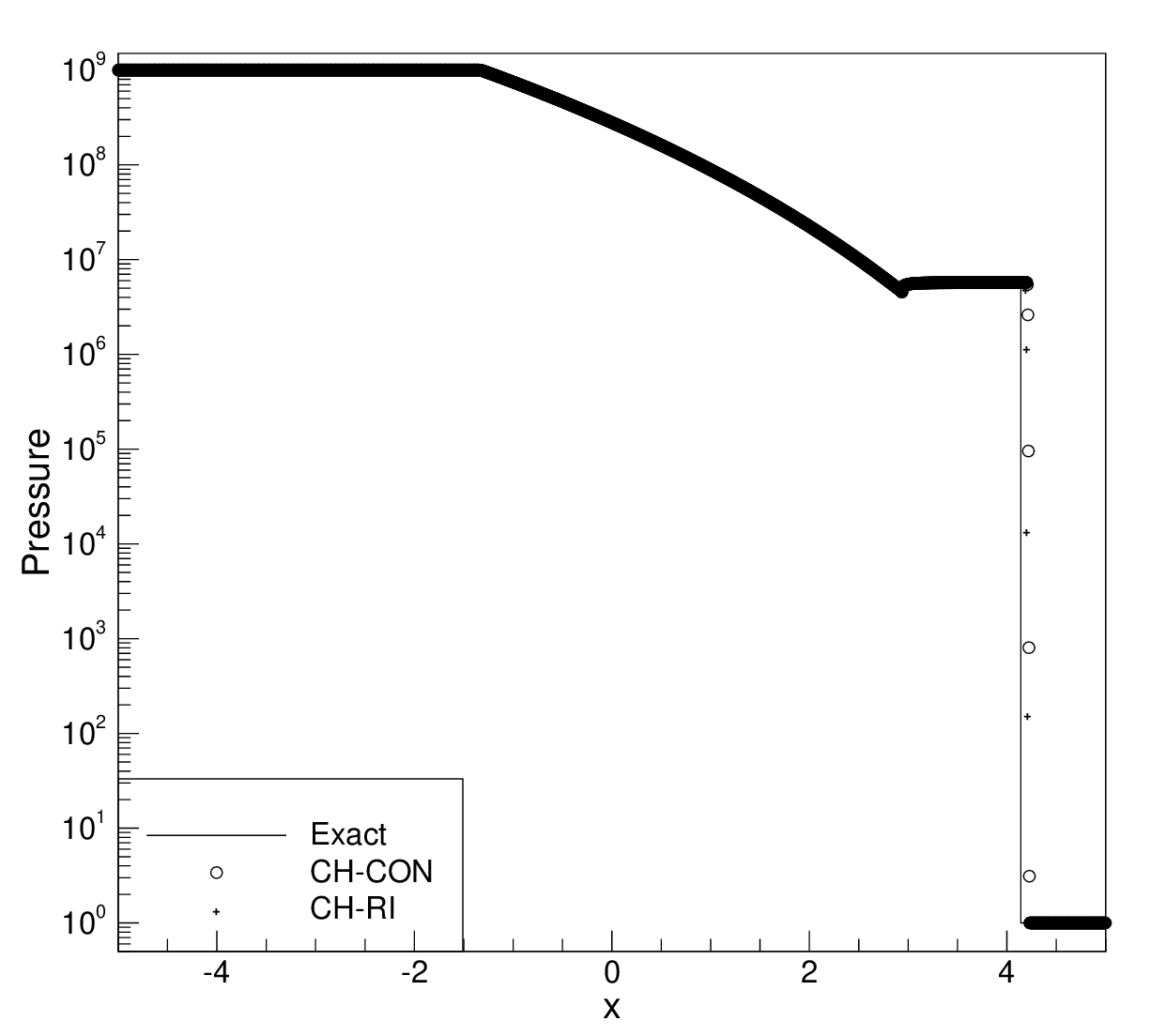}
        \caption{AWENO9}
    \end{subfigure}

    \caption{The LeBlanc shock tube problem at \(T = 5 \times 10^{-5}\) with \(N=2000\) grids.}
    \label{fig:LeBlanc}
\end{figure}

The initial condition is \((\rho, u, p) = (2, 0, 10^{9})\) for \(x < 0\) and \((\rho, u, p) = (10^{-3}, 0, 1)\) for \(x > 0\). The computation domain is \([-5, 5]\). We solve the problem to \(T = 5 \times 10^{-5}\). The results are shown in \Cref{fig:LeBlanc}. The extreme pressure ratio makes the problem very stiff and challenges the PP property of the scheme. With our PP limiters, both CH-RI and CH-CON perform quite well. 

\subsubsection{The critical double rarefaction problem}

\begin{figure}[htbp]
    \centering

    \begin{subfigure}{0.32\linewidth}
        \centering
        \includegraphics[width=\linewidth]{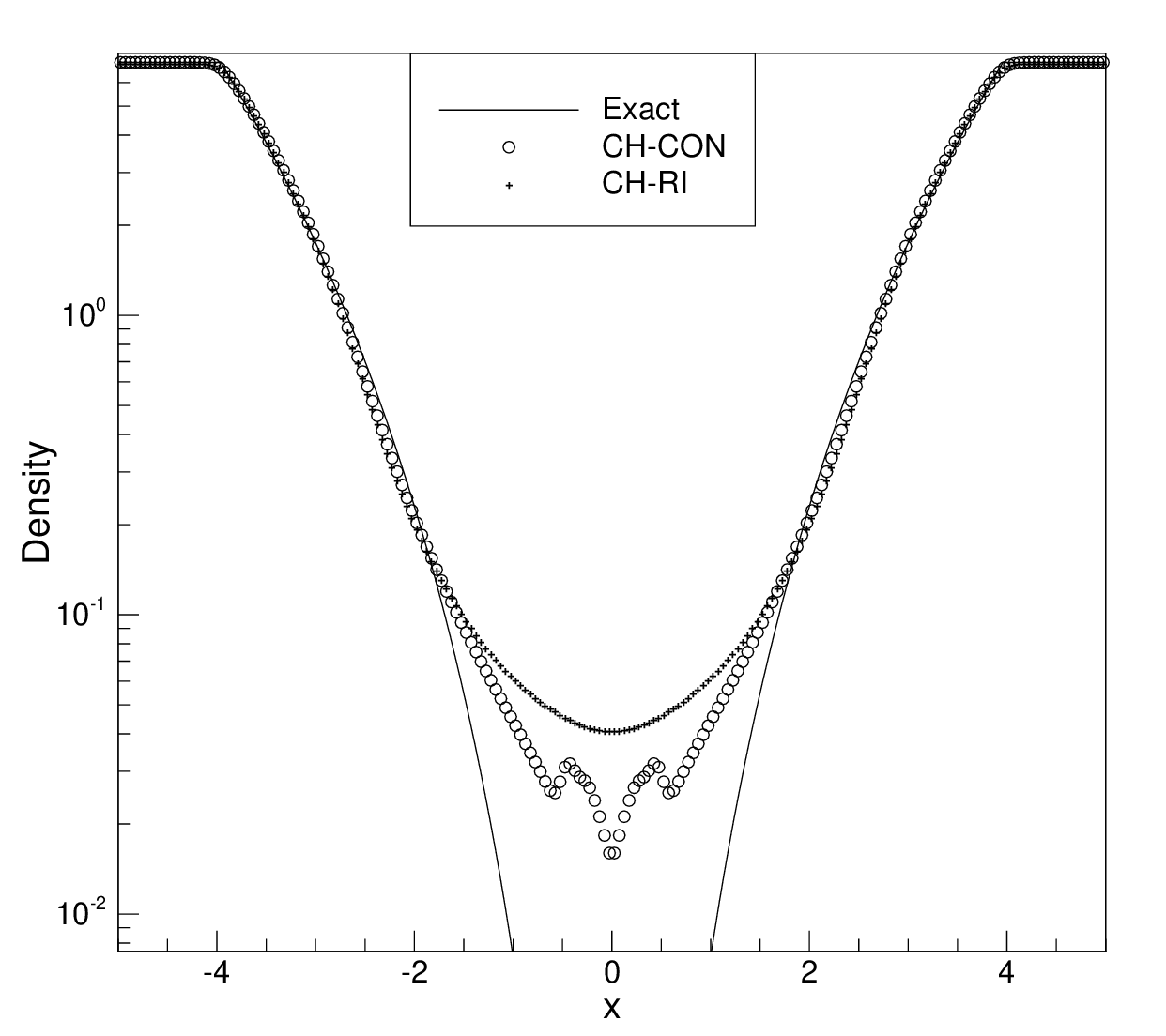}
        \caption{AWENO5}
    \end{subfigure}
    \begin{subfigure}{0.32\linewidth}
        \centering
        \includegraphics[width=\linewidth]{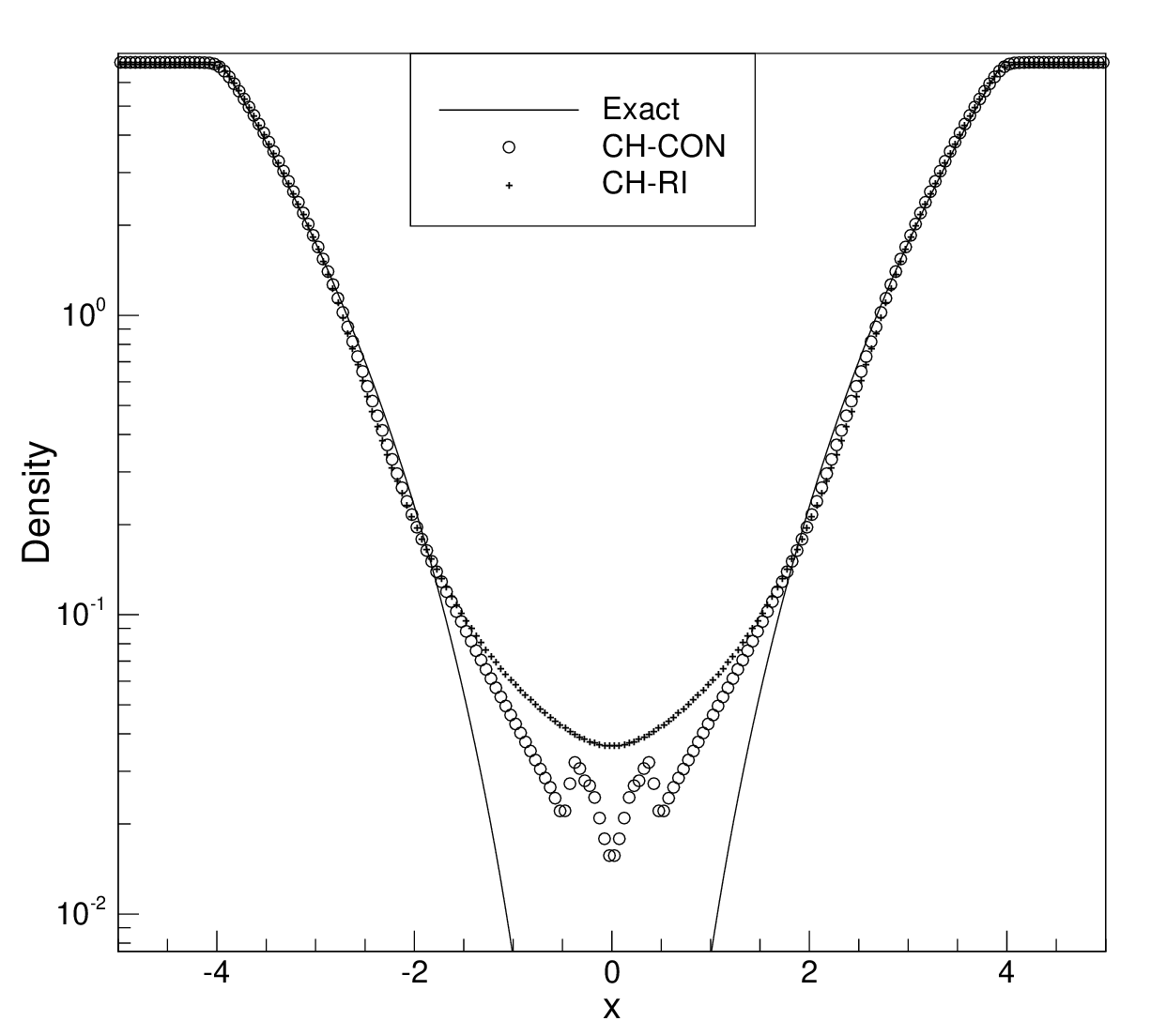}
        \caption{AWENO7}
    \end{subfigure}
    \begin{subfigure}{0.32\linewidth}
        \centering
        \includegraphics[width=\linewidth]{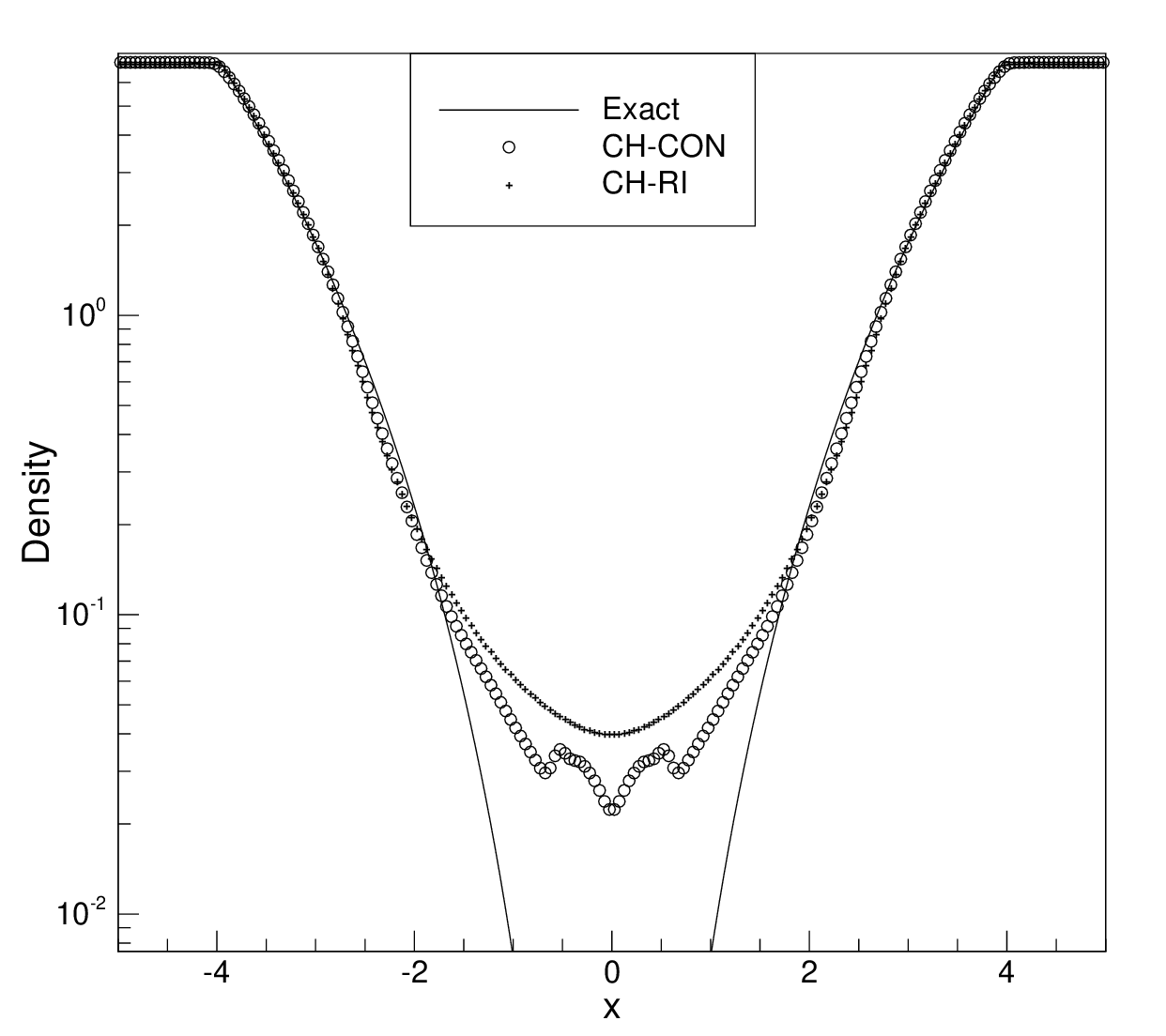}
        \caption{AWENO9}
    \end{subfigure}

    \begin{subfigure}{0.32\linewidth}
        \centering
        \includegraphics[width=\linewidth]{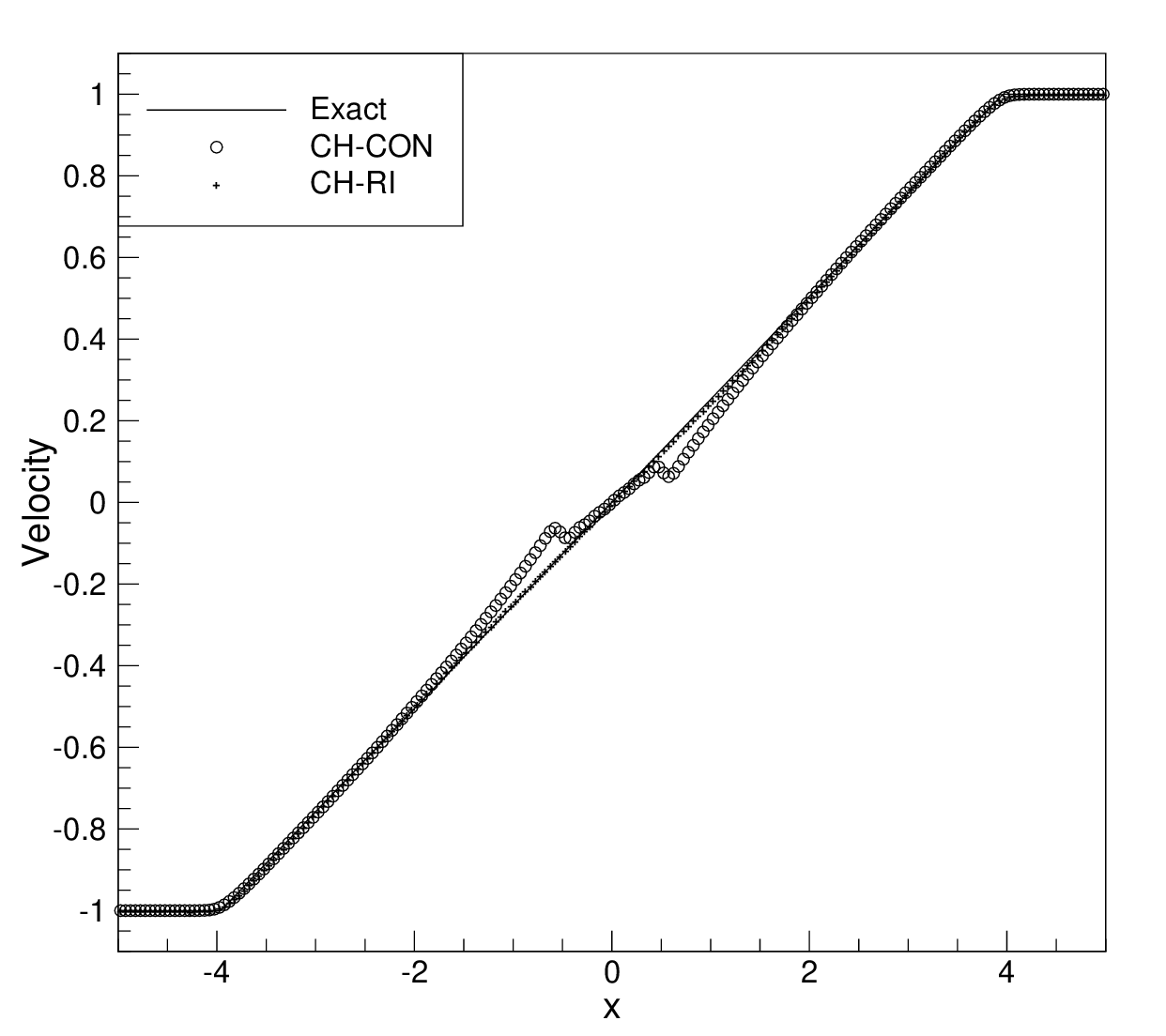}
        \caption{AWENO5}
    \end{subfigure}
    \begin{subfigure}{0.32\linewidth}
        \centering
        \includegraphics[width=\linewidth]{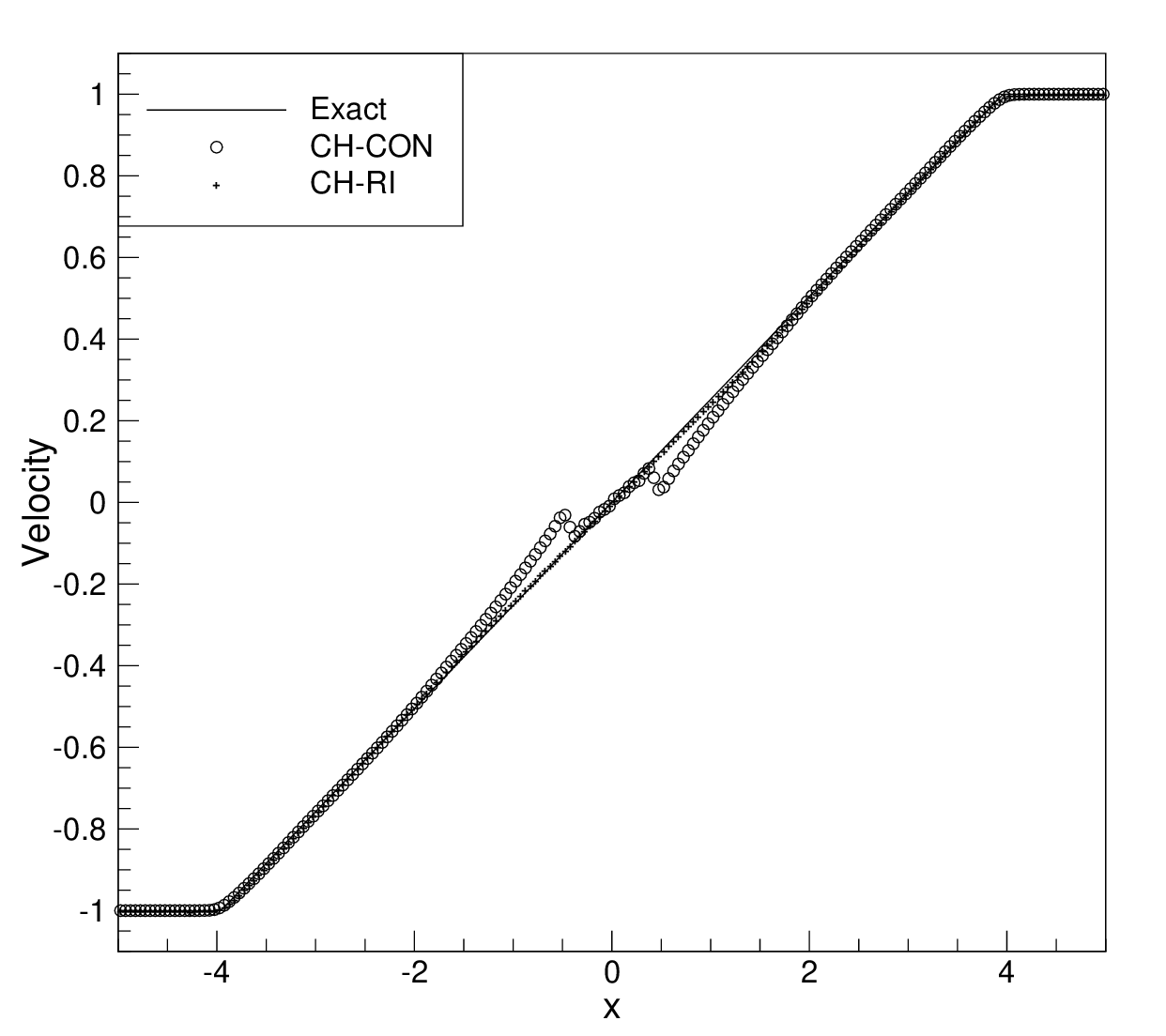}
        \caption{AWENO7}
    \end{subfigure}
    \begin{subfigure}{0.32\linewidth}
        \centering
        \includegraphics[width=\linewidth]{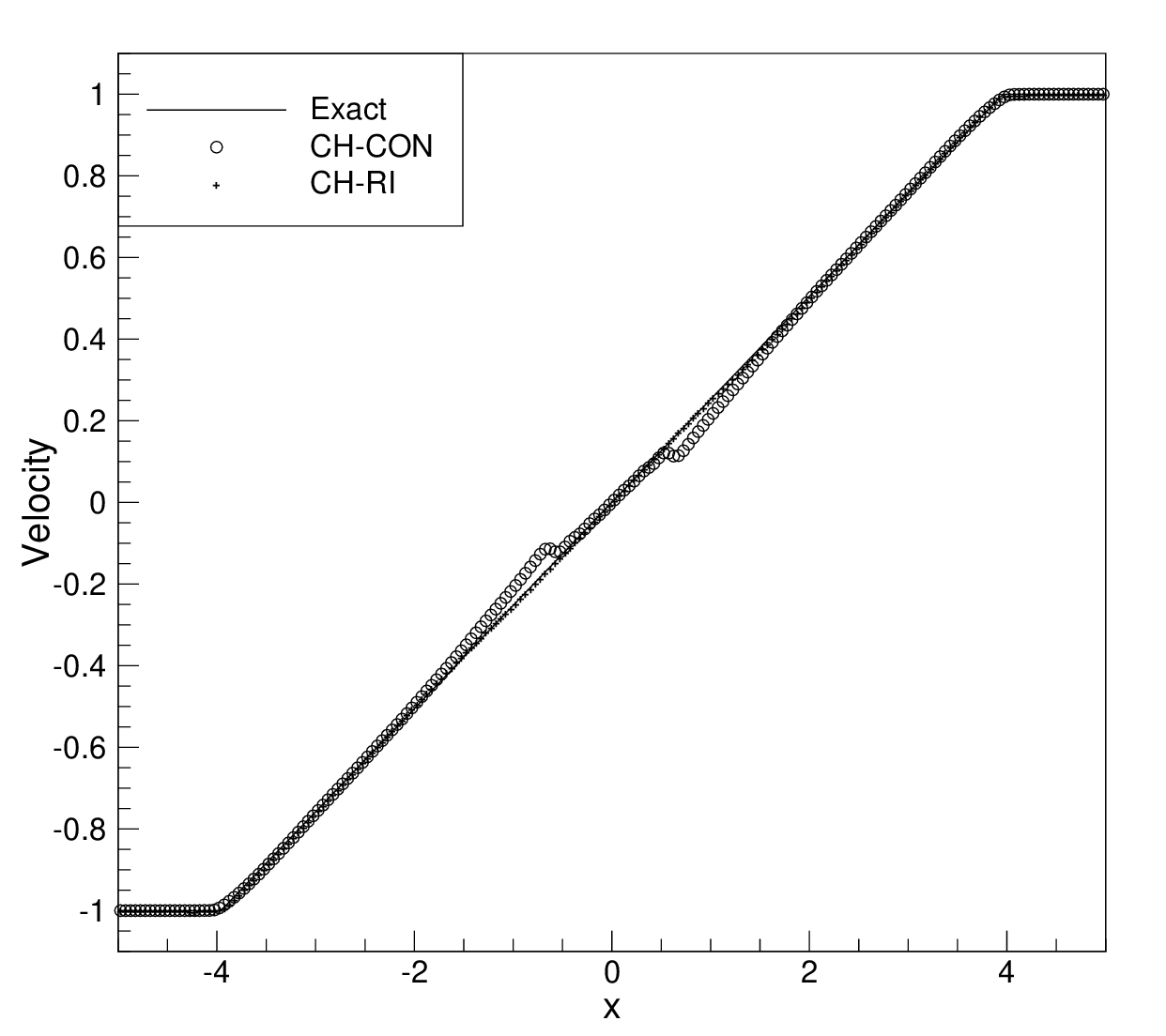}
        \caption{AWENO9}
    \end{subfigure}

    \begin{subfigure}{0.32\linewidth}
        \centering
        \includegraphics[width=\linewidth]{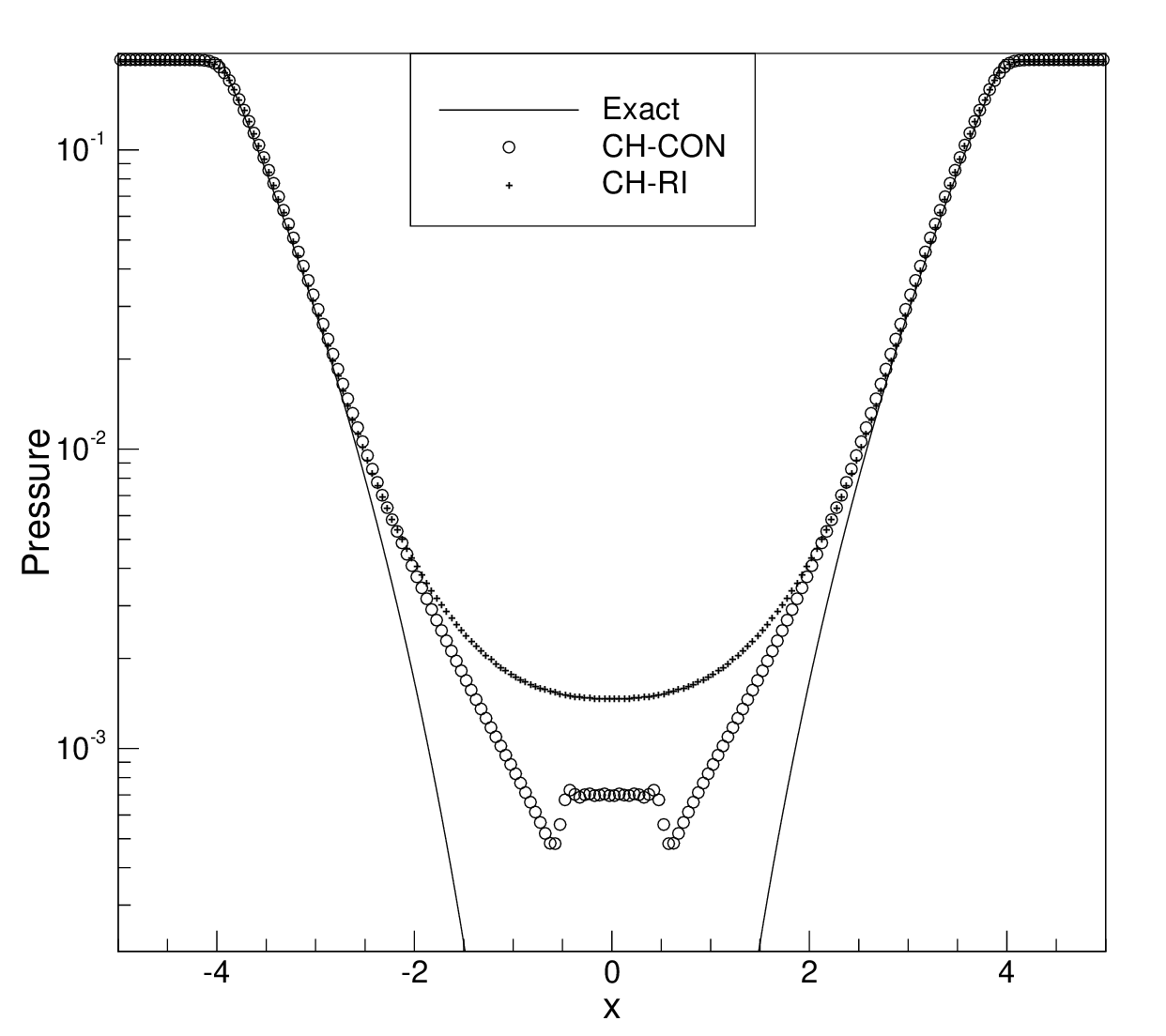}
        \caption{AWENO5}
    \end{subfigure}
    \begin{subfigure}{0.32\linewidth}
        \centering
        \includegraphics[width=\linewidth]{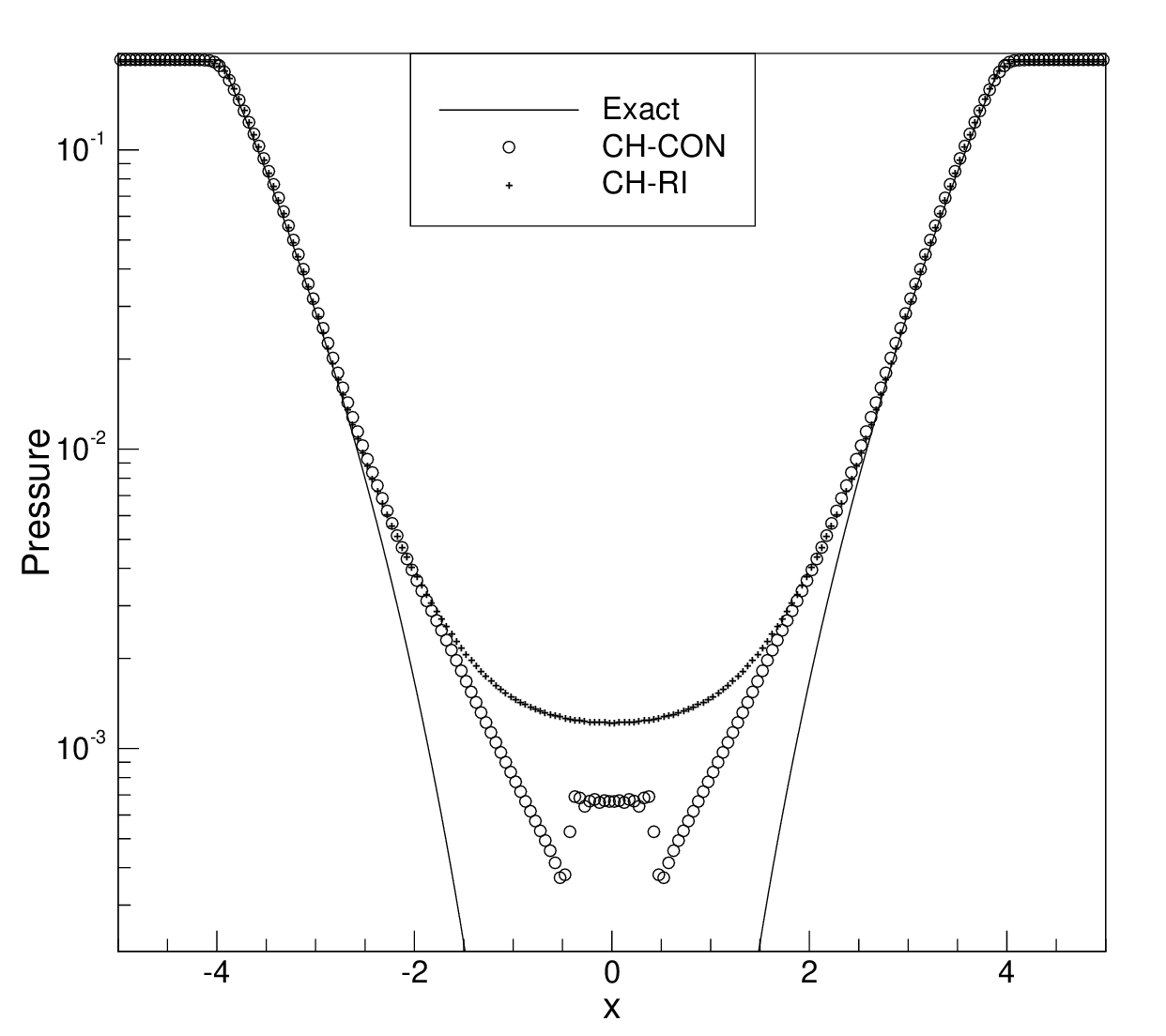}
        \caption{AWENO7}
    \end{subfigure}
    \begin{subfigure}{0.32\linewidth}
        \centering
        \includegraphics[width=\linewidth]{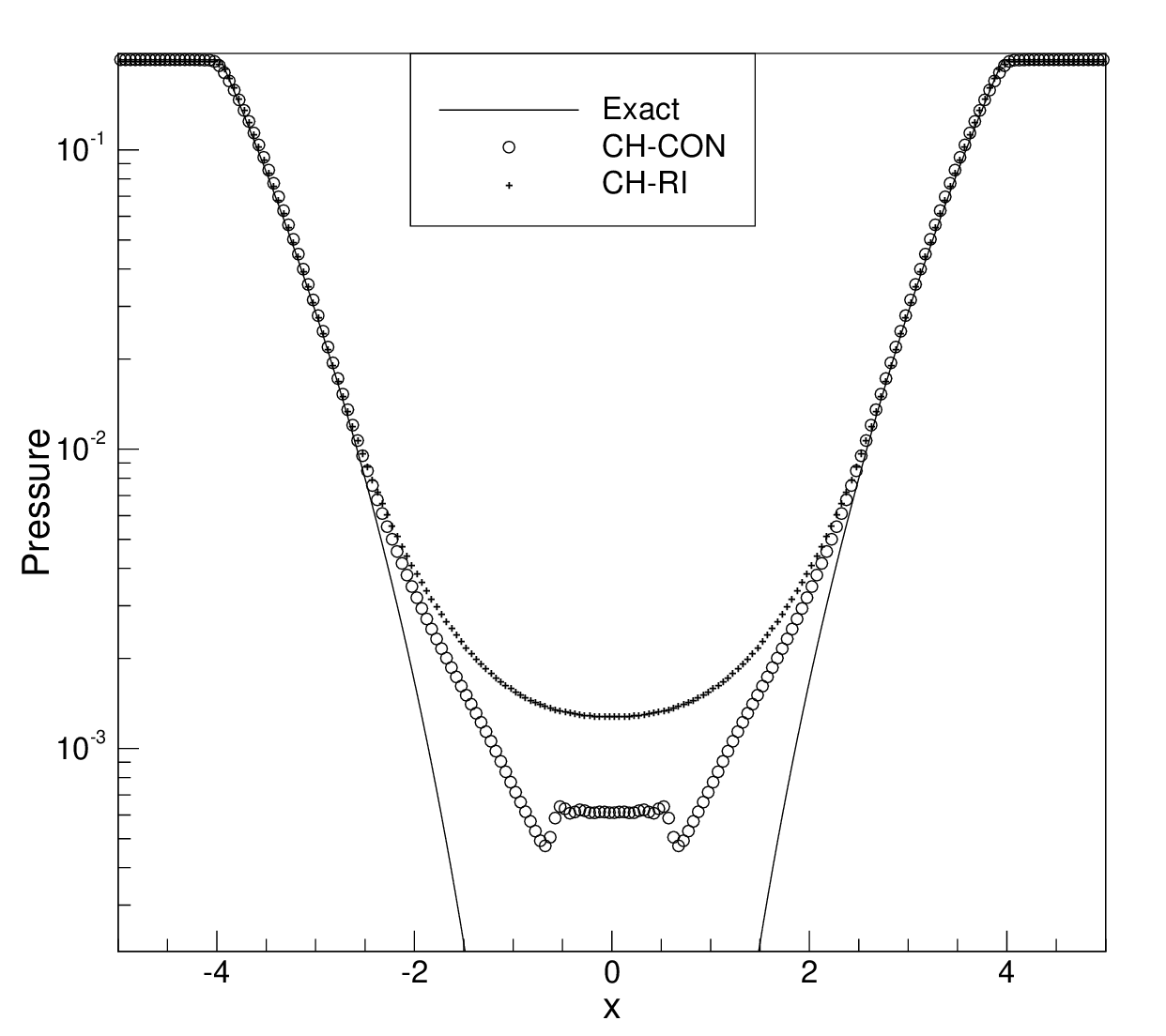}
        \caption{AWENO9}
    \end{subfigure}

    \caption{The critical double rarefaction problem at \(T = 3.3\) with \(N=200\) grids.}
    \label{fig:double-rarefaction}
\end{figure}

The initial condition is \((\rho, u, p) = (7, -1, 0.2)\) for \(x < 0\) and \((\rho, u, p) = (7, 1, 0.2)\) for \(x > 0\). The computation domain is \([-5, 5]\). We solve the problem to \(T = 3.3\). The results are shown in \Cref{fig:double-rarefaction}. The point \(x=0\) is the only vacuum point in the exact solution for \(t > 0\), so it causes difficulty for the scheme. Because the PP interpolation limiters for the two versions are different, there seems to be some difference near the vacuum point in the numerical solution. 

\subsubsection{The blast wave interaction problem}

\begin{figure}[htbp]
    \centering

    \begin{subfigure}{0.32\linewidth}
        \centering
        \includegraphics[width=\linewidth]{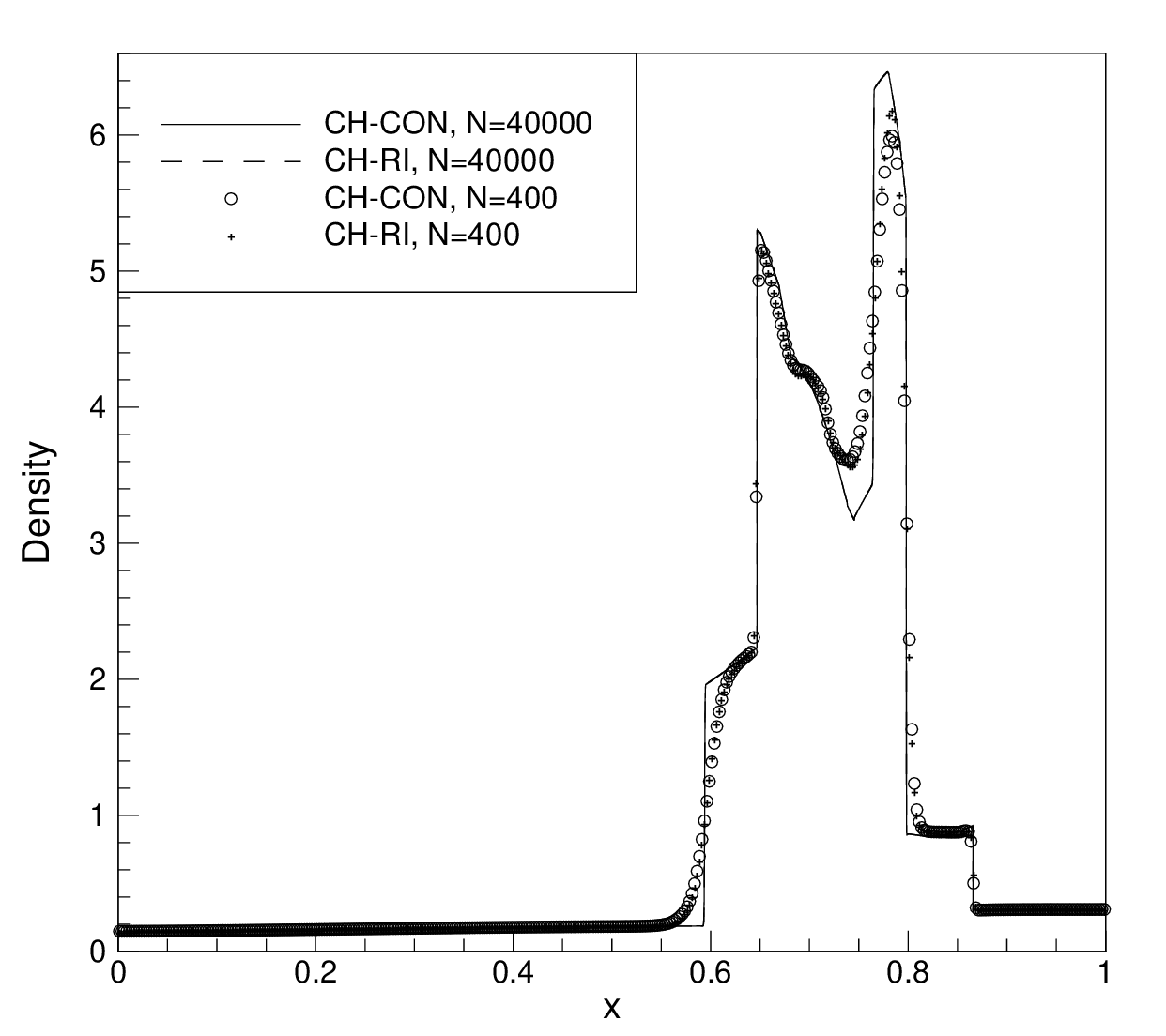}
        \caption{AWENO5}
    \end{subfigure}
    \begin{subfigure}{0.32\linewidth}
        \centering
        \includegraphics[width=\linewidth]{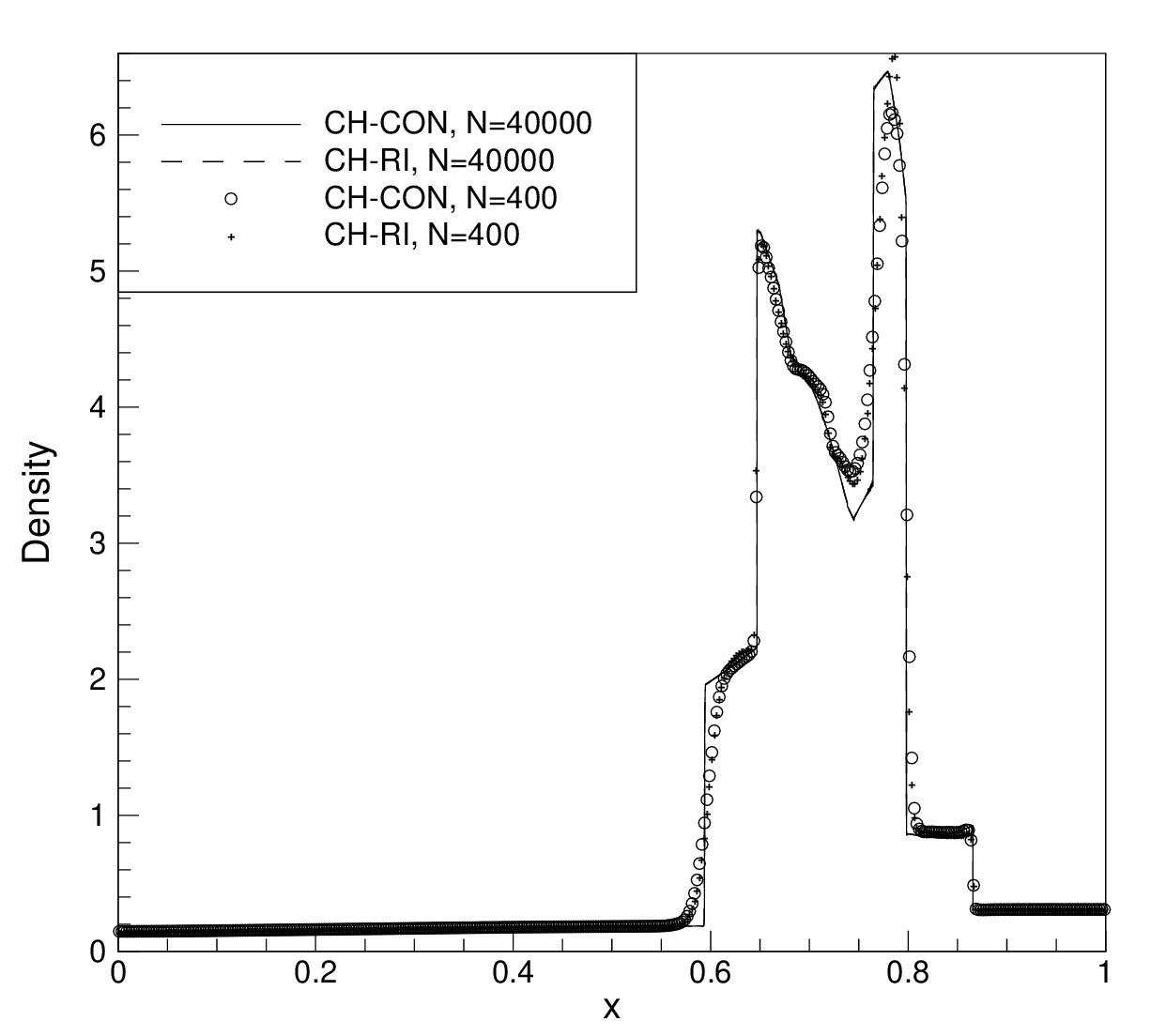}
        \caption{AWENO7}
    \end{subfigure}
    \begin{subfigure}{0.32\linewidth}
        \centering
        \includegraphics[width=\linewidth]{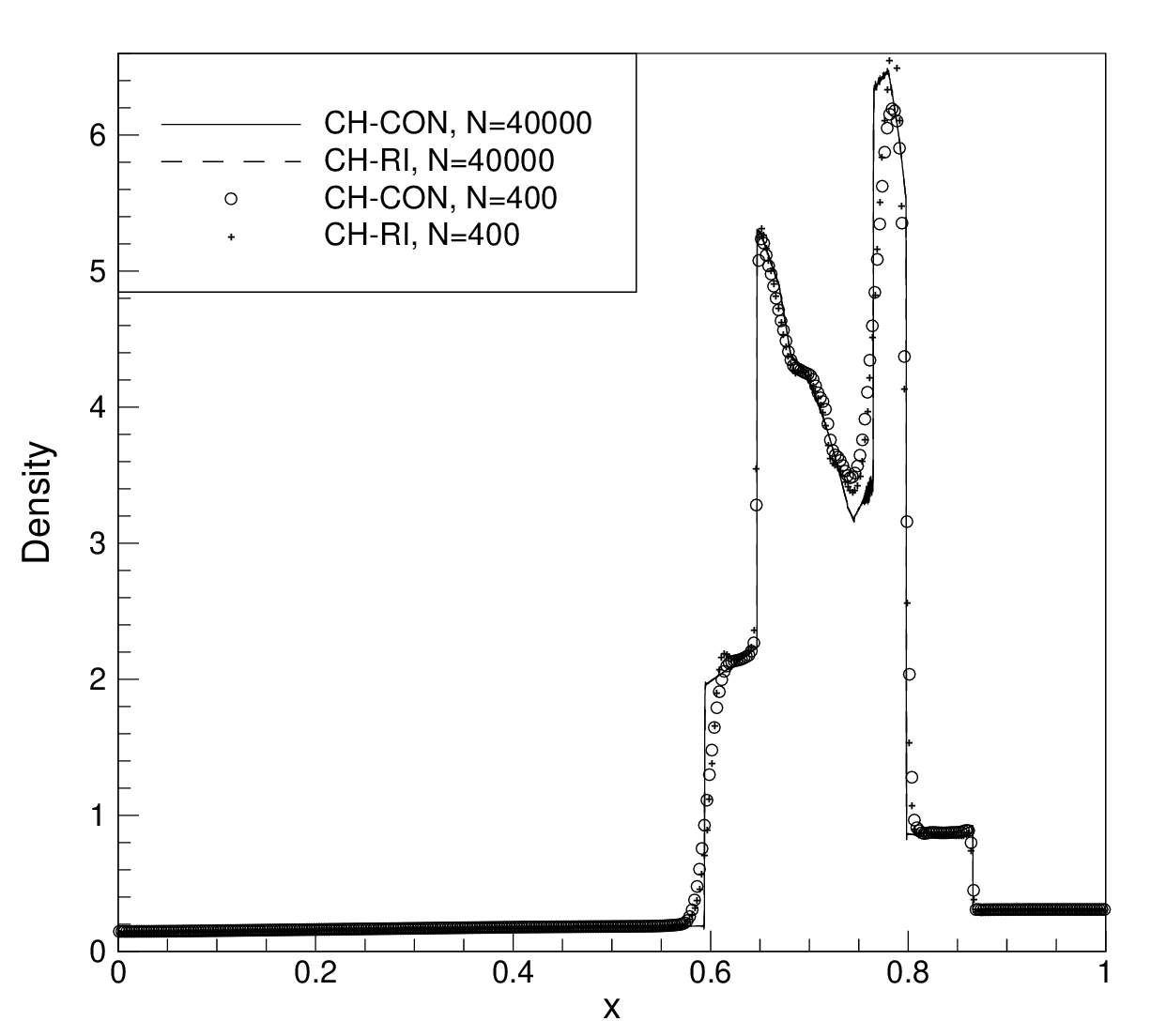}
        \caption{AWENO9}
    \end{subfigure}

    \begin{subfigure}{0.32\linewidth}
        \centering
        \includegraphics[width=\linewidth]{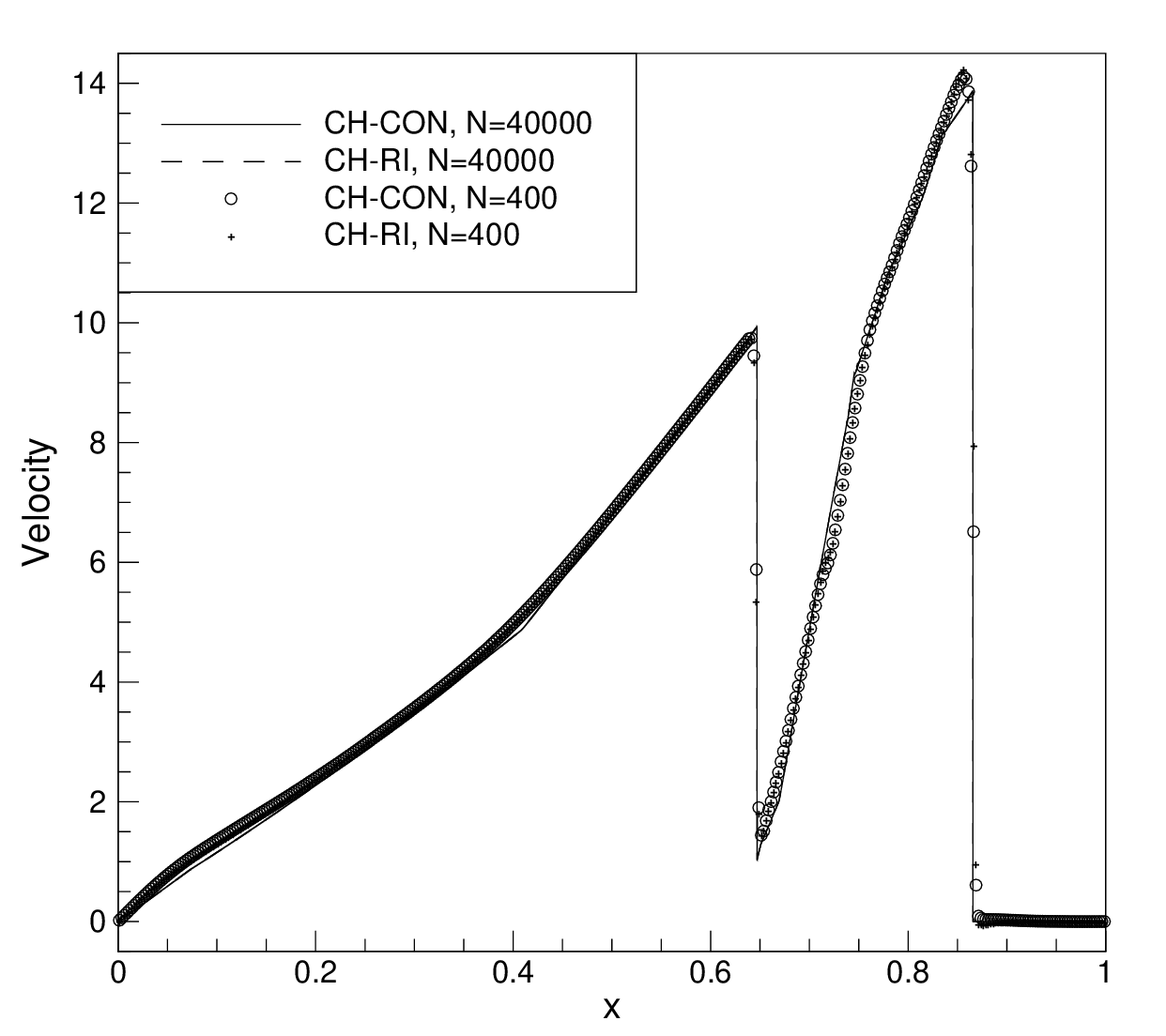}
        \caption{AWENO5}
    \end{subfigure}
    \begin{subfigure}{0.32\linewidth}
        \centering
        \includegraphics[width=\linewidth]{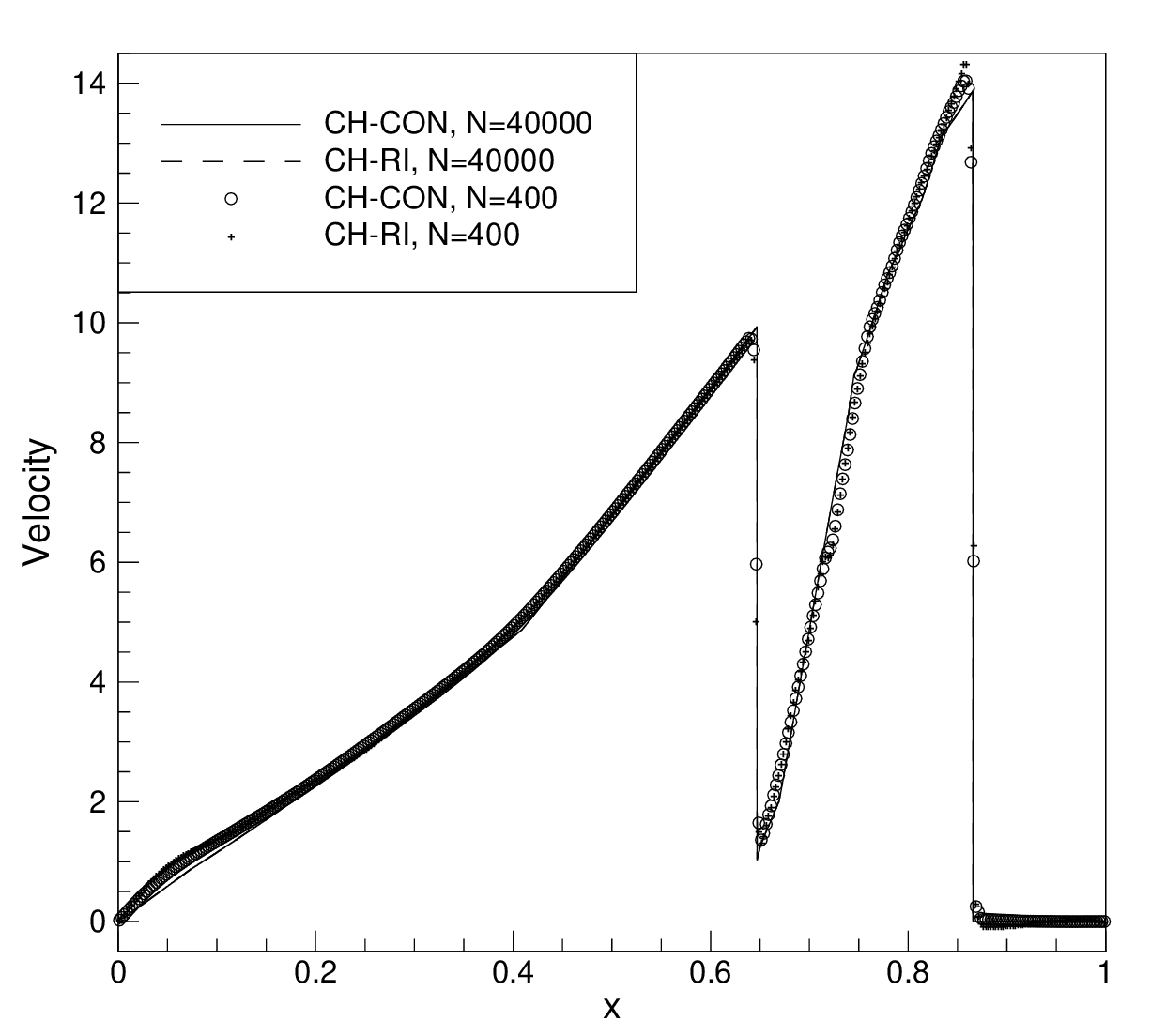}
        \caption{AWENO7}
    \end{subfigure}
    \begin{subfigure}{0.32\linewidth}
        \centering
        \includegraphics[width=\linewidth]{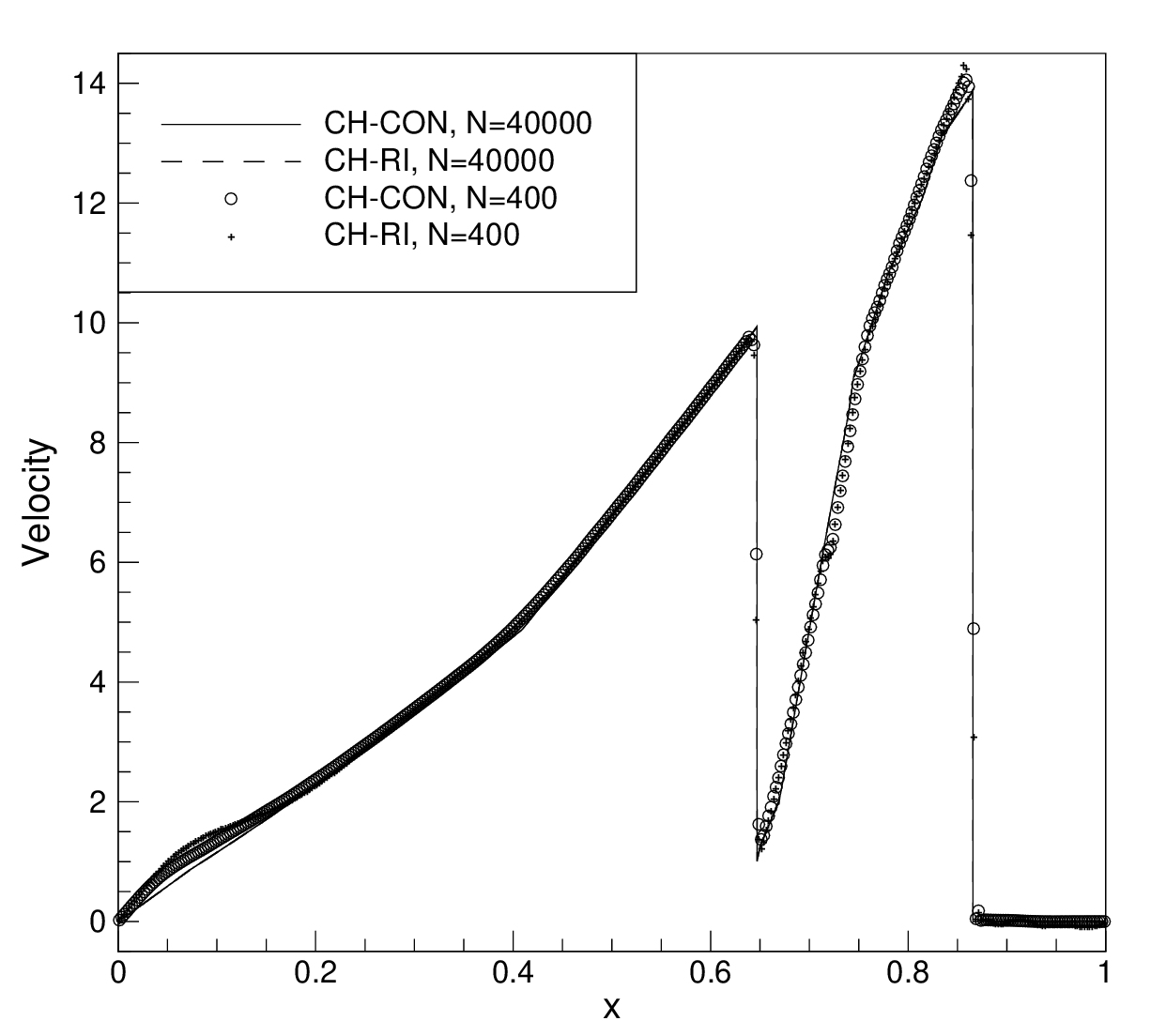}
        \caption{AWENO9}
    \end{subfigure}

    \begin{subfigure}{0.32\linewidth}
        \centering
        \includegraphics[width=\linewidth]{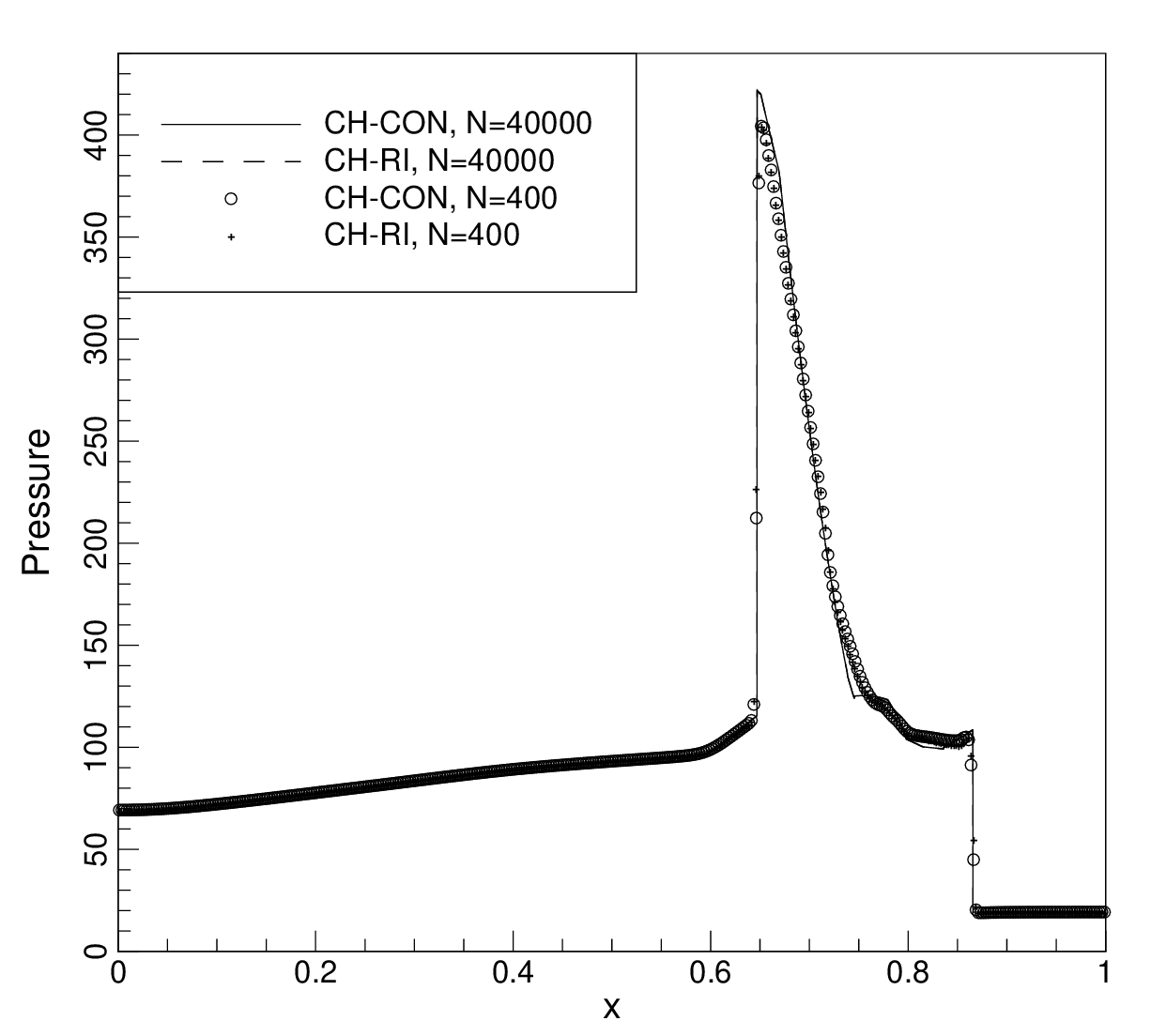}
        \caption{AWENO5}
    \end{subfigure}
    \begin{subfigure}{0.32\linewidth}
        \centering
        \includegraphics[width=\linewidth]{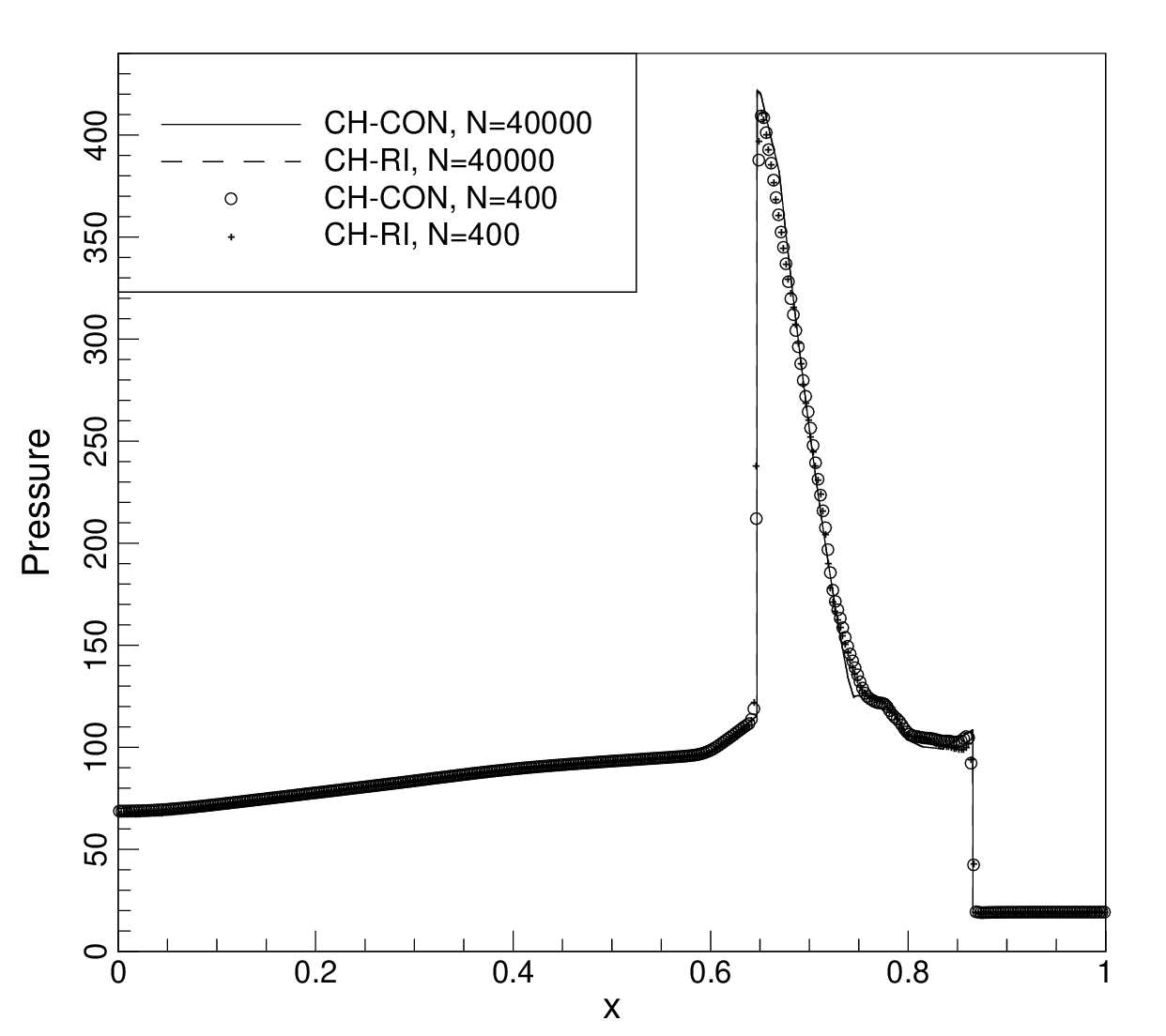}
        \caption{AWENO7}
    \end{subfigure}
    \begin{subfigure}{0.32\linewidth}
        \centering
        \includegraphics[width=\linewidth]{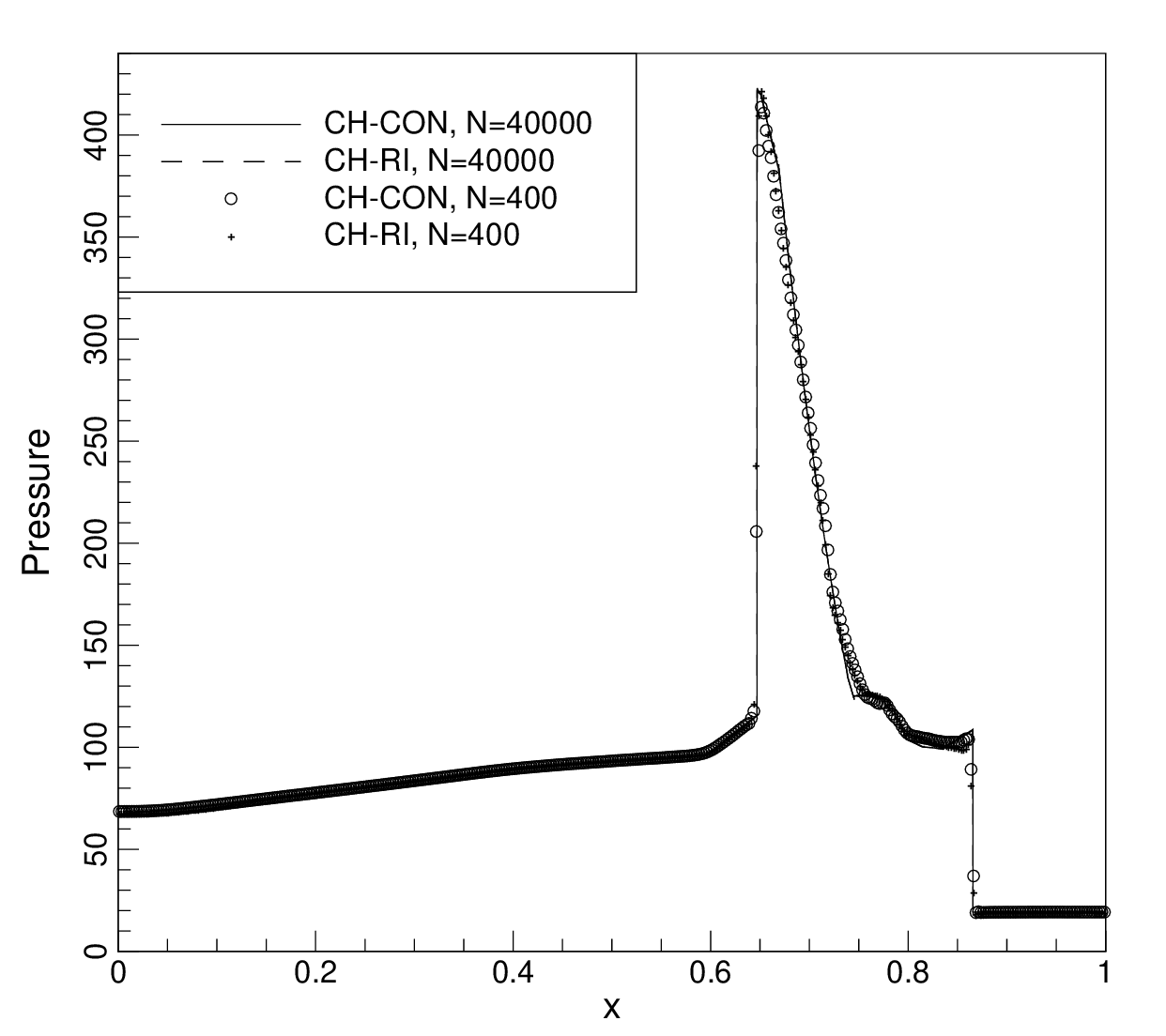}
        \caption{AWENO9}
    \end{subfigure}

    \caption{The blast wave interaction problem at \(T = 0.038\).}
    \label{fig:blast-wave}
\end{figure}

The initial condition is \(\rho = 1\) and \(u = 0\) with \(p = 1000\) for \(x \in (0, 0.1)\), \(p = 0.01\) for \(x \in (0.1, 0.9)\) and \(p = 100\) for \(x \in (0.9, 1)\). The computation domain is \([0, 1]\), and two reflective walls are set at both ends. We solve the problem to \(T = 0.038\). This model problem describes the interaction of two blast waves initialized near the boundaries. The results are shown in \Cref{fig:blast-wave}. Because there is no exact solutions for this problem, we use the solution with \(N = 40000\) grid points as the reference solution. Both versions perform well. The CH-RI seems less dissipative at the density peak than the CH-CON. 

\subsubsection{The Sedov blast wave problem}

\begin{figure}[htbp]
    \centering

    \begin{subfigure}{0.32\linewidth}
        \centering
        \includegraphics[width=\linewidth]{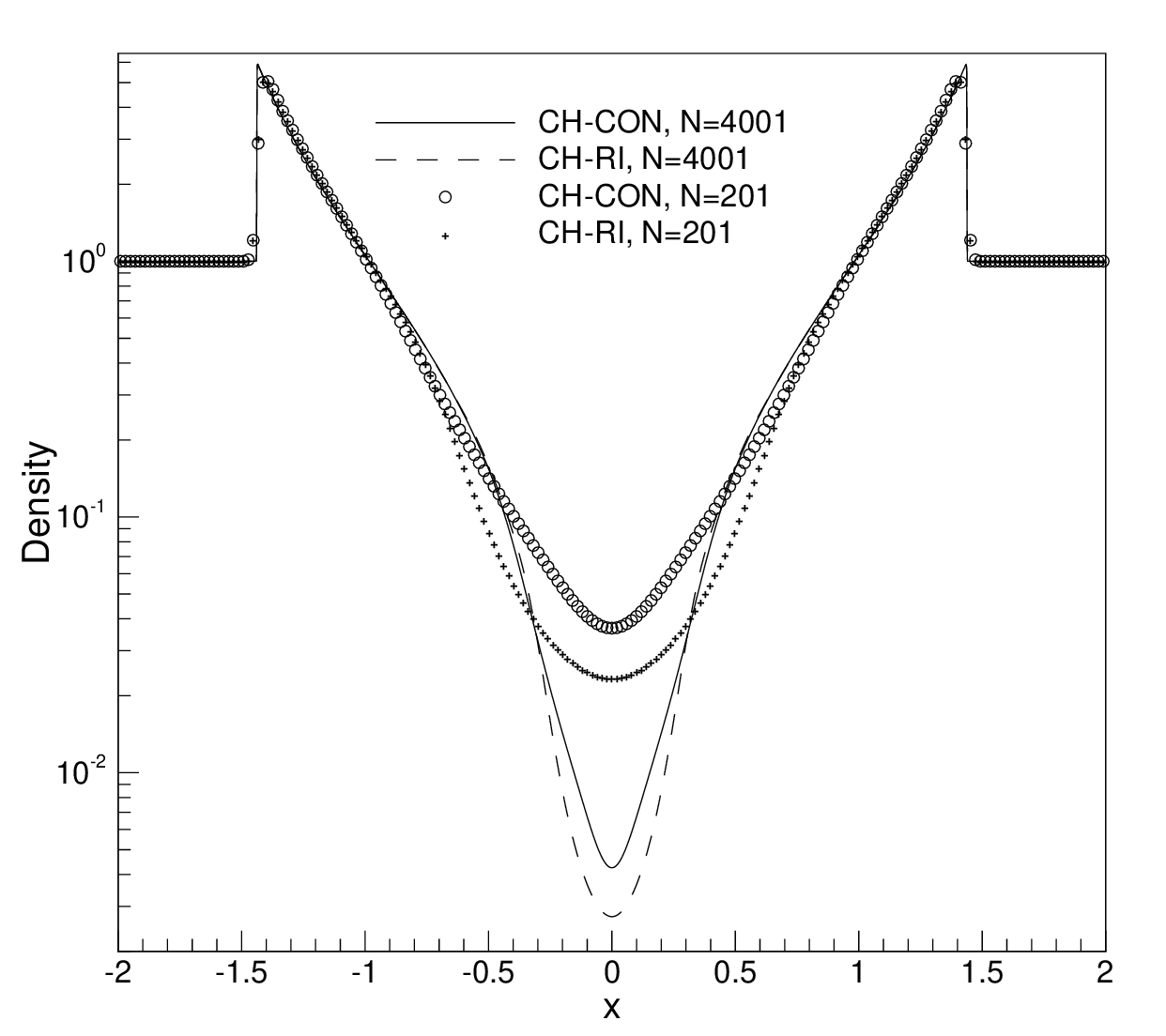}
        \caption{AWENO5}
    \end{subfigure}
    \begin{subfigure}{0.32\linewidth}
        \centering
        \includegraphics[width=\linewidth]{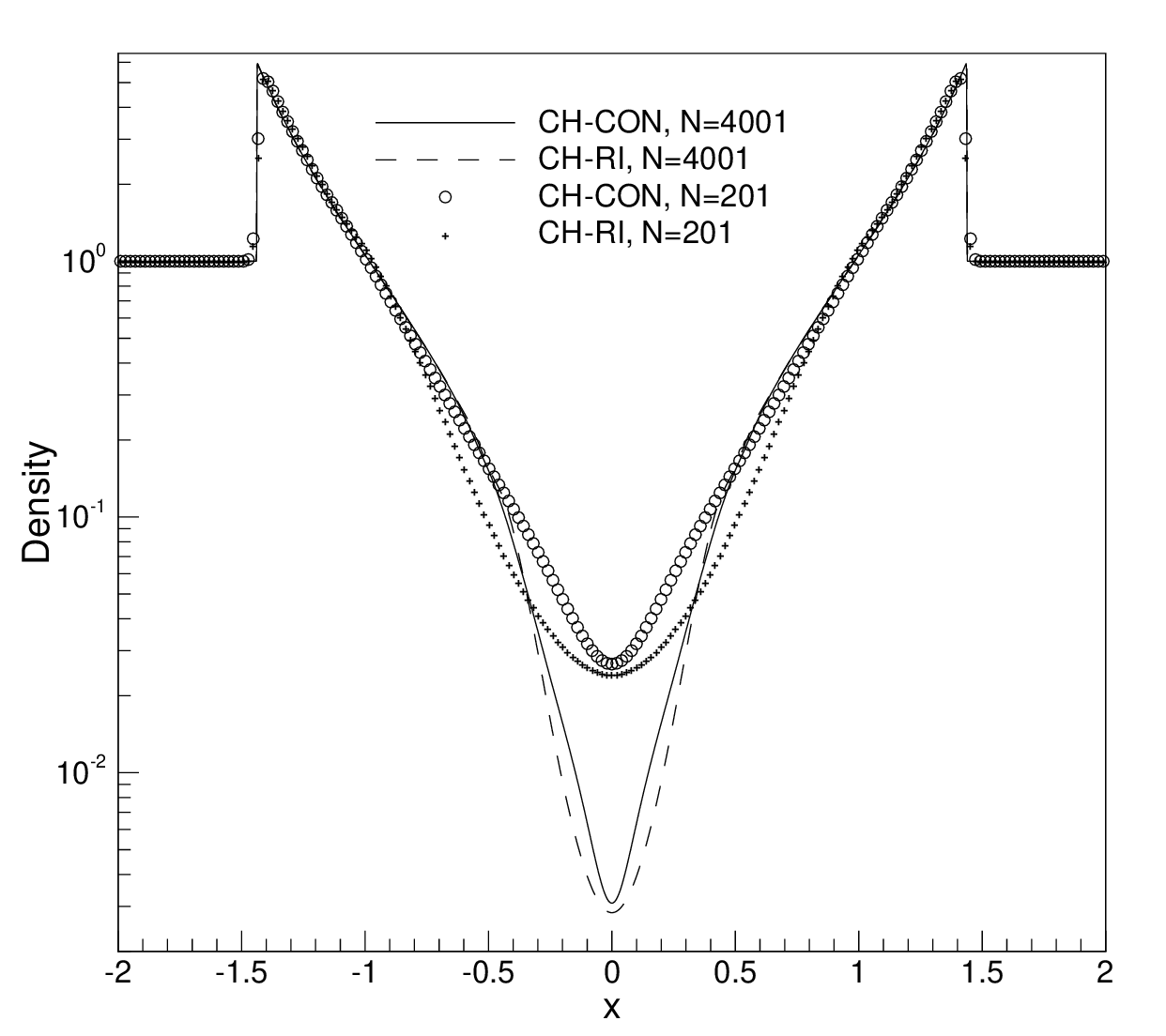}
        \caption{AWENO7}
    \end{subfigure}
    \begin{subfigure}{0.32\linewidth}
        \centering
        \includegraphics[width=\linewidth]{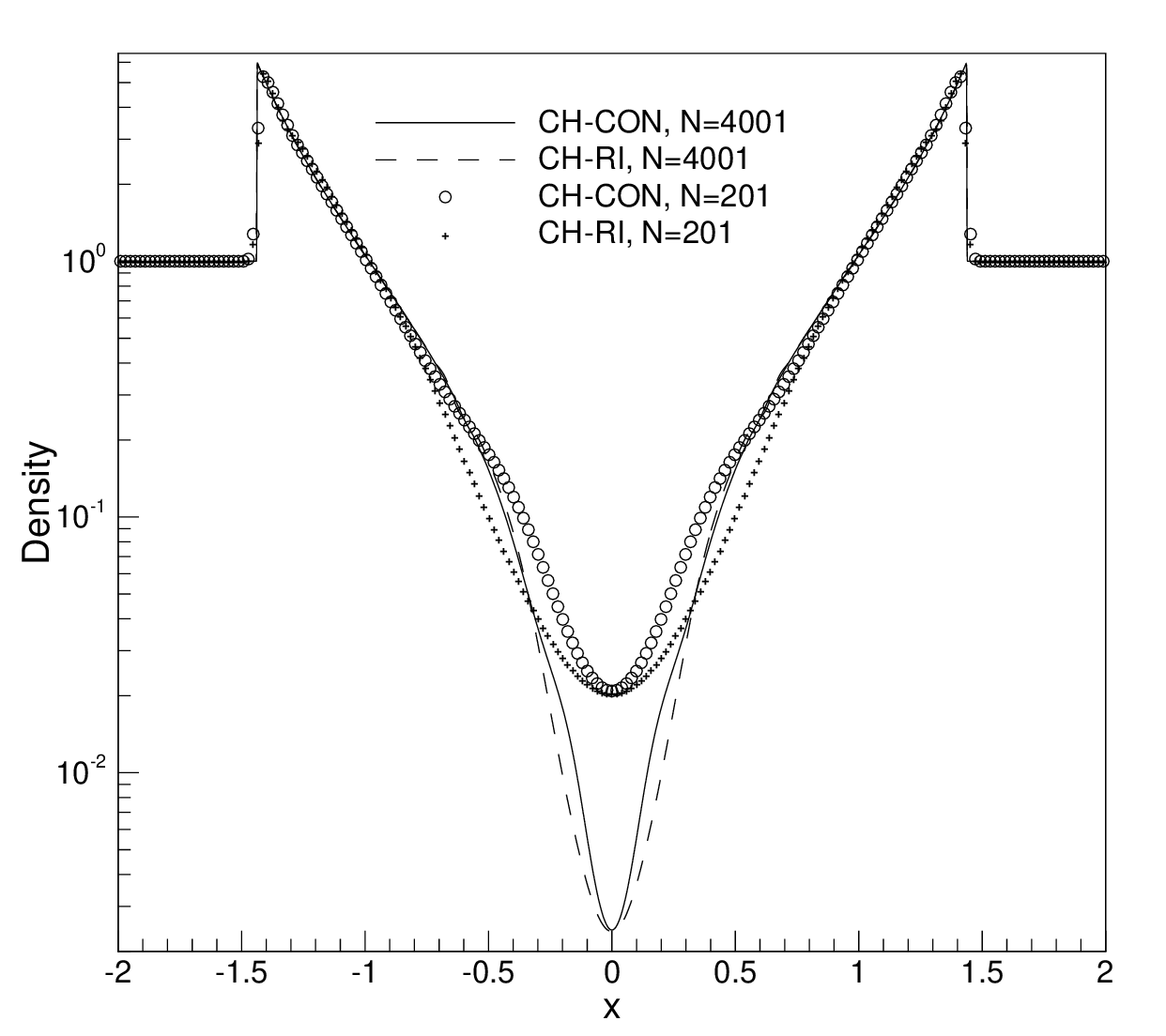}
        \caption{AWENO9}
    \end{subfigure}

    \begin{subfigure}{0.32\linewidth}
        \centering
        \includegraphics[width=\linewidth]{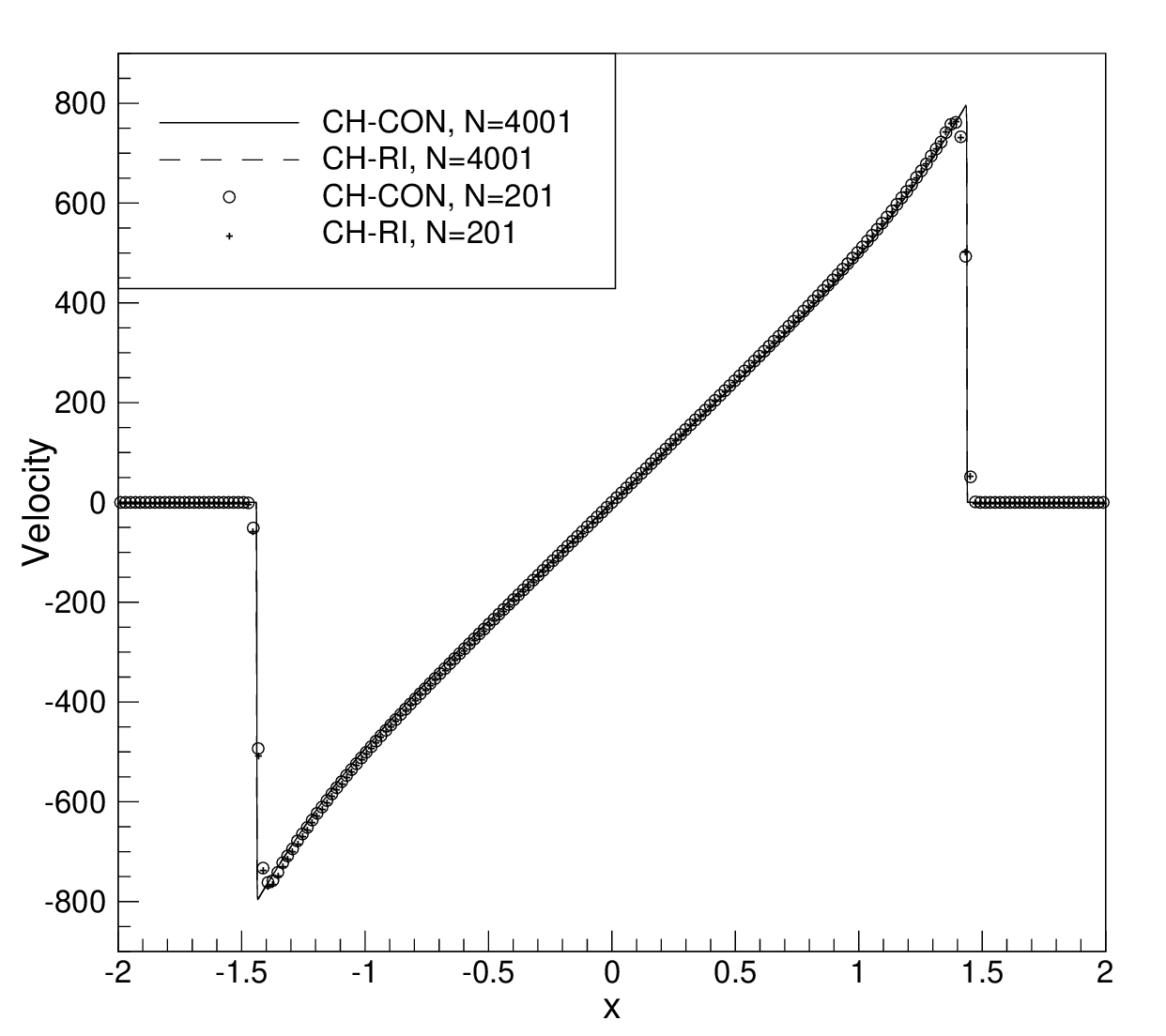}
        \caption{AWENO5}
    \end{subfigure}
    \begin{subfigure}{0.32\linewidth}
        \centering
        \includegraphics[width=\linewidth]{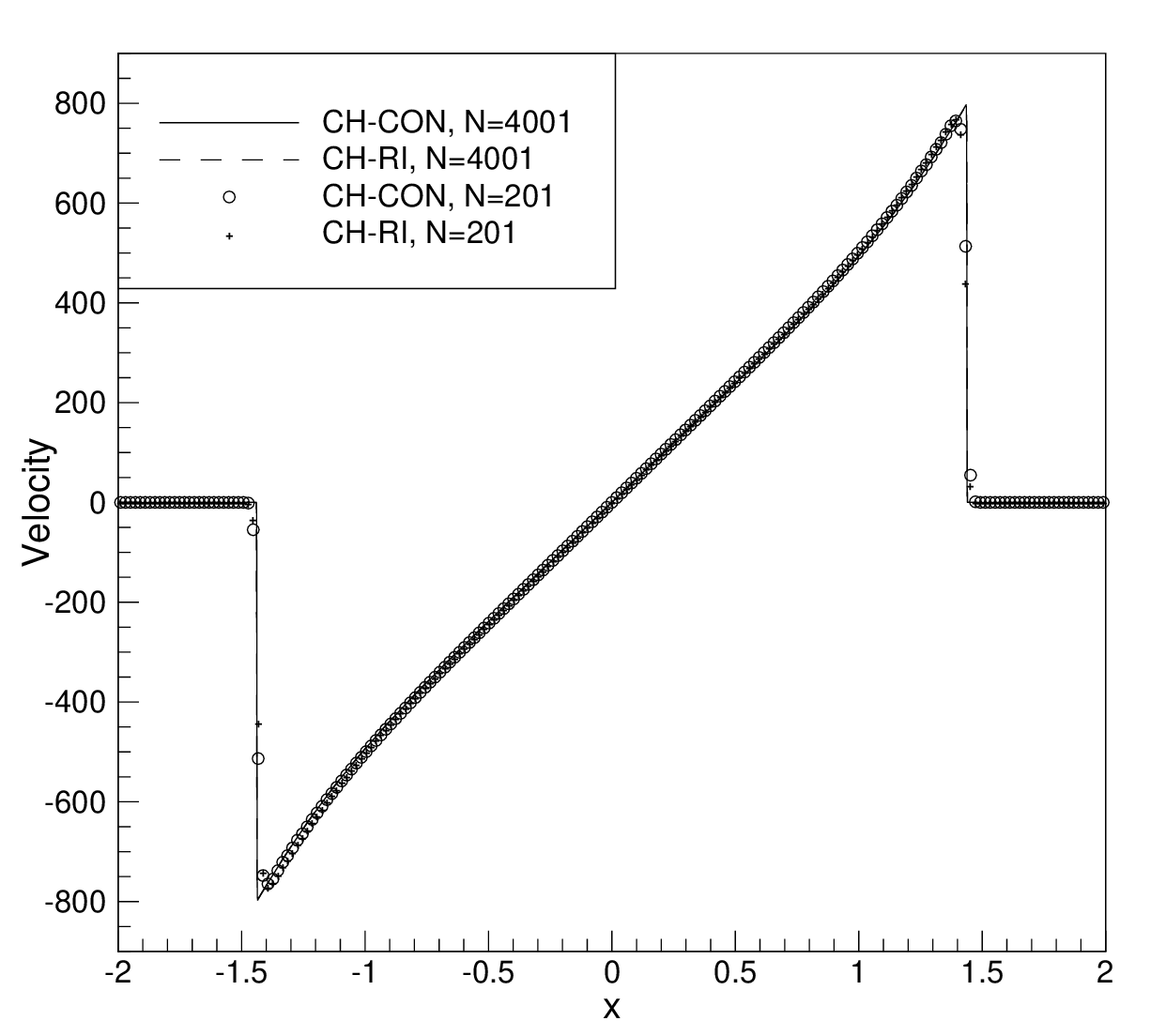}
        \caption{AWENO7}
    \end{subfigure}
    \begin{subfigure}{0.32\linewidth}
        \centering
        \includegraphics[width=\linewidth]{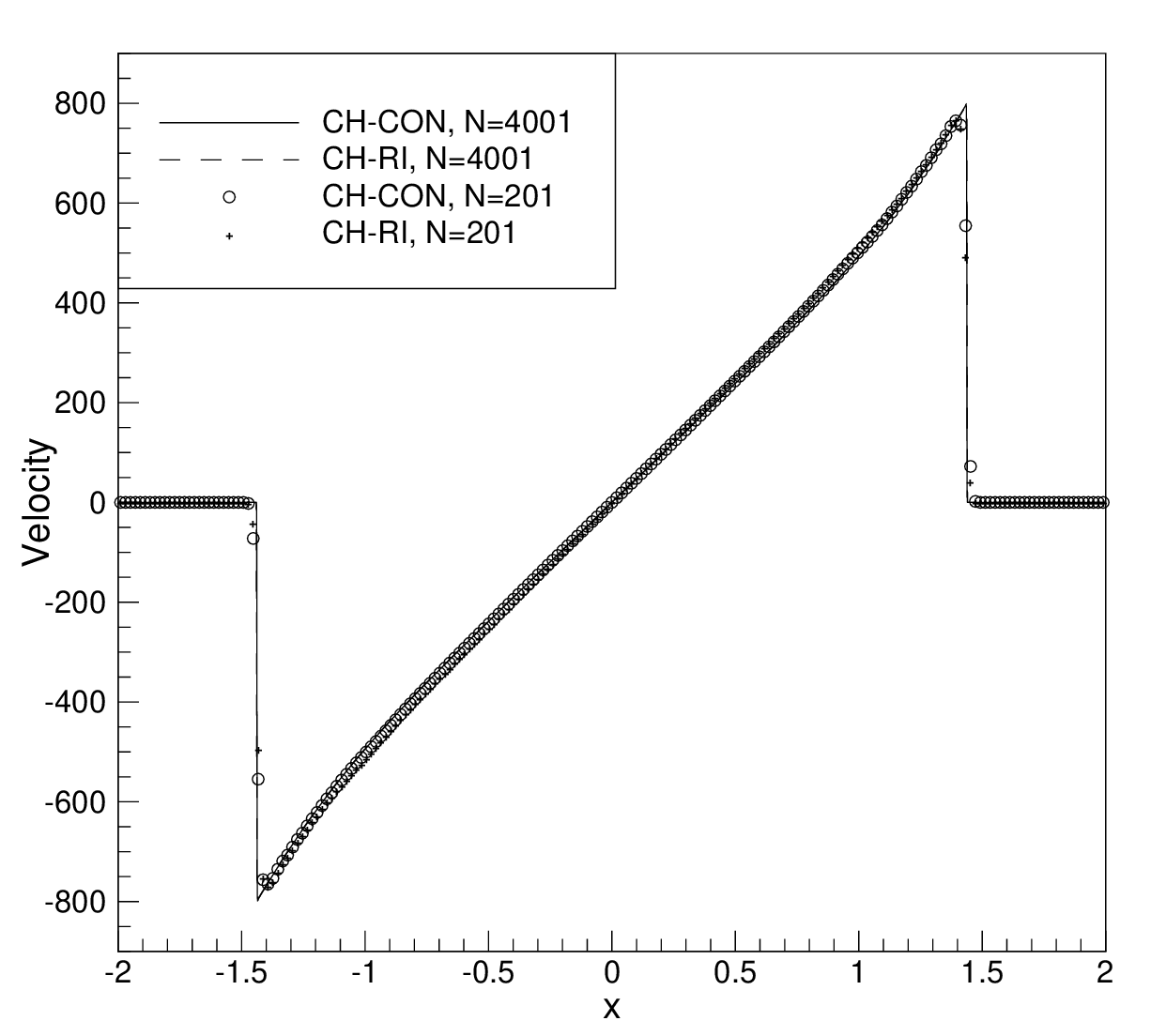}
        \caption{AWENO9}
    \end{subfigure}

    \begin{subfigure}{0.32\linewidth}
        \centering
        \includegraphics[width=\linewidth]{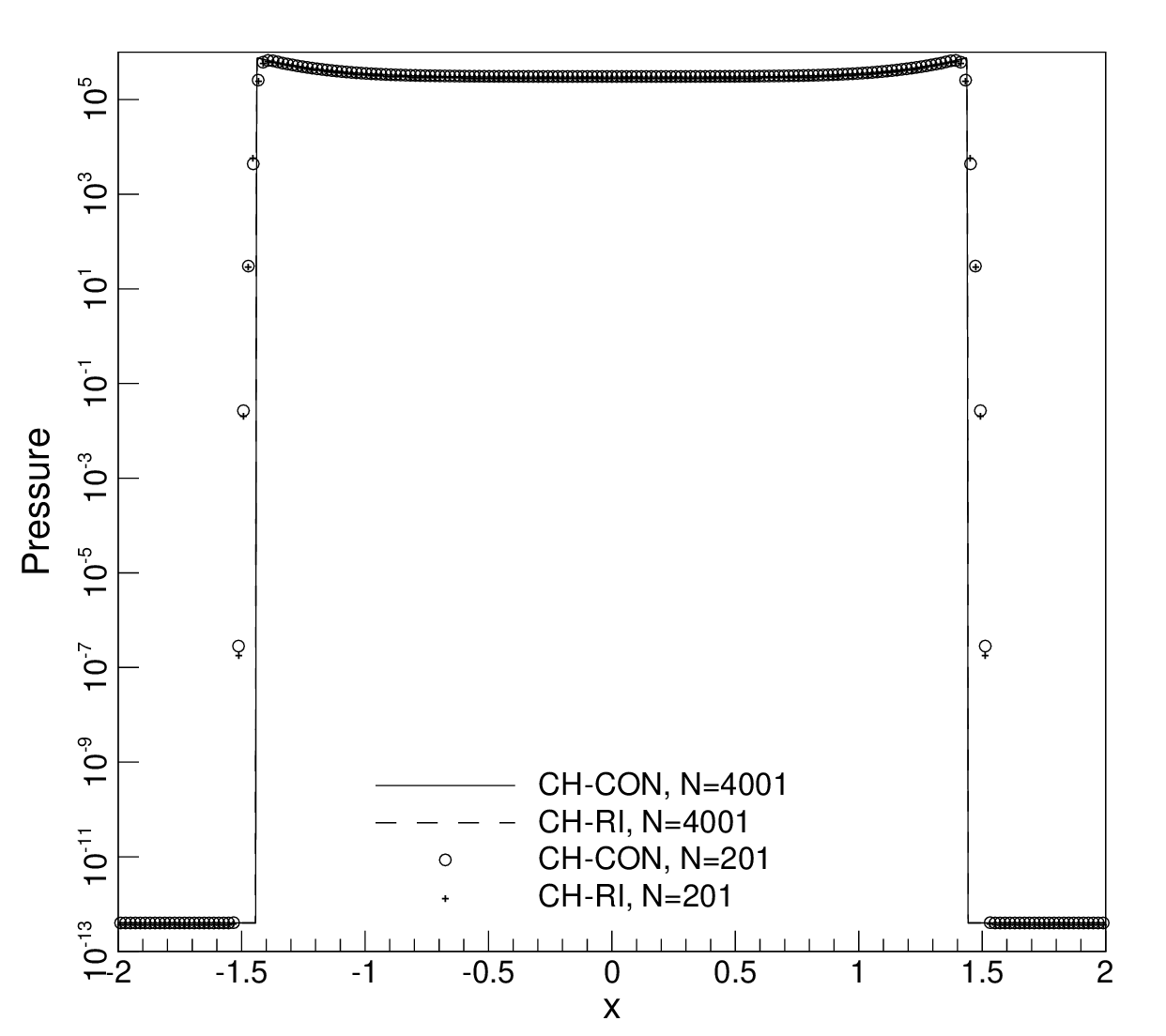}
        \caption{AWENO5}
    \end{subfigure}
    \begin{subfigure}{0.32\linewidth}
        \centering
        \includegraphics[width=\linewidth]{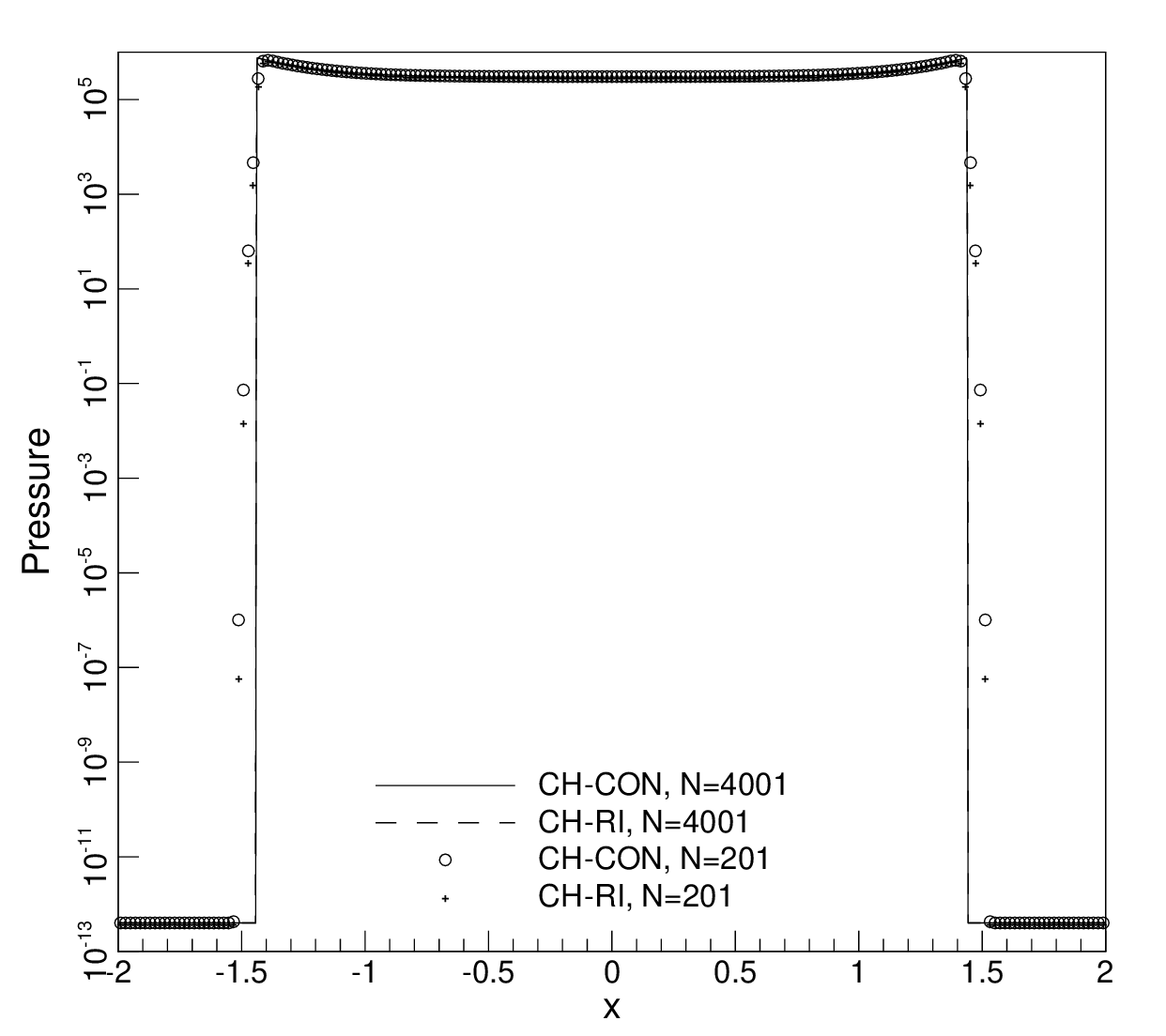}
        \caption{AWENO7}
    \end{subfigure}
    \begin{subfigure}{0.32\linewidth}
        \centering
        \includegraphics[width=\linewidth]{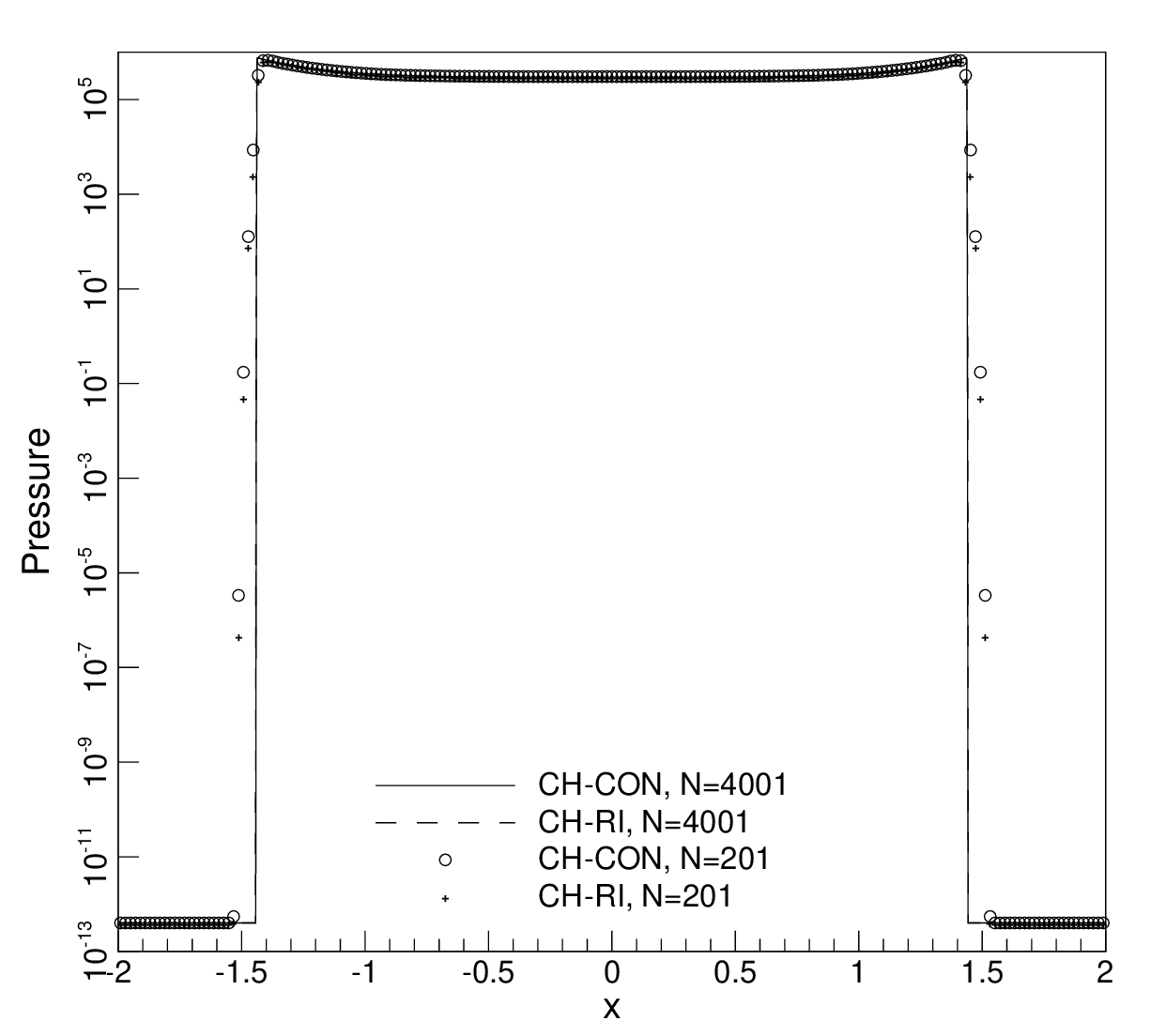}
        \caption{AWENO9}
    \end{subfigure}

    \caption{The Sedov blast wave problem at \(T = 0.001\).}
    \label{fig:Sedov}
\end{figure}

For the initial condition, the background state is \((\rho, u, p) = (1, 0, 10^{-12})\), and a singular delta total energy source of strength \(E_0\) is added at the origin. The computation domain is \([-2, 2]\). We solve the problem to \(T = 0.001\). This problem is a simple model of point explosion. It is an extreme test for scheme stability and has been used to test PP limiters of DG schemes \cite{10.1016/j.jcp.2010.08.016}. We choose \(E_0 = 3.2 \times 10^{6}\). For initialization, we use odd number of cells and set the additional total energy at the central node to be \(\frac{E_0}{h}\). The results are shown in \Cref{fig:Sedov}. The \(N = 4001\) solution is used for reference. The CH-RI seems less dissipative near the central point in the sense that it produces a smaller minimum density. Besides this region, both methods agree well. 

\subsubsection{The shock-density wave interaction problem}

\begin{figure}[htbp]
    \centering

    \begin{subfigure}{0.32\linewidth}
        \centering
        \includegraphics[width=\linewidth]{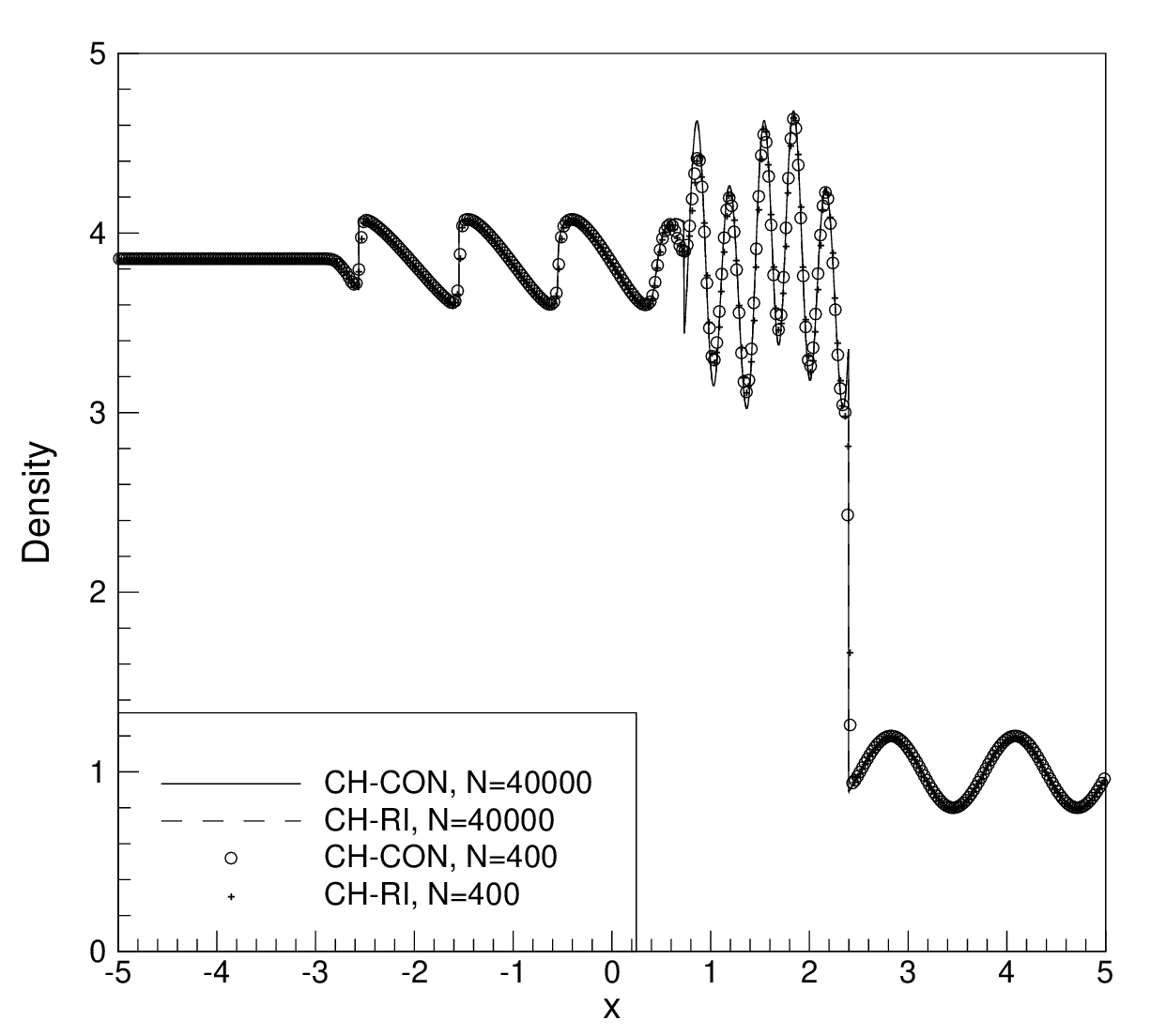}
        \caption{AWENO5}
    \end{subfigure}
    \begin{subfigure}{0.32\linewidth}
        \centering
        \includegraphics[width=\linewidth]{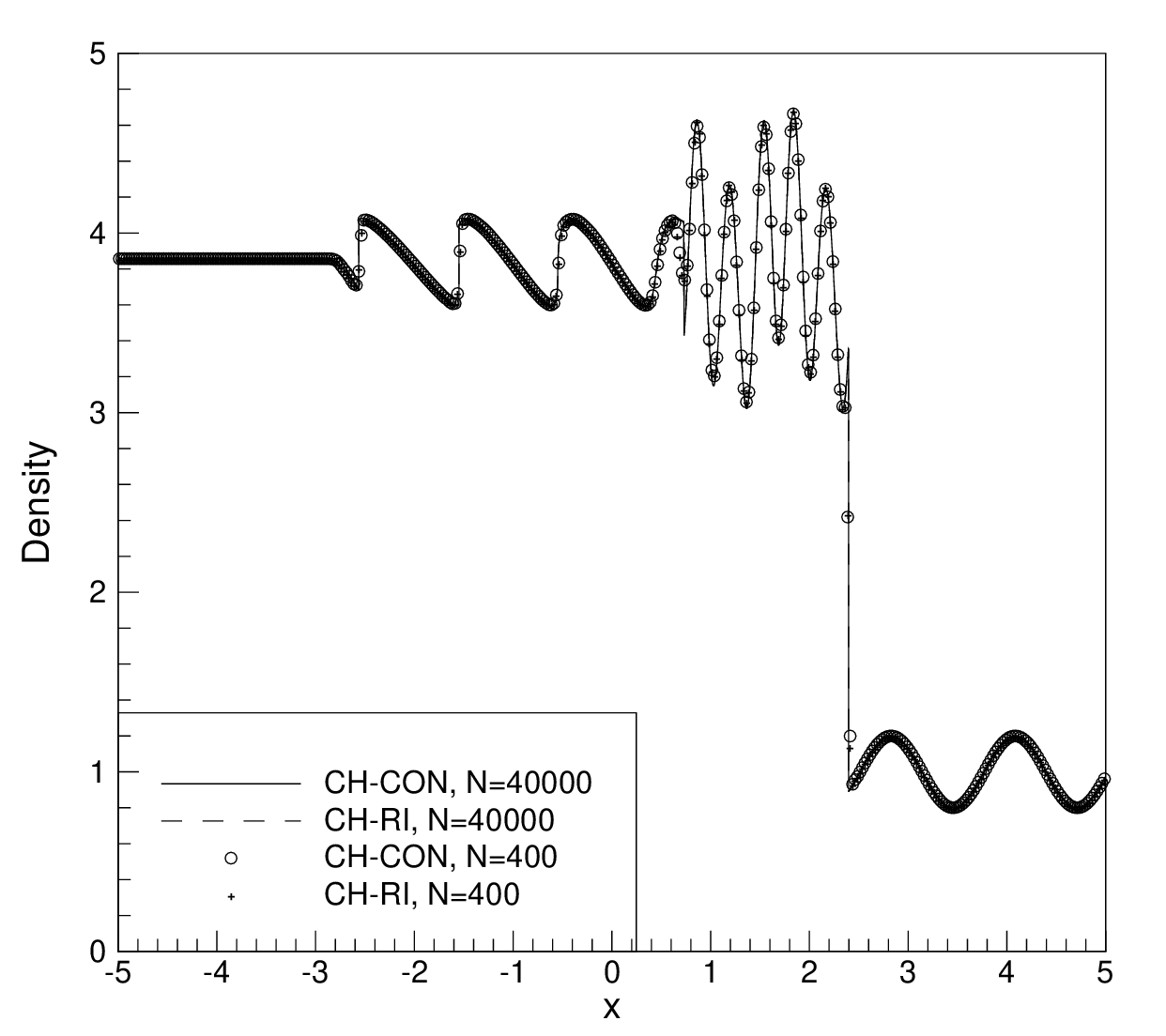}
        \caption{AWENO7}
    \end{subfigure}
    \begin{subfigure}{0.32\linewidth}
        \centering
        \includegraphics[width=\linewidth]{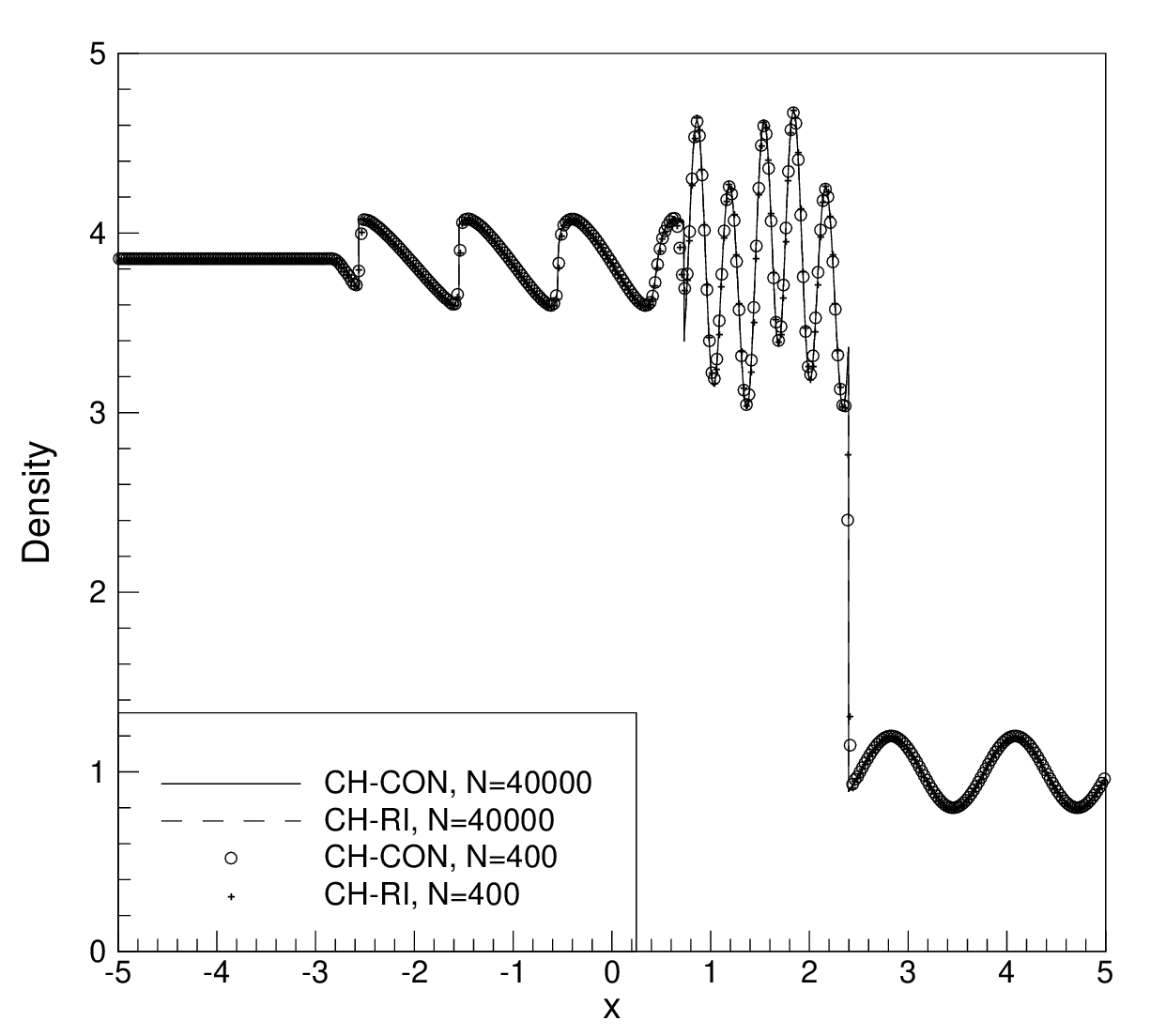}
        \caption{AWENO9}
    \end{subfigure}

    \begin{subfigure}{0.32\linewidth}
        \centering
        \includegraphics[width=\linewidth]{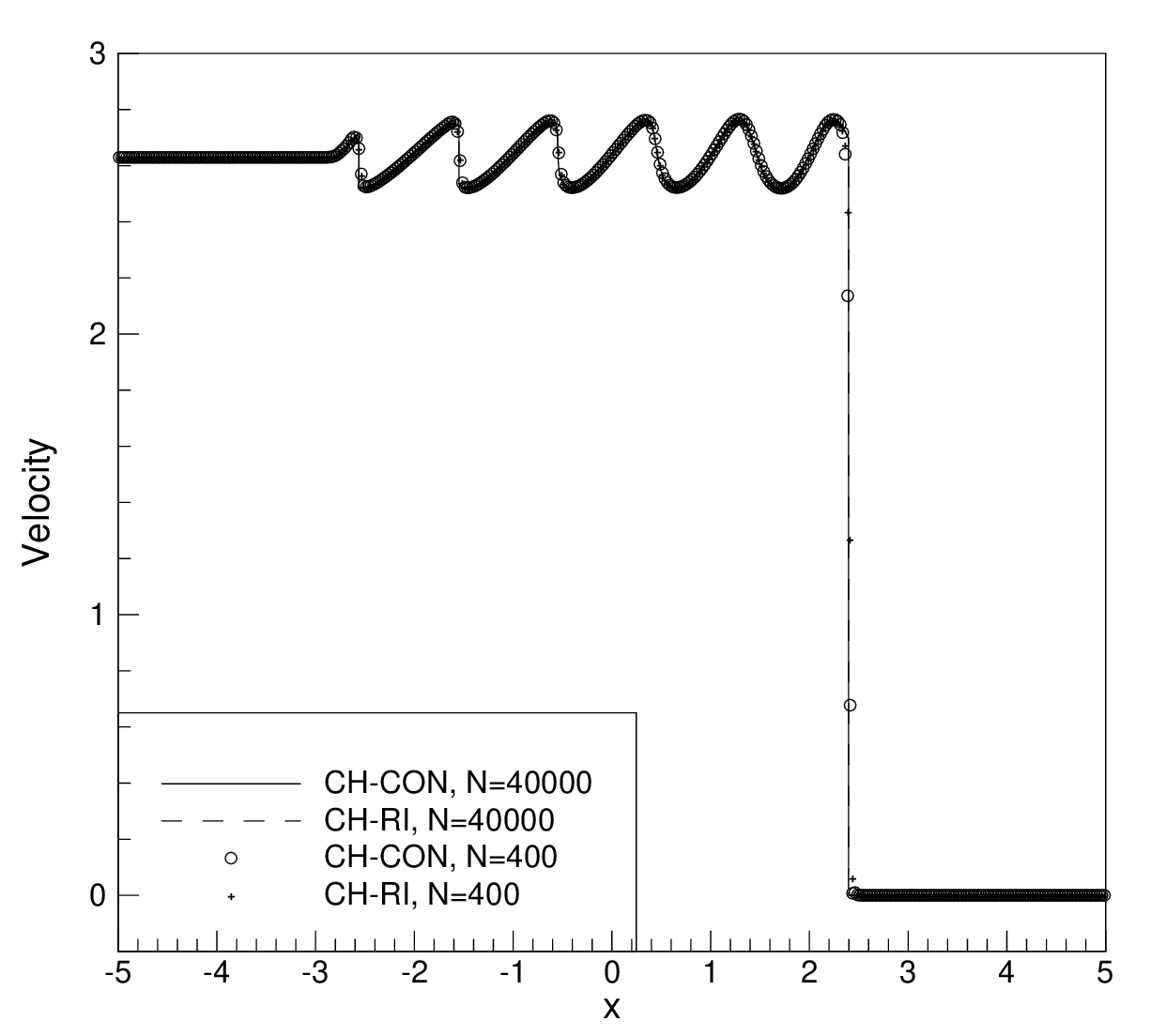}
        \caption{AWENO5}
    \end{subfigure}
    \begin{subfigure}{0.32\linewidth}
        \centering
        \includegraphics[width=\linewidth]{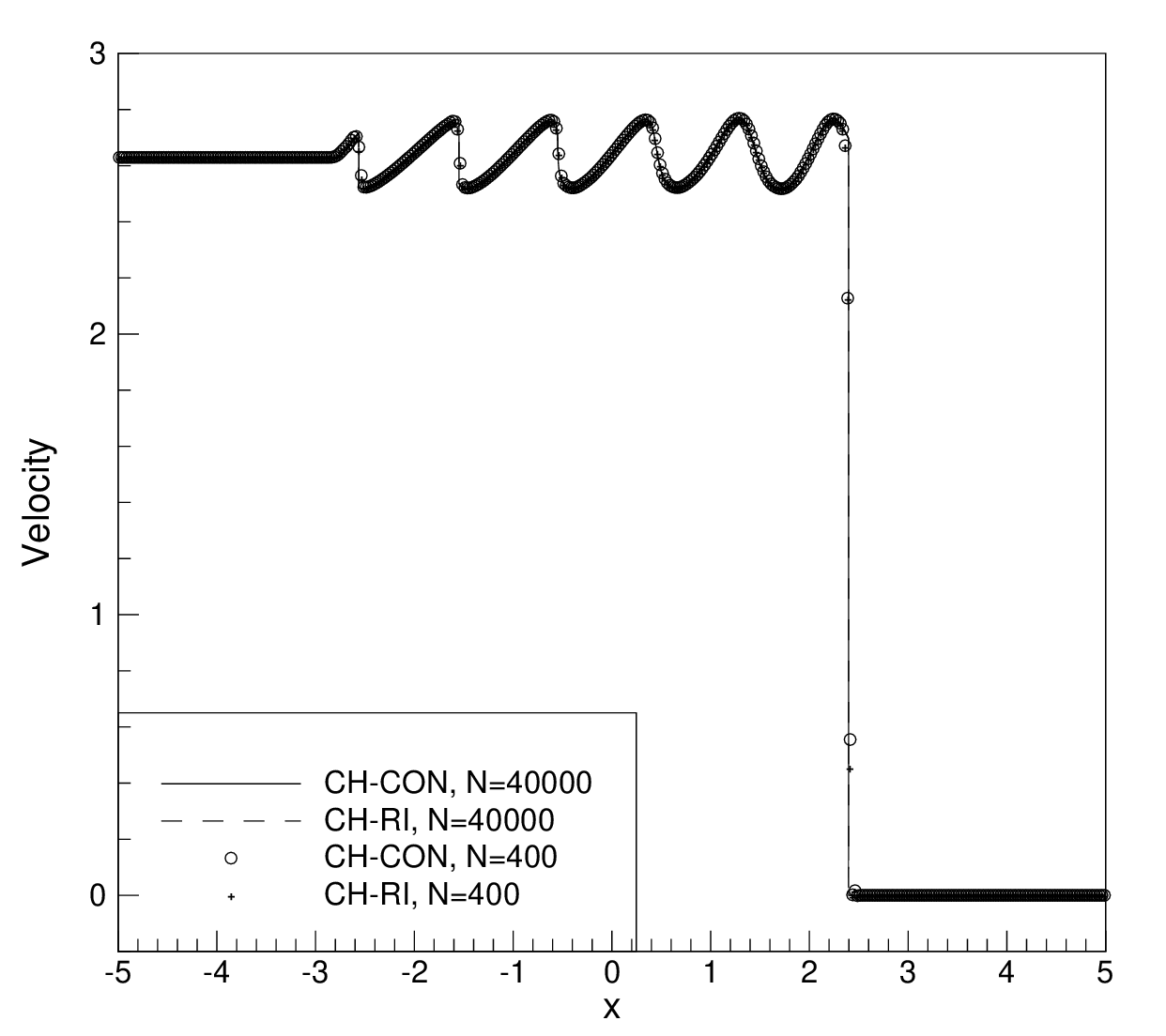}
        \caption{AWENO7}
    \end{subfigure}
    \begin{subfigure}{0.32\linewidth}
        \centering
        \includegraphics[width=\linewidth]{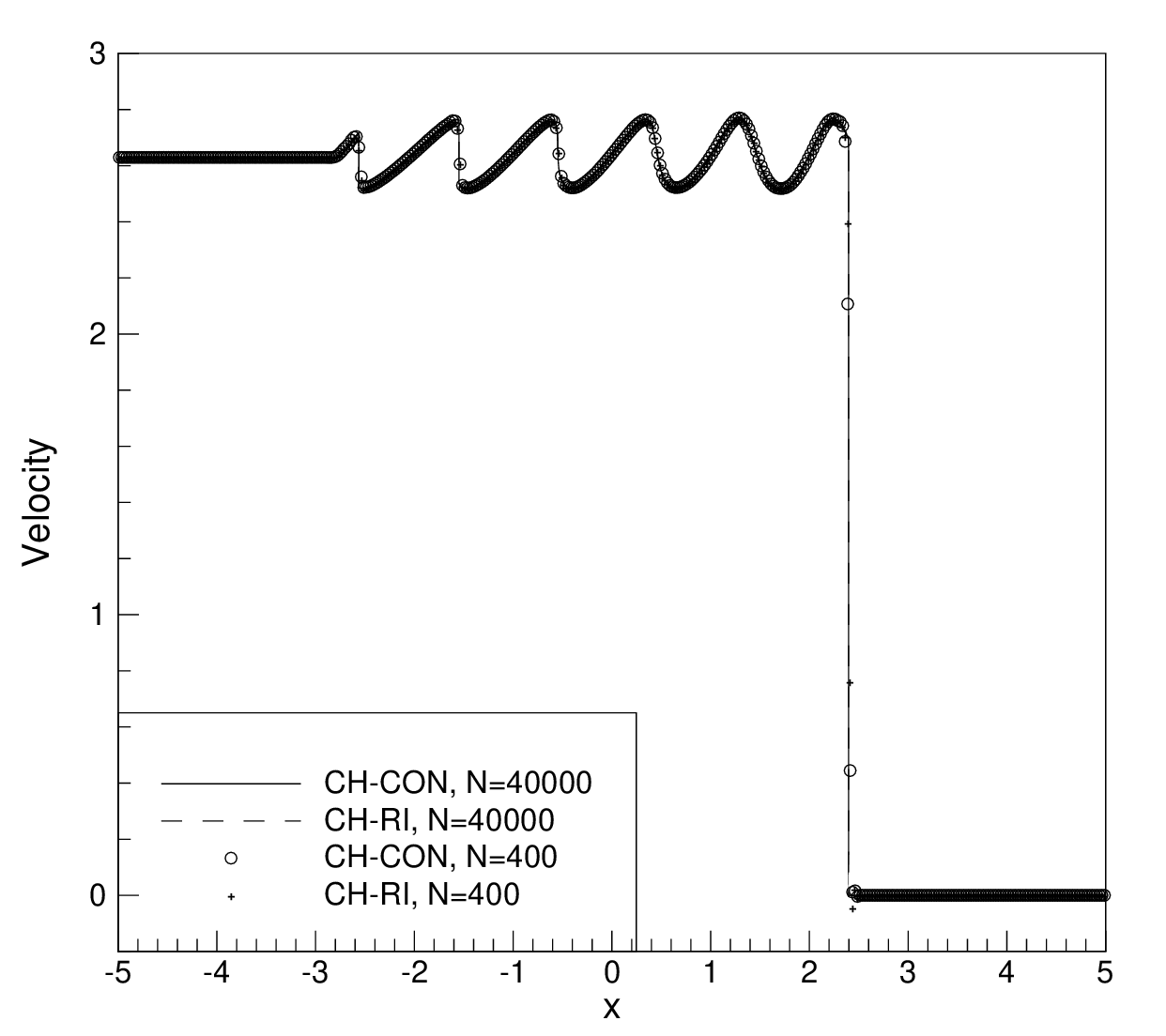}
        \caption{AWENO9}
    \end{subfigure}

    \begin{subfigure}{0.32\linewidth}
        \centering
        \includegraphics[width=\linewidth]{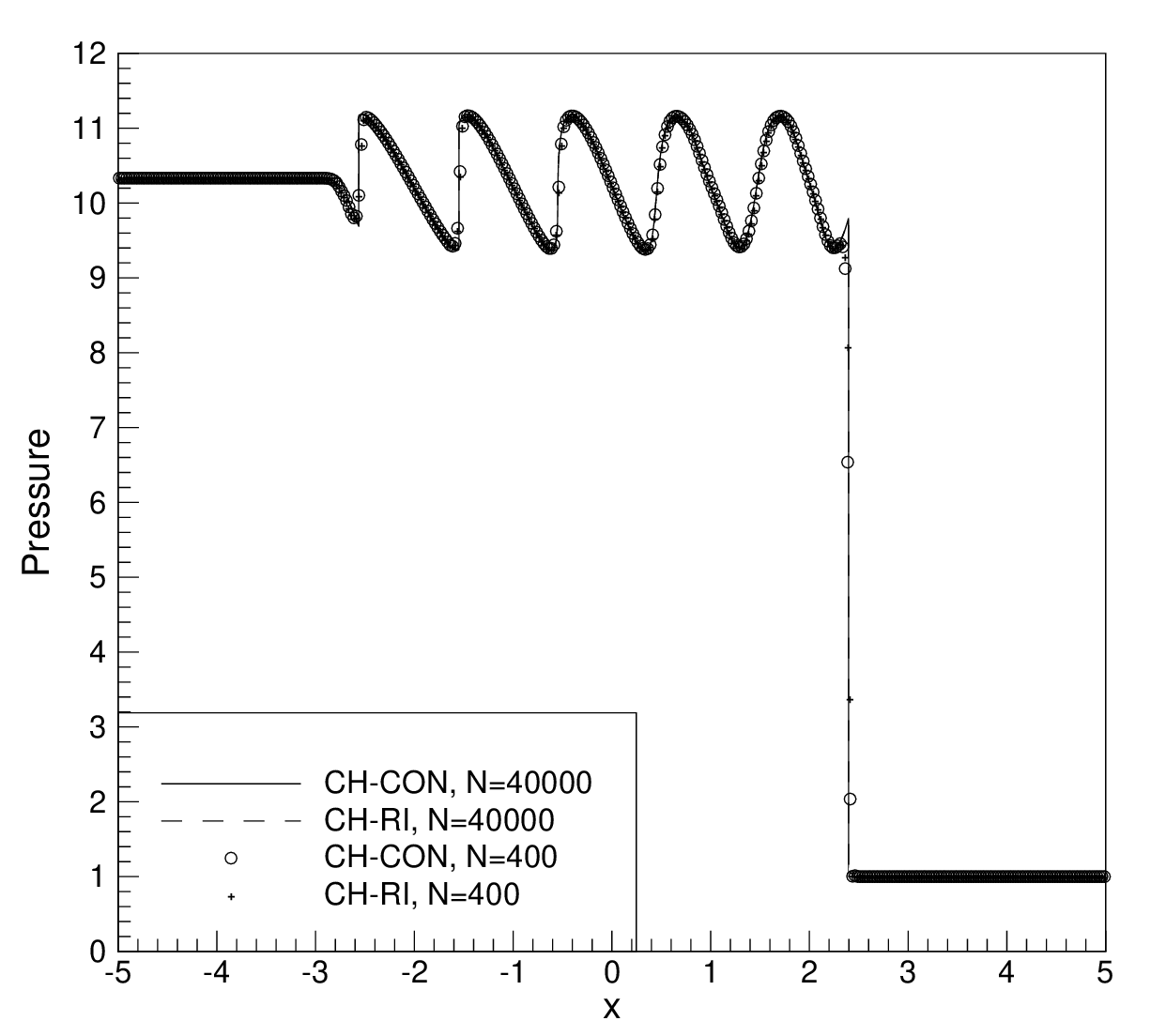}
        \caption{AWENO5}
    \end{subfigure}
    \begin{subfigure}{0.32\linewidth}
        \centering
        \includegraphics[width=\linewidth]{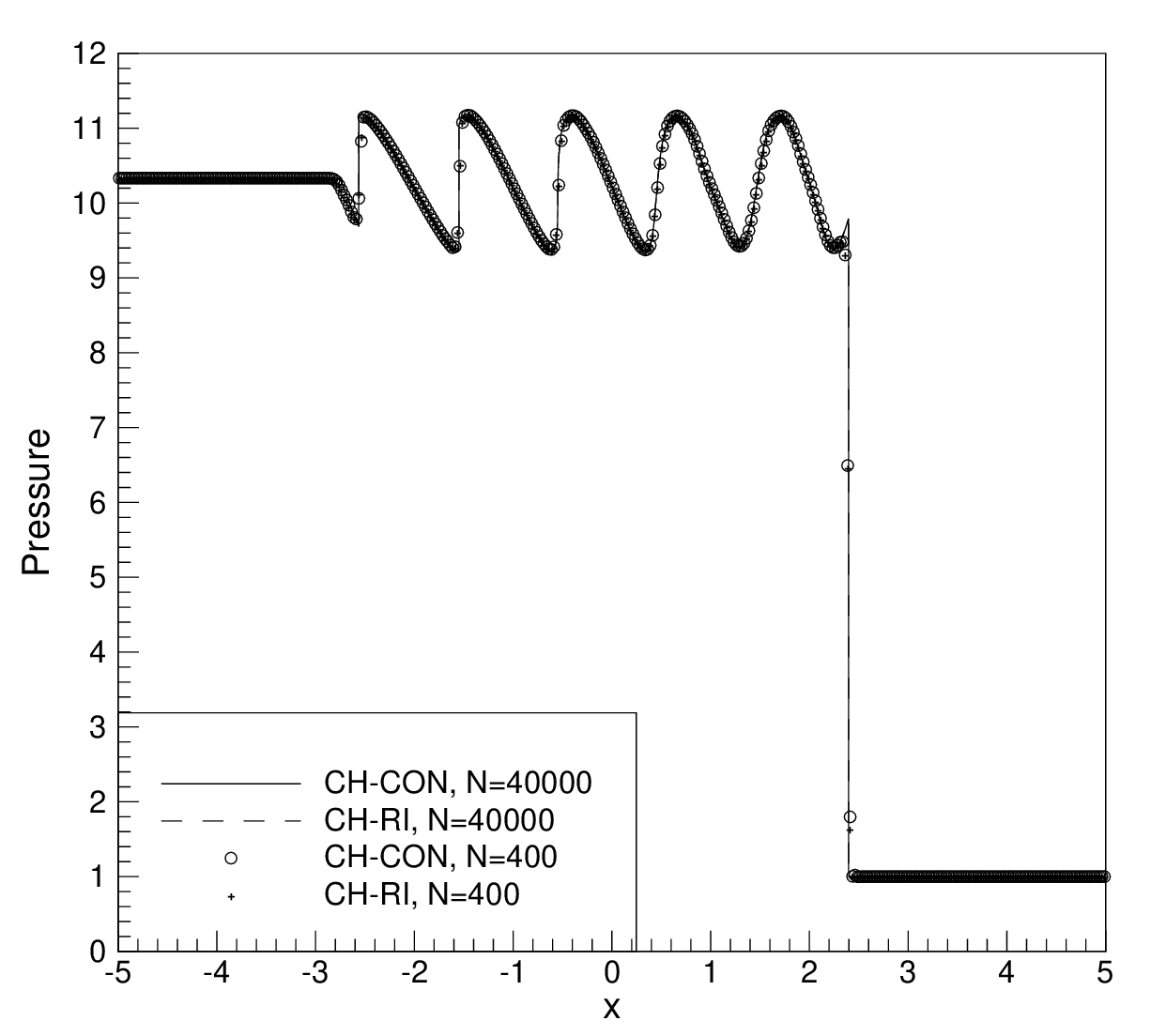}
        \caption{AWENO7}
    \end{subfigure}
    \begin{subfigure}{0.32\linewidth}
        \centering
        \includegraphics[width=\linewidth]{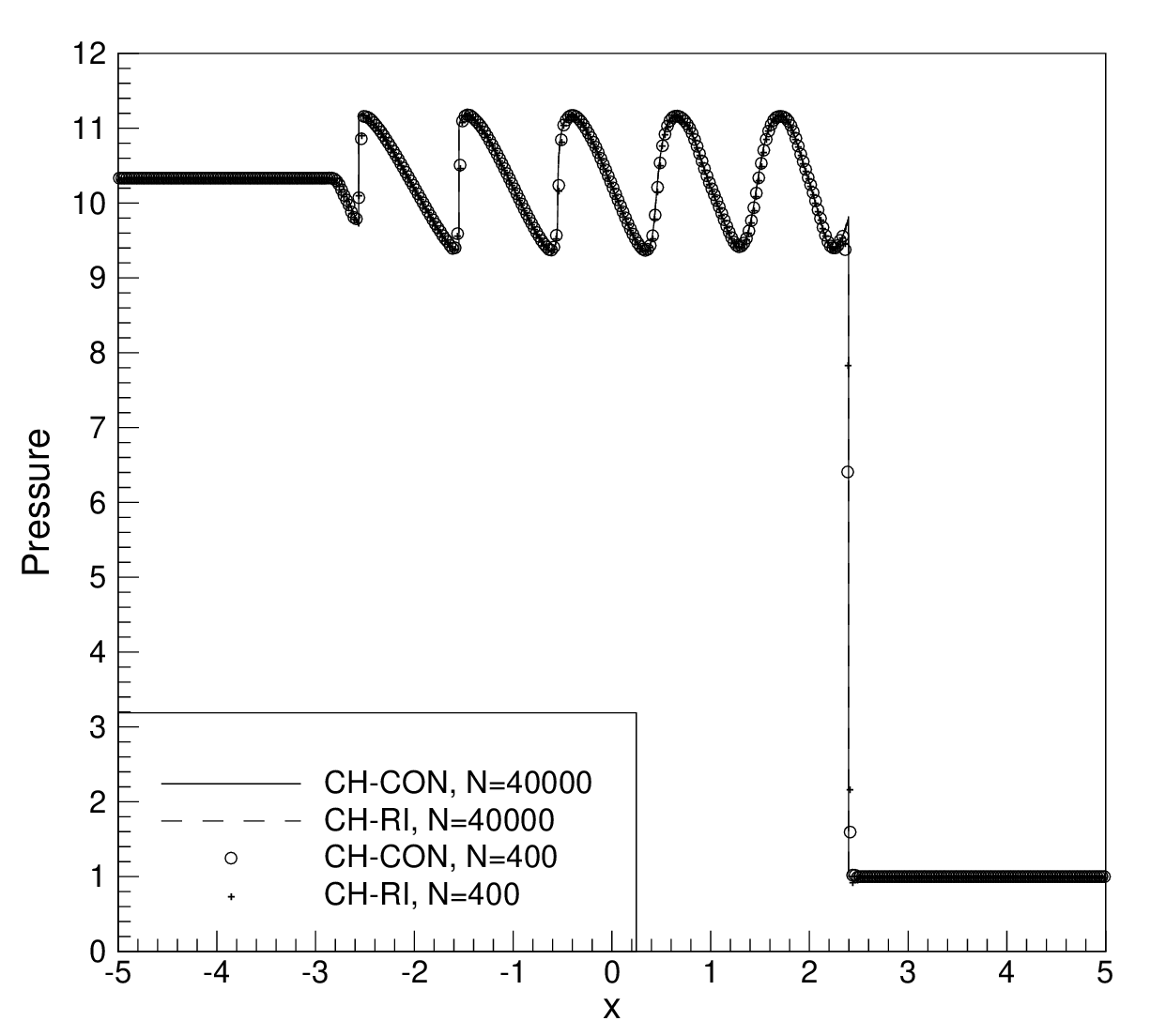}
        \caption{AWENO9}
    \end{subfigure}

    \caption{The shock-density wave interaction problem at \(T = 1.8\).}
    \label{fig:shock-density}
\end{figure}

The initial condition is \((\rho, u, p) = (\frac{27}{7}, \frac{4}{9}\sqrt{35}, \frac{31}{3})\) for \(x < -4\) and \((\rho, u, p) = (1 + \varepsilon\sin(5x), 0, 1)\) for \(x > -4\). The computation domain is \([-5, 5]\). We solve the problem to \(T = 1.8\). It describes the interaction between a right-moving Mach 3 shock and a stationary density wave \cite{10.1016/0021-99918990222-2}. We choose the parameter \(\epsilon = 0.2\). The results are shown in \Cref{fig:shock-density}. We use \(N=40000\) as references. The CH-RI and CH-CON versions agree well with each other and capture the high frequency part after the shock well. 

\subsection{2D tests}

\subsubsection{Accuracy test}

We consider the rotated simple density-transport problem, whose exact solution is \((\rho, u, v, p) = (1 + 0.2\sin(\pi(x+y)), 1, 1, 1)\). The computational domain is \([0, 2] \times [0, 2]\) with periodic boundaries. We solve up to \(T = 2\). We use \((N_{x,0}, N_{y,0}) = (20, 20)\) as the coarsest mesh and use \(h_{0} = \frac{2}{N_0}\) to denote the coarsest mesh size. Using SSP-RK(3,3), the time steps in the following tests with uniform mesh are taken to be \(\Delta t = \mathrm{CFL} \frac{\mathrm{CFL}}{\frac{\alpha}{\Delta x} + \frac{\beta}{\Delta y}} \left(\frac{h}{h_0}\right)^{\frac{k}{3} - 1}\), which makes the time discretization error \(k\)th-order while maintaining stability. Here, \(\alpha\) and \(\beta\) are the maximum wave speeds in the \(x\)- and \(y\)-direction respectively. The results are shown in \Cref{tab:2D-accuracy-test-1-k=5,tab:2D-accuracy-test-1-k=7,tab:2D-accuracy-test-1-k=9}. 

\begin{table}[htbp]
    \centering
    \caption{Accuracy test. Error table of the density. \(k = 5\).}
    \label{tab:2D-accuracy-test-1-k=5}
    \begin{tabular}{ccccccccc}
        \hline
        & \multicolumn{4}{c}{CH-RI} & \multicolumn{4}{c}{CH-CON} \\
        \(N\) & \(\ell_h^2\) error & order & \(\ell_h^\infty\) error & order & \(\ell_h^2\) error & order & \(\ell_h^\infty\) error & order \\ \hline
        20 & 2.08E-03 & - & 1.86E-03 & - & 1.20E-03 & - & 8.48E-04 & - \\ 
        40 & 9.55E-05 & 4.446  & 1.00E-04 & 4.216  & 3.88E-05 & 4.956  & 3.27E-05 & 4.696 \\ 
        60 & 1.37E-05 & 4.795  & 1.66E-05 & 4.442  & 5.03E-06 & 5.038  & 4.46E-06 & 4.917 \\ 
        80 & 3.31E-06 & 4.931  & 4.22E-06 & 4.750  & 1.18E-06 & 5.039  & 1.07E-06 & 4.955 \\ 
        100 & 1.08E-06 & 4.999  & 1.46E-06 & 4.754  & 3.84E-07 & 5.034  & 3.48E-07 & 5.038 \\ 
        120 & 4.32E-07 & 5.045  & 5.87E-07 & 5.005  & 1.53E-07 & 5.030  & 1.39E-07 & 5.043 \\ 
        140 & 1.97E-07 & 5.087  & 2.71E-07 & 5.013  & 7.07E-08 & 5.028  & 6.41E-08 & 5.008 \\ 
        160 & 9.95E-08 & 5.128  & 1.35E-07 & 5.235  & 3.61E-08 & 5.030  & 3.24E-08 & 5.123 \\ 
        \hline
    \end{tabular}
\end{table}

\begin{table}[htbp]
    \centering
    \caption{Accuracy test. Error table of the density. \(k = 7\).}
    \label{tab:2D-accuracy-test-1-k=7}
    \begin{tabular}{ccccccccc}
        \hline
        & \multicolumn{4}{c}{CH-RI} & \multicolumn{4}{c}{CH-CON} \\
        \(N\) & \(\ell_h^2\) error & order & \(\ell_h^\infty\) error & order & \(\ell_h^2\) error & order & \(\ell_h^\infty\) error & order \\ \hline
        20 & 3.01E-04 & - & 4.26E-04 & - & 7.07E-05 & - & 7.82E-05 & - \\ 
        40 & 5.29E-06 & 5.833  & 1.04E-05 & 5.361  & 1.21E-06 & 5.874  & 1.89E-06 & 5.368  \\ 
        60 & 4.36E-07 & 6.153  & 1.04E-06 & 5.667  & 1.15E-07 & 5.797  & 2.21E-07 & 5.292  \\ 
        80 & 6.72E-08 & 6.500  & 1.82E-07 & 6.068  & 2.03E-08 & 6.017  & 4.59E-08 & 5.472  \\ 
        100 & 1.40E-08 & 7.020  & 4.10E-08 & 6.671  & 4.87E-09 & 6.405  & 1.25E-08 & 5.838  \\ 
        120 & 3.57E-09 & 7.506  & 1.09E-08 & 7.279  & 1.40E-09 & 6.852  & 3.94E-09 & 6.319  \\ 
        140 & 1.06E-09 & 7.852  & 3.28E-09 & 7.782  & 4.55E-10 & 7.280  & 1.37E-09 & 6.846  \\ 
        160 & 3.63E-10 & 8.050  & 1.10E-09 & 8.161  & 1.64E-10 & 7.634  & 5.16E-10 & 7.322 \\ 
        \hline
    \end{tabular}
\end{table}

\begin{table}[htbp]
    \centering
    \caption{Accuracy test. Error table of the density. \(k = 9\).}
    \label{tab:2D-accuracy-test-1-k=9}
    \begin{tabular}{ccccccccc}
        \hline
        & \multicolumn{4}{c}{CH-RI} & \multicolumn{4}{c}{CH-CON} \\
        \(N\) & \(\ell_h^2\) error & order & \(\ell_h^\infty\) error & order & \(\ell_h^2\) error & order & \(\ell_h^\infty\) error & order \\ \hline
        20 & 2.72E-05 & - & 3.39E-05 & - & 6.30E-06 & - & 4.11E-06 & - \\ 
        40 & 9.20E-08 & 8.205  & 1.00E-07 & 8.402  & 1.26E-08 & 8.968  & 8.93E-09 & 8.847  \\ 
        60 & 3.00E-09 & 8.440  & 3.92E-09 & 7.996  & 3.29E-10 & 8.990  & 2.39E-10 & 8.929  \\ 
        80 & 2.45E-10 & 8.715  & 3.41E-10 & 8.484  & 2.48E-11 & 8.977  & 1.83E-11 & 8.937  \\ 
        100 & 3.40E-11 & 8.851  & 4.75E-11 & 8.836  & 4.20E-12 & 7.961  & 3.61E-12 & 7.271  \\ 
        \hline
    \end{tabular}
\end{table}

\subsubsection{Double Mach reflection}

\begin{figure}[htbp]
    \centering
    \begin{subfigure}{0.32\linewidth}
        \centering
        \includegraphics[width=\linewidth]{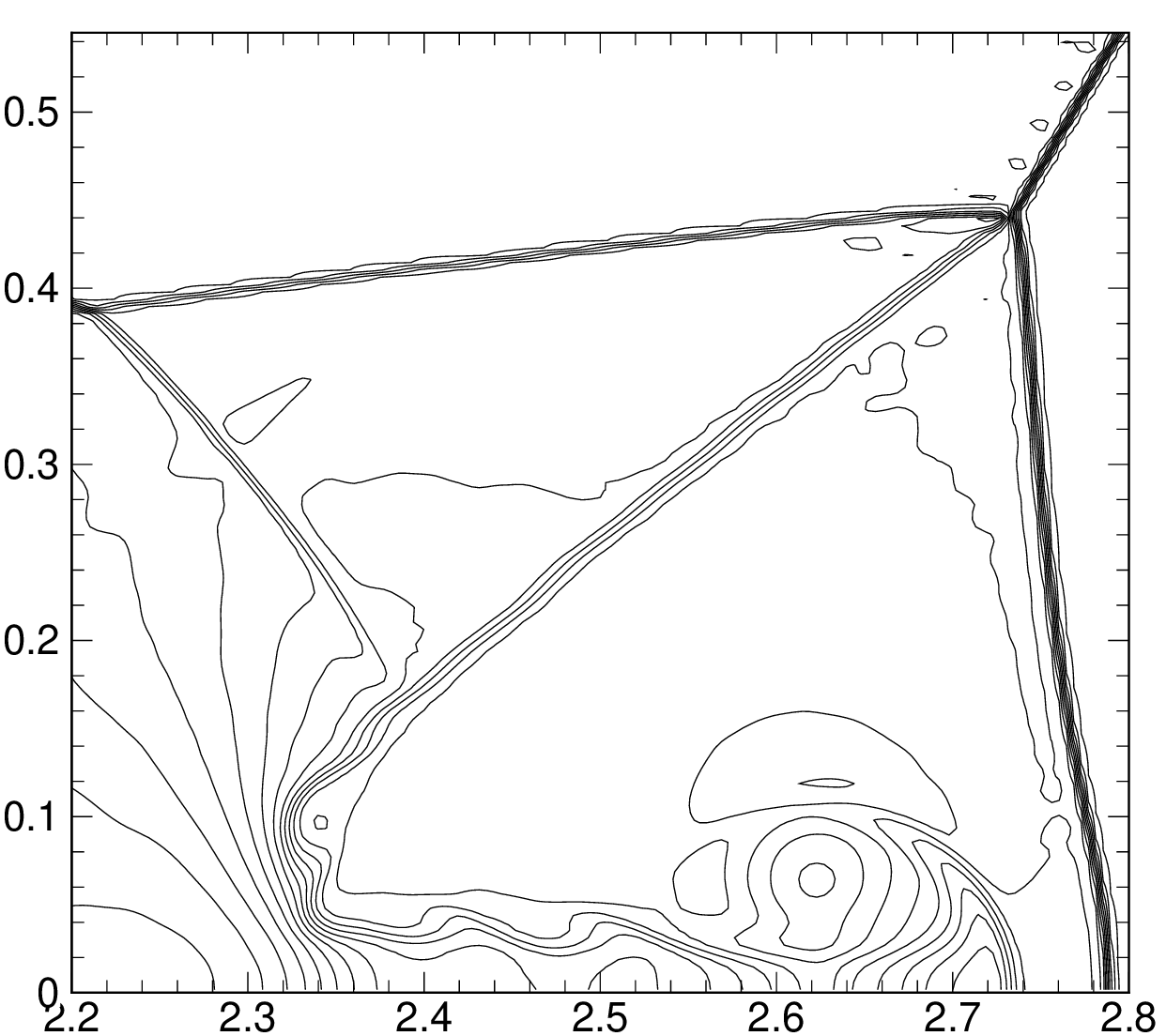}
        \caption{AWENO5, CH-RI}
    \end{subfigure}
    \begin{subfigure}{0.32\linewidth}
        \centering
        \includegraphics[width=\linewidth]{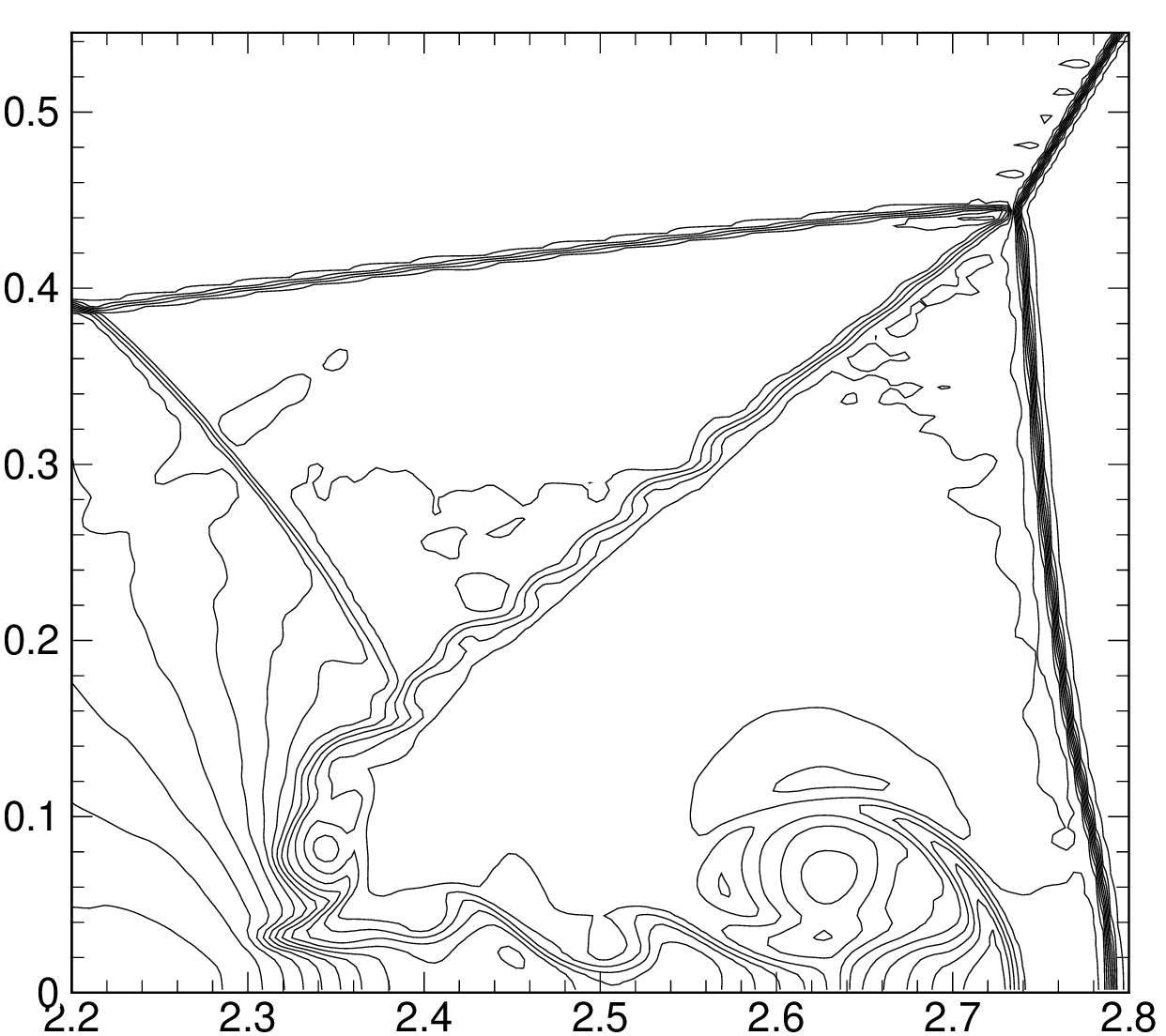}
        \caption{AWENO7, CH-RI}
    \end{subfigure}
    \begin{subfigure}{0.32\linewidth}
        \centering
        \includegraphics[width=\linewidth]{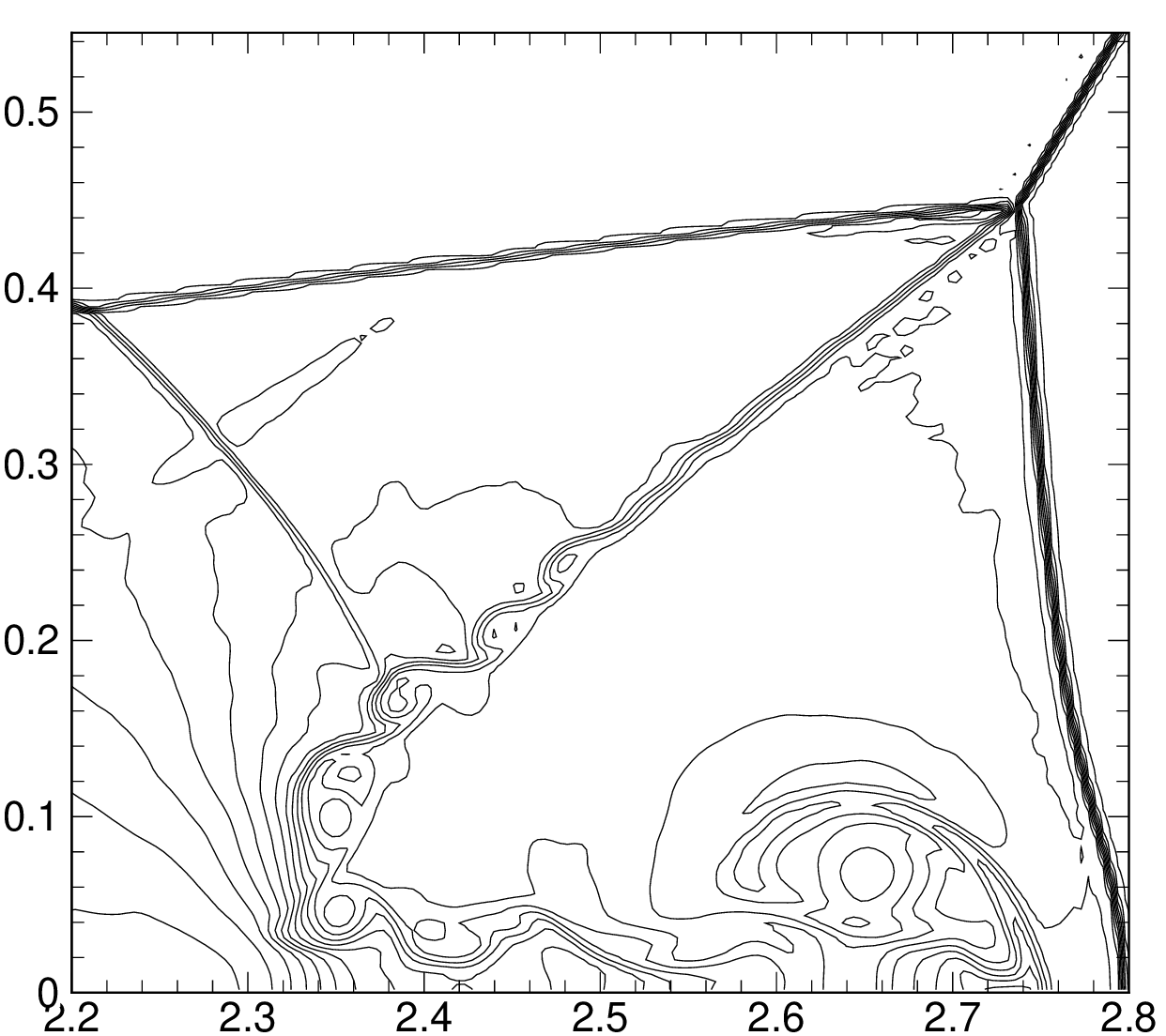}
        \caption{AWENO9, CH-RI}
    \end{subfigure}

    \begin{subfigure}{0.32\linewidth}
        \centering
        \includegraphics[width=\linewidth]{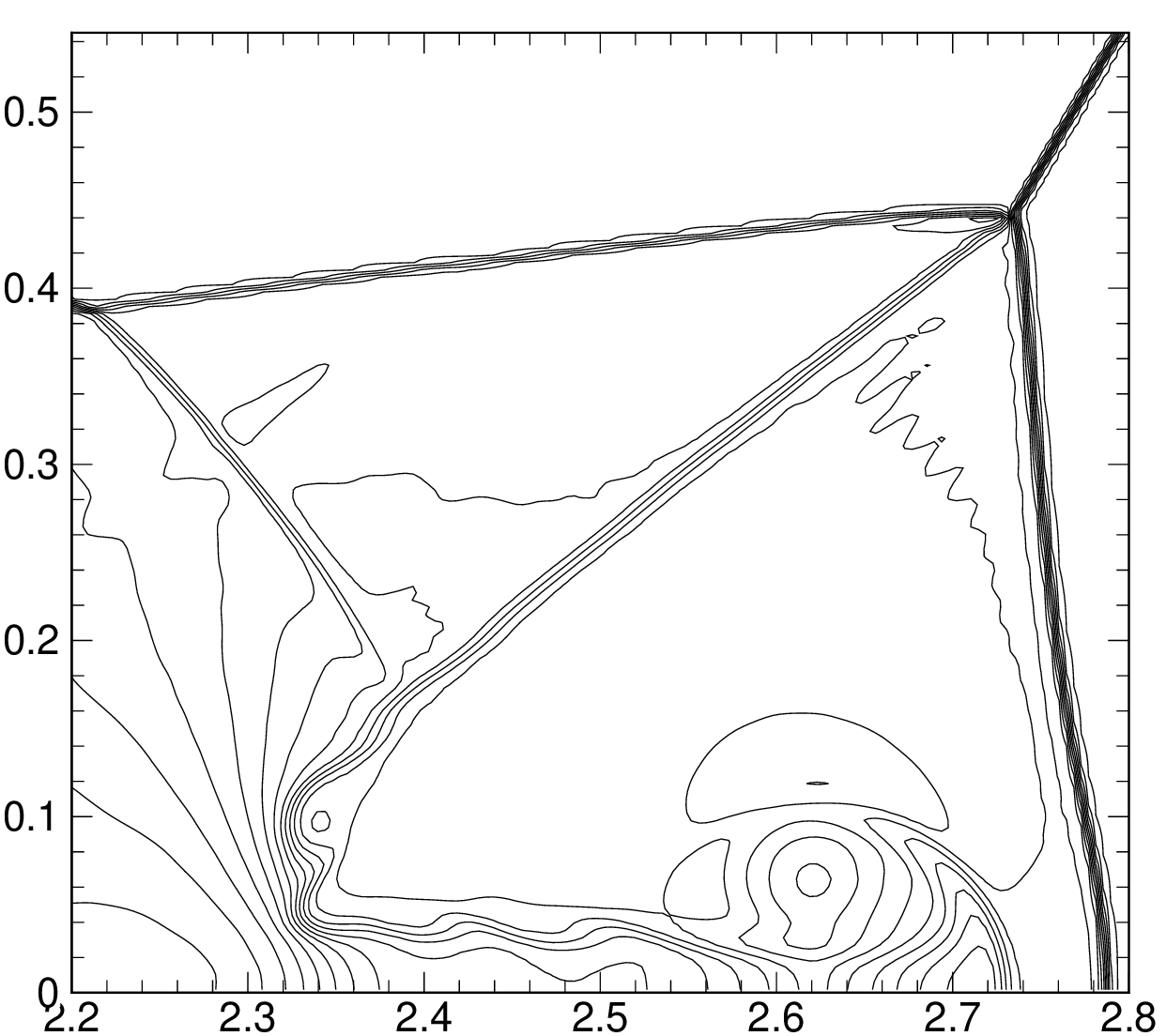}
        \caption{AWENO5, CH-CON}
    \end{subfigure}
    \begin{subfigure}{0.32\linewidth}
        \centering
        \includegraphics[width=\linewidth]{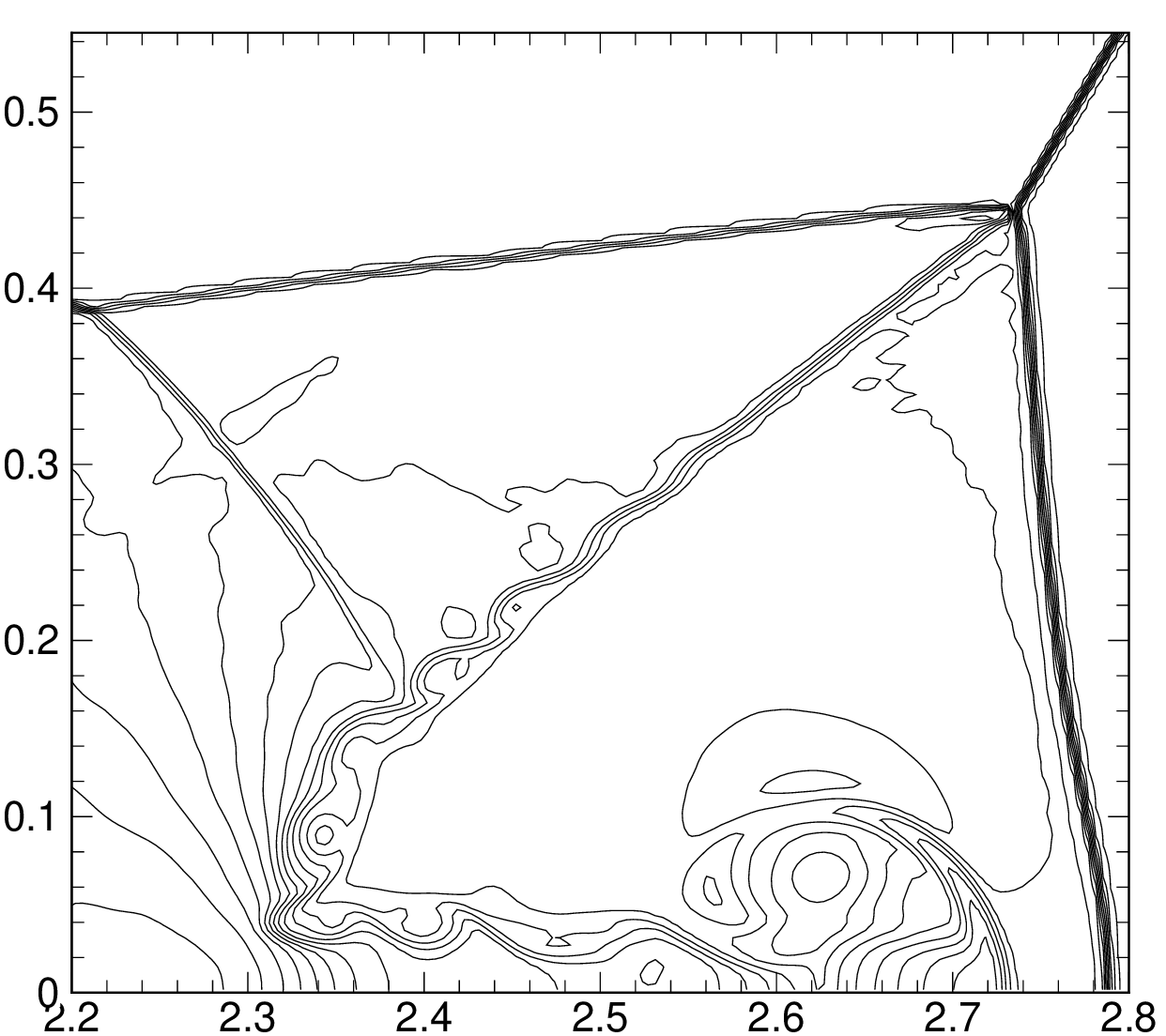}
        \caption{AWENO7, CH-CON}
    \end{subfigure}
    \begin{subfigure}{0.32\linewidth}
        \centering
        \includegraphics[width=\linewidth]{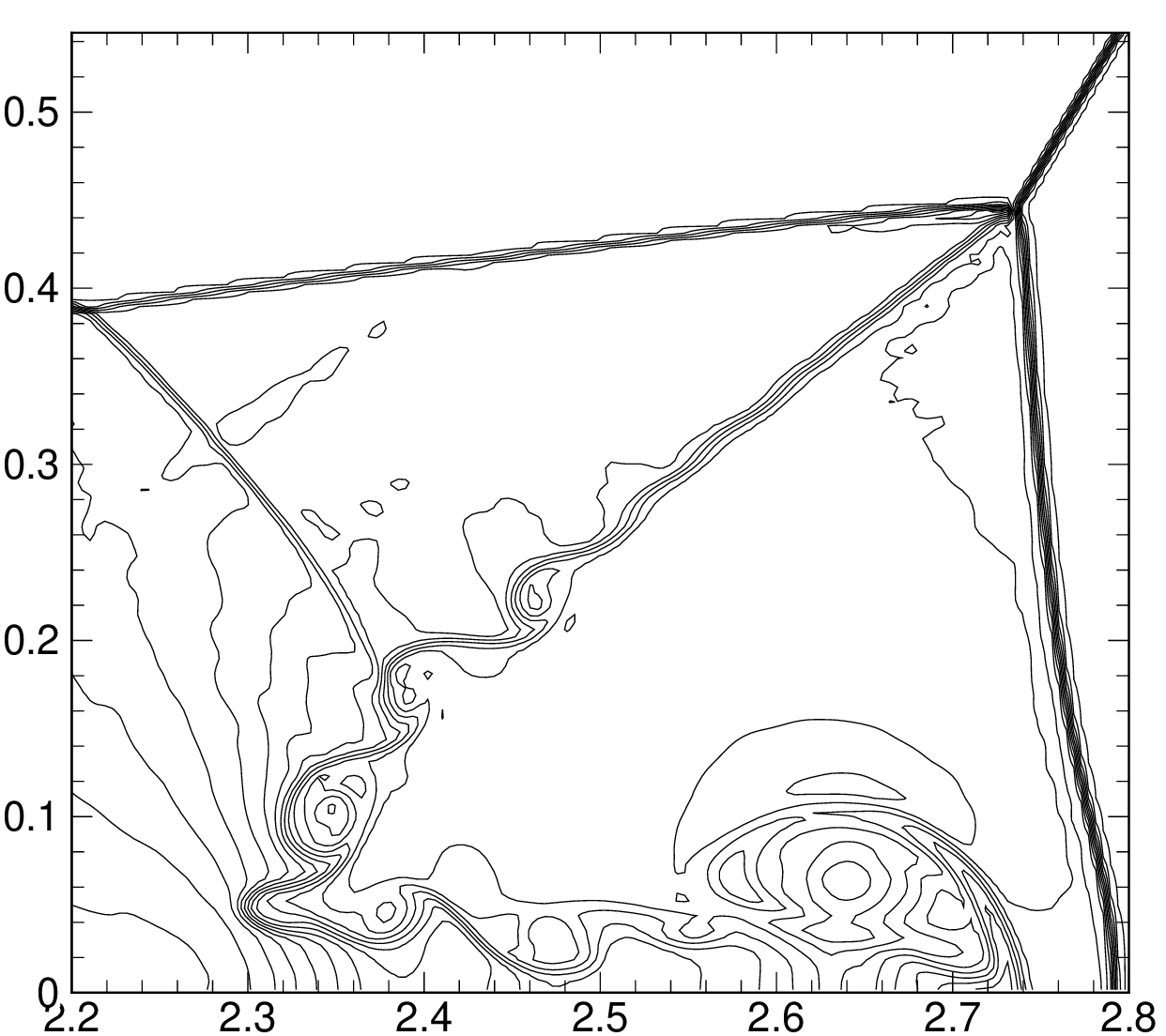}
        \caption{AWENO9, CH-CON}
    \end{subfigure}

    \caption{On the coarse mesh \((Nx, Ny) = (960, 240)\).}
    \label{fig:double-mach-reflection-1}
\end{figure}

\begin{figure}[htbp]
    \centering
    \begin{subfigure}{0.32\linewidth}
        \centering
        \includegraphics[width=\linewidth]{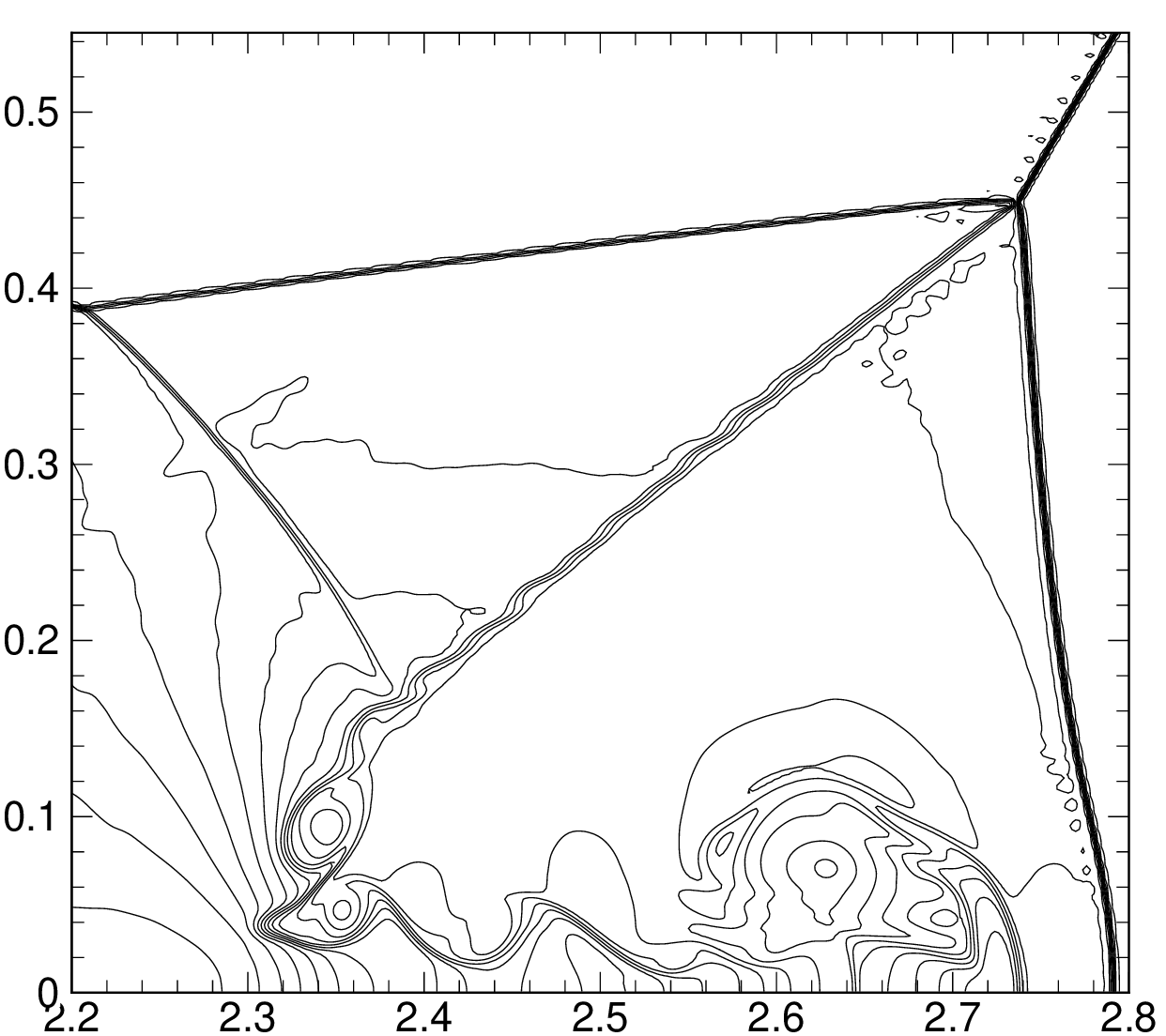}
        \caption{AWENO5, CH-RI}
    \end{subfigure}
    \begin{subfigure}{0.32\linewidth}
        \centering
        \includegraphics[width=\linewidth]{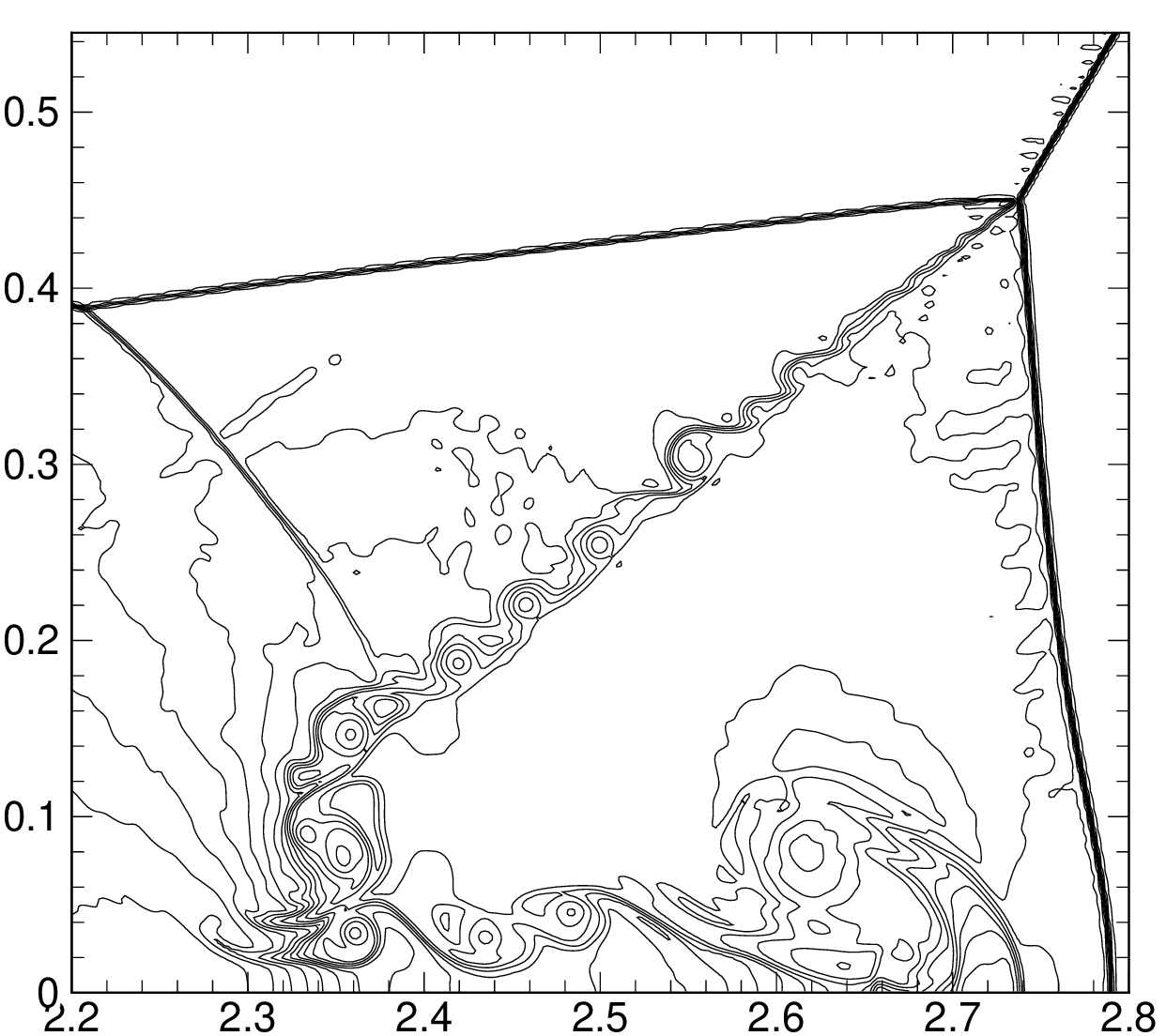}
        \caption{AWENO7, CH-RI}
    \end{subfigure}
    \begin{subfigure}{0.32\linewidth}
        \centering
        \includegraphics[width=\linewidth]{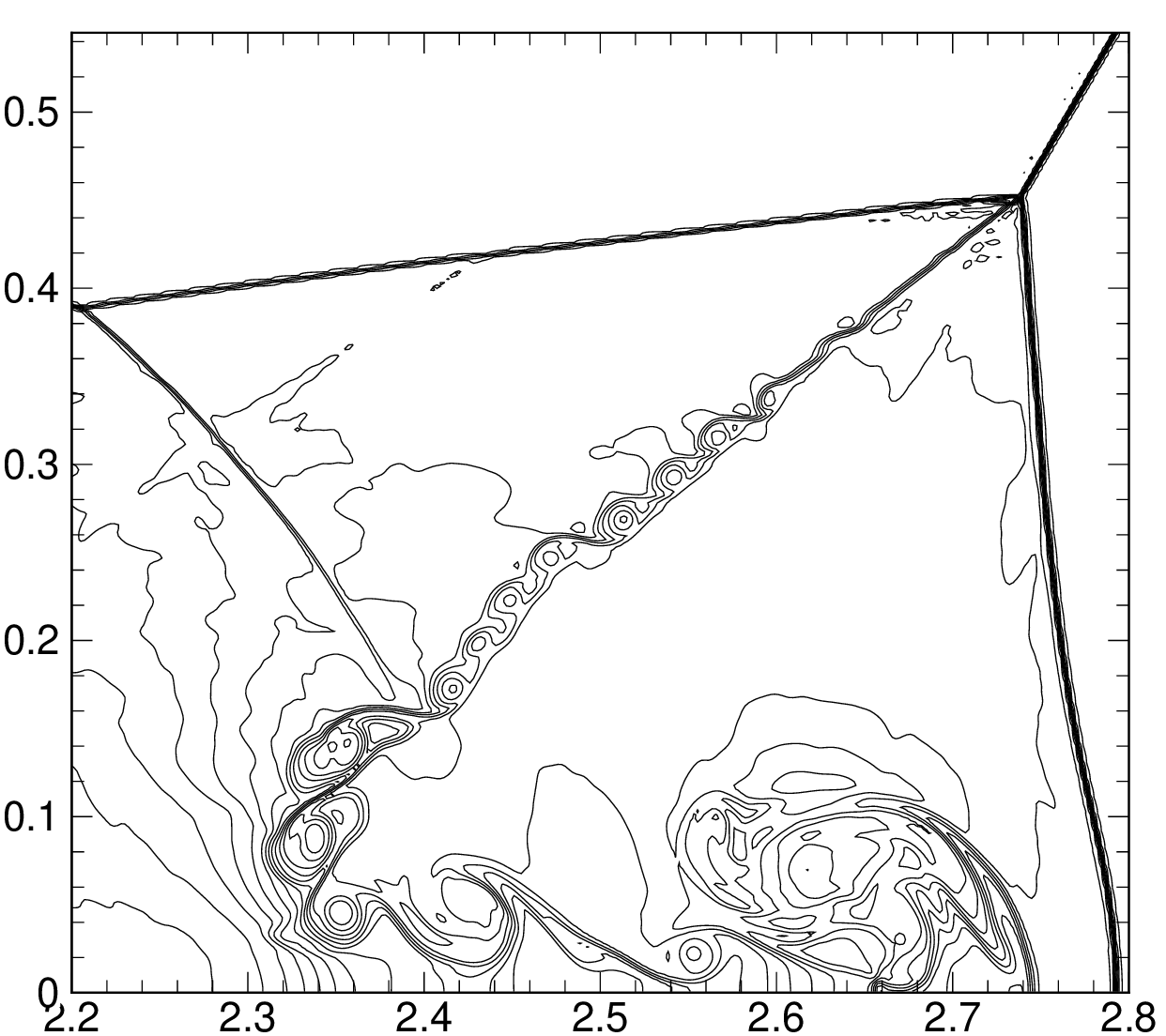}
        \caption{AWENO9, CH-RI}
    \end{subfigure}

    \begin{subfigure}{0.32\linewidth}
        \centering
        \includegraphics[width=\linewidth]{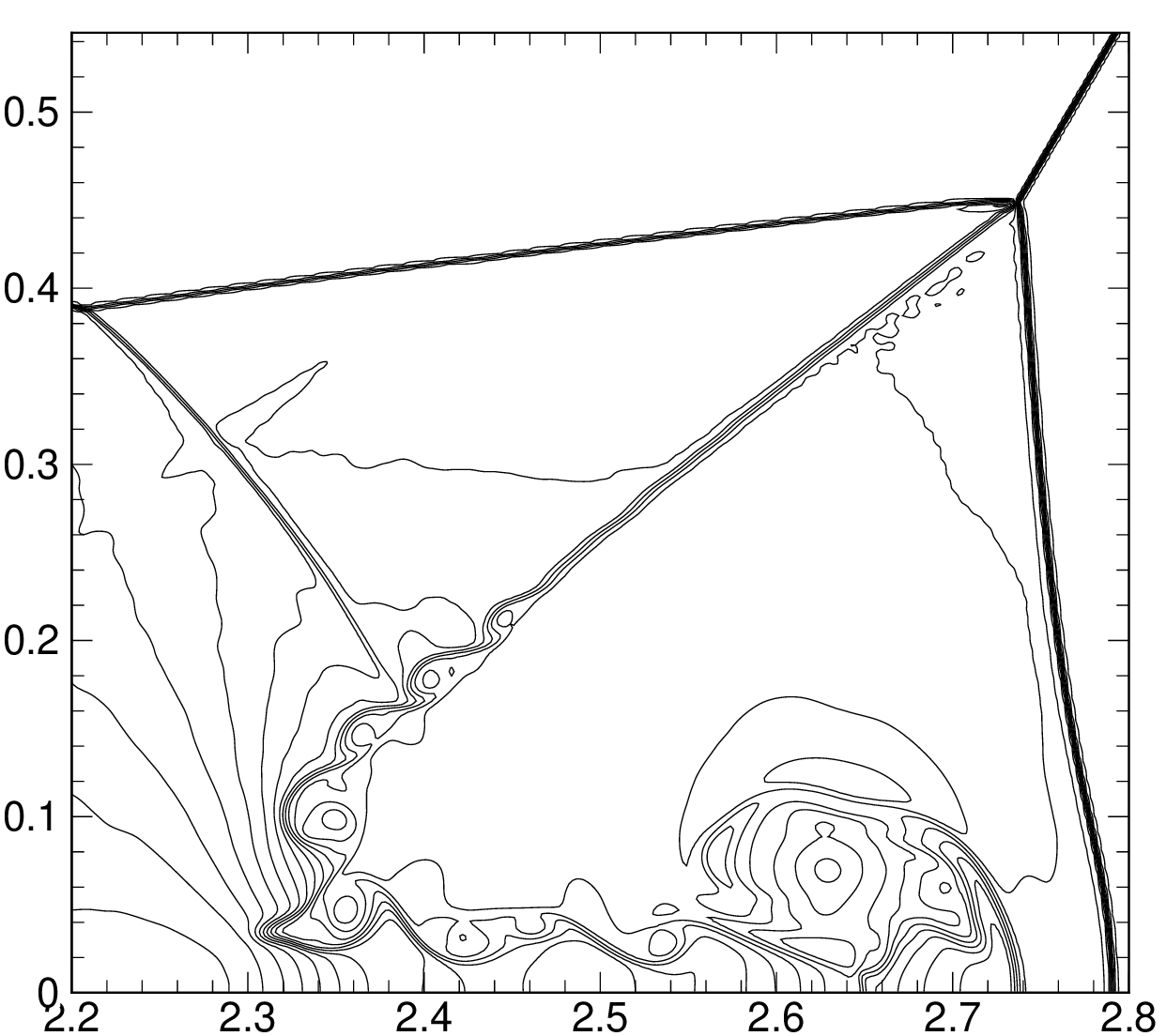}
        \caption{AWENO5, CH-CON}
    \end{subfigure}
    \begin{subfigure}{0.32\linewidth}
        \centering
        \includegraphics[width=\linewidth]{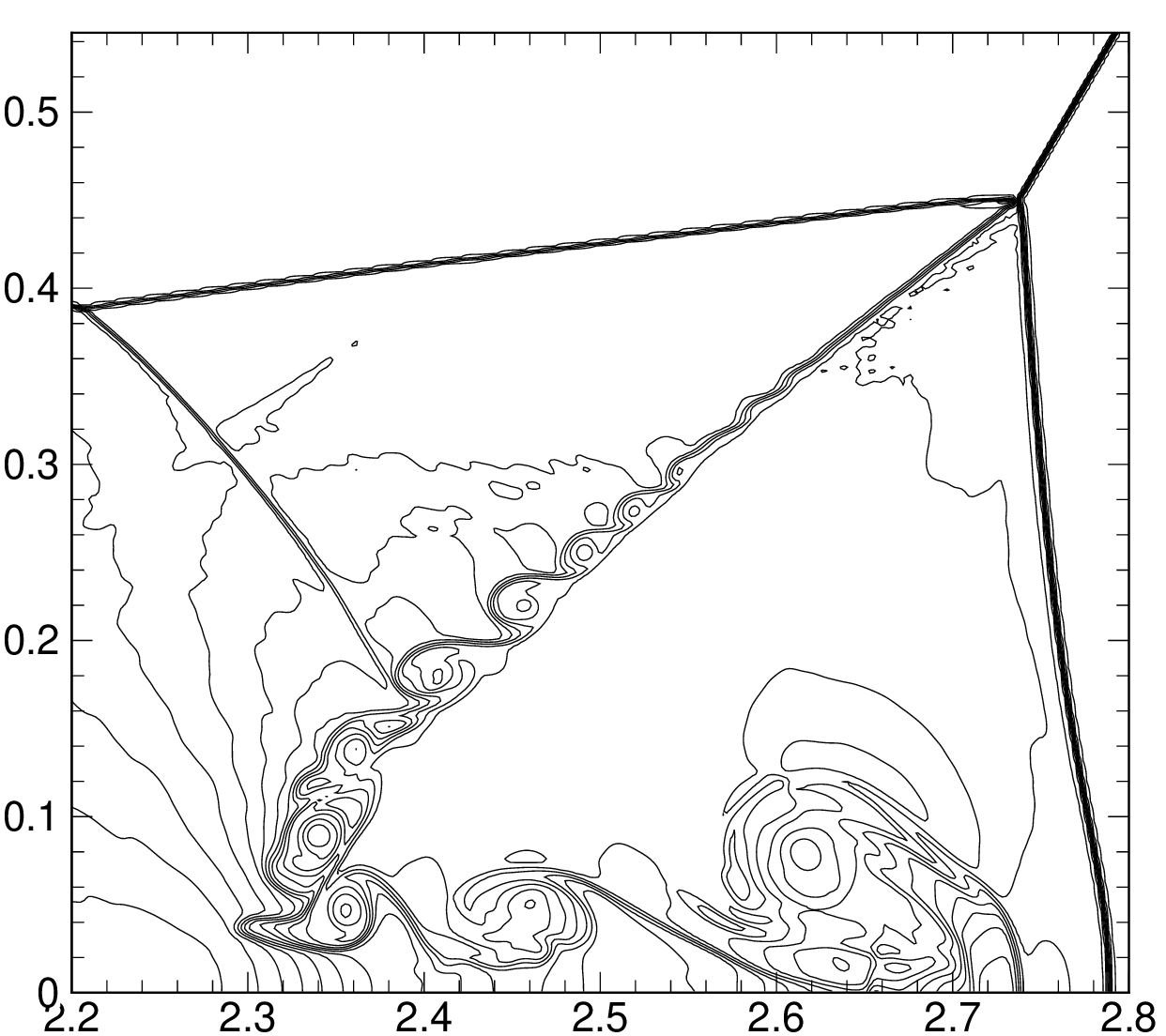}
        \caption{AWENO7, CH-CON}
    \end{subfigure}
    \begin{subfigure}{0.32\linewidth}
        \centering
        \includegraphics[width=\linewidth]{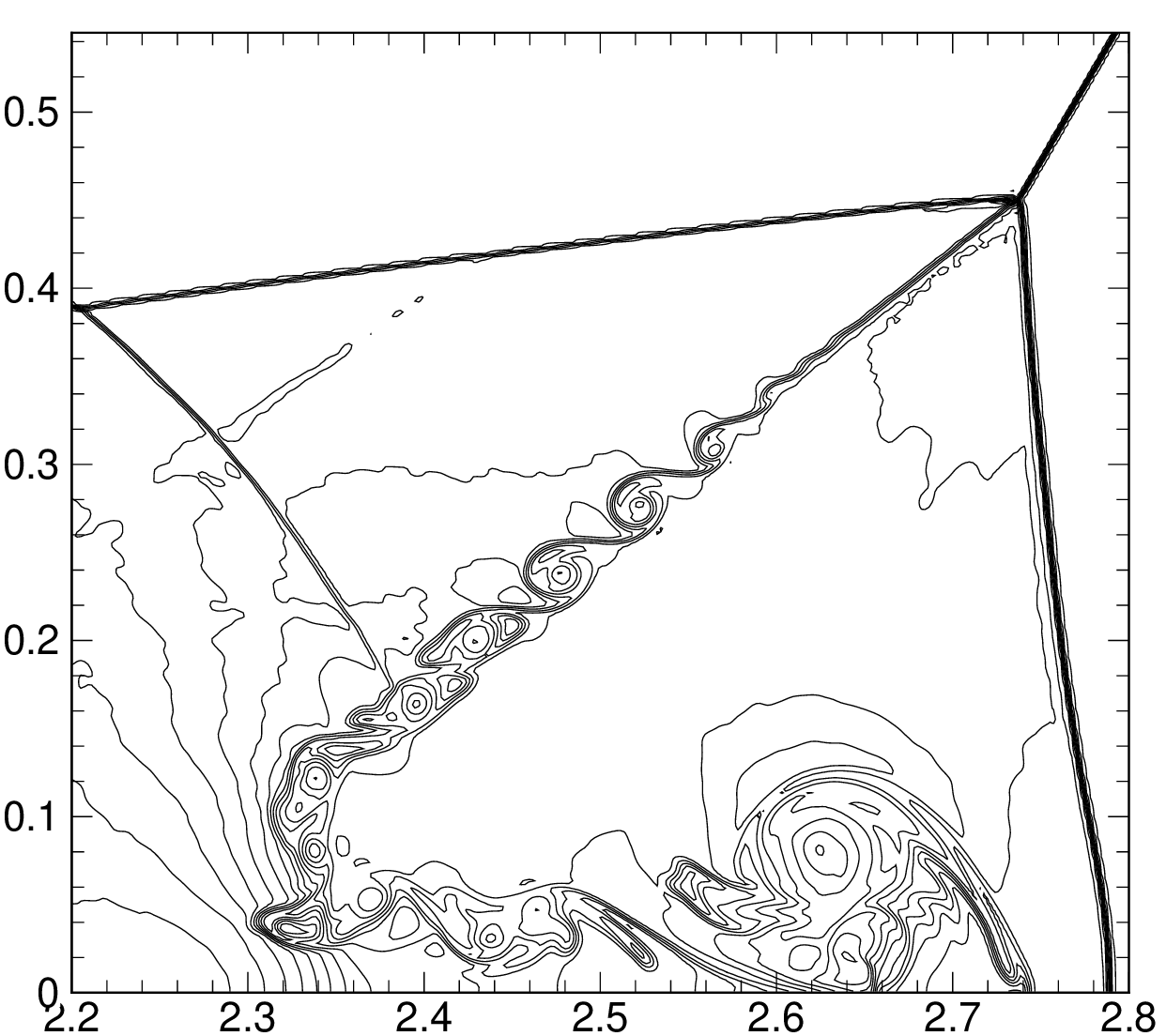}
        \caption{AWENO9, CH-CON}
    \end{subfigure}

    \caption{On the fine mesh \((Nx, Ny) = (1920, 480)\).}
    \label{fig:double-mach-reflection-2}
\end{figure}

We consider the double Mach reflection (DMR) problem \cite{10.1016/0021-99918490142-6}. The computational domain is \([0, 4] \times [0, 1]\). Let \(x_0 = \frac{1}{6}\) denote the start of the reflective wall and \(\theta = \frac{\pi}{6}\) denote the angle between the incoming flow and the reflective wall. The reflective wall lies at the bottom boundary starting from \(x_0\). Initially, a Mach 10 shock is positioned at \((x,y) = (x_0, 0)\), and it is perpendicular to the incoming flow direction. The post-shock (to the left of the shock) primitive variables are \((\rho, u, v, p) = (8, 8.5 \cos(-\theta), 8.5 \sin(-\theta), 116.5)\) and the pre-shock (to the right of the shock) primitive variables are \((\rho, u, v, p) = (1.4, 0, 0, 1)\). The velocity of the shock front is \((u_s, v_s) = (10 \cos(-\theta), 10 \sin(-\theta))\). Except for the reflective wall, all other boundaries are posed with the exact motion of the shock. We compute to \(T = 0.2\) and plot the density profile using 30 equally spaced contours from \(1.5\) to \(22.7\). The results are shown in \Cref{fig:double-mach-reflection-1,fig:double-mach-reflection-2}. We can see that both methods agree well. 

\subsubsection{Rayleigh--Taylor instability}

\begin{figure}[htbp]
    \centering
    \begin{subfigure}{0.28\linewidth}
        \centering
        \includegraphics[width=\linewidth]{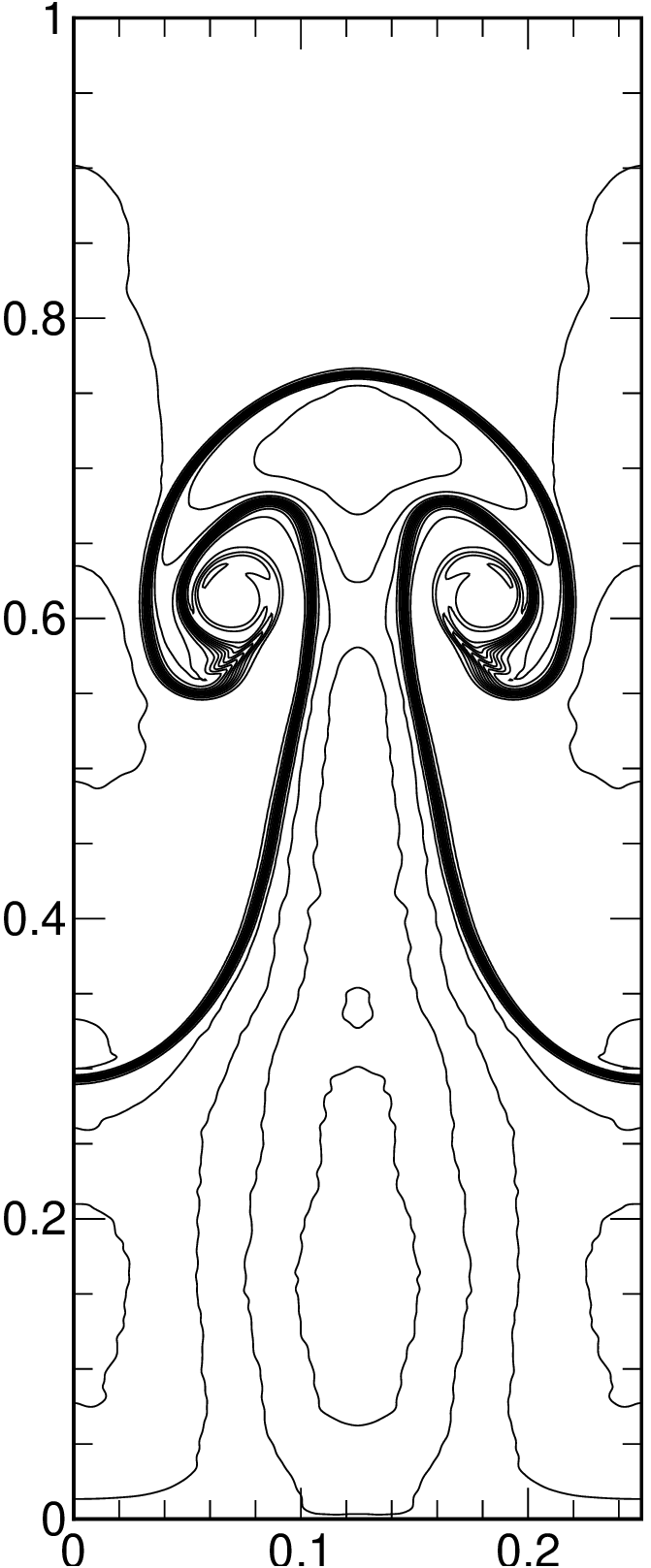}
        \caption{AWENO5, CH-RI}
    \end{subfigure}
    \begin{subfigure}{0.28\linewidth}
        \centering
        \includegraphics[width=\linewidth]{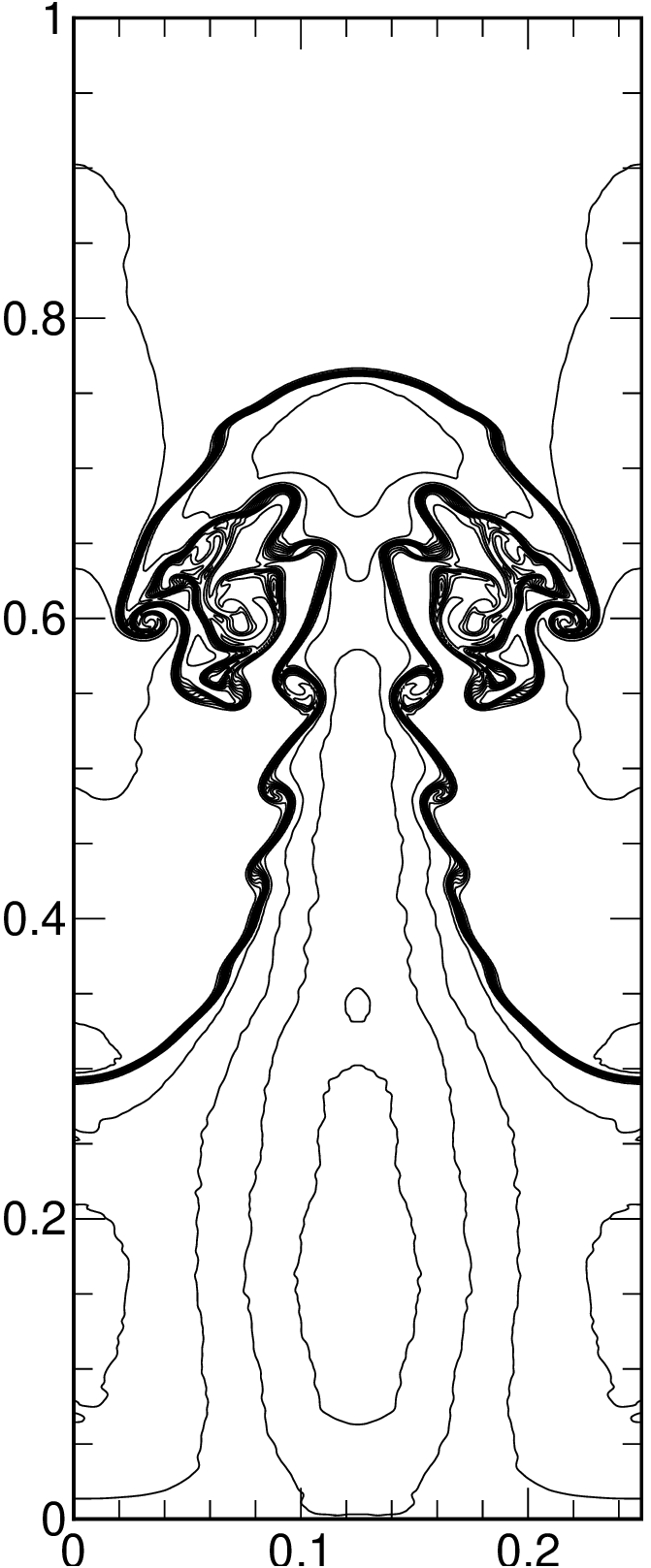}
        \caption{AWENO7, CH-RI}
    \end{subfigure}
    \begin{subfigure}{0.28\linewidth}
        \centering
        \includegraphics[width=\linewidth]{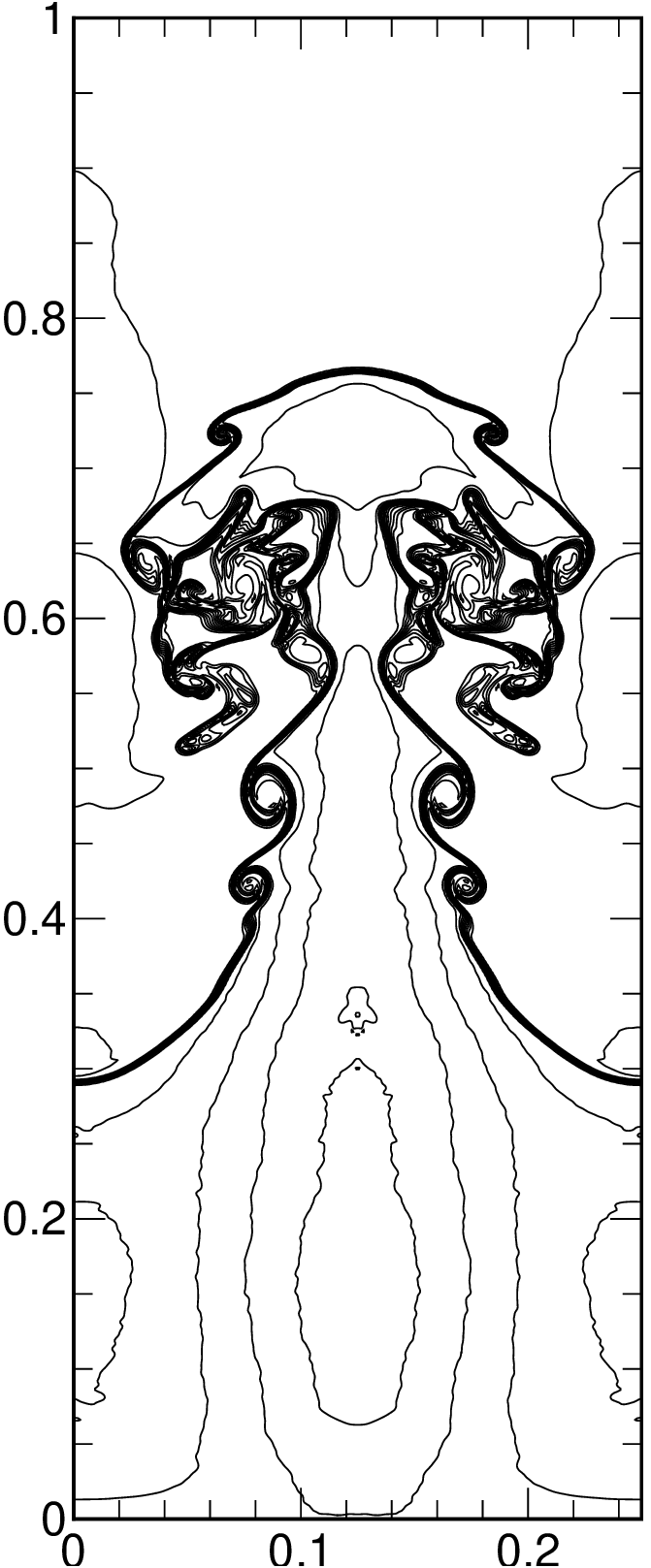}
        \caption{AWENO9, CH-RI}
    \end{subfigure}

    \begin{subfigure}{0.28\linewidth}
        \centering
        \includegraphics[width=\linewidth]{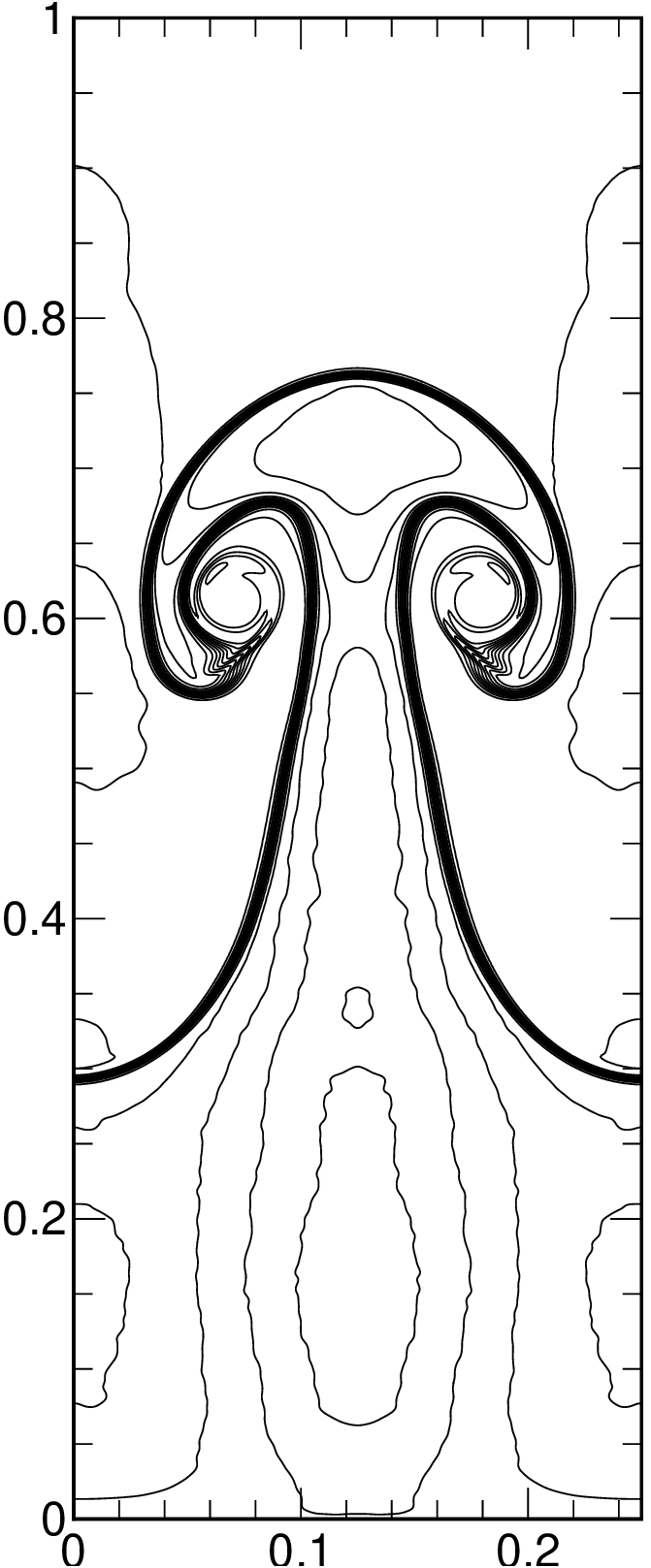}
        \caption{AWENO5, CH-CON}
    \end{subfigure}
    \begin{subfigure}{0.28\linewidth}
        \centering
        \includegraphics[width=\linewidth]{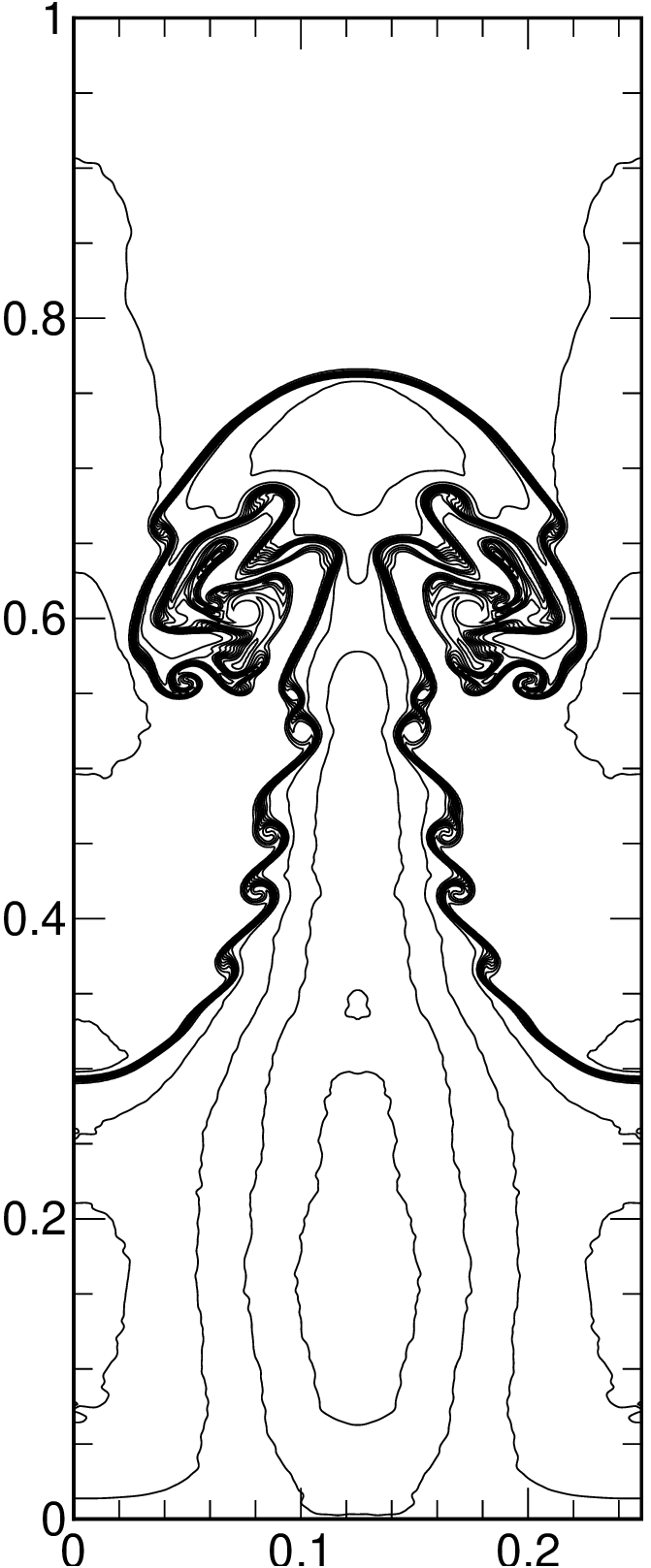}
        \caption{AWENO7, CH-CON}
    \end{subfigure}
    \begin{subfigure}{0.28\linewidth}
        \centering
        \includegraphics[width=\linewidth]{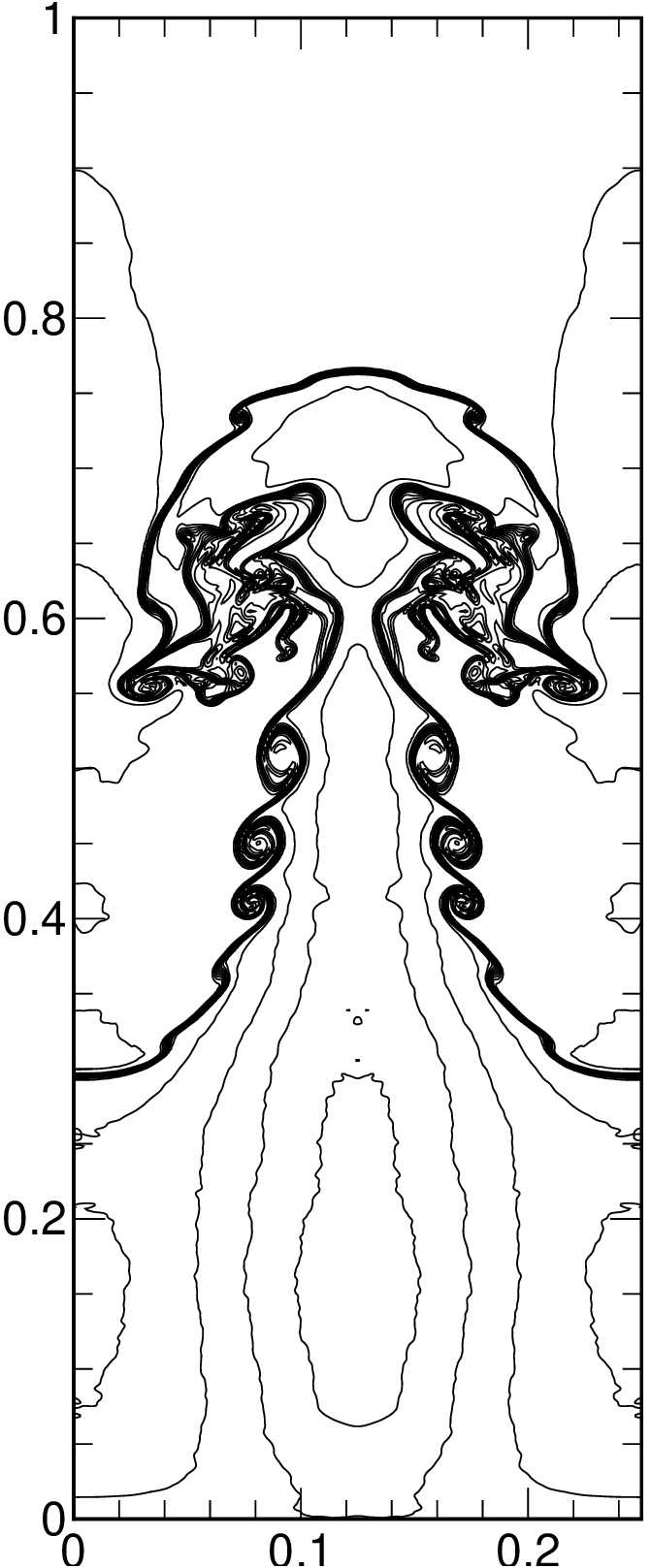}
        \caption{AWENO9, CH-CON}
    \end{subfigure}

    \caption{On the mesh \((Nx, Ny) = (240, 960)\).}
    \label{fig:RTI-3}
\end{figure}

We consider the Rayleigh--Taylor instability (RTI) problem \cite{10.1016/S0021-99910300094-9}. The computational domain is \([0, \frac{1}{4}] \times [0, 1]\). Initially, the interface is located at \(y = \frac{1}{2}\). The heavy fluid with density \(\rho = 2\) is below the interface, and the light fluid with density \(\rho = 1\) is above the interface, with acceleration in the positive \(y\)-direction, i.e. \((g_x, g_y) = (0, 1)\). The pressure \(p\) is continuous across the interface. A small perturbation is given to the \(y\)-direction fluid speed: for \(0 \leqslant y < \frac{1}{2}\), \((\rho, u, v, p) = (2, 0, -0.025 c \cdot \cos(8 \pi x), 2 y + 1)\), and for \(\frac{1}{2} \leqslant y \leqslant 1\), \((\rho, u, v, p) = (1, 0, -0.025 c \cdot \cos(8 \pi x), y + \frac{3}{2})\), where \(c = \sqrt{\frac{\gamma p}{\rho}}\) is the sonic speed. The specific heat capacity ratio is set to \(\gamma = \frac{5}{3}\). Reflective wall boundary conditions are imposed for both left and right boundaries. For the top boundary, the flow values are set to be \((\rho, u, v, p) = (1, 0, 0, 2.5)\). For the bottom boundary, the flow values are set to be \((\rho, u, v, p) = (2, 0, 0, 1)\). The source terms \(\rho g_y\) and \(\rho v g_y\) are added to the right hand side of the third and fourth equations. The final time is \(T = 1.95\). The results are shown in \Cref{fig:RTI-3} using 15 equally spaced contours from \(0.9\) to \(2.2\). 

\subsubsection{Kelvin--Helmholtz instability}

\begin{figure}[htbp]
    \centering
    \begin{subfigure}{0.32\linewidth}
        \centering
        \includegraphics[width=\linewidth]{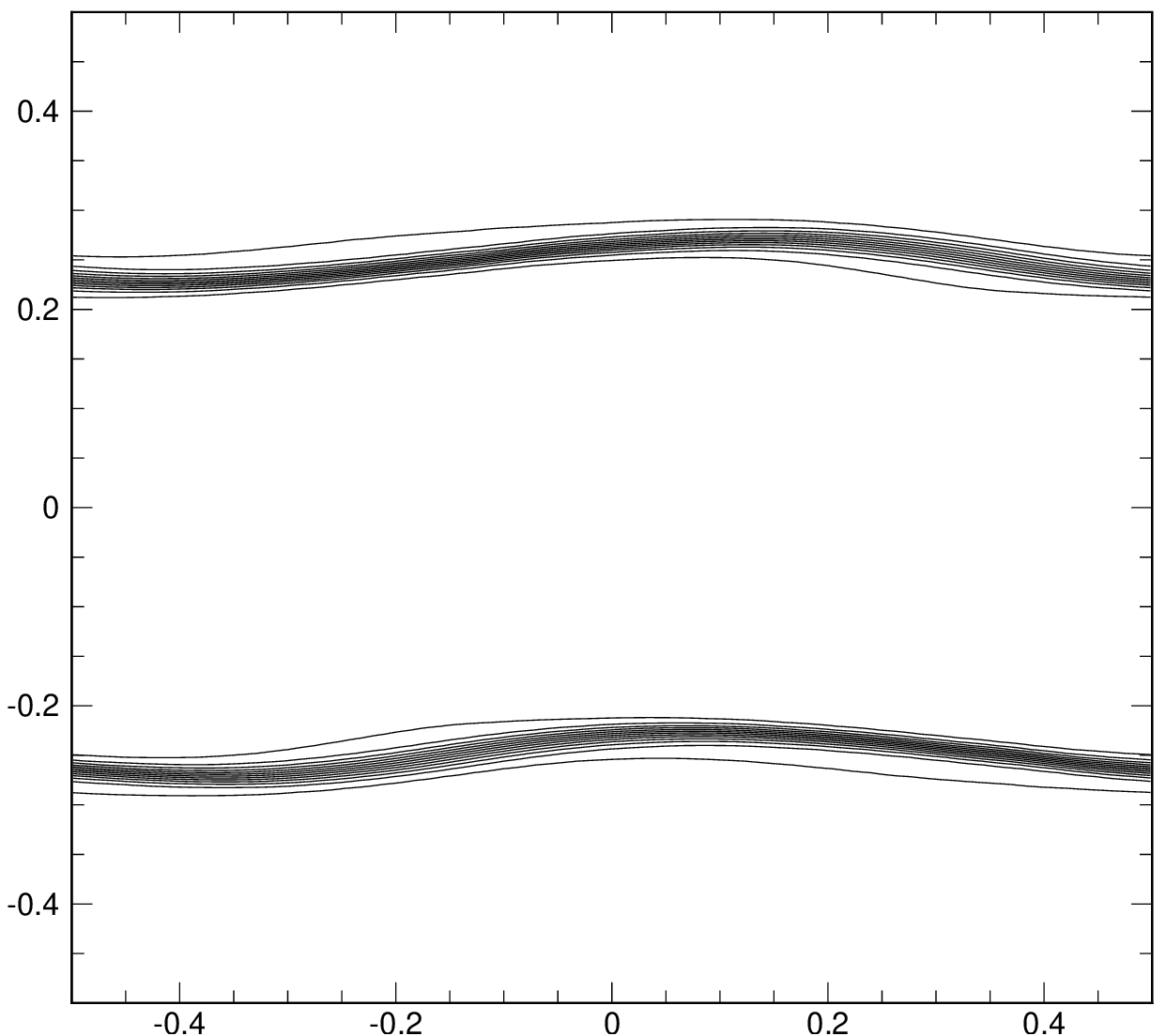}
        \caption{AWENO5, CH-RI}
    \end{subfigure}
    \begin{subfigure}{0.32\linewidth}
        \centering
        \includegraphics[width=\linewidth]{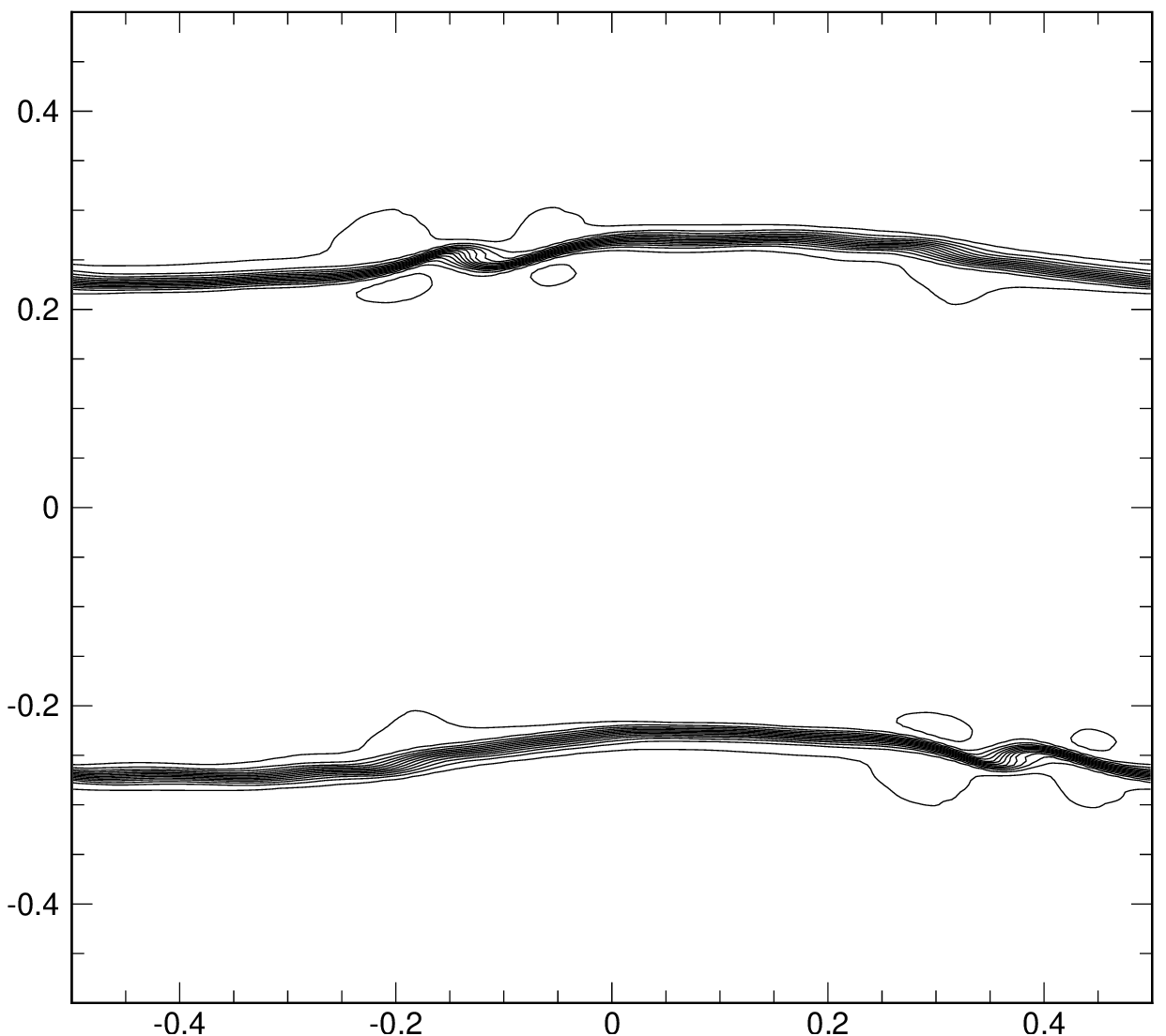}
        \caption{AWENO7, CH-RI}
    \end{subfigure}
    \begin{subfigure}{0.32\linewidth}
        \centering
        \includegraphics[width=\linewidth]{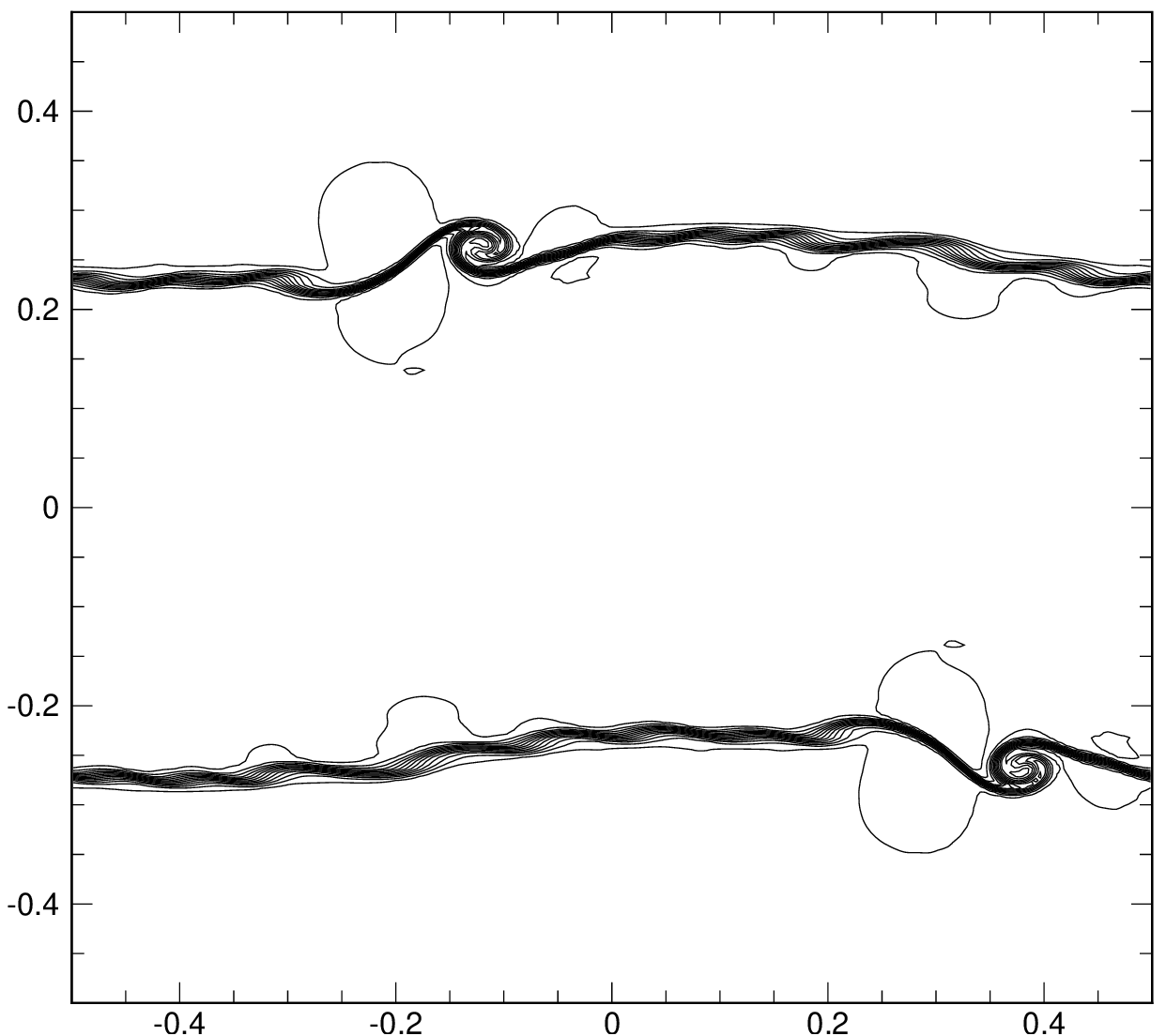}
        \caption{AWENO9, CH-RI}
    \end{subfigure}

    \begin{subfigure}{0.32\linewidth}
        \centering
        \includegraphics[width=\linewidth]{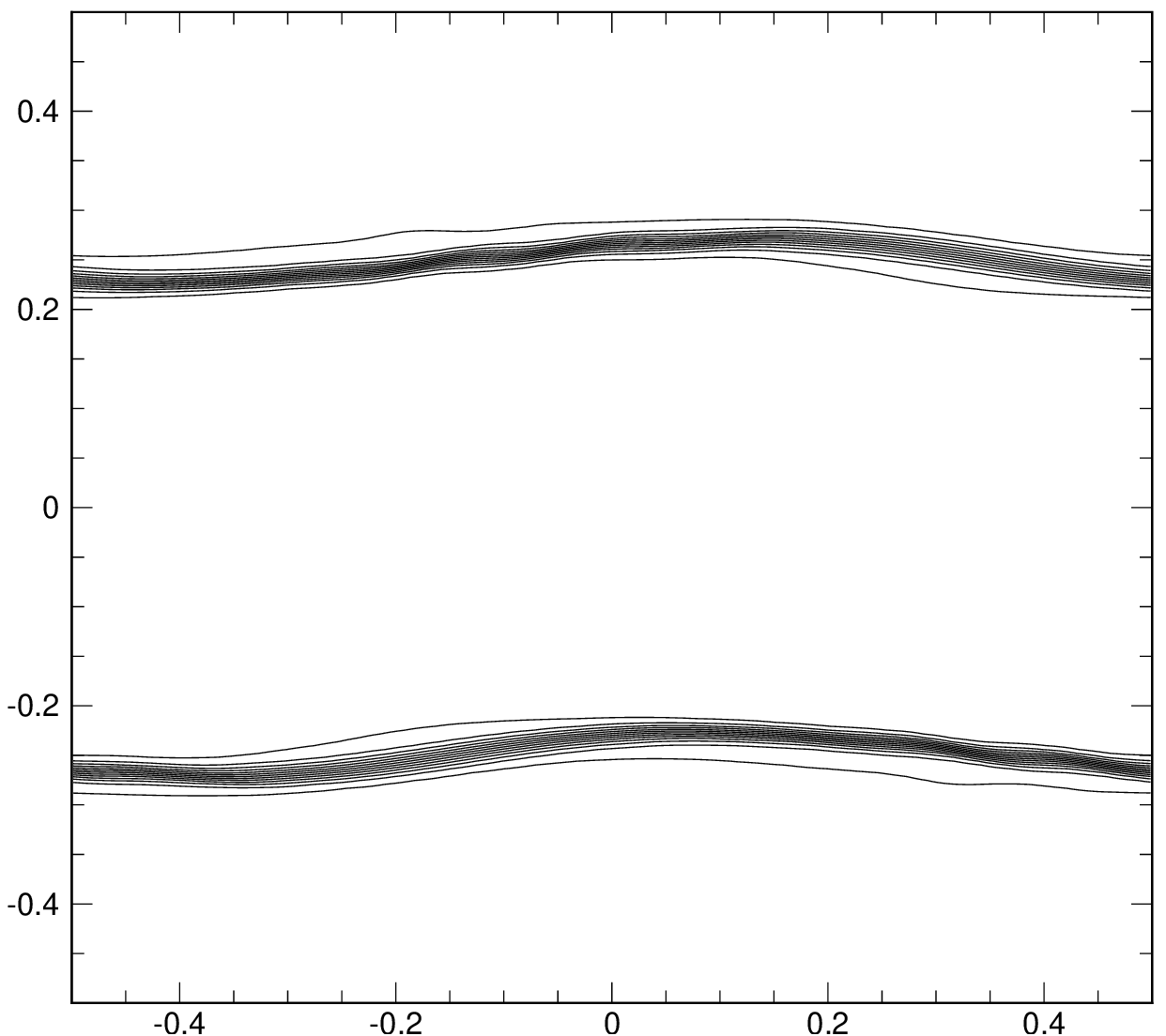}
        \caption{AWENO5, CH-CON}
    \end{subfigure}
    \begin{subfigure}{0.32\linewidth}
        \centering
        \includegraphics[width=\linewidth]{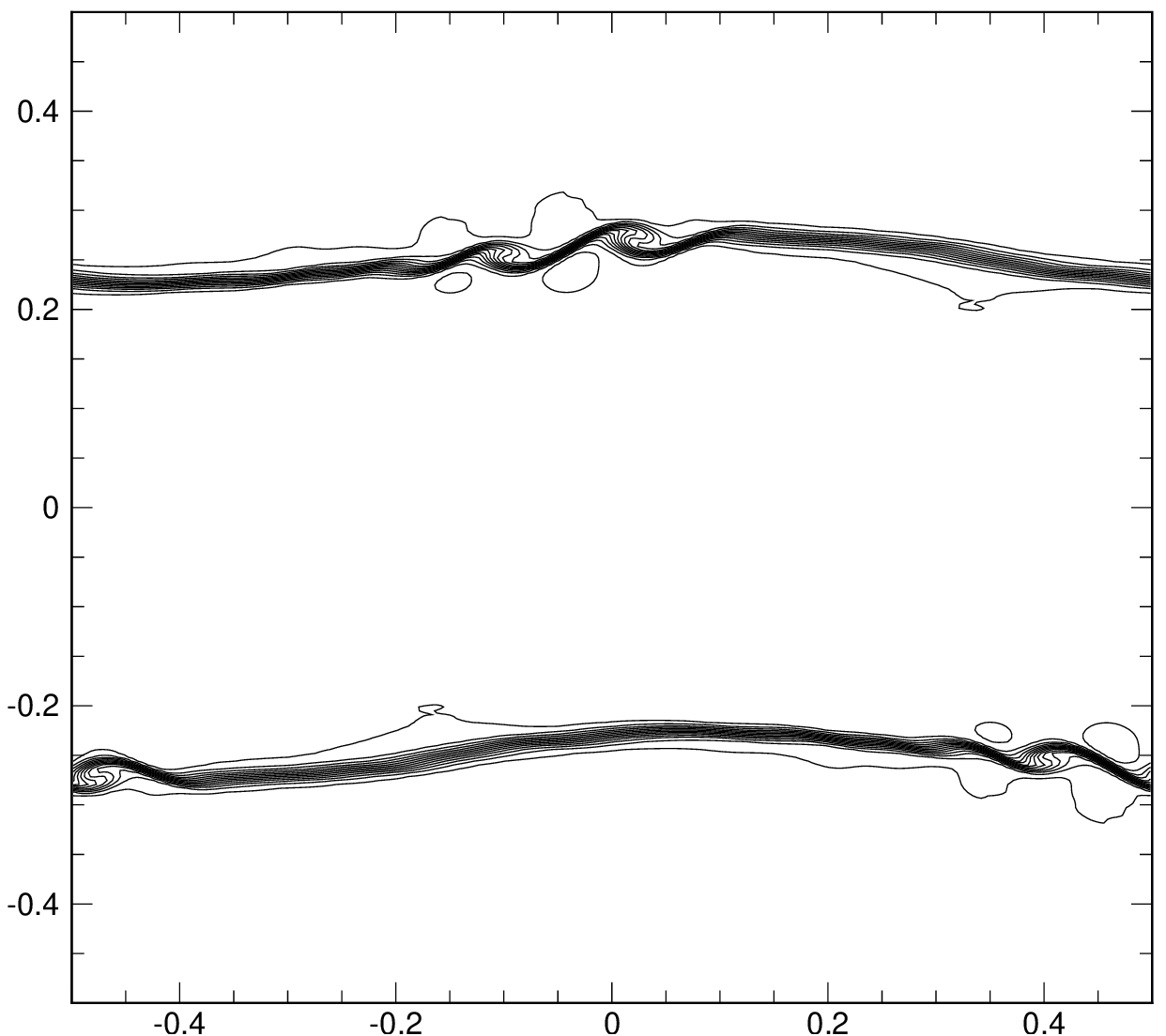}
        \caption{AWENO7, CH-CON}
    \end{subfigure}
    \begin{subfigure}{0.32\linewidth}
        \centering
        \includegraphics[width=\linewidth]{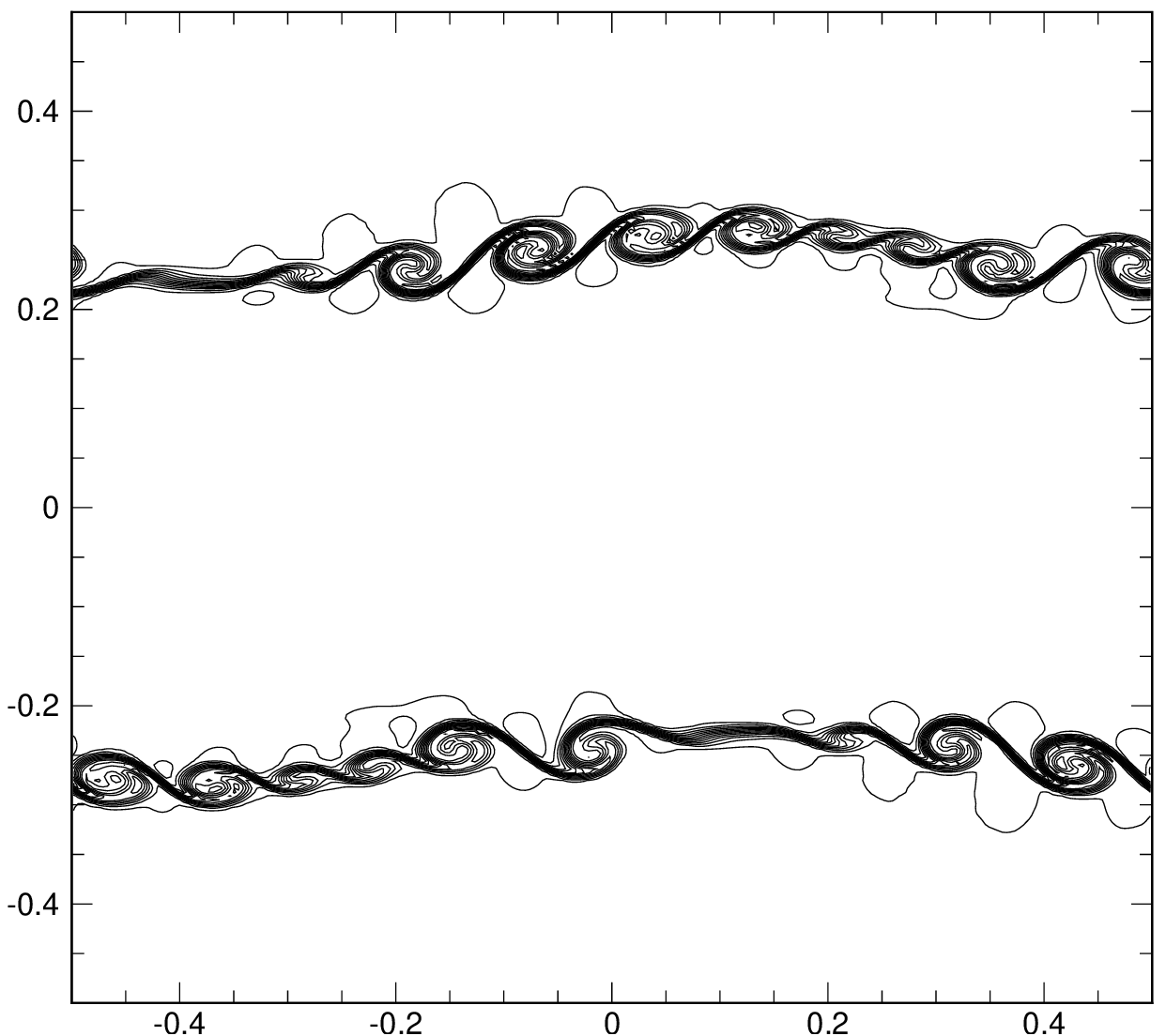}
        \caption{AWENO9, CH-CON}
    \end{subfigure}

    \caption{On the coarse mesh \((Nx, Ny) = (256, 256)\).}
    \label{fig:KHI-1}
\end{figure}

\begin{figure}[htbp]
    \centering
    \begin{subfigure}{0.32\linewidth}
        \centering
        \includegraphics[width=\linewidth]{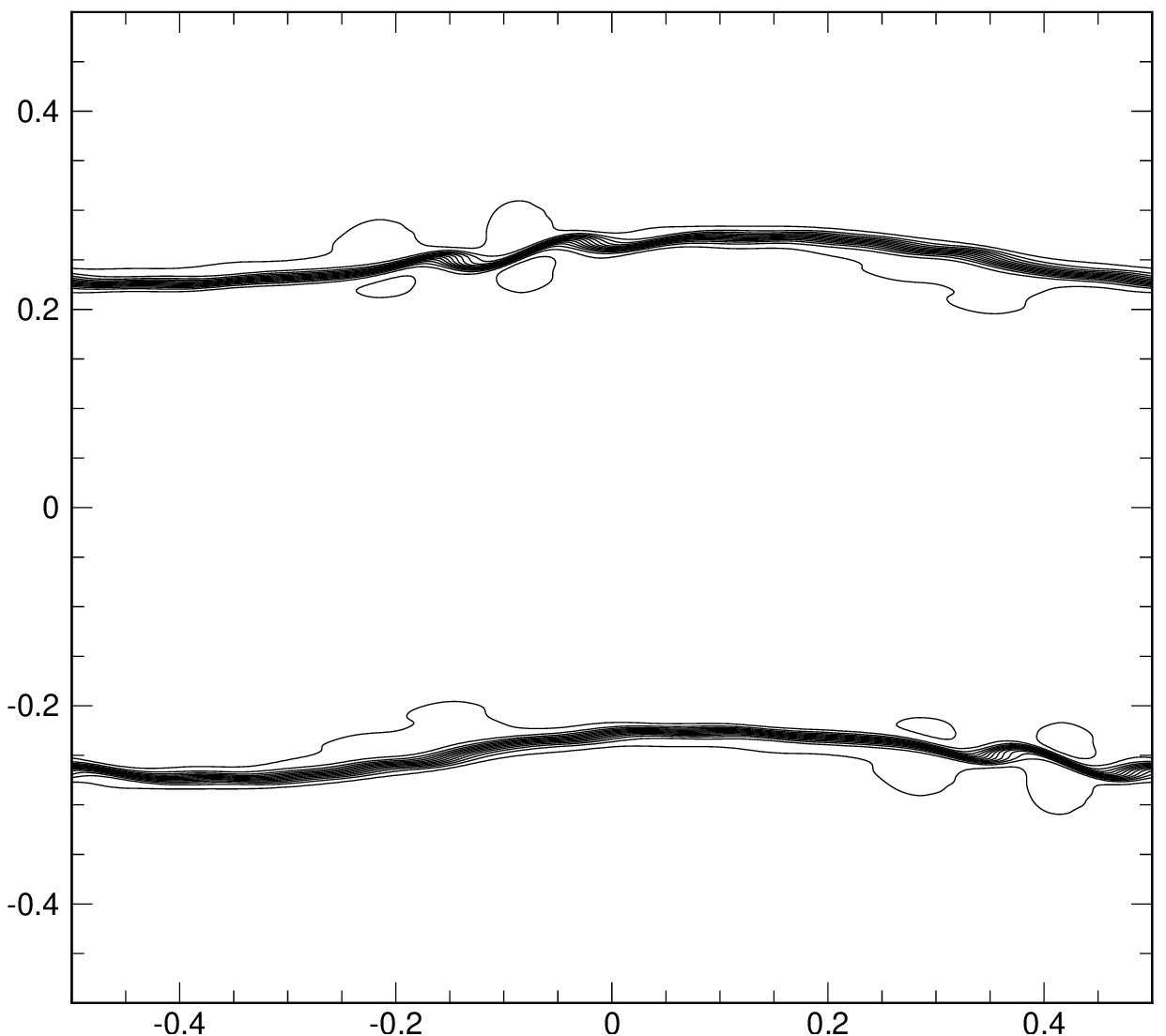}
        \caption{AWENO5, CH-RI}
    \end{subfigure}
    \begin{subfigure}{0.32\linewidth}
        \centering
        \includegraphics[width=\linewidth]{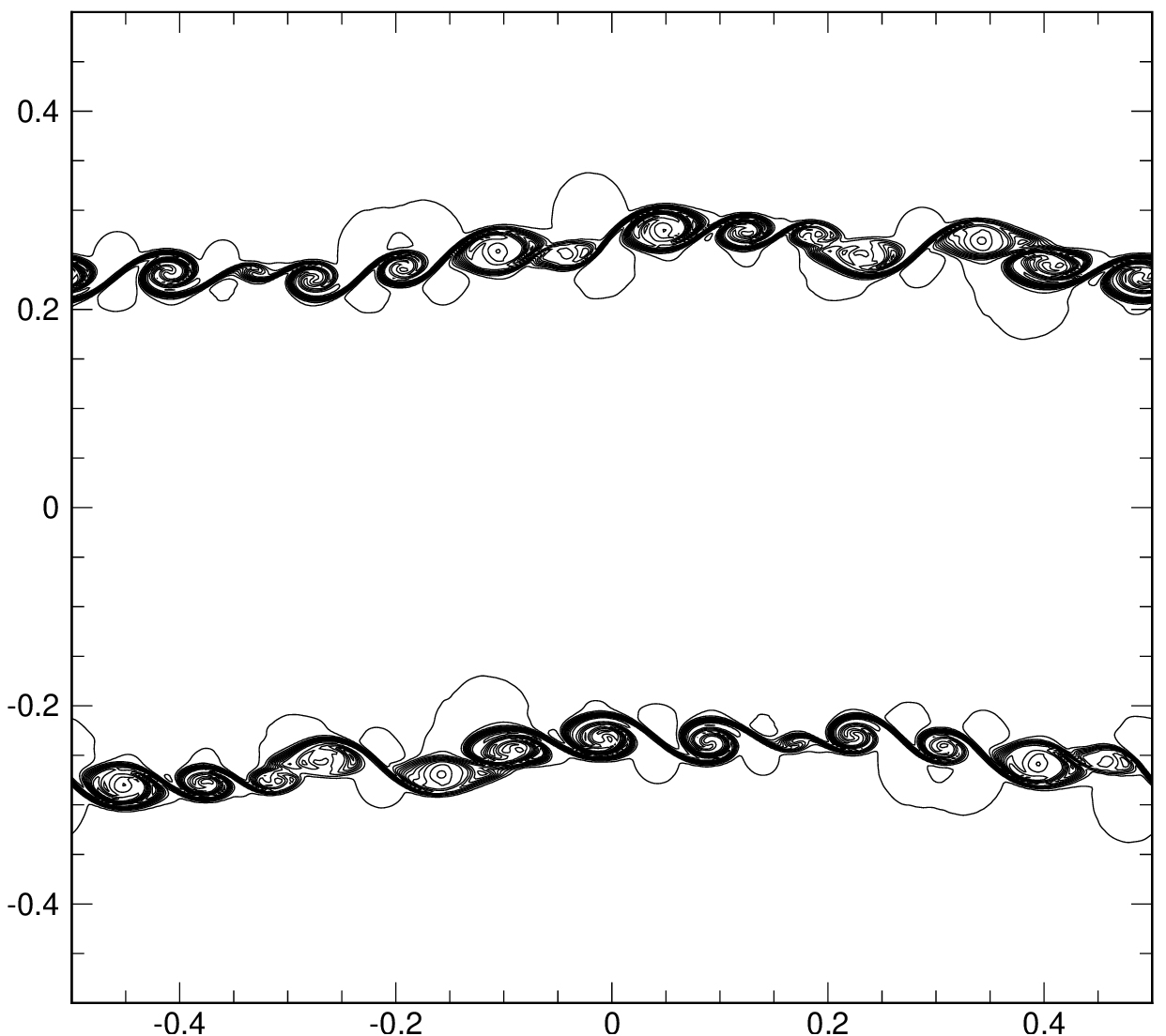}
        \caption{AWENO7, CH-RI}
    \end{subfigure}
    \begin{subfigure}{0.32\linewidth}
        \centering
        \includegraphics[width=\linewidth]{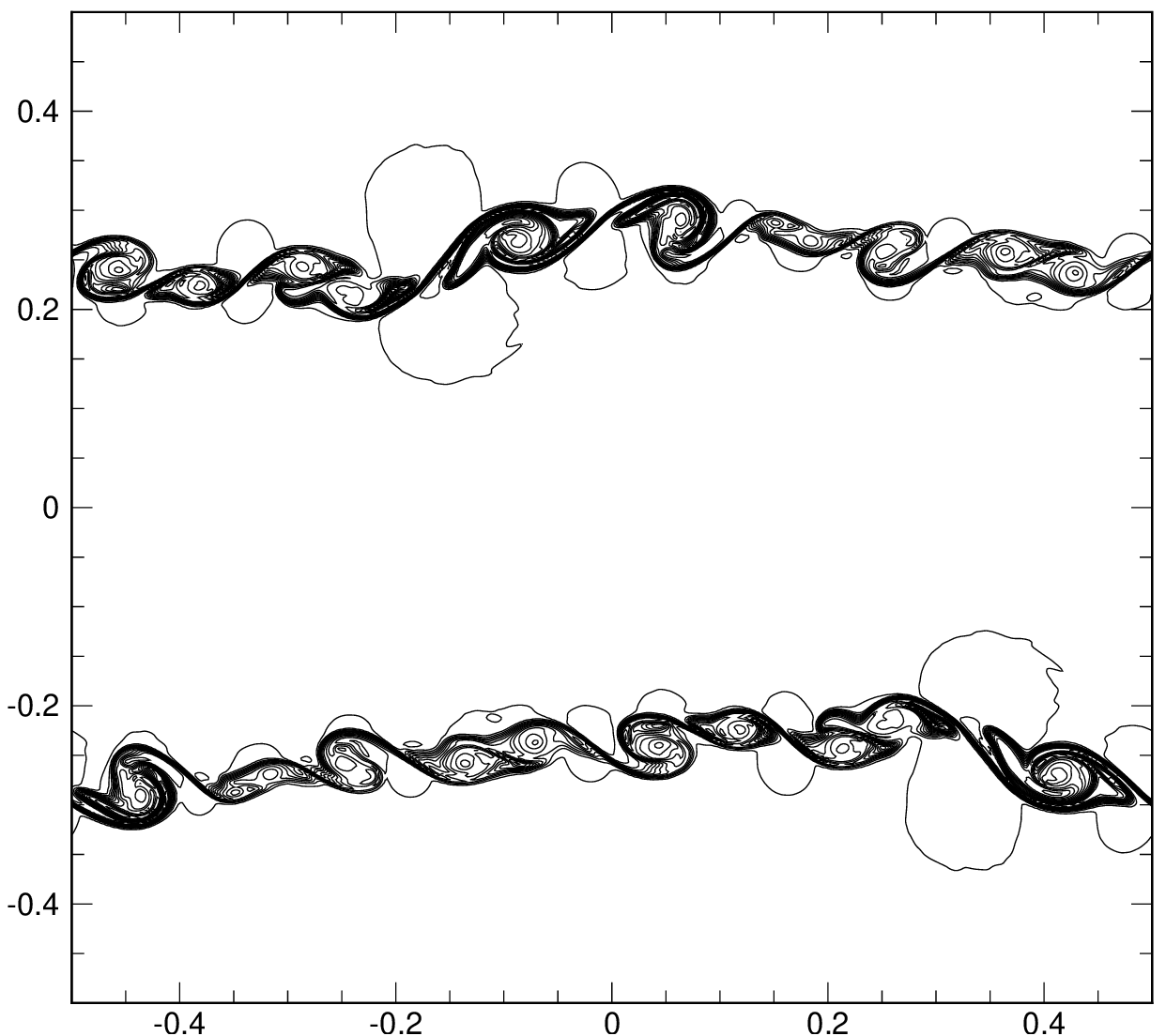}
        \caption{AWENO9, CH-RI}
    \end{subfigure}

    \begin{subfigure}{0.32\linewidth}
        \centering
        \includegraphics[width=\linewidth]{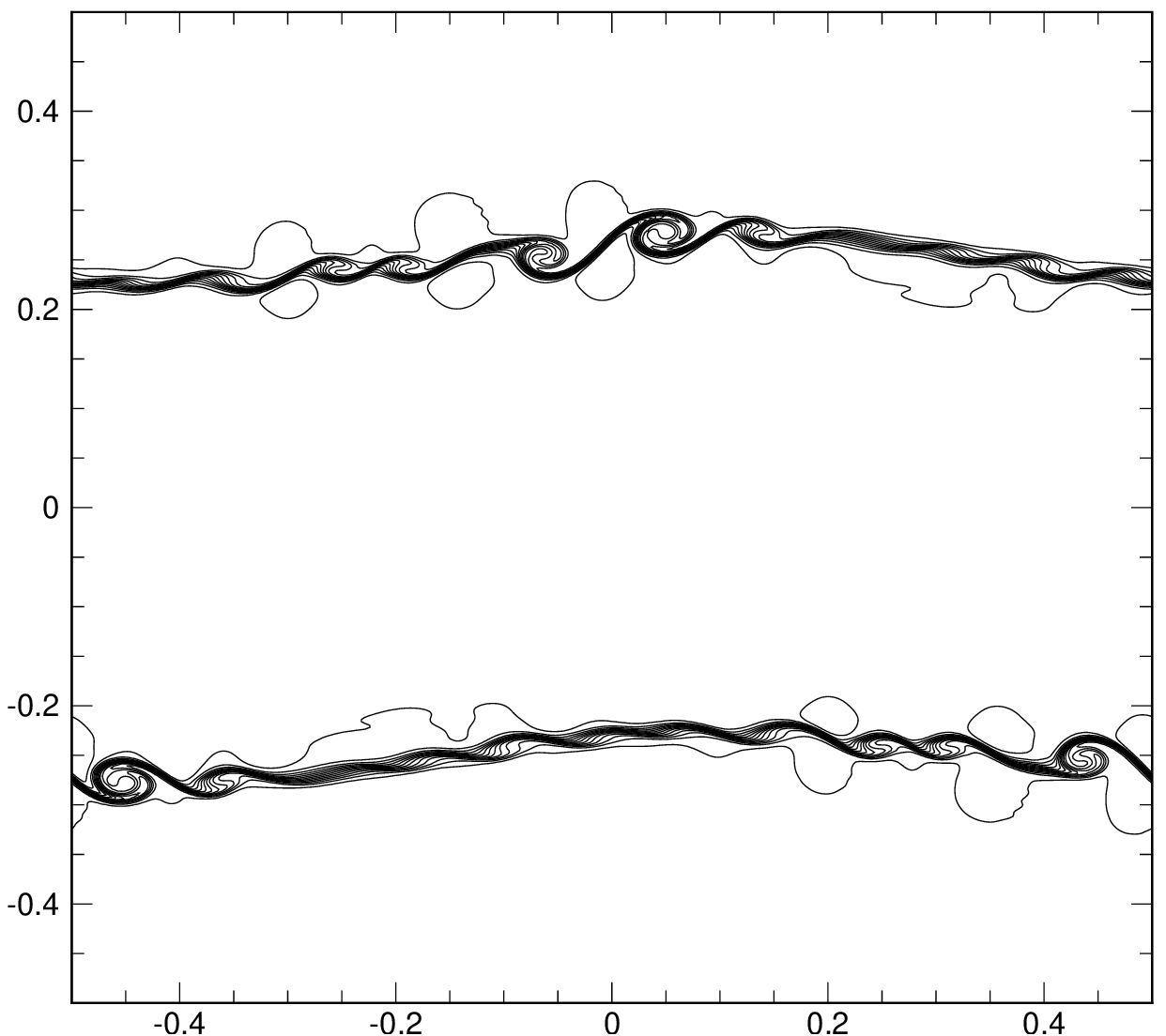}
        \caption{AWENO5, CH-CON}
    \end{subfigure}
    \begin{subfigure}{0.32\linewidth}
        \centering
        \includegraphics[width=\linewidth]{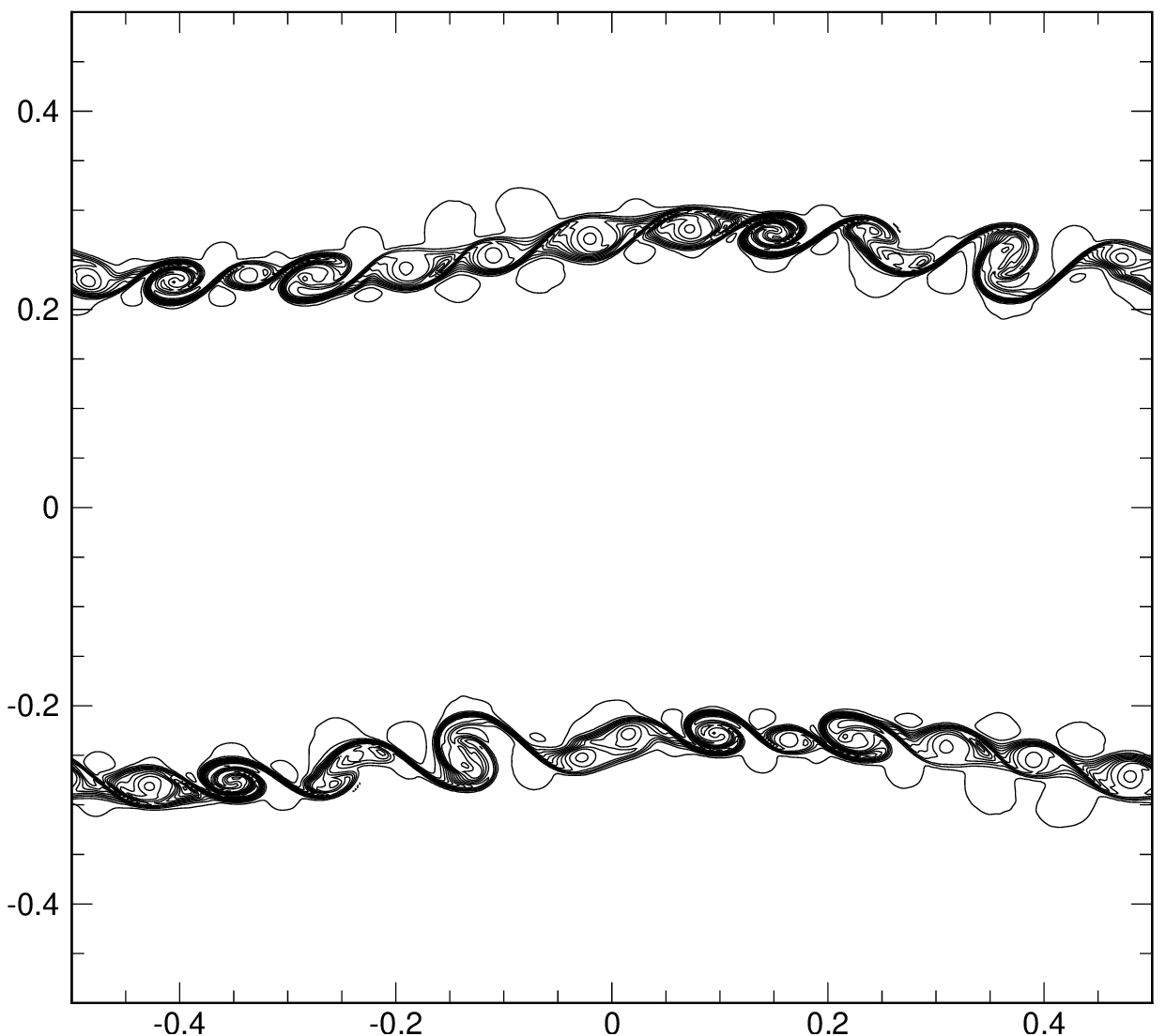}
        \caption{AWENO7, CH-CON}
    \end{subfigure}
    \begin{subfigure}{0.32\linewidth}
        \centering
        \includegraphics[width=\linewidth]{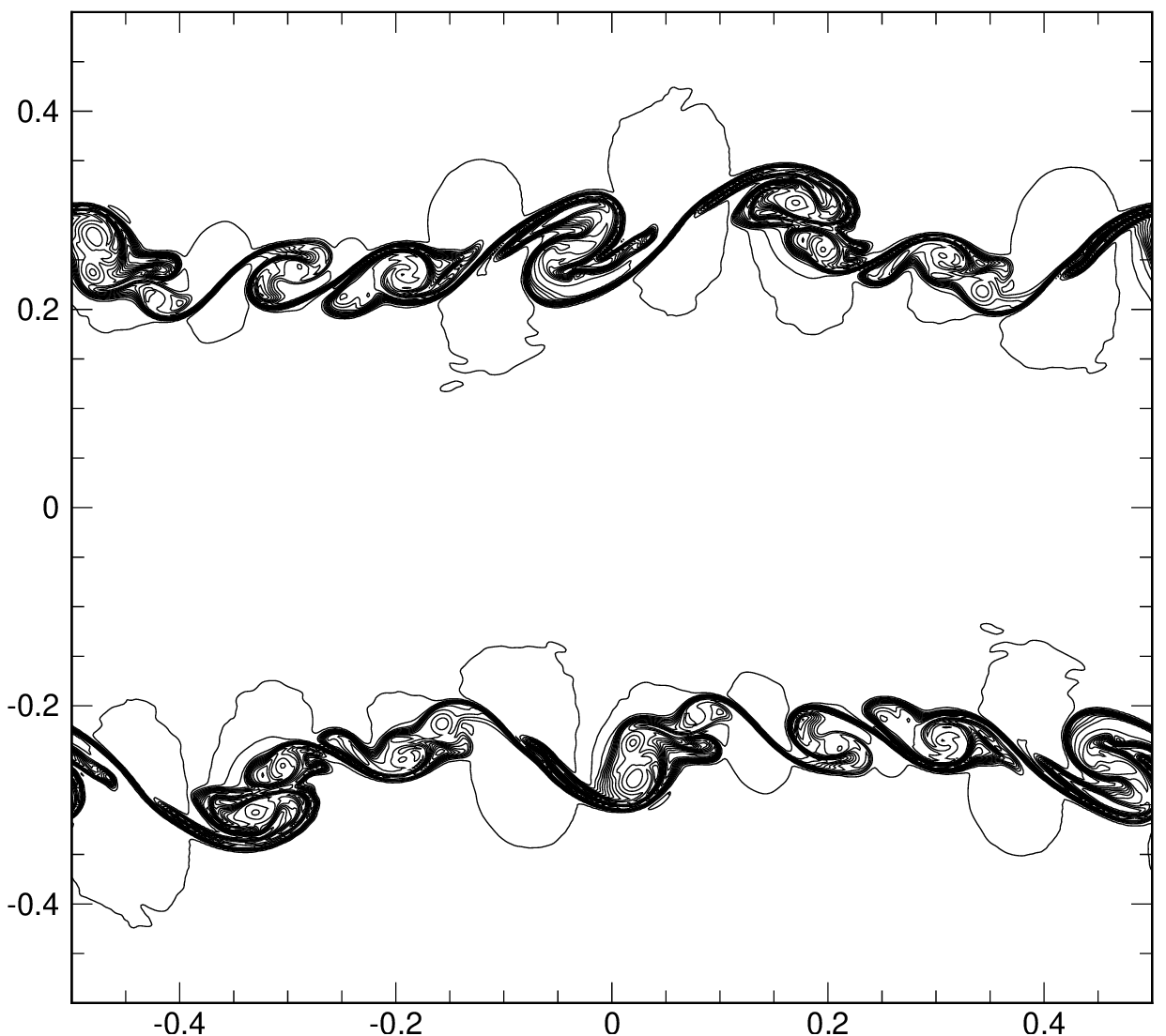}
        \caption{AWENO9, CH-CON}
    \end{subfigure}

    \caption{On the fine mesh \((Nx, Ny) = (512, 512)\).}
    \label{fig:KHI-2}
\end{figure}

We consider the Kelvin--Helmholtz instability (KHI) problem \cite{10.1016/j.compfluid.2015.04.026}. The computational domain is \([-0.5, 0.5] \times [-0.5, 0.5]\) with periodic boundary conditions in all directions. For \(\abs{y} \leqslant 0.25\), the initial condition is \((\rho, u, v, p) = (2, -0.5, 0.01 \sin(2 \pi x), 2.5)\), and for \(\abs{y} > 0.25\), the initial condition is \((\rho, u, v, p) = (1, 0.5, 0.01 \sin(2 \pi x), 2.5)\). The final simulation time is \(T = 1.0\). The results are shown in \Cref{fig:KHI-1,fig:KHI-2} using 20 equally spaced contours from \(0.7\) to \(2.2\). The CH-RI is more diffusive than the CH-CON in this example. 

\section{Conclusion}\label{se:con}

In this study, we have developed a finite difference A-WENO scheme with a Riemann invariant-based LCD procedure. We use the Riemann invariants as the transform variables, perform the LCD procedure on these variables and conduct WENO interpolation on the associated characteristic field. This change in the variable improves the sparsity of the eigenmatrix and thus greatly reduces the part of per-node computational cost that scales linearly with the stencil width or the order of accuracy increases. In the meanwhile, the scheme still maintains the convergence order of smooth solutions, satisfies the E-property and is free of spurious oscillations near discontinuities. Moreover, a combination of interpolation and flux PP limiters can make it as robust as those classical WENO schemes in the vicinity of low density or pressure. The scheme is tested on a variety of 1D and 2D problems to demonstrate its accuracy, efficiency and robustness. Our experiment shows that it can reduce \(23\%\) of the total time or \(88\%\) of the LCD time in comparison to the original characteristic-wise version at \(9\)th-order. With our cost analysis, greater acceleration can be expected for higher-order and 3D cases. 
This work merely focuses on the compressible Euler equations with the \(\gamma\)-law EOS. Generalizations to other EOS will probably change the Riemann invariants, but the sparsity pattern of the eigenmatrices will remain the same as long as the physical entropy exists, so such acceleration is still possible, and similar efficiency can be expected. Besides the compressible Euler equations, using Riemann invariants to reduce the LCD cost is also possible for other hyperbolic systems, such as the multi-component Euler equations, the magnetohydrodynamics (MHD) equations and the nonlinear elasticity equations.

\begin{appendices}

\section{Supplemental material for the compressible Euler equations}\label{ap:supplemental}

Following the notation in \Cref{def:hyperbolicity-in-1D}, the 2D compressible Euler equations \eqref{eq: Euler 2D} with the \(\gamma\)-law EOS \eqref{eq: EOS} is a hyperbolic system with the following eigenstructure in the \(x\)-direction. 
\begin{equation}\label{eq:Euler-eigen-values}
    \bm\Lambda = \diag\left(u-c, u, u, u+c\right), 
\end{equation}
\begin{subequations}\label{eq:Euler-right-eigen}
    \begin{eqnarray}
        && \bm r_1 = \frac{1}{2} \left(0, 1, 0, u\right)^T - \frac{1}{2 c} \left(1, u, v, H\right)^T, \\
        && \bm r_2 = -\frac{1}{c} \left(1, u, v, \frac{\abs{\vec{u}}^2}{2}\right)^T, \\
        && \bm r_3 = \left(0, 0, 1, v\right)^T, \\
        && \bm r_4 = \frac{1}{2} \left(0, 1, 0, u\right)^T + \frac{1}{2 c} \left(1, u, v, H\right)^T. 
    \end{eqnarray}
\end{subequations}
\begin{subequations}\label{eq:Euler-left-eigen}
    \begin{eqnarray}
        && \bm l_1 = \left(-u, 1, 0, 0\right) - \frac{\gamma-1}{c} \left(\frac{\abs{\vec u}^2}{2}, -u, -v, 1\right), \\
        && \bm l_2 = \frac{\gamma-1}{c} \left(\frac{\abs{\vec u}^2}{2}, -u, -v, 1\right) - \left(c, 0, 0, 0\right), \\
        && \bm l_3 = \left(-v, 0, 1, 0\right), \\
        && \bm l_4 = \left(-u, 1, 0, 0\right) + \frac{\gamma-1}{c} \left(\frac{\abs{\vec u}^2}{2}, -u, -v, 1\right). 
    \end{eqnarray}
\end{subequations}

Using the transform variables \(\bm v(\bm u) = (u - \frac{2}{\gamma-1}c, S^{\frac{1}{2\gamma}}, v, u + \frac{2}{\gamma-1}c)\), the transformed 2D compressible Euler equations \eqref{eq: Euler 2D} with the \(\gamma\)-law EOS \eqref{eq: EOS} have the following normalized eigenvectors in the \(x\)-direction: 
\begin{subequations}\label{eq:transform-right-eigen}
    \begin{eqnarray}
        && \tilde{\bm r_1} = \left(1, 0, 0, 0\right)^T, \\
        && \tilde{\bm r_2} = \left(-\frac{2\sqrt{\gamma} p^{\frac{\gamma-1}{2\gamma}}}{\gamma-1}, 1, 0, \frac{2\sqrt{\gamma} p^{\frac{\gamma-1}{2\gamma}}}{\gamma-1}\right)^T, \\
        && \tilde{\bm r_3} = \left(0, 0, 1, 0\right)^T, \\
        && \tilde{\bm r_4} = \left(0, 0, 0, 1\right)^T. 
    \end{eqnarray}
\end{subequations}
and 
\begin{subequations}\label{eq:transform-left-eigen}
    \begin{eqnarray}
        && \tilde{\bm l_1} = \left(1, \frac{2\sqrt{\gamma} p^{\frac{\gamma-1}{2\gamma}}}{\gamma-1}, 0, 0\right), \\
        && \tilde{\bm l_2} = \left(0, 1, 0, 0\right), \\
        && \tilde{\bm l_3} = \left(0, 0, 1, 0\right), \\
        && \tilde{\bm l_4} = \left(0, -\frac{2\sqrt{\gamma} p^{\frac{\gamma-1}{2\gamma}}}{\gamma-1}, 0, 1\right). 
    \end{eqnarray}
\end{subequations}

\begin{proposition}\label{prop:no-coordinate-Riemann-invariant}
    The compressible Euler equations \eqref{eq: Euler 2D} with the \(\gamma\)-law EOS \eqref{eq: EOS} is \emph{not} endowed with any coordinate system of Riemann invariants. 
    \begin{proof}
        Using the primitive variable \(\bm v = (\rho, u, v, p)\) as the transform variable, the normalized left eigenvectors are 
        \begin{eqnarray*}
            && \tilde{\bm l_1} = \left(0, 1, 0, -\frac{1}{\rho c}\right), \\
            && \tilde{\bm l_2} = \left(-\frac{\gamma}{\rho}, 0, 0, \frac{\gamma}{\rho c^2}\right), \\
            && \tilde{\bm l_3} = \left(0, 0, 1, 0\right), \\
            && \tilde{\bm l_4} = \left(0, 1, 0, \frac{1}{\rho c}\right). 
        \end{eqnarray*}
        Suppose there exists a coordinate system of Riemann invariants, then there has to be \(2,3,4\)-Riemann invariant \(w(\bm v)\) that satisfies \(\frac{\partial w}{\partial \bm v} \cdot \tilde{\bm r_i} = 0\) for \(i = 2, 3, 4\). This yields that there is a non-zero function \(\varphi(\bm v)\) such that \(\frac{\partial w}{\partial \bm v} = \varphi(\bm v) \tilde{\bm l_1}\). By the expression we deduce that \(w(\bm v)\) is independent of the density \(\rho\) and so is \(\varphi(\bm v)\) and \(\frac{\partial w}{\partial \bm v} = \varphi(\bm v) \tilde{\bm l_1}\). This is a contradiction since the last component of \(\varphi(\bm v) \tilde{\bm l_1}\) will depend on \(\rho\) in this case. 
    \end{proof}
\end{proposition}

\section{WENO-JS interpolation}\label{ap:WENO interp}

The \(k = (2r-1)\)th-order WENO-JS interpolation procedure for \(u_{j+\frac{1}{2}}^{-}\) is 
\begin{equation}
    u_{j+\frac{1}{2}}^{-} = \sum_{s=0}^{r-1} \omega_s u_{j+\frac{1}{2}}^{-,(s)},
\end{equation}
where 
\begin{equation}
    \omega_s = \frac{\alpha_s}{\sum_{\tau=0}^{r-1} \alpha_\tau}, 
\end{equation}
and 
\begin{equation}
    \alpha_s = \frac{d_s}{\left(\epsilon + \beta_s\right)^2}. 
\end{equation}
Here, the small constant \(\epsilon\) is usually taken to be \(10^{-6}\). Other parameters are given below. 

For \(k = 3\) (\(r = 2\)): 
\begin{subequations}
    \begin{eqnarray}
        & d_0 = \frac{3}{4}, & u_{j+\frac{1}{2}}^{-, (0)} = \frac{1}{2} u_{j} + \frac{1}{2} u_{j+1}, \\
        & d_1 = \frac{1}{4}, & u_{j+\frac{1}{2}}^{-, (1)} = -\frac{1}{2} u_{j-1} + \frac{3}{2} u_{j}. 
    \end{eqnarray}
\end{subequations}
\begin{subequations}
    \begin{eqnarray}
        && \beta_0 = \left(u_j - u_{j+1}\right)^2, \\
        && \beta_1 = \left(u_{j-1} - u_{j}\right)^2. 
    \end{eqnarray}
\end{subequations}

For \(k = 5\) (\(r = 3\)): 
\begin{subequations}
    \begin{eqnarray}
        & d_0 = \frac{5}{16}, & u_{j+\frac{1}{2}}^{-, (0)} = \frac{3}{8} u_{j} + \frac{3}{4} u_{j+1} - \frac{1}{8} u_{j+2}, \\
        & d_1 = \frac{5}{8}, & u_{j+\frac{1}{2}}^{-, (1)} = -\frac{1}{8} u_{j-1} + \frac{3}{4} u_{j} + \frac{3}{8} u_{j+1}, \\
        & d_2 = \frac{1}{16}, & u_{j+\frac{1}{2}}^{-, (2)} = \frac{3}{8} u_{j-2} - \frac{5}{4} u_{j-1} + \frac{15}{8} u_{j}. 
    \end{eqnarray}
\end{subequations}
\begin{subequations}
    \begin{eqnarray}
        && \beta_0 = \frac{13}{12} \left(u_{j} - 2 u_{j+1} + u_{j+2}\right)^2 + \frac{1}{4} \left(3 u_{j} - 4 u_{j+1} + u_{j+2}\right)^2, \\
        && \beta_1 = \frac{13}{12} \left(u_{j-1} - 2 u_{j} + u_{j+1}\right)^2 + \frac{1}{4} \left(u_{j-1} - u_{j+2}\right)^2, \\
        && \beta_2 = \frac{13}{12} \left(u_{j-2} - 2 u_{j-1} + u_{j}\right)^2 + \frac{1}{4} \left(u_{j-2} - 4 u_{j-1} + 3 u_{j}\right)^2. 
    \end{eqnarray}
\end{subequations}

For \(k = 7\) (\(r = 4\)): 
\begin{subequations}
    \begin{eqnarray}
        & d_0 = \frac{7}{64}, & u_{j+\frac{1}{2}}^{-, (0)} = \frac{5}{16} u_{j} + \frac{15}{16} u_{j+1} - \frac{5}{16} u_{j+2} + \frac{1}{16} u_{j+3}, \\
        & d_1 = \frac{35}{64}, & u_{j+\frac{1}{2}}^{-, (1)} = -\frac{1}{16} u_{j-1} + \frac{9}{16} u_{j} + \frac{9}{16} u_{j+1} - \frac{1}{16} u_{j+2}, \\
        & d_2 = \frac{21}{64}, & u_{j+\frac{1}{2}}^{-, (2)} = \frac{1}{16} u_{j-2} - \frac{5}{16} u_{j-1} + \frac{15}{16} u_{j} + \frac{5}{16} u_{j+1}, \\
        & d_3 = \frac{1}{64}, & u_{j+\frac{1}{2}}^{-, (3)} = -\frac{5}{16} u_{j-3} + \frac{21}{16} u_{j-2} - \frac{35}{16} u_{j-1} + \frac{35}{16} u_{j}. 
    \end{eqnarray}
\end{subequations}
\begin{subequations}
    \begin{eqnarray}
        && \begin{aligned}
            \beta_0 = & \frac{\left(-15 u_{j} + 25 u_{j+1} - 13 u_{j+2} + 3 u_{j+3}\right)^2}{64} + \frac{13}{12} \left(2 u_{j} - 5 u_{j+1} + 4 u_{j+2} - u_{j+3}\right)^2 \\& + \frac{781}{720} \left(-u_{j} + 3 u_{j+1} - 3 u_{j+2} + u_{j+3}\right)^2, 
        \end{aligned} \\
        && \begin{aligned}
            \beta_1 = & \frac{\left(-3 u_{j-1} - 3 u_{j} + 7 u_{j+1} - u_{j+2}\right)^2}{64} + \frac{13}{12} \left(u_{j-1} - 2 u_{j} + u_{j+1}\right)^2 \\& + \frac{781}{720} \left(-u_{j-1} + 3 u_{j} - 3 u_{j+1} + u_{j+2}\right)^2, 
        \end{aligned} \\
        && \begin{aligned}
            \beta_2 = & \frac{\left(u_{j-2} - 7 u_{j-1} + 3 u_{j} + 3 u_{j+1}\right)^2}{64} + \frac{13}{12} \left(u_{j-1} - 2 u_{j} + u_{j+1}\right)^2 \\& + \frac{781}{720} \left(-u_{j-2} + 3 u_{j-1} - 3 u_{j} + u_{j+1}\right)^2, 
        \end{aligned} \\
        && \begin{aligned}
            \beta_3 = & \frac{\left(-3 u_{j-3} + 13 u_{j-2} - 25 u_{j-1} + 15 u_{j}\right)^2}{64} + \frac{13}{12} \left(-u_{j-3} + 4 u_{j-2} - 5 u_{j-1} + 2 u_{j}\right)^2 \\& + \frac{781}{720} \left(-u_{j-3} + 3 u_{j-2} - 3 u_{j-1} + u_{j}\right)^2. 
        \end{aligned}
    \end{eqnarray}
\end{subequations}

For \(k = 9\) (\(r = 5\)): 
\begin{subequations}
    \begin{eqnarray}
        & d_0 = \frac{9}{256}, & u_{j+\frac{1}{2}}^{-, (0)} = \frac{35}{128} u_{j} + \frac{35}{32} u_{j+1} - \frac{35}{64} u_{j+2} + \frac{7}{32} u_{j+3} - \frac{5}{128} u_{j+4}, \\
        & d_1 = \frac{21}{64}, & u_{j+\frac{1}{2}}^{-, (1)} = -\frac{5}{128} u_{j-1} + \frac{15}{32} u_{j} + \frac{45}{64} u_{j+1} - \frac{5}{32} u_{j+2} + \frac{3}{128} u_{j+3}, \\
        & d_2 = \frac{63}{128}, & u_{j+\frac{1}{2}}^{-, (2)} = \frac{3}{128} u_{j-2} - \frac{5}{32} u_{j-1} + \frac{45}{64} u_{j} + \frac{15}{32} u_{j+1} - \frac{5}{128} u_{j+2}, \\
        & d_3 = \frac{9}{64}, & u_{j+\frac{1}{2}}^{-, (3)} = -\frac{5}{128} u_{j-3} + \frac{7}{32} u_{j-2} - \frac{35}{64} u_{j-1} + \frac{35}{32} u_{j} + \frac{35}{128} u_{j+1}, \\
        & d_4 = \frac{1}{256}, & u_{j+\frac{1}{2}}^{-,(4)} = \frac{35}{128} u_{j-4} - \frac{45}{32} u_{j-3} + \frac{189}{64} u_{j-2} - \frac{105}{32} u_{j-1} + \frac{315}{128} u_{j}. 
    \end{eqnarray}
\end{subequations}
\begin{subequations}
    \begin{eqnarray}
        && \begin{aligned}
            \beta_0 = & \frac{\left(-35 u_j + 70 u_{j+1} - 56 u_{j+2} + 26 u_{j+3} - 5 u_{j+4}\right)^2}{256} \\& + \frac{\left(4613 u_{j} - 13772 u_{j+1} + 15198 u_{j+2} - 7532 u_{j+3} + 1493 u_{j+4}\right)^2}{2246400} \\& + \frac{781}{2880} \left(-5 u_{j} + 18 u_{j+1} - 24 u_{j+2} + 14 u_{j+3} - 3 u_{j+4}\right)^2 \\& + \frac{1421461}{1310400} \left(u_{j} - 4 u_{j+1} + 6 u_{j+2} - 4 u_{j+3} + u_{j+4}\right)^2, 
        \end{aligned} \\
        && \begin{aligned}
            \beta_1 = & \frac{\left(-5 u_{j-1} - 10 u_{j} + 20 u_{j+1} - 6 u_{j+2} + u_{j+3}\right)^2}{256} \\& + \frac{\left(1493 u_{j-1} - 2852 u_{j} + 1158 u_{j+1} + 268 u_{j+2} - 67 u_{j+3}\right)^2}{2246400} \\& + \frac{781}{2880} \left(-3 u_{j-1} + 10 u_{j} - 12 u_{j+1} + 6 u_{j+2} - u_{j+3}\right)^2 \\& + \frac{1421461}{1310400} \left(u_{j-1} - 4 u_{j} + 6 u_{j+1} - 4 u_{j+2} + u_{j+3}\right)^2, 
        \end{aligned} \\
        && \begin{aligned}
            \beta_2 = & \frac{\left(u_{j-2} - 10 u_{j-1} + 10 u_{j+1} - u_{j+2}\right)^2}{256} \\& + \frac{\left(-67 u_{j-2} + 1828 u_{j-1} - 3522 u_{j} + 1828 u_{j+1} - 67 u_{j+2}\right)^2}{2246400} \\& + \frac{781}{2880} \left(-u_{j-2} + 2 u_{j-1} - 2 u_{j+1} + u_{j+2}\right)^2 \\& + \frac{1421461}{1310400} \left(u_{j-2} - 4 u_{j-1} + 6 u_{j} - 4 u_{j+1} + u_{j+2}\right)^2, 
        \end{aligned} \\
        && \begin{aligned}
            \beta_3 = & \frac{\left(-u_{j-3} + 6 u_{j-2} - 20 u_{j-1} + 10 u_{j} + 5 u_{j+1}\right)^2}{256} \\& + \frac{\left(-67 u_{j-3} + 268 u_{j-2} + 1158 u_{j-1} - 2852 u_{j} + 1493 u_{j+1}\right)^2}{2246400} \\& + \frac{781}{2880} \left(u_{j-3} - 6 u_{j-2} + 12 u_{j-1} - 10 u_{j} + 3 u_{j+1}\right)^2 \\& + \frac{1421461}{1310400} \left(u_{j-3} - 4 u_{j-2} + 6 u_{j-1} - 4 u_{j} + u_{j+1}\right)^2, 
        \end{aligned} \\
        && \begin{aligned}
            \beta_4 = & \frac{\left(5 u_{j-4} - 26 u_{j-3} + 56 u_{j-2} - 70 u_{j-1} + 35 u_{j}\right)^2}{256} \\& + \frac{\left(1493 u_{j-4} - 7532 u_{j-3} + 15198 u_{j-2} - 13772 u_{j-1} + 4613 u_{j}\right)^2}{2246400} \\& + \frac{781}{2880} \left(3 u_{j-4} - 14 u_{j-3} + 24 u_{j-2} - 18 u_{j-1} + 5 u_{j}\right)^2 \\& + \frac{1421461}{1310400} \left(u_{j-4} - 4 u_{j-3} + 6 u_{j-2} - 4 u_{j-1} + u_{j}\right)^2. 
        \end{aligned}
    \end{eqnarray}
\end{subequations}

We remark here that to derive the classical smoothness indicator up to \(r = 5\), we employ the following sum-of-squares (SOS) decomposition that is similar to the one given in \cite[Appendix B]{10.1016/j.compfluid.2020.104519}. For any \(p(x) \in \mathcal{P}^{4}([-\frac{h}{2}, \frac{h}{2}])\), 
\begin{equation}
    \begin{aligned}
        & \int_{-\frac{h}{2}}^{\frac{h}{2}} \sum_{l=1}^{\infty} h^{2l-1} \left(p^{(l)}(x)\right)^2 \d x \\
        = & \left(h p^{(1)}(0) + \frac{1}{24} h^3 p^{(3)}(0)\right)^2 + \frac{\left(520 h^2 p^{(2)}(0) + 21 h^4 p^{(4)}(0)\right)^2}{249600} + \frac{781}{720} \left(h^3 p^{(3)}(0)\right)^2 + \frac{1421461}{1310400} \left(h^4 p^{(4)}(0)\right)^2. 
    \end{aligned}
\end{equation}
Then, to compute the smoothness indicator \(\beta\)'s in any cell \([x_{j-\frac{1}{2}}, x_{j+\frac{1}{2}}]\), it only remains to provide the interpolation-based derivative values at the cell center \(x_j\).

\section{High-order flux correction terms}\label{ap:high-order corr}

We give the \(2 r = (k + 1)\)th-order central difference implementation here. All the coefficients are computed by \texttt{Mathematica}, and are verified to be consistent with those in \cite{10.1016/j.compfluid.2020.104519} and \cite[Appendix A]{10.1007/s42967-023-00360-z}, to ensure that 
\begin{equation}
    \frac{1}{\Delta x} \left(\left(\bm f_{j+\frac{1}{2}} + \bm f_{j+\frac{1}{2}}^{\mathrm{cor}, r}\right) - \left(\bm f_{j-\frac{1}{2}} + \bm f_{j-\frac{1}{2}}^{\mathrm{cor}, r}\right) \right) = \frac{\d}{\d x}\bm f|_{j} + \mathcal{O}(\Delta x^{2r}). 
\end{equation}

For \(k = 1\) (\(r = 1\)):  
\begin{equation}
    \bm f_{j+\frac{1}{2}}^{\mathrm{cor}, 1} = 0. 
\end{equation}

For \(k = 3\) (\(r = 2\)):  
\begin{equation}
    \bm f_{j+\frac{1}{2}}^{\mathrm{cor}, 2} 
    = -\frac{1}{48} \left(\bm f_{j-1} + \bm f_{j+2}\right) + \frac{1}{48} \left(\bm f_{j} + \bm f_{j+1}\right). 
\end{equation}

For \(k = 5\) (\(r = 3\)):  
\begin{equation}
    \bm f_{j+\frac{1}{2}}^{\mathrm{cor}, 3} 
    = \frac{19}{3840} \left(\bm f_{j-2} + \bm f_{j+3}\right) - \frac{137}{3840} \left(\bm f_{j-1} + \bm f_{j+2}\right) + \frac{59}{1920} \left(\bm f_{j} + \bm f_{j+1}\right). 
\end{equation}

For \(k = 7\) (\(r = 4\)): 
\begin{equation}
    \bm f_{j+\frac{1}{2}}^{\mathrm{cor}, 4} 
    = -\frac{81}{71680} \left(\bm f_{j-3} + \bm f_{j+4}\right) + \frac{2279}{215040} \left(\bm f_{j-2} + \bm f_{j+3}\right) - \frac{9859}{215040} \left(\bm f_{j-1} + \bm f_{j+2}\right) + \frac{7823}{215040} \left(\bm f_{j} + \bm f_{j+1}\right). 
\end{equation}

For \(k = 9\) (\(r = 5\)):  
\begin{equation}
    \begin{aligned}
        \bm f_{j+\frac{1}{2}}^{\mathrm{cor}, 5} 
        = & \frac{5359}{20643840} \left(\bm f_{j-4} + \bm f_{j+5}\right) - \frac{60841}{20643840} \left(\bm f_{j-3} + \bm f_{j+4}\right) + \frac{81491}{5160960} \left(\bm f_{j-2} + \bm f_{j+3}\right) \\& - \frac{274129}{5160960} \left(\bm f_{j-1} + \bm f_{j+2}\right) + \frac{413017}{10321920} \left(\bm f_{j} + \bm f_{j+1}\right). 
    \end{aligned}
\end{equation}

\end{appendices}



\bibliographystyle{amsplain}
\bibliography{refs}


\end{document}